\documentclass[a4paper, 11pt, openany]{report} % A4 paper size and default 11pt font size

\usepackage{amsmath}
\usepackage{amssymb}
\usepackage{enumerate}
\usepackage{graphicx}
\usepackage{float}
\usepackage{longtable}
\usepackage{tabularx}
\graphicspath{ {./images/} }

\newtheorem{theorem}{Theorem}
\newcolumntype{C}{>{\arraybackslash}X}

\begin{document}
\title{\sc{Research Problems in Function Theory\\ \vspace{1cm}Draft Copy - Replies Deadline Jan 6$^{\sc{th}}$ 2019}}

\author{W. K. Hayman \& E. F. Lingham}
\date{}
\maketitle

\tableofcontents

\chapter*{A Note for the Reader}

\noindent
Thank you for your interest in this book. 

\vspace{0.4cm}
\noindent
In order for this work to be as useful to the complex analysis community as possible, we want it to be as comprehensive and as accurate as can be. We also aim to include new problems in all nine chapters.

\vspace{0.4cm}
\noindent
We therefore ask for your help in three ways:

\vspace{0.2cm}
\begin{enumerate}[{\bf (a)}]
\item \textbf{Let us know of any errors and omissions.}\\
While you are reading this work, you may come across some errors or omissions. Please let us know about these by emailing Eleanor Lingham (e.lingham@shu.ac.uk). 

\item \textbf{Send us some new problems.}\\
For any mathematics research area to thrive, there needs to be an abundance of interesting problems to work on. If you have some interesting unsolved problems, we ask you to please consider sharing them with the wider mathematics community by publishing them in this work. 

\noindent To do so, email  Eleanor Lingham (e.lingham@shu.ac.uk) with the following:
\begin{enumerate}[(i)]
\item A clear statement of the problem, as you would like it to appear.
\item Full details of any definitions needed (you could check this work to see if the definition is already included).
\item Full details of any references needed.
\item A direction to which chapter the problem should appear in.
\item An attribution for the problem, if one exists.
\end{enumerate}

\item \textbf{Please let other interested members of the complex analysis community know about this project.}\\
\end{enumerate}

\noindent\textbf{\underline{Please send all replies by January 6$^{th}$ 2019 - Walter K. Hayman's 93$^{rd}$ birthday.}}\\

\vspace{0.4cm}
\noindent
With best regards and thanks,

\vspace{0.4cm}
\noindent
Walter Hayman and Eleanor Lingham \hspace{2.75cm}\textbf{Contact:} e.lingham@shu.ac.uk

\chapter*{Introduction}
By Walter K. Hayman and Eleanor F. Lingham (to appear)

\chapter{Meromorphic Functions}

\underline{Preface by A. Eremenko (to appear)}\\

\noindent\textbf{Notation} We use the usual notations of Nevanlinna theory (see for example Nevanlinna \cite{MR0417418}, \cite{MR0344426} and Hayman \cite{M.F.}).\\ 
If $f(z)$ is meromorphic in $|z|<R$, and $0<r<R$ we write 
\[n(r,a)=n(r,a,f)\]
for the number of roots of the equation $f(z)=a$ in $|z|\leq r$, when \textit{multiple roots are counted according to multiplicity}, and $\overline{n}(r,a)$ when \textit{multiple roots are counted only once}. We also define
\[N(r,a)=\int^r_0\frac{[ n(t,a)-n(0,a)] \,dt}{t}+n(0,a)\log r,\]
\[\overline{N}(r,a)=\int^r_a\frac{[\overline{n}(t,a)-\overline{n}(0,a)]\,dt}{t}+\overline{n}(0,a)\log r,\]
\[m(r,f)=m(r, \infty,f)=\frac{1}{2\pi}\int^{2\pi}_0\log^+|f(re^{i\theta})|\,d\theta,\]
where $\log^+x=\max\big\{\log x,0\big\}$,
\[m(r,a,f)=m\left(r,\infty,\frac{1}{f-a}\right), a\neq\infty.\]
and 
\[T(r,f)=m(r,\infty,f)+N(r,\infty,f).\]
Then for every finite $a$, we have by the first fundamental theorem (see Hayman \cite[p. 5]{M.F.}),
\begin{equation}
\label{1.1}
T(r,f)=m(r,a,f)+N(r,a,f)+O(1), \hspace{1cm}\text{ as }r\to R.
\end{equation}
We further define the deficiency,
\[\delta(a,f)=\liminf_{r\to R}\frac{m(r,a,f)}{T(r,f)}=1-\limsup_{r\to R}\frac{N(r,a,f)}{T(r,f)},\]
the Valiron deficiency
\[\Delta(a,f)=\limsup_{r\to R}\frac{m(r,a,f)}{T(r,f)},\]
and further
\[\Theta(a,f)=1-\limsup_{r\to R}\frac{\overline{N}(r,a,f)}{T(r,f)}.\]
We then have the \textit{``defect relation''} (see Hayman \cite[p. 43]{M.F.}),
\begin{equation}
\label{1.2}
\sum_a \delta(a,f)\leq\sum_a \Theta(a,f)\leq 2,
\end{equation}
provided that either $R=\infty$ and $f(z)$ is not constant, or $R<+\infty$ and
\[\limsup_{r\to R}\frac{T(r,f)}{\log\big(1/(R-r)\big)}=+\infty.\]
If $R=+\infty$ we also define the \textit{lower order} $\lambda$ and \textit{order} $\rho$,
\[\lambda=\liminf_{r\to R}\frac{\log T(r,f)}{\log r},\hspace{1cm}\rho=\limsup_{r\to R}\frac{\log T(r,f)}{\log r}.\]
If $\delta(a,f)>0$ the value $a$ is called \textit{deficient}. If follows from (\ref{1.2}) that there are at most countably many deficient values if the conditions for (\ref{1.2}) are satisfied.\\

\noindent\textbf{Problem 1.1} Is (\ref{1.2}) all that is true in general? In other words, can we construct a meromorphic function $f(z)$ such that $f(z)$ has an arbitrary sequence $a_n$ of deficient values and no others, and further that $\delta(a_n,f)=\delta_n$, where $\delta_n$ is an arbitrary sequence subject to $\sum\delta_n\leq 2$? (If $f(z)$ is an entire function $\delta(f,\infty)=1$, so that $\sum_{a\neq\infty}\delta(a,f)\leq1$. For a solution of the problem in this case see \cite[p. 80]{M.F.}).\\

\noindent\textbf{Update 1.1}
This problem has been completely settled by Drasin \cite{zbMATH03455609}. He constructs a meromorphic function $f(z)$ with arbitrary deficiencies and branching indices on a presassigned sequence $a_n$ of complex numbers with $f(z)$ growing arbitrarily slowly, subject to having infinite order.\\

\noindent\textbf{Problem 1.2} How big can the set of Valiron deficiencies be for functions in the plane? It is known that
\begin{equation}
\label{1.3}
N(r,a)=T(r,f)+O\big(T(r,f)^{\frac{1}{2}+\varepsilon}\big)
\end{equation}
as $r\to\infty$, for all $a$ outside a set of capacity zero (see Nevanlinna \cite[pp. 260-264]{MR0344426}.\\
In the case $R<+\infty$ this is more or less best-possible, but in the plane we only know from an example of Valiron \cite{Valiron1925} that the corresponding set of $a$ can be non-countably infinite. It is also not known whether (\ref{1.3}) can be sharpened.\\

\noindent\textbf{Update 1.2}
Hyllengren \cite{Hyllengren1970} has shown that all values of a set $E$ can have Valiron deficiency greater than a positive constant for a function of finite order in the plane, if and only if there exists a sequence of complex numbers $a_n$ and a $k>0$ such that each point of $E$ lies in infinitely many of the discs $\{z: |z-a_n|<e^{-e^{kn}}\}$.\\
Hayman \cite{Hayman1972} proved that all values of any $F_\sigma$ set of capacity zero can be Valiron deficiencies for an entire function of infinite order and a little more. \\

\noindent\textbf{Problem 1.3} If $f(z)$ is meromorphic of finite order $\rho$ and $\sum\delta(a,f)=2$, it is conjectured that $\rho=n/2$, where $n$ is an integer and $n\geq 2$, and all the deficiencies are rational. F. Nevanlinna \cite{FNevanlinna1929} has proved this result on the condition that $f(z)$ has no multiple values, so that $n(r,a)=\overline{n}(r,a)$ for every $a$ (see also R. Nevanlinna \cite{Newvanlinna1932}).\\

\noindent\textbf{Update 1.3}
Weitsmann proved that the number of deficiencies is at most twice the order in this case \cite{Weitsman1969}. The conjecture was completely proved by Drasin \cite{Drasin1987}.  Eremenko \cite{MR1266090} gave a simpler
proof of a stronger result, see Problem 1.33. \\

\noindent\textbf{Problem 1.4} Let $f(z)$ be an entire function of finite order $\rho$, and let $n_1(r,a)$ denote the number of simple zeros of the equation $f(z)=a$. If 
\[n_1(r,a)=O(r^c),\hspace{1cm} n_1(r,b)=O(r^c),\hspace{1cm} \text{ as } r\to\infty,\]
where $a\neq b, c<\rho$, is it true that $\rho$ is an integral multiple of $\frac{1}{2}$? More strongly, is this result true if $\Theta(a)=\frac{1}{2}=\Theta(b)$? (For a somewhat weaker result in this direction, see Gol'dberg and Tairova \cite{GoldbergTairova}.)\\
  
\noindent\textbf{Update 1.4}
The answer is `no', even in a very weak sense. Gol'dberg \cite{Gol'dberg1973} has constructed an example of an entire function for which 
\[n_1(r,a)=O((\log r)^{2+\varepsilon}), \hspace{1cm} n_1(r,b)=O((\log r)^{2+\varepsilon}) \hspace{1cm}\text{ as }r\to\infty,\]
but the order is not a multiple of $\frac{1}{2}$.\\

\noindent\textbf{Problem 1.5} Under what conditions can $\sum\delta(a,f)$ be nearly $2$ for an entire function of finite order $\rho$? Pfluger \cite{Pfluger1946} proved that if $\sum\delta(a,f)=2$, then (see Hayman \cite[p.115]{M.F.})
$\rho$ is a positive integer $q$, the lower order $\lambda$ is such that $\lambda=\rho$ and all the deficiencies are integral multiplicities of $1/q$. If further
\[\sum\delta(a,f)>2-\varepsilon(\lambda),\]
where $\varepsilon(\lambda)$  is a positive quantity depending on $\lambda$, then Edrei and Fuchs (\cite{10.2307/1993455}, \cite{EdreiAlbert;FuchsWolfgangH.J.1959Vdev}) proved that these results remain true `nearly', in the sense that there exist `large' deficiencies which are nearly positive integral multiplicities of $1/q$, and whose sum of deficiencies is `nearly' $2$. Can there be a finite or infinite number of small deficiences as well in this case? \\

\noindent\textbf{Update 1.5}
No progress on this problem has been reported to us. Hayman suspects that the answer is `no'.\\

\noindent\textbf{Problem 1.6} Arakelyan \cite{MR0206286} has proved that, given $\rho>\frac{1}{2}$ and a countable set $E$, there exists an entire function $f(z)$ of order $\rho$, for which all the points of $E$ are deficient. Can $E$ be the precise set of deficiencies of $f$ in the sense that $f$ has no other deficient values? It is also conjectured that if the $a_n$ are deficient values for an entire function of finite order, then
\[\sum\big(\log[1/\delta(a_n,f)]\big)^{-1}<+\infty.\]
\textit{(N. U. Arakelyan)}\\

\noindent\textbf{Update 1.6}
Eremenko \cite{Eremenko1980} has proved the first conjecture. He also proved \cite{MR1145577} that the second conjecture is false: given $\rho>1/2$ and a sequence of complex numbers $(a_k)$, there is an entire function $f$ of order $\rho$ with the property $\delta(a_k,f)>c^k, k=1,2,\ldots,$ for some $c\in(0,1)$. On the other hand, Lewis and Wu \cite{MR942832} proved $\sum \delta(a_k, f)^\alpha<\infty$ for entire functions of finite order with an absolute constant $\alpha<1/3 - 2^{-264}$. The exact rate of decrease of deficiencies of an entire function of finite order remains unknown.\\  

\noindent\textbf{Problem 1.7} If $f(z)$ is an entire function of finite order $\rho$ which is not an integer, it is known that (see Pfluger \cite{Pfluger1946} and Hayman \cite[p.104]{M.F.}),
\[\sum \delta(a,f)\leq 2-K(\rho)\]
where $K(\rho)$ is a positive quantity depending on $\rho$. What is the best possible value for $K(\rho)$? Edrei and Fuchs \cite{EdreiAlbert;FuchsWolfgangH.J.1959Vdev} conjectured (see also Hayman \cite[p. 104]{M.F.}) that if $q$ is the integral part of $\rho$, and if $q\geq1$, then
\[K(\rho)=\frac{|\sin(\pi\rho)|}{q+|\sin(\pi\rho)|},\hspace{1cm}q\leq\rho<q+\frac{1}{2},\]
\[K(\rho)=\frac{|\sin(\pi\rho)|}{q+1}, \hspace{1cm}q+\frac{1}{2}\leq\rho<q+1.\]
This result would be sharp.\\
If $\rho\leq\frac{1}{2}$, there are no deficient values, so that $K(\rho)=1$. If $\frac{1}{2}<\rho<1$, Pfluger \cite{Pfluger1946} proved  that $K(\rho)=\sin(\pi\rho)$. See also Hayman \cite[p. 104]{M.F.}.\\

\noindent\textbf{Update 1.7} 
For Problem 1.7 and Problem 1.8 a better lower bound was found by Miles and Shea \cite{MR0324041}, who also obtained the exact lower bound for any order $\rho$ of 
\[\limsup_{r\to\infty}\frac{N(r,0)+N(r,\infty)}{m_2(r)},\]
where 
\[m_2(r)=\Big(\frac{1}{2\pi}\int^{2\pi}_0(\log|f(re^{i\theta})|)^2\,d\theta\Big)^{\frac{1}{2}}.\]
Hellerstein and Williamson \cite{MR0247087} have solved the problems completely for entire functions with zeros on a ray.\\

\noindent\textbf{Problem 1.8} Following the notation in Problem 1.7, if $f(z)$ is meromorphic in the plane of order $\rho$, it is conjectured by Pfluger \cite{Pfluger1946}, that for $a\neq b$
\[\limsup_{r\to \infty}\frac{N(r,a)+N(r,b)}{T(r,f)}\geq K(\rho).\]
This is known to be true for $0<\rho\leq1$. If equality holds in the above inequality, it is conjectured that $f(z)$ has regular growth, i.e. $\rho=\lambda$.\\

\noindent\textbf{Update 1.8}
See Update 1.7.\\

\noindent\textbf{Problem 1.9} If $f(z)$ is an entire function of finite order $\rho$ which has a finite deficient value, find the best possible lower bound for the lower order $\lambda$ of $f(z)$. (Edrei and Fuchs \cite{EdreiAlbert;FuchsWolfgangH.J.1959Vdev} showed that $\lambda>0$.)\\
Gol'dberg \cite{MR0170006} showed that for every $\rho>1, \lambda \geq1$ is possible.\\

\noindent\textbf{Update 1.9}
This had been settled by Gol'dberg in \cite{MR0170006}.\\

\noindent\textbf{Problem 1.10} If $f(z)$ is a meromorphic function of finite order with more than two deficient values, is it true that if $\sigma>1$, then
\[\limsup_{r\to\infty}\frac{T(\sigma r)}{T(r)}<+\infty.\]\\

\noindent\textbf{Update 1.10}
No progress on this problem has been reported to us. (The update in \cite{Bprevious} has been withdrawn).\\

\noindent\textbf{Problem 1.11} If $f(z)$ is a meromorphic function of finite order with at least one finite deficient value, does the conclusion of Problem 1.10 hold?\\

\noindent\textbf{Update 1.11}
Drasin writes that Kotman \cite{MR553047} has shown that the answer to this question is `no'.\\

\noindent\textbf{Problem 1.12} Edrei, Fuchs and Hellerstein \cite{MR0122999} ask if $f(z)$ is an entire function of infinite order with real zeros, is $\delta(0,f)>0$? More generally, is $\delta(0,f)=1$?\\

\noindent\textbf{Update 1.12}
This has been disproved by Miles \cite{MR543739}, who showed that $\delta(0,f)=0$ is possible. However, Miles also showed 
\[\frac{N(r,0)}{T(r,f)}\to 0\]
as $r\to\infty$ outside a fairly small set in this case.\\

\noindent\textbf{Problem 1.13} If $f(z)$ is an entire function of finite order $\rho$ and lower order $\lambda$ with real zeros, find the best possible bound $B=B(\rho,\lambda)$ such that $\delta(0,f)\geq B$. From Edrei, Fuchs and Hellerstein \cite{MR0122999} it is known that $B>0$ if $2<\rho<\infty$, and it is conjectured that $B\to1$ as $\rho\to+\infty$.\\

\noindent\textbf{Update 1.13}
An affirmative answer with the exact value of $B(\lambda)$ was given by Hellerstein and Shea \cite{MR0393486}. \\

\noindent\textbf{Problem 1.14} If $f(z)$ is a meromorphic function of finite order, then it is known (see Hayman \cite[pp. 90, 98]{M.F.}) that $\sum\delta(a,f)^\alpha$ converges if $\alpha>\frac{1}{3}$, but may diverge if  $\alpha<\frac{1}{3}$. What happens when $\alpha=\frac{1}{3}$?\\

\noindent\textbf{Update 1.14}
This has been completely settled by Weitsman \cite{MR0387597} who proved that $\sum\big(\delta(a,f)\big)^{1/3}$ does indeed converge for any meromorphic function of finite order.\\

\noindent\textbf{Problem 1.15} \textbf{(Edrei's spread conjecture)} If $f(z)$ is meromorphic in the plane and of lower order $\lambda$, and if $\delta=\delta(a,f)>0$, is it true that, for a sequence $r=r_\nu\to\infty$, $f(z)$ is close to $a$ on a part of the circle $|z|=r_\nu$ having angular measure at least
\[\frac{4}{\lambda}\sin^{-1}\sqrt{\left(\frac{\delta}{2}\right)}+o(1)?\]
(For a result in this direction, see Edrei \cite{MR0180678}.)\\

\noindent\textbf{Update 1.15} 
This result was proved by Baernstein \cite{MR0374429} by means of the function $T^*(r,\theta)$, where
\[T^*(r,\theta)=\sup_E\frac{1}{2\pi}\int_E\log|f(re^{i\varphi})|\,d\varphi+N(r,f),\]
where $E$ runs over all sets of measure exactly $2\theta$. See also Baernstein \cite{MR0294640}.\\

\noindent\textbf{Problem 1.16} For any function $f(z)$ in the plane, let $n(r)=\sup_a n(r,a)$ be the maximum number of roots of the equation $f(z)=a$ in $|z|<r$, and 
\[A(r) = \frac{1}{\pi}\int\int_{|z|<r}\frac{|f'(z)|^2}{\{1+|f(z)|^2\}^2}\,dx \,dy =  \frac{1}{\pi}\int\int_{|a|<\infty}\frac{n(r,a)\,|da|^2}{(1+|a|^2)^2}.\]
Then $\pi A(r)$ is the area, with due count of multiplicity, of the image on the Riemann sphere of the disc $|z|<r$ under $f$, and $A(r)$ is the average value of $n(r,a)$ as $a$ moves over the Riemann sphere. It is known (see Hayman \cite[p. 14]{M.F.}) that
\[1\leq\liminf_{r\to\infty}\frac{n(r)}{A(r)}\leq e.\]
Can $e$ be replaced by any smaller quantity, and in particular, by $1$?\\

\noindent\textbf{Update 1.16} 
Toppila \cite{MR0507887} has shown that $e$ cannot be replaced by $1$. He has constructed an example of a meromorphic function for which
\[\liminf_{r\to\infty}\frac{n(r)}{A(r)}\geq\frac{80}{79}\]
for every sufficiently large $r$. The question remains open for entire functions. Among other examples, Toppila shows that for an entire function the following can occur
\[\limsup\frac{n(r,0)}{A(7r/6)}\geq\frac{9}{5}\]
and
\[\limsup\frac{n(r)}{A(Kr)}=\infty\]
for every $K$, $K\geq1$.\\
Miles \cite{doi:10.1080/17476939808815143} gave a positive answer, by showing that for every meromorphic
function
\[\liminf_{r\to\infty}\frac{\max_a n(r,a)}{A(r)}\leq e-10^{-28}.\]\\

\noindent\textbf{Problem 1.17} \textbf{(Paley's conjecture)} For any entire function $f(z)$ of finite order $\rho$ in the plane, we have 
\[1\leq\liminf_{r\to\infty}\frac{\log M(r,f)}{T(r,f)}\leq C(\rho),\]
where $C(\rho)$ depends on $\rho$ only. This follows very simply from Hayman \cite[Theorem 1.6, p.18]{M.F.}. It is known by Wahlund \cite{Wahlund1929} that the best possible value of $C(\rho)$ is $\pi\rho/\sin(\pi\rho)$ for $0<\rho<\frac{1}{2}$, and it is conjectured that $C(\rho)=\pi\rho$ is the corresponding result for $\rho>\frac{1}{2}$.\\

\noindent\textbf{Update 1.17}
This inequality has been proved by Govorov \cite{MR0247079} for entire functions, and by Petrenko \cite{MR0257363} for meromorphic functions.\\

\noindent\textbf{Problem 1.18} Suppose that $f(z)$ is meromorphic in the plane, and that $f(z)$ and $f^{(l)}(z)$ have no zeros, for some $l\geq2$. Prove that $f(z)=e^{az+b}$ or $(Az+B)^{-n}$.\\
The result is known if $f(z)$ has only a finite number of poles, (see Clunie \cite{MR0143906} and Hayman \cite[p.67]{M.F.}), or if $f(z)$ has finite order and $f\neq0, f'\neq0, f''\neq0$, and 
\[\liminf_{r\to\infty}\frac{\log n(r,f)}{\log r}<+\infty,\]
(see Hayman \cite{MR0110807}), or if none of the derivatives of $f(z)$ have any zeros and $f(z)$ has unrestricted growth (see P{\-o}lya \cite{MR1544505}, Hayman \cite[p. 63]{M.F.}).\\

\noindent\textbf{Update 1.18}
The conjecture was proved by Mues \cite{MR0276471} if $f$ has finite order and $ff''\neq0$ (instead of $ff'f''\neq0$). For $l>2$, the conjecture was proved by Frank \cite{MR0422615}. Since then, Frank, Polloczek and Hennekemper \cite{MR0430250} have obtained various extensions. Thus, if $f$, $f^{(l)}$ have only a finite number of zeros, and $l>2$, then 
\[f(z)=\frac{p_1}{p_2}e^{p_3},\]
 where $p_1, p_2, p_3$ are polynomials. However, the paper \cite{MR0430250} contains mistakes, including a gap in the proof of the case $l=2$.\\
The last case which remained unsolved, $l=2$, was settled by Langley \cite{MR1241784} who proved that the only meromorphic functions $f$ for which $ff''$ is zero-free, are $f(z)=e^{az+b}$ and $f(z)=(az+b)^{-n}$.\\

\noindent\textbf{Problem 1.19} Suppose that $f(z)$ is meromorphic in the plane and \mbox{$f'(z)f(z)^n\neq1$}, where $n\geq1$. Prove that $f(z)$ is constant. Hayman \cite{MR0110807} has shown this to be true for $n\geq3$.\\

\noindent\textbf{Update 1.19}
The case $n=2$ has been settled by Mues \cite{MR516609}. The last case which remained unsolved, $n=1$, was settled by Bergweiler and Eremenko \cite{MR1344897}: for every non-constant meromorphic function $f$, the equation $f'(z)f(z)=c$ has solutions for every $c$, $c\neq0,\infty$. This was first proved by Bergweiler and Eremenko for functions of finite order; then, Bergweiler and Eremenko \cite{MR1344897}, Chen and Fang \cite{MR1360682}, and Zalcman \cite{MR3539851} independently noticed that a general result of Pang \cite{MR1057999} permits an extension to arbitrary meromorphic functions. The proof actually applies whenever $n\geq1$.\\
There were many extensions and generalizations of this result
of Bergweiler and Eremenko. The strongest result so far is due to
Jianming Chang \cite{doi:10.1112/blms/bds003}:
{\em Let $f$ be a transcendental
meromorphic function whose derivative is bounded on the set
of zeros of f.
Then the equation $f(z) = c$ has infinitely many solutions for
every $c\in\mathbb{C}\setminus\{0\}$.}\\

\noindent\textbf{Problem 1.20} If $f(z)$ is non-constant and meromorphic in the plane, and $n=3$ or $4$, prove that $\phi(z)=f'(z)-f(z)^n$ assumes all finite complex values. This is known to be true if $f(z)$ is an entire function, or if $n\geq5$ in the case where $f(z)$ is meromorphic, see Hayman \cite{MR0110807}.\\
In connection to this, it would be most interesting to have general conditions under which a polynomial in $f(z)$ and its derivatives, can fail to take some complex value. When $f(z)$ is a meromorphic rather than an entire function, rather little is known, see however Clunie (\cite{MR0143906}, \cite{MR0214769}) and Hayman \cite[Chp. 3]{M.F.}. \\

\noindent\textbf{Update 1.20}
This question is closely related to Problem 1.19. Mues \cite{MR516609} proved that $\phi(z)$ may omit a finite non-zero value
when $n=3$ or $4$. He also showed that $\phi$ must have infinitely
many zeros for $n=4$. The remaining case of zeros for $n=3$
was settled by Bergweiler and  Eremenko \cite{MR1344897}. \\

\noindent\textbf{Problem 1.21} If $f(z)$ is non-constant in the plane, it is known (see Hayman \cite[pp.55-56]{M.F.}) that
\[
\alpha_f=\limsup_{r\to\infty}\frac{T(r,f)}{T(r,f')}\geq
\begin{cases}
\frac{1}{2}  & \text{if } f(z) \text{ is meromorphic},\\
1 & \text{if } f(z) \text{ is an entire function}.\\
\end{cases}
\]
These inequalities are sharp. It is not known whether
\[\beta_f=\liminf_{r\to\infty}\frac{T(r,f)}{T(r,f')}\]
can be greater than one, or even infinite. It is known that $\beta_f$ is finite if $f(z)$ has finite order. Examples show that $\alpha_f$ may be infinite for entire functions of any order $\rho$, i.e. $0\leq\rho\leq\infty$, and that given any positive constants $K$, $\rho$ there exists an entire function of order at most $\rho$, such that
\[\frac{T(r,f)}{T(r,f')}>K\]
on a set of $r$ having positive lower logarithmic density. For this and related results, see Hayman \cite{MR0180679}.\\

\noindent\textbf{Update 1.21}
Let $f$ be meromorphic in the plane. The relation between $T(r,f')$ and $T(r,f)$ constitutes an old problem of Nevanlinna theory. It is classical that
\[m(r,f')\leq m(r,f)+m\Big(r,\frac{f'}{f}\Big)\leq m(r,f)+O(\log T(r,f))\]
outside an exceptional set. In particular, if $f$ is entire so that $m(r,f)=T(r,f)$, we deduce that $T(r,f')<(1+o(1))T(r,f)$ outside an exceptional set.\\
The question of a corresponding result in the opposite direction, had been open until fairly recently. Hayman \cite{MR0168763} has shown that there exist entire functions of finite order $\rho$ for which $T(r,f)>KT(r,f')$ on a set having positive lower logarithmic density, for every positive $\rho$ and $K>1$. Toppila \cite{MR0585685} has given a simple example for which
\[\beta_f=\limsup_{\to\infty}\frac{T(r,f)}{T(r,f')}\geq1+\frac{7}{10^7}.\]
for all sufficiently large $r$. In this example, he takes for $f'$ the square of the sine product, having permuted the zeros in successive annuli, to the positive or negative axis. The result is that $f'$ is sometimes large on  each half-axis, and so $f$, the integral of $f'$, is always large on and near the real axis.\\
On the other hand, Hayman and Miles \cite{MR1040924} have proved that $\beta_f\leq 3e$ if $f$ is meromorphic, and $\beta_f\leq 2e$ if $f$ is entire. Density estimates are also given to show that previous examples are fairly sharp.\\

\noindent\textbf{Problem 1.22} The defect relation (\ref{1.2}) is a consequence of the inequality (see Hayman \cite[formula (2.9), p. 43]{M.F.}), which is called the \textit{``second fundamental theorem''},
\begin{equation}
\label{1.4}
\sum^k_{\nu=1}\overline{N}(r,a_\nu, f)\geq \big(q-2+o(1)\big)T(r,f) 
\end{equation}
which holds for any distinct numbers $a_\nu$ and $q\geq3$, as $r\to\infty$ outside a set $E$ of finite measure, if $f(z)$ is meromorphic in the plane. The exceptional set $E$ is known to be unnecessary if $f(z)$ has finite order. Does (\ref{1.4}) also hold as $r\to\infty$ without restriction if $f(z)$ has infinite order?\\

\noindent\textbf{Update 1.22}
A negative answer to this question is provided by the examples discussed in connection with Problem 1.2. These show that the second fundamental theorem fails to hold on the sequence $r=r_\nu$.\\

\noindent\textbf{Problem 1.23} Under what circumstances does $f(z_0+z)$ have the same deficiencies as $f(z)$? It was shown by Dugu{\'e} \cite{MR0021975} that this need not be the case for meromorphic functions, and by Hayman \cite{MR0054705} that it is not necessarily true for entire functions of infinite order. The case of functions of finite order remains open. Valiron \cite{MR0021976} notes that a sufficient condition is
\[\frac{T(r+1,f)}{T(r,f)}\to 1, \hspace{1cm}\text{ as }r\to\infty,\]
and this is the case in particular if $\rho-\lambda<1$. Since for entire functions of lower order $\lambda$, $\lambda\leq\frac{1}{2}$ there are no deficiencies anyway, it follows that the result is true at any rate, for entire functions of order $\rho<\frac{3}{2}$ and, since $\lambda\geq0$ always, for meromorphic functions of order less than one.\\

\noindent\textbf{Update 1.23}
Gol'dberg and Ostrovskii \cite{MR0280720} give examples of meromorphic functions of finite order for which the deficiency is not invariant under change of origin. See also Gol'dberg and Ostrovskii \cite{MR2435270}; and Wittich (\cite{MR0077620}, \cite{MR0244490}) for details.\\
Miles \cite{MR698931} provided a counter-example of an entire $f$ of large finite order. Gol'dberg, Eremenko and Sodin \cite{MR916442} have constructed such $f$ with preassigned order $\rho$, such that $5<\rho<\infty$.\\
  
\noindent\textbf{Problem 1.24} If $f$ is meromorphic in the plane, can $n(r,a)$ be compared in general with its average value
\[A(r)=\frac{1}{\pi}\int\int_{|z|<r}\frac{|f'(z)|^2}{(1+|f(z)|^2)^2}\,dx\, dy\]
in the same sort of way that $N(r,a)$ can be compared with $T(r)$? In particular, is it true that $n(r,a)\tilde A(r)$ as $r\to\infty$, outside an exceptional set of $r$, independent of $a$, and possibly an exceptional set of $a$? (Compare Problem 1.16.) \\
\textit{(P. Erd\"os)}\\

\noindent\textbf{Update 1.24}
Miles \cite{MR0254245} had shown that 
\[\lim_{r\to\infty, r\notin E}\frac{n(r,a)}{A(r)}=1,\]
for all $a$ not in $A$, a set of inner capacity zero, and all $r$ not in $E$ where $E$ is a set of finite logarithmic measure.\\

\noindent\textbf{Problem 1.25} In the opposite direction to Problem 1.24, does there exist a meromorphic function such that for every pair of distinct values $a, b$, we have
\[\limsup_{r\to\infty}\frac{n(r,a)}{n(r,b)}=\infty\hspace{1cm}\text{ and }\hspace{1cm}\liminf_{r\to\infty}\frac{n(r,a)}{n(r,b)}=0.\]
Note, of course, that either of the above limits for all distinct $a, b$ implies the other.\\
(Compare the result (1.3) quoted in Problem 1.2, which shows that this certainly cannot occur for the $N$-function.)
The above question can also be asked for entire functions.\\
\textit{(P. Erd\"os)}\\

\noindent\textbf{Update 1.25}
The exceptional sets are necessary. Both Gol'dberg \cite{MR0466551} and Toppila \cite{MR0507887} have produced examples of entire functions for which
\[\limsup_{r\to\infty}\frac{n(r,a)}{n(r,b)}=\infty,\]
for every finite unequal pair $(a,b)$. A corresponding example for meromorphic functions has also been given by Toppila \cite{MR0507887}.\\

\noindent\textbf{Problem 1.26} The analogue of Problem 1.7 may be asked for meromorphic functions. The proposers conjecture that in this case
\[\sum\delta(a,f)\leq\max\{\Lambda_1(\rho),\Lambda_2(\rho)\},\]
where for $\rho\geq1, q=[2\rho]$ we have
\[\Lambda_1(\rho)=2-\frac{2\sin\big(\frac{1}{2}\pi(2\rho-q)\big)}{q+2\sin\big(\frac{1}{2}\pi(2\rho-q)\big)},\]
\[\Lambda_2(\rho)=2-\frac{2\cos\big(\frac{1}{2}\pi(2\rho-q)\big)}{q+1}.\]
Weitsman \cite{MR0387598} shows that this result would be sharp. The correct bound is known for $0\leq\rho\leq1$.\\
\textit{(D. Drasin and A. Weitsman)}\\

\noindent\textbf{Update 1.26} 
No progress on this problem has been reported to us.\\

\noindent\textbf{Problem 1.27} Let $E$ be the set for which $m(r,a)\to\infty$ as $r\to\infty$. How large can $E$ be if:
\begin{enumerate}[(a)]
\item $f$ is entire and of order $\frac{1}{2}$ mean type,
\item $f$ is meromorphic of order $\rho$, where $0\leq\rho\leq\frac{1}{2}$.
\end{enumerate}
The proposers settled this problem in all other cases (see Update 2.1 for more details).\\
\textit{(D. Drasin and A. Weitsman)}\\

\noindent\textbf{Update 1.27}
For part (b) above, Damodaran \cite{MR0450558} proved the existence of meromorphic functions of growth $T(r,f)=O(\rho(r)(\log r)^3)$, where $\rho(r)\to\infty$ arbitrarily slowly, such that $m(r,a)\to\infty$ for all $a$ in an arbitrarily prescribed set of capacity zero. Lewis \cite{MR0499158} and Eremenko \cite{MR501765} independently improved this to $(\log r)^2$ in place of $(\log r)^3$. This is best-possible, following from an old result of Tumura \cite{MR0016488}.\\

\noindent\textbf{Problem 1.28}  Are there upper bounds of any kind on the set of asymptotic values of a meromorphic function of finite order?\\
\textit{(D. Drasin and A. Weitsman)}\\

\noindent\textbf{Update 1.28} 
A negative answer to this question has been given by Eremenko \cite{MR522410} who has constructed meromorphic functions of positive and of zero order, having every value in the closed plane as an asymptotic value. This has been improved by Canton, Drasin and Granados \cite{10.2307/24903359} who proved that for every $\phi(r)\to+\infty$ and every
analytic (Suslin) set $A$, there exists a meromorphic function
$f$ with the property $T(r,f)=O(\phi(r)\log^2r)$ and
whose set of asymptotic values coincides with $A$.\\

\noindent\textbf{Problem 1.29} Under what circumstances does there exist a meromorphic function $f(z)$ of finite order $\rho$ with preassigned deficiencies $\delta_n=\delta(a_n,f)$ at a preassigned sequence of complex numbers? Weitsman has solved this problem (see Update 1.14) by showing that it is necessary that
\begin{equation}
\label{1.5}
\sum\delta_n^{\frac{1}{3}}<\infty,
\end{equation} 
but the bound of the sum of the series depends on the largest term $\delta_1$. On the other hand, Hayman \cite[p.98]{M.F.} showed that the condition 
\begin{equation}
\label{1.6}
\sum\delta_n^{\frac{1}{3}}<A,
\end{equation} 
with $A=9^{-\frac{1}{3}}$ is sufficient to yield a meromorphic function of order $1$ mean type, such that $\delta(a_n,f)\geq\delta_n$. Possibly (\ref{1.6}) with a constant $A$, $A=A_1(\rho)$ is sufficient, and with a larger constant, $A$, $A=A_2(\rho)$ is necessary. If $\rho$ is allowed to be arbitrary but finite, the problem may be a little easier.\\

\noindent\textbf{Update 1.29}
Eremenko \cite{MR853888} has constructed an example of a function of finite order, having preassigned deficiencies $\delta_n=\delta(a_n,f)$, subject to $0<\delta_n<1$, $\sum\delta_n<2$ and $\sum\delta_n^{\frac{1}{3}}<\infty$,  and no other conditions. In view of the results reported in Update 1.3 and Update 1.33, this result
is a complete solution of the Inverse Problem in the class of functions of
(unspecified) finite  order.\\

\noindent\textbf{Problem 1.30} Can one establish an upper bound on the number of finite asymptotic values of a meromorphic function $f(z)$ in $\mathbb{C}$, taking into account both the order of $f$, and the angular measure of its tracts?\\
\textit{(W. Al-Katifi)}\\

\noindent\textbf{Update 1.30}
No progress on this problem has been reported to us.\\

\noindent\textbf{Problem 1.31} Let the function $f$ be meromorphic in the plane, and not rational, and satisfy the condition
\begin{equation}
\label{star}
\frac{T(r,f)}{(\log r)^3}\to\infty,\hspace{1cm}\text{ as }r\to\infty,
\end{equation}
where $T(r,f)$ is the Nevanlinna characteristic. A theorem of Yang Lo \cite{MR653387} states that then there exists a direction $\theta_0\in[0,2\pi)$ such that for every positive $\varepsilon$, either $f$ attains every finite value infinitely often in $D_\varepsilon=\{z:|\arg z - \theta_0|<\varepsilon\}$, or else $f^{(k)}$ attains every value, except possibly zero, infinitely often in $D_\varepsilon$ for all positive integers $k$. Can the condition (\ref{star}) be dropped completely? Or, possibly, can it be replaced be the `more usual' condition
\[\frac{T(r,f)}{(\log r)^2}\to\infty,\hspace{1cm}\text{ as }r\to\infty\,?\]
One cannot expect any more from Yang Lo's method of finding $\theta_0$ through the use of `filling discs'. Rossi \cite{MR1304116} has shown that (\ref{star}) cannot be improved if $\theta_0$ is sought in this way.\\
\textit{(D. Drasin; communicated by J. Rossi)}\\

\noindent\textbf{Update 1.31}
Rossi writes that there is an incorrect paper of Zhu \cite{MR1317004} where he purports to use filling discs to solve this problem. However, Fenton and Rossi \cite{MR1651013} remark that Zhu's approach is wrong, and points to the example in Rossi \cite{MR1304116}. Some work on this problem has been produced by Sauer \cite{MR2478263}.\\

\noindent\textbf{Problem 1.32} Let $f$ be meromorphic in $\mathbb{C}$, and let $f^{-1}$ denote any element of the inverse function that is analytic in a neighbourhood of a point $w$. A well-known theorem of Gross \cite{MR1548975}, states that $f^{-1}$ may be continued analytically along almost all rays beginning at $w$. Is it possible to refine the exceptional set in this theorem?\\
\textit{(A. Eremenko)}\\

\noindent\textbf{Update 1.32}
No progress on this problem has been reported to us. \\

\noindent\textbf{Problem 1.33} Let $f$ be a meromorphic function of finite order $\rho$. Does the condition 
\[N(r,1/f')+2N(r,f)-N(r,f')=o(T(r,f)),\hspace{1cm}\text{ as }r\to\infty,\]
imply that $2\rho$ is an integer?\\
\textit{(A. Eremenko)}\\

\noindent\textbf{Update 1.33}
This is a slightly more precise conjecture than Problem 2.25. Both problems are solved completely by the following theorem of Eremenko \cite{MR1266090}:\\
suppose that $f$ is a meromorphic function of finite lower order $\lambda$, and that
\[N_1(r,f):=N(r,1/f')+2N(r,f)-N(r,f')=o(T(r,f)).\]
Then
\begin{enumerate}[(a)]
\item $2\lambda$ is an integer greater than or equal to $2$.
\item $T(r,f)=r^\lambda l(r)$, where $l$ is a slowly varying function in the sense of Karamata,
\item $\sum\delta(a,f)=2$, all deficient values are asymptotoic, and all deficiencies are multiples of $1/\lambda$.\\
\end{enumerate}

\noindent\textbf{Problem 1.34} Let $n_1(r,a,f)$ denote the number of simple zeros of $f(z)-a$ in $\{|z|\leq r\}$. Selberg \cite{Selberg1928} has shown that if: 
\begin{enumerate}[(a)]
\item $f$ is a meromorphic function of finite order $\rho$, and
\item $n_1(r,a,f)=O(1)$, as $r\to\infty$ for four distinct values of $a$, then 
\end{enumerate}
$\rho$ is an integral multiple of $\frac{1}{2}$ or $\frac{1}{3}$.\\ 
Does this conclusion remain true if $(b)$ is replaced by:
\begin{enumerate}[(c)]
\item  $n_1(r,a,f)=o(T(r,f))$, as $r\to\infty$, for four distinct values of $a$?
\end{enumerate}
Gol'dberg \cite{MR0335805} has constucted an entire function of arbitrary prescribed order which satisfies the condition $(c)$ for \text{two} distinct values of $a$.\\
\textit{(A. Eremenko)}\\

\noindent\textbf{Update 1.34}
The answer is `no'. K\"unzi \cite{MR0069893} has shown that $\rho$ can be arbitrary, subject to $1<\rho<\infty$, and Gol'dberg \cite{MR0335805} has a counter-example for arbitrary positive $\rho$.\\

\noindent\textbf{Problem 1.35} Determine the upper and lower estimates for the growth of entire and meromorphic solutions of algebraic ordinary differential equations (AODE). (This is a classical problem.)\\
For AODEs of first order, it is known that the meromorphic solutions $f$ must have finite order (see Gol'dberg \cite{MR0085396}) and that $(\log r)^2=O(T(r,f))$ (see Eremenko \cite{MR667974}, \cite{Eremenko1984}). (The latter two references contain a general account of first order AODEs, including modern proofs of some classical results.) The order of entire solutions of first order AODEs is an integral multiple of $\frac{1}{2}$ (see Malmquist \cite{MR1555170}). For AODEs of second order, it is known that the order of entire solutions is positive, see Zimogljad \cite{MR0291462}. There is no upper estimate valid for all entire or meromorphic solutions of AODEs of order grreater than one, but there is an old conjecture that $\log |f(z)|\leq exp_n(|z|)$ for entire solutions $f$ of an AODE of order $n$.\\
\textit{(A. Eremenko)}\\

\noindent\textbf{Update 1.35}
Steinmetz (\cite{MR595424}, \cite{MR686638}) proved that every meromorphic solution of a homogeneous algebraic differential equation of second order has the form $f=(g_1/g_2)exp(g_3)$, where the $g_i$ are entire functions of finite order. Thus $T(r,f)=O(exp(r^k))$ for some positive $k$. Due to Wiman and Valiron (see, for example \cite{MR0061658}), it is known that `most' algebraic differntial equations do not have entire solutions of infinite order. A precise statement of this sort is contained in Hayman \cite{MR1438606}.\\

\noindent\textbf{Problem 1.36} Let $F$ be a polynomial in two variables, and let $y$ be a meromorphic solution of the algebraic ordinary differential equation $F(y^{(n)},y)=0$. Is it true that $y$ must be an elliptic function, or a rational function of exponentials, or a rational function? This is known in the following cases:
\begin{enumerate}[(a)]
\item $n=1$: a classical result, probably due to Abel;
\item $n=2$: an old result of Picard \cite{Picard1880}, and independently, Bank and Kaufman \cite{MR624232};
\item $n$ is odd and $y$ has at least one pole, Eremenko \cite{MR686077}; and
\item the genus of the curve $F(x_1,x_2)=0$ is at least equal to one, Eremenko \cite{MR686077}.
\end{enumerate}
\textit{(A. Eremenko)}\\

\noindent\textbf{Update 1.36} This has been solved by Eremenko, Liao and Tuen-Wai Ng \cite{EremLiaoNg} who prove (c) for all $n$, and give an example of an entire solution which is neither rational, nor a rational function of an exponential.\\

\noindent\textbf{Problem 1.37} Find criteria for and/or give explicit methods for the construction of meromorphic functions $f$ in $\mathbb{C}$ with the following properties:
\begin{enumerate}[(a)]
\item all poles of $f$ are of odd multiplicity;
\item all zeros of $f$ are of even multiplicity.
\end{enumerate}
(Here `explicit methods' means that all computations must be practicable.) The background of this problem lies in the question of meromorphic solutions of the differential equation $y''+A(z)y=0$ in the whole plane.\\
\textit{(J. Winkler)}\\

\noindent\textbf{Update 1.37}
We mention a result of Schmieder \cite{MR630449},which may be relevant: on every open Riemann surface, there exists an analytic function with prescribed divisors of zeros and critical points, subject to the trivial restrictions.\\

\noindent\textbf{Problem 1.38}
\begin{enumerate}[(a)]
\item Let $f$ be non-constant and meromorphic in the open unit disc $\mathbb{D}$, with $\alpha<+\infty$, and define 
\begin{equation}
 \label{alphadef}
\alpha = \limsup_{r\to1}\frac{T(r, f)}{-\log (1-r)},
\end{equation}
and
\[\Psi=(f)^{m_0}(f')^{m_1}\ldots(f^{(k)})^{m_k}.\]
It is known that $\Psi$ assumes all finite values, except possibly zero, infinitely often, provided that $m_0\geq3$ and $\alpha>2/(m_0-2)$, (or $m_0\geq2$ and \mbox{$\alpha>2/(m_0-1)$}, if $f$ is analytic). For which smaller values of $\alpha$ does the same conclusion hold?
\item Let $f$ be non-constant and meromorphic in $\mathbb{D}$, with $\alpha<+\infty$
in (\ref{alphadef}); assume also that $f$ has only finitely many zeros and poles in $\mathbb{D}$. Let $l$ be a positive integer, and write $\Psi=\sum^l_{\nu=0}a_\nu f^{(\nu)}$, where the $a_\nu$ are functions in $\mathbb{D}$ for which $T(r,a_\nu)=O(T(r,f))$ as $r\to1$ (for each $\nu$). It is known that if $\Psi$ is non-constant, then $\Psi$ assumes every finite value, except possibly zero, infinitely often, provided that $\alpha>\frac{1}{2}l(l+1)+1$. For which smaller values of $\alpha$ does the same conclusion hold?
\end{enumerate}
\textit{(L. R. Sons)}\\

\noindent\textbf{Update 1.38}
For part (b), Gunsul \cite{zbMATH05900730} has put forth a condition that enables smaller values of $\alpha$ for which the same conclusion holds.\\

\noindent\textbf{Problem 1.39} Let $f$ be a function meromorphic in $\mathbb{D}$, for which $\alpha<+\infty$ in (\ref{alphadef}). 
\begin{enumerate}[(a)]
\item Shea and Sons \cite[Theorem 5]{MR862041} have shown that if $f(z)\neq0,\infty$ and $f'(z)\neq1$ in $\mathbb{D}$, then $\alpha\leq 2$. Is $2$ best possible?
\item Shea and Sons \cite{MR862041} have shown that, if $f(z)\neq0$ and $f'(z)\neq1$ in $\mathbb{D}$, then $\alpha\leq7$. What is the best possible $\alpha$ in this case?
\end{enumerate}
\textit{(L. R. Sons)}\\

\noindent\textbf{Update 1.39}
No progress on this problem has been reported to us.\\

\noindent\textbf{Problem 1.40} Let $f$ be a function meromorphic in $\mathbb{D}$ of finite order $\rho$. Shea and Sons \cite{MR862041} have shown that
\[\sum_{a\neq\infty}\delta(a,f)\leq\delta(0,f')\big(1+k(f)\big)+\frac{2}{\lambda}(\rho+1),\]
where
\[k(f)=\limsup_{r\to1}\frac{\overline{N}(r,\infty,f)}{T(r,f)+1}\hspace{1cm}\text{ and }\hspace{1cm}\lambda(f)=\liminf_{r\to\infty}\frac{T(r,f)}{\log(1/(1-r))}.\]
Can the factor $2$ be eliminated? (If so, the result is then best possible).\\
\textit{(L. R. Sons)}\\

\noindent\textbf{Update 1.40}
No progress on this problem has been reported to us.\\

\noindent\textbf{Problem 1.41} Let $f$ be a function meromorphic in $\mathbb{D}$, for which $\alpha=+\infty$ in (\ref{alphadef}).  Then it is known that
\[\sum_{a\in\mathbb{C}\cup\{\infty\}}\delta(a,f)\leq2.\]
Are there functions which have an `arbitrary' assignment of deficiencies at an arbitrary sequence of complex numbers, subject only to these conditions?\\
For analytic functions, Girynk \cite{MR0437767} has a result; whereas for arbitrary meromorphic functions, there is a result of Krutin \cite{MR0335807}.\\
\textit{(L. R. Sons)}\\

\noindent\textbf{Update 1.41}
No progress on this problem has been reported to us.\\

\noindent\textbf{Problem 1.42} Let $f$ be meromorphic in $\mathbb{C}$, and suppose that the function
\[F(z)=f^{(k)}(z)+\sum^{k-2}_{j=0}a_j(z)f^{(j)}(z)\]
is non-constant, where $k\geq3$ and the coefficients $a_j$ are polynomials. Characterise those functions $f$ for which $f$ and $F$ have no zeros.\\
The case where $f$ is entire has been settled by Frank and Hellerstein \cite{MR868452}. If all the $a_j$ are constant, then the problem has also been solved by Steinmetz \cite{MR928649} using results from Frank and Hellerstein \cite{MR868452}. It seems possible that if the $a_j$ are not all constants, then the only solutions with infinitely many poles are of the form $f=(H')^{\frac{1}{2}(k-1)}H^{-l}$, where $l$ is a positive integer, and $H''/H'$ is a polynomial.\\
\textit{(G. Frank and J. K. Langley)}\\

\noindent\textbf{Update 1.42}
This was  solved by Br\"uggemann \cite{MR1159367}, who proved the following: let a linear differential operator
\[L(f)=f^{(n)}+\sum^{n-2}_{j=0}a_j f^j\]
with polynomial coefficients $a_j$ be given, with at least one non-constant $a_j$. Then the only meromorphic functions $f$ with infinitely many poles, satisfying $fL(f)\neq0$, are of the form $f=(H')^{-(n-1)/2}H^{-l}$, where $l$ is a positive integer, and $H''/H'$ is a polynomial.\\

An extension to rational coefficients has been given by Langley:\\
Linear differential polynomials in zero-free meromorphic functions, Ann. Acad. Sci. Fenn., to appear. 

\noindent\textbf{Problem 1.43} Let $f$ be a meromorphic function of lower order $\lambda$. Let
\[m_0(r,f)=\inf\{|f(z)|:|z|=r\}\]
and
\[M(r,f)=\sup\{|f(z)|:|z|=r\}\]
and suppose that
\[\log r=o(\log M(r,f)),\hspace{1cm}\text{ as }r\to\infty.\]
Gol'dberg and Ostrovskii \cite{MR0280720} proved that if $0<\lambda<\frac{1}{2}$, then 
\[\limsup_{r\to\infty}\frac{\log m_0(r,f)}{\log M(r,f)}+\pi\lambda\sin(\pi\lambda)\limsup_{r\to\infty}\frac{N(r,f)}{\log M(r,f)}\geq\cos(\pi\lambda).\]
Does this inequality remain valid for $\frac{1}{2}\leq\lambda<1$?  See also Gol'dberg and Ostrovskii \cite{MR2435270}.\\
\textit{(A. A. Gol'dberg and I. V. Ostrovskii)}\\

\noindent\textbf{Update 1.43}
No progress on this problem has been reported to us.

\section*{New Problems}

\noindent To appear.

% ==================== CHAPTER 2 ==========================================

\chapter{Entire Functions}

\underline{Preface by P. Rippon (to appear)}\\

\noindent\textbf{Notation} Let $f(z)$ be an entire function. We say that $a$ is an \textit{asymptotic value} of $f(z)$, if 
\[f(z)\to a,\] 
as $z\to\infty$ along a path $\Gamma$, called a corresponding \textit{asymptotic path}. Some of the most interesting open problems concerning entire functions centre on these asymptotic values and paths. It follows from a famous result of Ahlfors \cite{Ahlfors1930}, that an entire function of finite order $\rho$ can have at most $2\rho$ distinct finite asymptotic values. On the other hand, by a theorem of Iversen \cite{Iversen}, $\infty$ is an asymptotic value of every entire function. Some of the following problems are concerned with generalisations arising out of the above two theorems. \\
Throughout this section
\[M(r)=M(r,f)=\max_{|z|=r}|f(z)|\]
denotes the \textit{maximum modulus} of $f(z)$.\\

\noindent\textbf{Problem 2.1} Suppose that $f(z)$ is an entire function of finite order. What can we say about the set $E$ of values $w$ such that
\begin{enumerate}[(a)]
\item\[ \lambda(r,f-w)=\min_{|z|=r}|f(z)-w|\to0,\hspace{1cm}\text{ as }r\to\infty;\]
\end{enumerate}
or
\begin{enumerate}[(b)]
\item\[m\Big(r,\frac{1}{f-w}\Big)=\frac{1}{2\pi}\int^{2\pi}_0\log^+\Big|\frac{1}{f(re^{i\theta})-w}\Big|\,d\theta, \hspace{1cm}\text{ as }r\to\infty\,?\]
\end{enumerate}
Clearly $(b)$ implies $(a)$. By the result of Arakelyan \cite{MR0206286} (see Problem 1.6 also), the set of deficient values, which is clearly contained in $E$, can include any countable set. Can $E$ be non-countably infinite in case $(b)$, or contain interior points in case $(a)$?\\

\noindent\textbf{Update 2.1} Part$(a)$ was settled by the examples of Arakelyan \cite{MR0206286}. Let $E$ be a dense countable set in the place, every value of which is deficient. Then clearly $\lambda(r,f-w)\to0$ as $r\to\infty$, for every $w$ in the plane. \\
Part $(b)$ was settled by Drasin and Weitsman \cite{zbMATH03353695}. The set of $w$ for which $m(r,\frac{1}{f-w})\to\infty$ as $r\to\infty$ must have capacity zero, and an arbitrary set of capacity zero may occur.\\

\noindent\textbf{Problem 2.2} Produce a general method for constructing an entire function of finite order, and in fact, minimal growth, which tends to different asymptotic values $w_1, w_2, \ldots, w_k$ as $z\to\infty$, along preassigned asymptotic paths $C_1, C_2, \ldots, C_k$. (Known methods by Kennedy \cite{MR0083033} and Al-Katifi \cite{MR0199388} only seem to work if the $w_\nu$ are all equal, unless the $C_\nu$ are straight lines.)\\

\noindent\textbf{Update 2.2}
Such a construction has been given by Hayman \cite{MR0244487}.\\

\noindent\textbf{Problem 2.3} If $\phi(z)$ is an entire function growing slowly compared with the function $f(z)$, we can consider $\phi(z)$ to be an asymptotic function of $f(z)$, if $f(z)-\phi(z)\to0$ as $z\to\infty$ along a path $\Gamma$. Is it true that an entire function of order $\rho$ can have almost $2\rho$ distinct asymptotic functions of order less than $\frac{1}{2}$?\\ 
(If $\phi=\phi_1(z)-\phi_2(z)$ and the minimum modulus of $\phi$ tends to zero, then $\phi$ has lower order at least $\frac{1}{2}$ mean type (See Hayman \cite[p. 288]{zbMATH00194706}).\\
A positive result in this direction is due to Denjoy \cite{MR0075294}, but only when the paths are straight lines. The result when the $\phi_\lambda(z)$ are polynomials is true (and is a trivial consequence of Ahlfors' theorem for asymptotic values \cite{Ahlfors1930}).\\

\noindent\textbf{Update 2.3}
An answer has been given by Somorjai \cite{MR575388} with $\frac{1}{30}$ instead of $\frac{1}{2}$, and Fenton \cite{MR715310} has obtained the same conclusion, if the orders are less than $\frac{1}{4}$. \\

\noindent\textbf{Problem 2.4} Suppose that $f(z)$ is a meromorphic function in the plane, and that for some $\theta$, $0\leq\theta<2\pi$, $f(z)$ assumes every value infinitely often, with at most two exceptions, in every angle $\theta-\varepsilon<\arg z<\theta + \varepsilon$, when $\varepsilon>0$. Then the ray $\arg z=\theta$ is called a \textit{Julia line}. Lehto \cite{MR0107003} has shown that if $f(z)$ is an entire function, or if $f(z)$ is meromorphic and 
\[\limsup_{r\to\infty}\frac{T(r,f)}{(\log r)^2}=+\infty,\]
(but not necessarily otherwise), at least one direction of Julia exists. What can we say about the exceptional values at different Julia lines? In particular, can an entire function $f(z)$ have one exceptional (finite) value $a$ at one Julia line $\Gamma_a$, and a different exceptional value $b$ at a different Julia line $\Gamma_b$?\\
\textit{(C. R\'enyi)}\\

\noindent\textbf{Update 2.4}
An example has been given by Toppila \cite{MR0283202} of a function having different exceptional values at each of $n$ Julia lines. See also Gol'dberg \cite{MR0229827}.\\

\noindent\textbf{Problem 2.5} What can we say about the set $E$ of values $a$ which an entire function $f(z)$ assumes infinitely often in every angle? Simple examples show that $E$ may be the whole open plane, e.g. if 
\[f(z)=\sigma(z)=z\prod_{(m,n)\neq (0,0)}\left(1-\frac{z}{z_{m,n}}\right)
\exp\left\{\frac{z}{z_{m,n}}+\frac{1}{2}\left(\frac{z}{z_{m,n}}\right)^2\right\}.\]
where $z_{m,n}=m+ni$, or the whole plane except one point, if e.g. $f(z)=e^{\sigma(z)}$.
If $z_m=2^me^{im}$, and 
\[f(z)=e^z \prod^\infty_{m=1}\Big(1-\frac{z}{z_m}\Big),\]
then $E$ consists of the value $0$ only, since clearly $f(z)\to0$ as $z\to\infty$, uniformly for $\pi/2+\varepsilon<\arg z<3\pi/2-\varepsilon$, if $\varepsilon>0$. Can $E$ consist of exactly two values?\\
\textit{(C. R\'enyi)}\\

\noindent\textbf{Update 2.5} 
Gol'dberg \cite{goldbergunknown} has answered several of the questions posed in this problem. In particular, he showed that, given any countable set $A$, there exists $E$ such that $A\subset E\subset \overline{A}$, where $\overline{A}$ is the closure of $A$. \\

\noindent\textbf{Problem 2.6} Let $f(z)$ be an entire function. Then Boas (unpublished) proved that there exists a path $\Gamma_\infty$ such that, for every $n$, 
\begin{equation}
\label{2.1}
\left|\frac{f(z)}{z^n}\right|\to\infty, \hspace{1cm}\text{ as }z\to\infty \text{ along } \Gamma.
\end{equation}
Can we improve this result if something is known about the lower growth of $M(r,f)$?
Hayman \cite{MR0125223} has shown that there exist functions of infinite order and, in fact, growing arbitrarily rapidly, such that, on every path $\Gamma$ on which $f(z)\to\infty$, we have
\[\log \log |f(z)|=O(\log |z|),\]
i.e. $f(z)$ has finite order on $\Gamma$.\\

\noindent\textbf{Update 2.6}
Talpur \cite{MR0393514} has shown that if $f$ has order $\rho$ and $\alpha<\rho<\frac{1}{2}$, $\varepsilon>0$, then we can find a path $\Gamma$ going to $\infty$ on which
\[\log|f(z)|>\log M(|z|^{(1-\rho/\alpha)/(1+\varepsilon)})\cos(\pi\alpha),\]
where $M(r,f)$ is the maximum modulus.
Eremenko \cite{MR592215} has proved the following result: let $f$ be an entire function of order $\rho$ and lower order $\lambda$. Then there exists an asymptotic path $\Gamma$, such that
\[\log |f(z)|>(A(\rho, \lambda)+o(1))\log|z|,\hspace{1cm}\text{ as }z\to\infty, z\in \Gamma,\]
where $A(\rho, \lambda)$ is some explicitly written function, with the property $A(\rho, \lambda)>0$, for $0<\lambda\leq\rho<\infty$. When $\lambda<1/2$, we have $A(\rho,\lambda)\geq\lambda$. See also Update 2.8.\\

\noindent\textbf{Problem 2.7} If $f(z)$  of finite order, can anything be asserted about the length of $\Gamma_\infty$, which is the path on which $f(z)$ tends to $\infty$, or the part of it in $|z|\leq r$? \\

\noindent\textbf{Update 2.7} 
If $f$ is entire, then a classical theorem of Iversen \cite{Iversen} asserts that \mbox{$f(z)\to\infty$} as $z\to\infty$ along some path $\Gamma_\infty$. It was asked how short $\Gamma_\infty$ can be. Let $\ell(r)$ be the length of the arc of $\Gamma_\infty$ to the first intersection with $|z|=r$. If 
\[T(r,f)=O((\log r)^2)\hspace{1cm}\text{ as }r\to\infty,\]
then Hayman \cite{MR0111839} showed that $\Gamma_\infty$ can be taken to be a straight line.\\ 
Eremenko and Gol'dberg \cite{GoldbergErem} have constructed examples for which   
$T(r,f)/(\log r)^2$  tends to $\infty$ arbitrarily slowly but $\ell(r)=O(r)$ fails to hold. An independent proof has been given by Toppila \cite{MR595173}. \\
On the other hand, Chang Kuan Heo \cite{MR0585642} has proved that if $f$ has finite order $\rho$, then for any positive $\varepsilon$,
\[\ell(r)=O(r^{(1+\frac{1}{2}\rho+\varepsilon)})\]
can always hold.\\
It is also possible that $f(z)\to a$ as $z\to\infty$ along a path $\Gamma_a$. In this case. Gol'dberg and Eremenko \cite{GoldbergErem} have constructed examples with $f$ having order arbitrarily close to $\frac{1}{2}$, while $\ell(r)\neq O(r)$. 
See also Update 2.10, and Lewis, Rossi and Weitsman \cite{MR735882}.\\

\noindent\textbf{Problem 2.8} Does (\ref{2.1}) remain true if the number $n(r)$ of poles of $f(z)$ in $|z|<r$ satisfies $n(r)=O(r^k)$, where $k<\frac{1}{2}<\lambda$, and $\lambda$ is the lower order of $f(z)$? Gol'dberg and Ostrovskii \cite{MR0280720} have shown that (\ref{2.1}) can be false if $\frac{1}{2}<k<\lambda$.\\

\noindent\textbf{Update 2.8}
The original version of this problem had order $\rho$ instead of lower order $\lambda$. However, Gol'dberg pointed out to Hayman orally that this would give a negative answer.\\
One may ask for corresponding results if $f$ is meromorphic with sufficiently few poles. If $\infty$ is Nevanlinna deficient, and $T(r,f)=O(\log r)^2$, then Anderson and Clunie \cite{MR0199397} showed that $f\to\infty$ along almost all straight lines. The result fails for functions of larger growth, according to an example of Hayman \cite{MR0492265}, even if the deficiency is one. He also proves, that if
\[\limsup_{r\to\infty}\frac{r^{\frac{1}{2}}}{T(r,f)}\int^\infty_r\frac{N(t,f)\,dt}{t^\frac{3}{2}}<2,\]
then $\infty$ is an asymptotic value of $f$. This is true, for instance, if the order of the poles of $f$ is smaller than $\frac{1}{2}$, and smaller than the lower order of $f$.\\
See also Gol'dberg and Ostrovskii (\cite{MR0280720}, \cite{MR2435270}).\\

\noindent\textbf{Problem 2.9} We ask the analogues of Problems 2.6, 2.7 and 2.8 if, in addition, $f(z)$ has another finite Picard value, e.g. $f(z)\neq0$. In this case, if $\infty$ has deficiency one, in the sense of Nevanlinna, (\ref{2.1}) remains true for functions of finite order (see Edrei and Fuchs \cite{EdreiAlbert;FuchsWolfgangH.J.1959Vdev}), but not necessarily for functions of infinite order, see Gol'dberg and Ostrovskii (\cite{MR0280720}, \cite{MR2435270}).\\

\noindent\textbf{Update 2.9} 
If $f(z)\neq0$, then for every $K$, the level set $|f(z)|=K$ contains a curve tending to infinity. Under this condition, Rossi and Weitsman \cite{MR787658} proved that there is an asymptotic curve $\Gamma$ with the following properties:
\begin{equation}
\label{Erem1}
\log |f(z)|>|z|^{1/2-\varepsilon(z)},\hspace{1cm}\text{where }\varepsilon(z)\to0,\hspace{1cm}\text{and}
\end{equation}
\[\int_\Gamma(\log|f|)^{-(2+\alpha)}\,|dz|<\infty,\hspace{1cm}\text{for all positive }\alpha.\]
On the other hand, Barth, Brannan and Hayman \cite{MR507611} constructed a zero-free entire function, having no asymptotic curve, satisfying (\ref{Erem1}), with $\varepsilon=0$. Furthermore, Brannan pointed out that for their example, every asymptotic curve $\Gamma$ satisfies
\[\int_\Gamma(\log|f(z)|)^{-2}|dz|=\infty.\]\\

\noindent\textbf{Problem 2.10} Huber \cite{MR0094452} proved that, for every positive $\mu$, there exists a path $C_\mu$ tending to infinity, such that
\begin{equation}
\label{2.2}
\int_{C_\mu}|f(z)|^{-\mu}|dz|<\infty,
\end{equation}
provided that $f(z)$ does not reduce to a polynomial.
Does there exist a $C_\infty$ such that (\ref{2.2}) holds for every positive $\mu$ with $C_\mu=C_\infty$?\\

\noindent\textbf{Update 2.10}
This question has been settled by Lewis, Rossi and Weitsman \cite{MR735882}, who have proved that there is a path $C_\infty$, suct that (\ref{2.2}) holds for every positive $\mu$, with $C_\infty$ instead of $C_\mu$, thus answering the question affirmatively. Further, they prove that $C_\mu$ is asymptotic, i.e. $f(z)\to\infty$ as $z\to\infty$ along $C_\mu$; and they also obtain estimates for the length $\ell(r)$ in Problem 2.7. The case of finite order was dealt with by Chang Kuan Heo \cite{MR0585642}.\\

%\section* *{Power series with gaps}

\noindent\textbf{Problem 2.11} If $f(z)=\sum a_nz^{\lambda_n}$ is an entire function, and $\sum(1/\lambda_n)$ converges, is it true that:
\begin{enumerate}[(a)]
\item $f(z)$ has no finite asymptotic value,
\item $\limsup_{r\to\infty}\frac{\log m_0(r,f)}{\log M(r,f)}=1$,
\end{enumerate}
where $m_0(r,f)=\inf_{|z|=r}|f(z)|$ is the minimum modulus of $f(z)$?\\
$(a)$ is known for $\lambda_n>n(\log n)^{1+\varepsilon}$ (see K\"ovari \cite{MR0188438}); and $(b)$ is known for $\lambda_n>n(\log n)^2$ (see K\"ovari \cite{MR0176076}). It is also known that $f(z)$ has no finite radial asymptotic value if $\sum(1/\lambda_n)$ converges, and that here this hypothesis cannot be replaced by any weaker condition (see Macintyre \cite{MR0050000}).\\

\noindent\textbf{Update 2.11}
Nazarov \cite{Nazarov} proved that each of the following conditions 
\begin{enumerate}[(a)]
\item $\lambda_{k+1}+\lambda_{k-1}\geq2\lambda_k$ and $\sum(1/\lambda_k)<\infty$,
\item $\sum(\log \log k)/\lambda_k<\infty$
\end{enumerate}
imply
\[{\lim \sup}_{z\to\infty,z\in\Gamma}\frac{\log|f(z)|}{\log M(|z|,f)}=1\]
for every curve $\Gamma$ tending to infinity.\\

\noindent\textbf{Problem 2.12} If the entire function $f(z)$ has finite order $\rho$, and the maximal density of non-zero coefficients is $\Delta$, is it true that if $\rho\Delta<\frac{1}{2}$, $f(z)$ cannot have a finite deficient value with deficiency one?\\

\noindent\textbf{Update 2.12}
No progress on this problem has been reported to us.\\

\noindent\textbf{Problem 2.12a} Under the same conditions as in Problem 2.12, is it true that if $\rho\Delta<1$, $f(z)$ cannot have a finite asymptotic value?
This is known if $\rho\Delta<1/\pi^2$, see K\"ovari \cite{MR0185122}. Is it true that, if $m_0(r,f)$ is the minimum modulus
\[\limsup_{r\to\infty}\frac{\log m_0(r,f)}{\log M(r,f)}\geq\cos(\pi\rho\Delta)?\]
See  K\"ovari \cite{MR0193235}.\\

\noindent\textbf{Update 2.12a}
Fryntov \cite{MR1687884} proved the following partial result: suppose that $f$ is an entire function of lower order $\lambda$, with density of non-zero exponents $\Delta$. If $\rho\Delta<1/3$, and $\Gamma$ is a curve which intersects each circle $|z|=r$ at most once, then
\[{\limsup}_{z\to\infty,z\in\Gamma}\frac{\log|f(z)|}{\log M(|z|,f)}\geq2\cos(\pi\rho\Delta)-1.\]\\

\noindent\textbf{Problem 2.13} If $f(z)=\sum a_n z^{\lambda_n}$ is an entire function, and $\lambda_n/n\to\infty$, is it true that $f(z)$ has
\begin{enumerate}[(a)]
\item no Picard value,
\item no Borel exceptional value,
\item no deficient value?
\end{enumerate}
All this is known for functions of finite order (see Fuchs \cite{MR0159933}). If the answer is `no', are $(b)$ and $(c)$ true for $\sum(1/\lambda_n)<\infty$? Certainly by a theorem of Biernacki \cite{biernacki}, $(a)$ is true in this case.\\
\textit{(T.K\"ovari)}\\

\noindent\textbf{Update 2.13}
Murai \cite{MR723947} has shown that $\sum(1/\lambda_n)<\infty$ does indeed imply $(b)$ and $(c)$, but that $\lambda_n/n\to\infty$ as $n\to\infty$ does not. He constructed an example to show that $\lambda_n/n\to\infty$ as $n\to\infty$, is consistent with $\delta(0,f)=1$. \\
Only the question whether $\sum(1/\lambda_n)<\infty$ implies $(a)$ remains open.\\

%\section* *{JKL: Further problems on entire functions}

\noindent\textbf{Problem 2.14} 
\begin{enumerate} [(a)]
\item  Let $f(z)=\sum a_n z^n$ be entire and $m(r)=\max_n |a_n|r^n$. If $C>\frac{1}{2}$ then does there exist an entire $f$ with 
\[m(r)/M(r,f)\to C ?\]
Any value of $C$ such that $0<C\leq\frac{1}{2}$ is possible. 

\item If $f(z)\neq0$, then
\[\liminf_{r\to\infty}\frac{m(r)}{M(r,f)}=0;\]
Is 
\[\lim_{r\to\infty}\frac{m(r)}{M(r,f)}=0\,?\]
\item What is the exact upper bound $\beta$ of
\[\beta_f=\liminf_{r\to\infty}\frac{m(r)}{M(r,f)}\,?\]
\end{enumerate}
It is known that $\frac{4}{7}<\beta<2/\pi$. See Clunie and Hayman \cite{MR0170005}.\\

\noindent\textbf{Update 2.14}
Davies \cite{MR509941} has shown that the upper limit in (b) can be positive.\\

\noindent\textbf{Problem 2.15} \textbf{(Blumenthal's conjecture)} Let $w=f_1(z), f_2(z)$ be entire functions. Is it true that if
\[M(r,f_1)=M(r,f_2),\hspace{1cm}0<r<\infty,\]
then $f_1(z), f_2(z)$ are equivalent, apart from rotations and reflections in the $z$ and $w$ planes? The corresponding problem for polynomials (of degree higher than about $6$) is also open.\\
The functions $(1-z)e^z$ and $e^{\frac{1}{2}z^2}$ have the same value of $M(r,f)$ for $0<r<2$, and $e^{z-z^2}$ and $e^{z^2+\frac{1}{2}}$ for $r\geq\frac{1}{4}$.\\

\noindent\textbf{Update 2.15}
Hayman, Tyler and White \cite{MR3155684} establish Blumenthal's conjecture for polynomials $f$ and $g$ with at most four non-zero terms (and so, in particular, for all quadratic and cubic polynomials). Short examples are also given to show that, in the general case, it is not sufficient only to consider arbitrarily large or arbitrarily small positive values of $r$. However, the following result is proved: if entire functions $f$ and $g$ are real on the real axis, $f(0)f'(0)\neq0$ and $M(r,f)=M(r,g)$ in some range $0<r<r_0$, then $f$ and $g$ are equivalent.\\
The corresponding problem for polynomials of degree higher than $3$ remains open.\\

\noindent\textbf{Problem 2.16} Let $\nu(r)$ be the number of points on $|z|=r$, such that \mbox{$|f(z)|=M(r,f)$}. Can we have
\begin{enumerate}[(a)]
\item $\limsup_{r\to\infty}\nu(r)=\infty$\,?
\item $\liminf_{r\to\infty}\nu(r)=\infty$\,?
\end{enumerate}
\textit{(P. Erd\"os)}\\

\noindent\textbf{Update 2.16}
Herzog and Piranian \cite{MR0257359} have shown that $(a)$ is indeed possible, however the answer to $(b)$ is still unknown. These authors also provided an example of a univalent function in $\mathbb{D}$ for which the analogue of $(a)$ holds. \\

\noindent\textbf{Problem 2.17} If $f(z)$ is a non-constant entire function and 
\[b(r)=\left(r\frac{d}{dr}\right)^2\log M(r,f),\]
then
\begin{equation}
\label{2.3}
\limsup_{r\to\infty} b(r)\geq A
\end{equation}
where $A$ is an absolute constant, such that $0.18<A\leq\frac{1}{4}$, see Hayman \cite{MR0252639}. What is the best value of $A$?\\
It seems fair to conjecture that the correct constant in (\ref{2.3}) is in fact $\frac{1}{4}$.\\

\noindent\textbf{Update 2.17}
Kjellberg \cite{Kjellberg} proved that this conjecture is false, and that \mbox{$0.24<A<0.248$}. In this problem, and also in Problem 2.18, the monomials $az^n$ should be excluded. Bo\u\i\v cuk and Gol'dberg \cite{MR0344465} have shown that the result is true if $f$ has positive coefficients; and in fact
\[\limsup_{r\to\infty}b(r)\geq\frac{1}{4}A^2,\hspace{1cm} \text{where }A=\limsup_{k\to\infty}(n_{k+1}-n_k).\]
Tyler (unpublished) has some numerical evidence that $A$ is just under $0.247$.\\

\noindent\textbf{Problem 2.18} Consider the function $b(r)$ of Problem 2.17. Since $\log M(r,f)$ is an analytic function of $r$, except for isolated points, $b(r)$ exists except at isolated points where the right and left limits $b(r\mp0)$ still exist, but may be different. It follows from Hadamard's convexity theorem (see Hayman \cite{MR0252639}) that $b(r)\geq0$. Is equality possible here for an entire function, or more generally, a function analytic on $|z|=r$, in the sense that
\[b(r+0)=b(r-0)=0?\]
Clunie notes that if $f(z)=(z-1)e^z$, then
\[M(r,f)=(r-1)e^r,\hspace{1cm}r>2, \hspace{1cm}\text{ and }b(2+0)=0.\]\\

\noindent\textbf{Update 2.18}
Let 
\[b(r)=\Big(\frac{d}{d\log r}\Big)^2\log M(r,f).\]
The quantity exists as a left or right limit everywhere and is non-negative. Examples had been constructed by Clunie previously to show that $b(r+0)$ or $b(r-0)$ may be zero. London \cite{MR0444948} has given a positive answer to this question by constructing an example of the form $f(z)=(1-z)e^{\alpha z+\beta z^2+\gamma z^3}$ for which $b(2+0)=b(2-0)=0$. The powers $z^n$ are the only functions for which $b(r)=0$ for a whole interval $r_1\leq r\leq r_2$ of values of $r$. \\
See also Update 2.17.\\

\noindent\textbf{Problem 2.19} If $f(z)$ is an entire function of exponential type, i.e. satisfying \mbox{$|f(z)|\leq Me^{K|z|}$} for some constants $M$, $K$, and if, further, 
$|f(x)|\leq A$ for negative $x$, and $|f(x)|\leq B$ for positive $x$, what is the sharp bound for $|f(z)|$? If $A=B$, it is known that 
\[|f(z)|\leq Ae^{Ky}\]
is true and sharp.\\

\noindent\textbf{Update 2.19}
An analogous problem for subharmonic functions was solved by Gol'dberg and
Levin \cite{MR0165105}. This gives an upper estimate but it is not exact for
entire functions. A sharp bound for $|f(x)|$, for entire functions was found by Eremenko \cite{MR948815}.\\

%\section* *{Iteration Theory}

\noindent\textbf{Problem 2.20} If $f(z)$ is an entire function, the iterates $f_n(z), n=1,2,\ldots$ are defined inductively by
\[f_{n+1}(z)=f(f_n(z)),\hspace{1cm}f_1(z)=f(z).\]
A point $z$ satisfying the equation $f_n(z)=z$, but such that $f_k(z)\neq z$ for $k<n$, is called a \textit{fixed point of exact order $n$}.
Prove that there always exist fixed points of exact order $n$ if $f(z)$ is transcendental, and $n\geq2$. 
For the case of polynomial or rational $f(z)$ see Baker \cite{MR0114007}. For a proof that fixed points of exact order $n$ exist, except for at most one value of $n$, see Baker \cite{MR0169989}.\\

\noindent\textbf{Update 2.20}
Bergweiler \cite{MR1123803} proved that: if $f$ is a transcendental entire function, and $n\geq2$, then $f$ has infinitely many fixed points of exact order $n$. This also follows from a result by Bergweiler \cite{MR1464835}.\\

\noindent\textbf{Problem 2.21} If, in the terminology of Problem 2.20, $z_0$ is a fixed point of exact order $n$ for $f(z)$, the fixed point is called \textit{repelling} if $|{f_n}'(z_0)|>1$. It is a problem of Fatou (see \cite{MR1504797}, \cite{MR1504792}) whether every entire transcendental function $f(z)$ has repelling fixed points. It is shown by Fatou (see \cite{MR1504797}, \cite{MR1504792}) that for rational $f(z)$ (including polynomials), all fixed points of sufficiently high (exact) order are repelling.\\

\noindent\textbf{Update 2.21}
The existence of repelling fixed points was proved for the first time by Baker \cite{MR0226009}, who used Ahlfors' theory of covering surfaces \cite{MR1509465}. Since then, the proof of this important result was generalised to meromorphic functions, and ultimately evolved into an elementary half-page argument of Berteloot and Duval \cite{MR1782673}.\\

\noindent\textbf{Problem 2.22} With the terminology of Problem 2.20, denote by $\mathcal{F}(f)$ the set of points where the sequence $\{f_n(z)\}$ is not normal. Fatou \cite{MR1555220} asks if there is an entire function $f(z)$ for which $\mathcal{F}(f)$ is the whole plane, and, in particular, if this is the case for $f=e^z$.
Since every point of $\mathcal{F}(f)$ is an accumulation point of fixed points of $f(z)$, this is equivalent to asking if the fixed points (of all orders) of $e^z$ are dense in the plane.\\

\noindent\textbf{Update 2.22}
Baker \cite{MR0264071} has shown that if $f(z)=kze^z$, where $k$ is a certain positive constant, then the set of non-normality does indeed occupy the whole plane. Misiurewicz \cite{MR627790} has proved that this is also the case for $f(z)=e^z$.\\
This answers the whole question positively.\\

\noindent\textbf{Problem 2.23} Baker \cite{MR0197725} has proved that if $f(z)$ is a transcendental entire function, then $\mathcal{F}(f)$ is not restricted to a straight line in the plane. This implies (see Problem 2.22) that, given a line $l$, there are fixed points (of sufficiently high order) not belonging to $l$. Is it already true that a transcendental $f(z)$ cannot have all its fixed points of order at most $2$ on $l$? This is indeed true for $f(z)$ of order less than $\frac{1}{2}$.\\

\noindent\textbf{Update 2.23}
Bergweiler, Clunie and Langley \cite{MR1325262} proved the conjecture by showing that, for every transcendental entire function $f$, and every line, infinitely many of the fixed points of every $n$-th iterate, $n\geq2$, do not lie on this line. Bergweiler \cite{MR1464835} improved this by showing that, for every line, there are infinitely many repelling fixed points of each $n$-th iterate, $n\geq2$, which do not lie on this line.  \\

\noindent\textbf{Problem 2.24} Can an entire function have all its zeros and ones on two distinct straight lines, having infinitely many on each line?\\
Edrei has proved (unpublished) that if $l, m$ are intersecting straight lines, then it is impossible for all the zeros of an entire function $f(z)$ to lie on $l$, and all the ones on $m$.\\
\textit{(A. Edrei)}\\

\noindent\textbf{Update 2.24} 
Al-Katifi has noted orally that $f(z)=\sin(z^2)$ has this property on the real and imaginary axes. Ozawa (by letter) made the same observation, and also proved some related uniqueness theorems.\\
When the lines in question are parallel, there are no entire functions
with zeros on one line and ones on another, except one explicitly listed
exceptional function. This was proved independently by
Baker \cite{Baker1964ae} and Kobayashi \cite{kobayashi1979}.\\
Bergweiler, Eremenko and Hinkkanen \cite{bergweiler_eremenko_hinkkanen_2018} proved a number of results on this subject, generalising Edrei's theorem.\\

\noindent\textbf{Problem 2.25} If $f, g$ are linearly independent entire functions of order $\rho$, which is not a positive multiple of $\frac{1}{2}$, can $fg'-gf'$ have order less than $\rho$? This is possible if $\rho=n/2, n\geq 2$. Clunie (unpublished) proved the result if $\rho<\frac{1}{3}$.\\
\textit{(A. Edrei)}\\

\noindent\textbf{Update 2.25} The original statement of this problem contained `distinct' instead of `linearly independent'. A negative answer follows from a stronger result by Eremenko \cite{MR1266090}, who writes that his solution of Problem 1.33 in \cite{MR1266090} also contains a solution of Problem 2.25. See Update 1.33.\\

\noindent\textbf{Problem 2.26} What is the least integer $k=k(N)$, such that every entire function $f(z)$ can be written as
\[f(z)=\sum^k_{\nu=1}[f_\nu(z)]^N,\]
where $f(z)$ and $f_\nu(z)$ are entire functions? It is enough to solve the problem for $f(z)=z$, since then one can substitute $f(z)$ for $z$. The equation
\[z=\frac{1}{N^2}\sum^N_{\nu=1}\frac{(1+\omega_\nu z)^N}{\omega_\nu},\]
where $\omega_\nu$ are the distinct $N$-th roots of unity, shows that $k(N)\leq N$. On the other hand, for $N=1, 2, 3$ we have $k(N)=N$.\\
To see e.g. that $k(3)\geq 3$, suppose that
\[z=f^3+g^3=(f+g)(f+\omega g)(f+\omega^2 g),\]
where $\omega=exp(2\pi i/3)$. It follows that the meromorphic function $\phi(z)=f(z)/g(z)$ satisfies $\phi(z)\neq-1, -\omega, -\omega^2$, except possibly at $z=0$. Thus, by Picard's Theorem, $\phi(z)$ must be rational, and so $\phi(z)$ assumes at least two of the three values $-1, -\omega, -\omega^2$. This gives a contradiction.\\
\textit{(H.A. Heilbronn)}\\

\noindent\textbf{Update 2.26}
Let $P$ be the class of polynomials, $E$ that of entire functions, $R$ that of rational functions and $M$ that of meromorphic functions in the plane. When $N\geq2$, and $C$ is one of the above classes, denote by $g_C(N)$ the least integer $k$ such that every $f$ in $C$ can be written as
\[f(z)=\sum^k_{\nu=1}f_\nu(z)^N,\hspace{1cm}\text{with }f_\nu(z)\in C.\]
An easy argument shows that it is enough to consider $f(z)=z$. The problem is to find or to estimate the numbers $g_C(N)$.\\ 
Since a number of authors have obtained partial results, it is worth noting that all known lower bounds for $k$ are immediate consequences of an old result of Cartan \cite{zbMATH03010843}. See also Gundersen and Hayman \cite[Theorem 5.1]{zbMATH02113438}.\\
This yields
\[k(k-1)>N\text{ in }P,\text{ if } N>2;\hspace{1cm}k(k-1)\geq N \text{ in }E;\]
\[k^2-1>N\text{ in }R; \hspace{1cm}k^2-1\geq N \text{ in }M.\]
These results are all sharp for $k=2$. (For the case $k=2$, $N=3$ in $M$, see Gross and Osgood \cite{MR722895}.) But for large $k$, no upper bounds substantially better than $k\leq N$ are known for any of the above classes.\\
For the related equation $1=\sum^k_{\nu=1}f_\nu(z)^N$, where the terms on the right hand side are supposed linearly independent, the above lower bounds for $k$ remain valid, and are proved in the same way; but an example due to Molluzzo \cite{Molluzzo}, quoted by Newman and Slater \cite{MR544895}, shows that they give the correct order of magnitude $\sqrt{N}$ for $k$ as a function of $N$.\\
It follows from this that the functional equation $f^n+g^n+h^n=1$ cannot have non-constant meromorphic solutions for $n\geq9$. Gundersen (\cite{MR1821651} ,\cite{MR1660942}) constructed examples of transcendental meromorphic solutions for $n=5$ and $n=6$. Thus, only the cases $n=7$ and $n=8$ remain unsolved.\\

\noindent\textbf{Problem 2.27} Let $\phi_1, \ldots, \phi_n$ denote entire functions of the form
\begin{equation}
\label{2.4}
\phi(z)=\sum e^{f_\nu(z)}/\sum e^{g_\nu(z)}
\end{equation}
where $f_\nu(z), g_\nu(z)$ are entire functions. Does there exist an entire function $f(z)$, not of the form $\phi(z)$, but satisfying an algebraic equation of the form \mbox{$f^n+\phi_1f^{n-1}+\ldots+\phi_n=0$}? The special cases $n=2$, or when the $f_\nu(z)$ are linear polynomials, may be easier to settle. Note that
\[f(z)=\frac{\sin \pi z^2}{\sin \pi z}\]
is not of the form $\sum e^{f_\nu(z)}$, although it is a ratio of such functions.\\

\noindent\textbf{Update 2.27}
No progress on this problem has been reported to us.\\

\noindent\textbf{Problem 2.28} A meromorphic function $f(z)$ in the plane, is said to be of \textit{bounded value distribution (b.v.d.)} if, for every positive $r$, there exists a fixed constant $C(r)$ such that the equation $f(z)=w$ never has more than $C(r)$ roots in any disc of radius $r$. (It is clearly enough to make the assumption for a single value of $r$.)
\begin{enumerate}[(a)]
\item If $f(z)$ is an entire function, prove that $C(r)=O(r)$ as $r\to\infty$, so that $f(z)$ has exponential type at most.\\
If a differential equation
\begin{equation}
\label{2.5}
y^{(n)}+f_1(z)y^{(n-1)}+\ldots+f_n(z)y=0,
\end{equation}
where the $f_n(z)$ are entire functions, has only b.v.d. solutions, Wittich \cite{MR0197828} proves that the $f_\nu$ are all constants. The converse is also true.
\item Is it sufficient to make the basic assumption, not for all values $w$, but for only three such values, to assure that $f(z)$ is of b.v.d.?
\end{enumerate}
\textit{(P. Tur\'an)}\\

\noindent\textbf{Update 2.28}
A positive answer to the first part was provided by Hayman \cite{MR0316693}, when he showed that b.v.d. functions are precisely those whose derivatives have bounded index. (Let $f(z)$ be an entire function and for each $z$ let $N(z)$ be the least integer such that 
\[\sup_{0\leq j<\infty}\Big|\frac{f^{(j)}(z)}{j!}\Big|=\frac{|f^{(N)}(z)|}{N!}.\]
If $N(z)$ is bounded above for varying $z$, then $f(z)$ is said to be of \textit{bounded index}, and the least upper bound $N$ of $N(z)$ is called the \textit{index} of $f(z)$.)\\ 
A negative answer to the second part was provided by Gol'dberg \cite{MR0320316}. Another such example is given by the sigma-function. If $E$ is any bounded set, and $\varepsilon<\frac{1}{2}$, then, when $\omega\in E$, the equation $\sigma(z)=\omega$ has exactly one root in $|z-m-in|<\varepsilon$ for integers $m, n$ with $m^2+n^2>r_0$ for $r_0$ depending on $\varepsilon$ and $E$. But $\sigma(z)$ has order $2$ and so cannot be a b.v.d. function.\\

\noindent\textbf{Problem 2.29} Is it possible to give an analogous characterisation of the solutions of (\ref{2.5}) in the case where the $f_\nu(z)$ are polynomials?\\
\textit{(P. Tur\'an)}\\

\noindent\textbf{Update 2.29}
Such a characterisation has been provided by Tijdeman (\cite{MR0259075}, \cite{MR0286986}). The coefficients $f_\nu(z)$ of the differential equation
\begin{equation}
\label{(3)}
y^{(n)}+f_1(z)y^{(n-1)}+\ldots+f_n(z)y=0
\end{equation}
are polynomials, if and only if there exist fixed numbers $p$ and $q$ such that each solution $g(z)$ of (\ref{(3)}) is $p$-valent in any disc $\{z:|z-z_0|<1/(1+r^q)\}$, where $r=|z_0|$.\\

\noindent\textbf{Problem 2.30} Let $S_k, k=1, 2, \ldots$ be sets which have no finite limit points. Does there exist a sequence $n_k$ and an entire function $f(z)$, so that whenever $z\in S_k$ we have $f^{(n_k)}(z)=0$?\\
\textit{(P. Erd\"os)} \\

\noindent\textbf{Update 2.30}
Functions satisfying the conditions of this problem have been constructed by Barth and Schneider \cite{MR0293089}.\\

\noindent\textbf{Problem 2.31} Let $A, B$ be two countable dense sets in the plane. Does there exist an entire function $f(z)$, so that $f(z)\in B$, if and only if $z\in A$? If the answer is negative, it would be desirable to have conditions on $A, B$ when this is so.\\
\textit{(P. Erd\"os)}\\

\noindent\textbf{Update 2.31}
Functions satisfying the conditions of this problem have been constructed by Barth and Schneider \cite{MR0301195}.\\

\noindent\textbf{Problem 2.32} Let $f(z)=\sum^\infty_{n=0} a_nz^n$ be a transcendental entire function where $a_n\geq0$ for $n\geq0$, and set
\[p_n(z)=\frac{a_nz^n}{f(z)}.\]
Then 
\[\sum^\infty_{n=0}p_n(z)=1.\]
In addition, if $f(z)=e^z=\sum^\infty_{n=0} z^n/n!$, we have 
\begin{equation}
\label{2.6}
\int^\infty_0 p_n(z) \,dz=1, \hspace{1cm} n=0\text{ to }\infty,
\end{equation} 
or, equivalently,
\begin{equation}
\label{2.6'}
\int^\infty_0 \frac{f(\rho z)}{f(z)}\,dz=\frac{1}{1-\rho}, \hspace{1cm} 0<\rho<1.
\end{equation} 
Does there exist any transcendental entire function $f(z)$ other than $e^z$ satisfying (\ref{2.6}) or (\ref{2.6'})?\\
\textit{(A. R\'enyi, St. Vincze)}\\

\noindent\textbf{Update 2.32}
Let $f(z)=\sum^\infty_{n=0}a_nz^n$ be a transcendental entire function with positive coefficients, and suppose that 
\[\int^\infty_0\frac{a_nz^n}{f(z)}\,dz=1,\hspace{1cm}n=0, 1, 2, \ldots .\]
R\'enyi and Vincze had asked whether these conditions imply that $f(z)=ce^z$, and this has been proved by Miles and Williamson \cite{MR829392}. Weaker results were proved earlier by Hall and Williamson \cite{MR0387585} and Hayman and Vincze \cite{MR513412}. \\

%\section* *{Minimum Modulus Problems}

\noindent\textbf{Notation} Let $f(z)$ be an entire function of order $\rho$, and lower order $\lambda$, and let 
\[m_0(r,f)=\inf_{|z|=r}|f(z)|,\hspace{1cm}M(r,f)=\sup_{|z|=r}|f(z)|.\]
It is a classical result that
\[\limsup_{r\to\infty}\frac{\log m_0(r,f)}{\log M(r,f)}\geq C(\lambda).\]
Here $C(\lambda)=\cos(\pi\lambda)$ for $0\leq\lambda\leq1$, and Hayman \cite{MR0056083} has shown that 
\[-2.19\log\lambda<C(\lambda)<-0.09\log\lambda\]
when $\lambda$ is large. We refer the reader to Barry \cite{PLMS:PLMS0445} for a general account
of the situation. For functions of infinite order the analogous result
is
\[\limsup_{r\to\infty}\frac{\log m_0(r,f)}{\log M(r,f)\log\log\log M(r,f)}\geq C_\infty,\]
where
\[-2.19 < C_\infty < -0. 09.\]

\noindent\textbf{Problem 2.33} Is it possible to obtain the exact value of $C_\infty$, or the asymptotic behaviour of $\frac{C_\lambda}{\log\lambda}$ as $\lambda\to\infty$? The question seems related to the number of zeros a function can have in a small disc centred on a point of $|z|=r$, see Hayman \cite{MR0056083}.\\

\noindent\textbf{Update 2.33}
No progress on this problem has been reported to us.\\

\noindent\textbf{Problem 2.34} Is it possible to say something more precise about $C(\lambda)$ when $\lambda$ is just greater than $1$? In particular, is it true that $C(\lambda)=-1$ for such $\lambda$, or alternatively, is $C(\lambda)$ a strictly decreasing function of $\lambda$?\\

\noindent\textbf{Update 2.34} 
Whenever $\rho>1$, Fryntov \cite{MR1687884} constructed an entire function $f$ of order $\rho$ with the property \[\limsup_{r\to\infty}(\log L(r))/(\log M(r,f))<-1.\]

\noindent\textbf{Problem 2.35} If $\Gamma$ is a continuum that recedes to $\infty$, it is known (see Hayman \cite{MR0056083}) that as $z\to\infty$ on $\Gamma$,
\[\limsup_{r\to\infty}\frac{\log |f(z)|}{\log M(|z|)}\geq -A,\]
where $A$ is an absolute constant. Is it true that $A=1$? This is certainly the case if $\Gamma$ is a ray through the origin, see Beurling \cite{MR0029980}. If $A>1$, is it possible to obtain a good numerical estimate for $A$?\\

\noindent\textbf{Update 2.35}
Hayman and Kjellberg \cite{MR820231} gave a positive answer by proving that, for any non-constant subharmonic function $u$, and $A>1$, the set \mbox{$\{z:u(z)+AB(z)\}$ }
where $B(z)=\max_{|\zeta|=|z|}u(\zeta)$ has no unbounded components. Furthermore, if the set $\{z:u(z)+B(|z|)<0\}$ has an unbounded component, then $u$ has infinite lower order; or else regular growth and mean or minimal type of order $\rho$, where $0<\rho<\infty$; or $u$ is linear.\\

\noindent\textbf{Problem 2.36} Suppose that $0<\rho<\alpha\leq1$, where $\rho$ is the order of an entire function $f$. Let $E_\alpha$ be the set of $r$ for which $\log m_0(r,f)>\cos (\pi\alpha) \log M(r,f)$. Besicovitch \cite{MR1512383} showed that the upper density of $E_\alpha$ is at least $1-\rho/\alpha$, and Barry \cite{MR0156993} proved the stronger result, that the same is true of the lower logarithmic density of $E_\alpha$. Examples given by Hayman \cite{MR0274761}, show that Barry's theorem is sharp; in these examples, the logarithmic density exists, but the upper density is larger. This suggests that Besicovitch's theorem may be sharpened.\\

\noindent\textbf{Update 2.36} 
No progress on this problem has been reported to us.\\

\noindent\textbf{Problem 2.37} Let $r_n$ be a sequence of P\'olya peaks (as defined by Edrei \cite{MR0267098}) of order $\rho$. Then Edrei \cite{MR0267098} showed that there exists $K=K(\alpha,\rho)$ such that $\log m_0(r,f)>\cos  (\pi\alpha) \log M(r,f)$ for some value $r$ in the interval $r_n\leq r\leq Kr_n$ and $n$ sufficiently large. Is $K(\alpha,\rho)$ independent of $\alpha$ for fixed $\rho$? Can it be taken arbitrarily near $1$?\\
\textit{(D. Drasin and A. Weitsman)}\\

\noindent\textbf{Update 2.37}
This has been proved by Eremenko, Shea and Sodin \cite{MR866467} for the P\'olya peaks of $N(r,0,f)$ and they have given an example showing that the answer is negative for the P\'olya peaks of $\log M(r,f)$ or $T(r,f)$.\\

\noindent\textbf{Problem 2.38} It was shown by Kjellberg \cite{MR0159942} that if $0<\alpha<1$ and 
\[\log m_0(r,f)<\cos (\phi\alpha)\log M(r,f)+O(1),\hspace{1cm}\text{ as }r\to\infty,\]
then
\[\lim_{r\to\infty}\frac{\log M(r,f)}{r^\alpha}=\beta,\]
where $0<\beta\leq\infty$. If $\alpha=1$, it was shown by Hayman \cite{Hayman1975} that unless $f(z)=Ae^{Bz}$, the corresponding result holds with $\beta=\infty$. Examples constructed by Hayman show that $m_0(r,f)M(r,f)\to\infty$ as $r\to\infty$ can occur for a function of order $1+\varepsilon$ for every positive $\varepsilon$. The case of functions of order $1$ and maximal type remains open.\\

\noindent\textbf{Update 2.38}
Drasin \cite{MR1405053} constructed  an entire function of order $1$ with the property $m_0(r,f)M(r,f)\to0$.\\

\noindent\textbf{Problem 2.39} We can also compare $m_0(r,f)$ with the characteristic $T(r)$. We have
\[\limsup_{r\to\infty}\frac{\log m_0(r,f)}{T(r)}\geq D(\lambda)\]
and ask for the best constant $D(\lambda)$. In view of Petrenko's solution of Problem 1.17, we certainly have $D(\lambda)\geq-\pi\lambda$ with $1\leq\lambda<\infty$. Also, Ess\'en and Shea \cite{MR0450584} show that $D(\lambda)\leq\frac{\pi\lambda}{1+|\sin(\pi\lambda)|}$ for $1<\lambda<\frac{3}{2}$, and $D(\lambda)\leq\frac{-\pi\lambda}{2}$ for $\frac{3}{2}<\lambda<\infty$. Further, it follows from results of Valiron \cite{MR1556886},  and Edrei and Fuchs (\cite{MR0477051} and \cite{MR0585641}) that
\[D(\lambda)=
\begin{cases}
\pi\lambda\cot(\pi\lambda) & \text{if }0\leq\lambda<\frac{1}{2}\\
\pi\lambda\cos(\pi\lambda) & \text{if }\frac{1}{2}\leq\lambda<1.
\end{cases}
\]
\textit{(D. Shea)}\\

\noindent\textbf{Update 2.39}
No progress on this problem has been reported to us.\\

\noindent\textbf{Problem 2.40} Let $f(z)$ be a non-constant entire function, and assume that for some constant $c$ the plane measure of the set $E(c)$ where $|f(z)|>c$ is finite. What is the minimum growth rate of $f(z)$? Hayman conjectures that
\[\int^\infty_0\frac{r\,dr}{\log\log M(r,f)}<\infty\]
is true and best possible.
If $E(c)$ has finite measure, is the same true for $E(c')$ for $c'<c$?\\
\textit{(P. Erd\"os)}\\

\noindent\textbf{Update 2.40}
Erd\"os had asked about the minimum growth of non-constant entire functions, bounded outside a set of finite area. Camera \cite{Camera} has established a conjecture of Hayman  according to which
\[\int^\infty_0\frac{r\,dr}{\log\log M(r,f)}<\infty\]
is true and best possible in the following sense: if $\phi(r)$ increases, and 
\[\int^\infty_0\frac{r\,dr}{\phi(r)}=\infty,\] 
then there exists $f(z)$ entire and such that $\log \log M(r,f)<\phi(r)$, and $f$ is bounded outside a set of finite area. Camera has also established the analogous result for subharmonic functions $u(x)$ in $\mathbb{R}^m$. If $B(r)=\sup_{|x|=r}u(x)$, then
\[\int^\infty_0\Big(\frac{r}{\log B(r)}\Big)^{m-1}dr<\infty\]
is true, and this is best possible in the same sense as above.\\
This was also proved independently Hansen \cite{MR0460634} and by Gol'dberg \cite{MR537357}. Gol'dberg also answered the second part of this problem by giving an example of a function $f$ for which $A(c)$ is finite for some $c$, but not for all $c$.\\

\noindent\textbf{Problem 2.41} Suppose that $f(z)$ has finite order, and that $\Gamma$ is a rectifiable path on which $f(z)\to\infty$. Let $\ell(r)$ be the length of $\Gamma$ in $|z|<r$. Find such a path for which $\ell(r)$ grows as slowly as possible, and estimate $\ell(r)$ in terms of $M(r,f)$. If $f(z)$ has zero order, or more generally, finite order, can a path be found for which $\ell(r)=O(r)$ as $r\to\infty$? If $\log M(r,f)=O(\log^2 r)$ as $r\to\infty$, but under no weaker growth condition, it is shown by Hayman \cite{MR0111839} and Piranian \cite{MR0111838} that we may choose a ray through the origin for $\Gamma$.\\
If $f(z)$ has a finite asymptotic value $a$, the corresponding question may be asked for paths on which $f(z)\to a$.\\
\textit{(P. Erd\"os)}\\

\noindent\textbf{Update 2.41}
This is a refined form of Problem 2.7. It was completely solved by Gol'dberg and Eremenko \cite{GoldbergErem} who showed 
that for every function $\phi(r)$ tending to infinity there exists an
entire function $f$ such that $\ell(r)\neq O(r)$ for every asymptotic curve.\\
Moreover, for every $\rho>1/2$ there exists an entire function of order $\rho$
with finite asymptotic value $a$ such that $\ell(r)\neq O(r)$ for every
asymptotic curve on which $f(z)\to a$.\\

\noindent\textbf{Problem 2.42} Let $f(z)$ be an entire function (of sufficiently high order) with $l$, $l\geq2$ different asymptotic values $a_k$, $k=1,\ldots, l$. Suppose that $\gamma_k$ is a path such that $f(z)\to a_k$ as $z\to\infty$, $z\in\gamma_k$. Let $n_1(r,a_k)$ be the number of zeros of $f(z)-a_k$ on $\gamma_k$, and in $|z|\leq r$. Can we find a function $f(z)$ such that
\[\frac{n_1(r,a_k)}{n(r,a_k)}\to b_k>0\]
as $r\to\infty$, for $k=1, 2, \ldots, n$? Can we take $b_k=1$?\\
\textit{(J. Winkler)}\\

\noindent\textbf{Update 2.42}
Winkler has supplied the following additional information: the problem deals with the question of value distribution in the tongues which are excluded in Ahlfors' island theorems (see for example Hayman \cite[Chp. 5]{M.F.}). The problem is connected with two previous papers of Winkler (\cite{MR0261006}, \cite{MR0288268}).
A complete answer was obtained by Barsegyan \cite{Barsegyan} who has shown that 
$\sum b_k\leq 1$ for entire functions, and $\sum b_k\leq 2$
for meromorphic functions.\\

\noindent\textbf{Problem 2.43} Let $f(z)$ be a transcendental entire function which permutes the integers, i.e. gives an injective mapping of the integers onto themselves. Is it true that $f(z)$ is at least of order $1$, type $\pi$? We can also ask the corresponding question for a function permuting the positive integers with the same conjectured answer. Note that $f(z)=z+\sin z$ satisfies both conditions and is of order $1$, type $\pi$.\\
If $f(z)$ assumes integer values on the positive integers, then Hardy and P\'olya proved (see \cite[Theorem 11, p. 55]{WhittakerJohnM1964Ift}) that $f(z)$ is at least of order $1$, type $2$; and if $f(z)$ assumes integer values on all the integers, then Buck \cite{MR0029984} has shown that $f(z)$ is at least of order $1$, type $\log(\frac{3+\sqrt{5}}{2})=0.962\ldots$.\\
\textit{(L. A. Rubel)}\\

\noindent\textbf{Update 2.43}
There was a mistake in the original statement of this problem, where $\log$ was omitted in the final sentence. 
Linden points out a contradiction with the original statement: that the function $f(z)=z+4\cos(\pi(z-2)/3)$ gives $(i)$ an injective mapping of the integers onto themselves, and $(ii)$ an injective mapping of the positive integers onto themselves and has type $\pi/3$.  There is no contradiction with the corrected statement.
Rubel adds that it would be interesting to find an entire function of reasonably small type that `really' permutes the positive integers, i.e. really scrambles them badly, as all the known examples are pretty tame.\\

\noindent\textbf{Problem 2.44} For $f(z)$ entire of order $\rho$, and non-constant, let $\nu(r)$ be the number of points on $|z|=r$ where $|f(z)|=1$. Is it true that
\[\limsup_{r\to\infty}\frac{\log \nu(r)}{\log r}=\rho\,?\]
If one replaces $\nu(r)$ by the number of points on $|z|=r$, where $f(z)$ is real, then Hellerstein and Korevaar \cite{MR0168761} have shown that the corresponding upper limit is always equal to $\rho$.\\
\textit{(J. Korevaar)}\\

\noindent\textbf{Update 2.44}
The example $f(z)=e^z$ show this to be false ($\nu(r)=2$), but Miles considers the above problem for the real part of meromorphic functions $F$. If $F=u+iv$, and $f=(F-1)(F+1)$, then $u=0 \Leftrightarrow |f|=1$. Thus the problem of the number $\phi(r)$ of real points for meromorphic functions is equivalent to that of points of modulus $1$ for meromorphic functions. For this problem, Miles and Townsend \cite{MR0492261} have shown that if
\[\Phi(r)=\int^r_0\frac{\phi(t)dt}{t},\]
then
\[\limsup_{r\to\infty}\frac{\log \Phi(r)}{\log r}\leq \rho\hspace{1cm}\text{ and }\hspace{1cm}\limsup_{r\to\infty}\frac{\log \phi(r)}{\log r}\leq \rho,\]
outside an exceptional set, which may exist if $f$ is meromorphic. However, for the original entire problem, there is no exceptional set.\\

\noindent\textbf{Problem 2.45} Let $J_0(z)$ be the Bessel function of order zero. Is it true that the equation $J_0(z)=1$ has at most one solution on each ray from the origin? An affirmative answer would show that the exceptional set in a theorem of Delsarte and Lions \cite{MR0102672} is, in fact, void. Asymptotic estimates show that there can be at most a finite number of solutions on any ray, and yield even stronger information.\\
\textit{(L. Zalcman)}\\

\noindent\textbf{Update 2.45} 
No progress on this problem has been reported to us.\\

\noindent\textbf{Problem 2.46} Let $\{f_\alpha(z)\}$ be a family of entire functions, and assume that for every $z_0$, there are only denumerably many distinct values of $f_\alpha(z_0)$. Then if $c=2^{\mathfrak{N}_0}>\mathfrak{N}_1$, the family $\{f_\alpha(z)\}$ is itself denumerable. The above result was proved by Erd\"os \cite{MR0168482}. If $c=\mathfrak{N}_1$, he constructed a counter-example.\\
Suppose now that $m$ is an infinite cardinal, $\mathfrak{N}_0<m<c$. Assume that for every $z_0$, there are at most $m$ distinct values $f_\alpha(z_0)$. Does it then follow that the family $\{f_\alpha(z)\}$ has at most power $m$? If $m^+<c$, where $m^+$ is the successor of $m$, it is easy to see that the answer is `yes'. However, if $c=m^+$, the counter-example fails. It is possible the question is undecidable.\\
\textit{(P. Erd\"os)} \\

\noindent\textbf{Update 2.46} 
No progress on this problem has been reported to us.\\

\noindent\textbf{Problem 2.47} Let $E_\rho$ be the linear space of entire functions $f$ such that \mbox{$|f(z)|\leq B\exp(A|z|^\rho)$} for some positive $A$ and $B$. Let $K_\rho$ be the family of functions $k(z)$ positive and continuous on $\mathbb{C}$, with $\exp(A|z|^\rho)=o(k(z))$ as $|z|\to\infty$, for all positive $A$. Let $S$ be a subset of $\mathbb{C}$, and $\|\cdot\|_{k,S},\|\cdot\|_k$ the semi-norms defined for $f\in E_\rho, k\in K_\rho$ by
\[\|f\|_{k,S}=\sup_S\Big\{\frac{|f(z)|}{k(z)}\Big\},\]
\[\|f\|_k=\sup_\mathbb{C}\Big\{\frac{|f(z)|}{k(z)}\Big\}.\]
We say that $S$ is a \textit{sufficient set} for $E_\rho$ if the topologies defined by the semi-norms $\{\|\cdot\|_k,k\in K_\rho\},\{\|\cdot\|_{k,S},k\in K_\rho\}$ coincide, see Ehrenpreis \cite{MR0285849}.
\begin{enumerate}[(a)]
\item Suppose that, whenever $f\in E_\rho$, $f$ is bounded on $S$ if and only if $f$ is bounded on $\mathbb{C}$. Does it follow that $S$ is a sufficient set for $E_\rho$?
\item Suppose that $S$ is a sufficient set for $E_\rho$. Then, if $f\in E_\rho$ and $f$ is bounded on $S$, does it follow that $f$ is bounded on $\mathbb{C}$, or (maybe) of small growth?
\end{enumerate}
Some recent work of Oliver \cite{Oliver} suggests the latter, at least, is likely.\\
\textit{(D. A. Brannan)}\\

\noindent\textbf{Update 2.47}
No progress on this problem has been reported to us.\\

\noindent\textbf{Problem 2.48} If $A, B$ are countable dense subsets of $\mathbb{R}$, $\mathbb{C}$ respectively, does there necessarily exist a transcendental entire function that maps $A$ onto $B$, and $\mathbb{R} \setminus A$ into $\mathbb{C} \setminus B$?\\
Suppose that $E, F$ are countable dense subsets of $\mathbb{R}$, and that there exists an entire function $f$, monotonic on $\mathbb{R}$, that maps $E$ onto $F$, and $\mathbb{R}\setminus E$ onto $\mathbb{R}\setminus F$. Find interesting conditions which imply that $f$ is trivial. For example, if $E, F$ are real rationals, under what conditions is $f$ necessarily linear with rational coefficients? One could also investigate this question in the case of  real-valued harmonic or subharmonic functions in $\mathbb{R}^n$,  $n\geq2$.\\
If $A, B$ are two countable dense subsets of $\mathbb{R}$, Barth and Schneider \cite{MR0269834} have shown that there exists a transcendental entire function, monotonic on $\mathbb{R}$, that maps $A$ onto $B$ and $\mathbb{R}\setminus A$ onto $\mathbb{R}\setminus B$; also, if $A,B$ are countable dense subsets of $\mathbb{C}$, see Barth and Schneider \cite{MR0301195}.\\
\textit{(K. F. Barth) }\\

\noindent\textbf{Update 2.48}
No progress on this problem has been reported to us.\\

\noindent\textbf{Problem 2.49} If $f(z)$ is a transcendental entire function, we define
\[M=\{z:|f(z)|=M(|z|,f)\}.\]
Tyler \cite{MR1662226} has shown that $M$ can have isolated points, and that, given any $N$, we can have $\nu(r)>N$ for infinitely many $r$, where $\nu(r)$ is the number of points in $M \cap \{|z|=r\}$. Herzog and Piranian \cite{MR0257359} have shown that $\limsup_{r\to\infty}\nu(r)$ can be infinite; is it true that $\liminf_{r\to\infty}\nu(r)<\infty$ for all entire $f$?\\
\textit{(J. G. Clunie)}\\

\noindent\textbf{Update 2.49}
No progress on this problem has been reported to us.\\

\noindent\textbf{Problem 2.50} Characterise those entire functions having at least one continuous maximum modulus path going from $0$ to $\infty$.\\
\textit{(W. Al-Katifi)}\\

\noindent\textbf{Update 2.50} 
No progress on this problem has been reported to us.\\

\noindent\textbf{Problem 2.51} Suppose that an entire function $f$ has exactly one curve $\Gamma$ of maximum  modulus (that is, $\Gamma$ is connected, joins $0$ to $\infty$, and $f$ has no other maximum modulus points). What can be said about the minimum rate of growth of $M(r,f)$, if one is given information about the geometry of the curve $\Gamma$, for example, that $\Gamma$ is a given infinitely-spiralling spiral? If $\Gamma$ is a radial line, clearly nothing much can be said. (In a sense, this is the opposite of a Phragm\'en-Lindel\"of principle).\\
\textit{(D. A. Brannan)}\\

\noindent\textbf{Update 2.51}
No progress on this problem has been reported to us.\\

\noindent\textbf{Problem 2.52} What is the best function $g(\sigma)$, $\sigma\geq 0$ such that, for a non-constant entire function $f(z)$ with maximum and minimum modulus $M(r,f)$ and $m_0(r,f)$ respectively, the assumption
\[\limsup_{r\to\infty}\frac{\log M(r,f)}{(\log r)^2}\leq\sigma\]
implies that
\[\limsup_{r\to\infty}\frac{m_0(r,f)}{M(r,f)}\geq g(\sigma)\,?\]
\textit{(P. D. Barry)}\\

\noindent\textbf{Update 2.52}
This is solved completely by Gol'dberg \cite{MR540239}, who proved the following stronger result: let $f$ be a meromorphic function of zero order, satisfying
\[\liminf_{r\to\infty}\frac{N(r,0,f)+N(r,\infty,f)}{\log^2 r}\leq\sigma<\infty,\]
then
\[\limsup_{r\to\infty}\frac{\min_\theta|f(re^{i\theta})|}{\max_\theta|f(re^{i\theta})|}\geq C(\sigma),\]
where
\[C(\sigma)=\Big(\prod^\infty_{n=1}\frac{1-q^{2n-1}}{1+q^{2n-1}}\Big)^2,\text{ where }q=\exp(-1/(4\sigma)),\]
and this estimate is best possible.\\
A different proof of Barry's original conjecture was given by Fenton \cite{MR601729}.\\

\noindent\textbf{Problem 2.53} For entire or, more generally, meromorphic functions $f$ and $g$, let `$f\leq g$' mean that, for any sequence $\{z_n\}^\infty_1$ for which $|f(z_n)|\to\infty$, then $|g(z_n)|\to\infty$. For entire functions, it can be proved that, if $f'\leq f$, then $f$ is of exponential type; what can be said if $f''\leq f$? Does $f''\leq f$ imply that $f$ is normal (that is, that $|f'(z)|(1+|f(z)|^2)^{-1}$ is bounded)? (An analogue in the case of the unit disc $\mathbb{D}$ has been proved by Pommerenke \cite{MR0325976}). In the above ordering, does there exist $f\wedge g$ and $f\vee g$ for any two entire functions $f$ and $g$? That is, given $f$ and $g$, does there exist an $h$ such that $h\leq f$ and $h\leq g$, and so that if $k\leq f$ and $k\leq g$, then $k\leq h$? Similarly, for $f\vee g$. Finally, if $f$ is meromorphic and $f'\leq f$, does this imply a growth restriction on $f$, e.g. is the order of $f$ at most two?\\
Note: Any constant $c$ satisfies $c\leq f$ for all entire $f$, and any non-constant polynomial $p$ satisfies $f\leq p$ for any entire $f$. Observe also that $e^{cz}$ are all equivalent if $c>0$.\\
\textit{(L. A. Rubel and J. M. Anderson)}\\

\noindent\textbf{Update 2.53}
No progress on this problem has been reported to us.\\

\noindent\textbf{Problem 2.54} Let $E$ be a closed set in $\mathbb{C}$, with the following properties: $(1)$ there exists a transcendental entire function $f(z)$ that is bounded on $E$; and $(2)$, there exists a transcendental entire function $g(z)$ that is bounded away from $0$ on the complement of $E$. For each such set $E$, must there exist one transcendental entire function that is simultaneously bounded on $E$ and bounded away from $0$ on the complement of $E$?\\
\textit{(L. A. Rubel)}\\

\noindent\textbf{Update 2.54}
No progress on this problem has been reported to us.\\

\noindent\textbf{Problem 2.55} Let $f_i(z)$, $i=1, 2, 3$ be non-constant entire functions of one complex variable, and
\[V=\{z:z=(z_1,z_2,z_3)\in\mathbb{C}^3,f_1(z_1)+f_2(z_2)+f_3(z_3)=0\}.\]
If $F$ is an entire function of $z=(z_1,z_2,z_3)$ that is bounded on $V$, is $F$ necessarily constant on $V$?\\
\textit{(J. M. Anderson and L. A. Rubel)}\\

\noindent\textbf{Update 2.55}
Rubel points out that there has been progress in one very special case by Demailly \cite{MR533896}.\\

\noindent\textbf{Problem 2.56} Prove or disprove the conjecture that an entire function $f$ of $n$ complex variables is an $L$-atom (where this is defined in a way analogous to the definition for $n=1$; see Problem 5.55 with $\mathbb{D}$ replaced by $\mathbb{C}$);  if and only if there are entire functions $f_2,f_3,\ldots,f_n$ of the $n$ variables such that $(f, f_2,f_3,\ldots,f_n)$ is an analytic automorphism of $\mathbb{C}^n$, that is, an injective biholomorphic map of $\mathbb{C}^n$ onto $\mathbb{C}^n$. Rubel can prove the result in the case $n=1$.\\
\textit{(L. A. Rubel)}\\

\noindent\textbf{Update 2.56}
No progress on this problem has been reported to us.\\

\noindent\textbf{Problem 2.57} If $f$ is an entire function such that $\log M(r,f)=O(\log r)^2$ as $r\to\infty$, then Hayman \cite{MR0111839} has shown that $\log |f(re^{i\theta})|\sim\log M(r,f)$, as $r\to\infty$, for $re^{i\theta}$ outside an exceptional set $E$ of circles subtending angles at the origin, whose sum is finite. In particular, 
\[\log |f(re^{i\theta})|\sim\log M(r,f),\text{ as }r\to\infty,\text{ for almost every }\theta.\]
Using this result, Anderson and Clunie \cite{MR0199397} showed that if $f$ is meromorphic and $T(r,f)=O(\log r)^2$ as $r\to\infty$, then a deficient value (there is at most one) must be asymptotic and, moreover, if $\delta(a,f)>0$, then 
\[f(re^{i\theta})\to a\text{ as }r\to\infty\text{ for almost every }\theta.\]
Now consider a new class of entire functions $I_\alpha$, $\alpha>1$, defined by
\begin{enumerate}[(a)]
\item $\log M(r,f)=O(\log r)^{1+\alpha}\hspace{1cm}\text{ as }r\to\infty$
\item $(\log r)^\alpha=o(\log M(r,f))\hspace{1cm}\text{ as }r\to\infty.$
\end{enumerate}
(Hayman's functions are $I_1$). Do the results of Hayman, and the corresponding results of Anderson and Clunie still hold? In other words, do these results depend on smallness of growth, or only on smoothness of growth?\\
\textit{(J. M. Anderson)}\\

\noindent\textbf{Update 2.57} 
This problem has been solved by Anderson \cite{MR536388} as follows: let $f(z)$ be meromorphic in $\mathbb{C}$ and suppose that, for some $a$, $\delta(a,f)>0$. Then under the hypothesis $T(r,f)\sim T(2r,f)$ as $r\to\infty$, $f(z)$ has as the asymptotic value $a$. More precisely, there is an asymptotic path $\Gamma$ such that
\[\liminf_{r\to\infty}\frac{1}{T(r,f)}\log\left|\frac{1}{f(z)-a}\delta(a)\right|\text{ as }z\to\infty\text{ along }\Gamma,\]
and the length $\ell(r)$ of $\Gamma\cap\{|z|<r\}$ satisfies $\ell(r)=r\{1+o(1)\}$ as $r\to\infty$. The hypothesis $T(r,f)\sim T(2r,f)$ forces $f(z)$ to be of order zero, but there are examples due to Gol'dberg and Eremenko of entire functions $g(z)$ satisfying $T(r,g)=O(\psi(r))(\log r)^2$ as $r\to\infty$, where $\psi(r)$ is any assigned real-valued function tending to $\infty$ with $r$ such that no asymptotic satisfies $\ell(r)=O(r)$.\\

\noindent\textbf{Problem 2.58} Suppose that $f$ is entire with a non-zero Picard exceptional value $\alpha$. Then $f$ has $\alpha$ as an asymptotic value. It can be shown that $f\to\alpha$ along a level curve $|f|=|\alpha|$, if $f$ is of finite order. This follows readily from the fact that $\arg f$ is monotone along such a curve, together with an application of the Denjoy-Carleman-Ahlfors theorem. We call such a level curve a \textit{natural asymptotic path} for $\alpha$. Does there exist an entire function of infinite order, with a non-zero Picard exceptional value $\alpha$, having no natural asymptotic paths?\\
\textit{(S. Hellerstein)}\\

\noindent\textbf{Update 2.58} An example of such a function has been given by Eremenko \cite{MR858915}, who constructed a suitable Riemann surface.\\

\noindent\textbf{Problem 2.59} \textbf{(A width conjecture)} Given a power series $\sum^\infty_{k=0}a_kz^k$, suppose that there is a non-negative $\rho$ such that all of the partial sums $S_n(z)=\sum^n_{k=0}a_kz^k$, $n=1, 2, 3, \ldots$ are non-zero in the region
\[S\rho=\{z=x+iy:|y|<Kx^{1-(\rho/2)},x>0\}.\]
We conjecture that $f(z)$ must be entire of order at most $\rho$. When $\rho=0$, $S_\rho$ is a sector with vertex at $z=0$, and the conjecture is a consequence of results of Carlson \cite{MR0025568}, which were later generalised by Rosenbloom \cite{Rosenbloom}, Ganelius \cite{MR0156980}, and Korevaar and McCoy \cite{MR0164005}.
Remark: If $f(z)=e^z$, which is of order $1$, Saff and Varga \cite{MR0414968} have shown that the partial sums $\sum^n_{k=0}z^k/k!$ are in fact zero-free in the parabolic region
\[\{z=z+iy:y^2\leq4(x+1), x>-1\}.\]
\textit{(E. B. Saff and R. S. Varga)}\\

\noindent\textbf{Update 2.59}
No progress on this problem has been reported to us.\\

\noindent\textbf{Problem 2.60} Let $\sum^\infty_{k=0}a_kz^k$ be a non-vanishing entire function, and let  \mbox{$S_n(z)=\sum^n_{k=0}a_kz^k$}. Given $\varepsilon>0$, must there exist a $z_0$ and an $n$ such that $S_n(z_0)=0$ and $|f(z_0)|<\varepsilon$? (The Hurwitz theorem shows that the result is true if $f$ has a zero.)\\
\textit{(D. J. Newman and A. Abian)}\\

\noindent\textbf{Update 2.60}
No progress on this problem has been reported to us.\\

\noindent\textbf{Problem 2.61} Let $\Gamma$ be a rectifiable curve. Suppose $f$ is a continuous function on the plane satisfying
\[\int_{\sigma(\Gamma)}f(z)dz=0\hspace{1cm}\text{ for all rigid motions }\sigma.\]
Does this imply that $f$ is an entire function? The answer is `yes' for some choices of $\Gamma$, and `no' for others. For example, the answer is `no' for the circle ($f(x+iy)=\sin(ax)$ is a counter-example for a suitable choice of $a$); `yes' for an ellipse; `yes' for any polygonal Jordan curve; and `yes' for the boundary of any convex set with at least one corner. Prove that the circle is the only closed rectifiable Jordan curve (or the only curve among the class of curves which are boundaries of bounded convex sets) for which the answer is `no'. This problem is related to the `Pompeiu Problem' discussed by Brown, Schreiber and Taylor \cite{Brown1973}.\\
\textit{(L. Brown)}\\

\noindent\textbf{Update 2.61}
No progress on this problem has been reported to us.\\

\noindent\textbf{Problem 2.62} Let $f$ denote a rational or entire function of a complex variable, and $f^n, n=1, 2, \ldots$, the $n$-th iterate of $f$, so that $f^1=f, f^{n+1}=f\circ f^n=f^n\circ f$. Provided that $f$ is not rational of degree $0$ or $1$, the set $C$ of those points where $\{f^n\}$ forms a normal family is a proper open subset of the plane, and is invariant under the map $z \mapsto f(z)$. A component $G$ of $C$ is a wandering domain of $f$ if $f^k(G)\cap f^n(G)=\emptyset$ for all $\{k,n\,|\,k\geq1,n\geq1,k\neq n\}$.\\
Jakobson has asked whether it is possible for a rational function $f$ to have a wandering domain.
Baker \cite{MR0419759} gave a transcendental entire function which does have such domains.\\
\textit{(I. N. Baker)}  \\

\noindent\textbf{Update 2.62}
Sullivan \cite{MR819553} has proved that rational functions have no wandering domains, so that the answer to Jakobson's question is `no'.\\

\noindent\textbf{Problem 2.63} Let $f$ be a rational function and $C$ be as in Problem 2.62. We say that $g$ is a \textit{limit function} for $f$ if $g$ is defined in some component $G$ of $C$ and is the limit of some subsequence of $(f^n)$ in $G$. In the simplest examples, the number of limit functions is finite, which implies that each has a constant $\alpha$, say, such that $f^k(\alpha)=\alpha$ for some positive integer $k$. If, in addition, $|(f^k)'(\alpha)|<1$ for each of the limit functions, we say that the function $f$ belongs to the class $N$.
\begin{enumerate}[(a)]
\item Does there exist a rational $f$ which has an infinity of constant limit functions?
\item Is the property of belonging to $N$ `generic' in some sense for rational functions?
\end{enumerate}
Some account of the older established results can be found in Fatou (\cite{MR1504797}, \cite{MR1504792}, \cite{MR1504787}) and a sketch from a more modern point of view is given by Guckenheimer \cite{Guckenheimer}.\\
\textit{(I. N. Baker)}\\

\noindent\textbf{Update 2.63}
Sullivan \cite{MR819553} proved that the answer to the first question is `no'.\\
On the other hand, Eremenko and Lyubich \cite{MR918638} have constructed an example of an entire function, iterates of which have an infinite set of constant limit functions in a component of normality. \\

\noindent\textbf{Problem 2.64} Let $f(z)$ be a real entire function (i.e. $f(z)$ is real for real $z$). It has been shown by Hellerstein and Williamson \cite{MR0435393} that if $f, f', f''$ have only real zeros, then $f$ is in the Laguerre-P\'olya class (i.e. $f(z)=e^{-az^2}g(z), a\geq0, g$ a real entire function of genus $0$ or $1$, with only real zeros), thus affirming a conjecture of P\'olya. Wiman  made the stronger conjecture, that the above is true, with no assumption on the zeros of $f'$. Even the simplest case, $f(z)=\exp(Q)$, where $Q$ is a real polynomial, has not been settled (see Problem 4.28). Wiman's conjecture has been proved by Levin and Ostrovskii \cite{MR0130979} for $f$ of infinite order, growing sufficiently fast.\\
\textit{(S. Hellerstein)} \\

\noindent\textbf{Update 2.64} 
Wiman's conjecture has been proved by Bergweiler, Eremenko and Langley \cite{MR2024413}. For functions of finite order, the conjecture was proved previously by Sheil-Small \cite{MR979605}. Both works use the methods developed by Levin and Ostrovskii \cite{MR0130979}.\\

\noindent\textbf{Problem 2.65} Since the knowledge of the zeros of an entire function $f$ leaves an unknown factor, $e^h$ say, in the Hadamard product for $f$, one can ask if $f$ is determined by the zeros of $f$, and of its first few derivatives. Does there exist an integer $k$, $k\geq2$ such that, if $f$ and $g$ are entire, and $f^{(n)}/g^{(n)}$ is entire and non-vanishing for $0\leq n\leq k$, then $f/g$ is constant, unless
\[f(z)=e^{az+b}, g(z)=e^{cz+d}\hspace{1cm}\text{ or }\hspace{1cm}f(z)=A(e^{az}-b), g(z)=B(e^{-az}-b^{-1})\,?\]
The proposer has shown (unpublished) that $k=2$ will do in certain cases; for example, when $f$ and $g$ have finite order. The example
\[f(z)=(e^{2z}-1)\exp(-ie^z), \hspace{1cm}g(z)=(1-e^{-2z})\exp(ie^{-z})\]
shows that one sometimes needs $k=3$.
One can ask a similar question for meromorphic functions, with the additional possibility that
\[f(z)=A(e^{h(z)}-1)^{-1}, \hspace{1cm}g(z)=B(1-e^{-h(z)})\]
for any non-constant entire function $h$.\\
\textit{(A. Hinkkanen)}\\

\noindent\textbf{Update 2.65}
K\"ohler \cite{MR998240} proved that the answer is `yes', for meromorphic functions and $k=6$. Namely, if $f$ and $g$ are meromorphic function, such that $f^{(n)}/g^{(n)}$ are entire and without zeros for $0\leq n\leq6$, then $f$ and $g$ satisfy one of the four relations suggested by Hinkkanen. If one makes additional assumptions about the growth of $f$ and $g$, one needs fewer derivatives to achieve the same conclusion (\cite{MR998240}, \cite{MR1047599}). See also Langley \cite{MR1374462} for related results.\\
Yang \cite{MR2140307} shows that the answer to Hinkkanen's problem is positive for meromorphic functions of finite order if $n = 1$, and an additional condition holds. Yang also generalises this result to meromorphic functions of hyper-order less than one. Examples show that the order restriction is sharp.\\

\noindent\textbf{Problem 2.66} Given a countable number of entire functions, one can find an entire function growing faster than any of these. Without making any assumption about the Continuum Hypothesis, can one associate with every countable ordinal number $\alpha$ an entire function $f_\alpha$ such that
\begin{enumerate}[(a)]
\item if $\alpha<\beta$, then $M(r,f_\alpha)/M(r,f_\beta)\to0$ as $r\to\infty$,
\end{enumerate}
and
\begin{enumerate}[(b)]
\item if $f$ is an entire function, then there exists $\gamma$ such that $M(r,f)/M(r,f_\gamma)\to0$ as $r\to\infty$\,?
\end{enumerate}
See also Problem 7.62.\\
\textit{(A. Hinkkanen)}\\

\noindent\textbf{Update 2.66}
In terms of the cardinality $A$ of the smallest set of entire functions which exhaust all orders of growth, Hinkkanen (unpublished) can show that $A$ can be anything between $\mathfrak{N}_1$ and $C$, having the same general properties as $C$, that is, $A\neq\sum_{i\in I}M_i$, where $\# M_i<A$ and $\#I<A$. The answer is independent of the axioms of set theory.\\

\noindent\textbf{Problem 2.67} Let $f$ be an entire function, and let $D$ be a component of the set in $\mathbb{C}$ where the family of iterates $\{f_n\}$ is normal. Can this family have an infinite bounded set of constant limit functions? Eremenko and Lyubich \cite{MR918638} have shown that the set of constant limit functions may be infinite.\\
\textit{(A. Eremenko)}\\  

\noindent\textbf{Update 2.67}
See Update 2.87 for details on the equivalence of Problem 2.77 and Problem 2.87, and the special case that Problem 2.67 is of these problems. There has been no progress reported to us on any of these problems.\\

\noindent\textbf{Problem 2.68} Let $f$ be an entire function satisfying the condition
\[\log M(r,f)\leq(1+o(1))r^\rho,\hspace{1cm}\text{ as }r\to\infty.\]
Suppose that there exists a curve $\Gamma$ tending to $\infty$ such that on $\Gamma$
\[\log|f(z)|\leq(\alpha+o(1))r^\rho,\hspace{1cm}\text{ as }r=|z|\to\infty,\]
for some $\alpha$ in $[-1,1)$; and denote by $E(r,\varepsilon)$ the angular measure of the set
\[\{re^{i\theta}:\log|f(re^{i\theta})|\leq(1-\varepsilon)r^\rho\}.\]
Eremenko conjectures that 
\[\limsup_{\varepsilon\to0,\, r\to\infty}E(r,\varepsilon)\geq\frac{2}{\rho}\arccos\alpha.\]
Jaenisch \cite{MR0182732} has proved some related results.\\
\textit{(A. Eremenko)}  \\

\noindent\textbf{Update 2.68}
No progress on this problem has been reported to us.\\

%\section* *{Entire functions (including iteration)}

\noindent\textbf{Problem 2.69} Hayman \cite{MR0180679} has shown that 
\begin{equation}
\label{asterix}
\liminf_{r\to\infty}\frac{T(r,f)}{T(r,f')}\leq1
\end{equation}
for transcendental entire functions $f$ of lower order zero. Toppila \cite{MR0585685} has shown that there exists an entire function of order one which does not satisfy (\ref{asterix}).
Does there exist a constant $d > 0$ such that (\ref{asterix}) holds for all transcendental entire
functions $f$ of order less than $d$? \\
\textit{(S. Toppila)}\\

\noindent\textbf{Update 2.69}
If $\rho>1/2$, Langley \cite{MR1200194} constructed an entire function of order $\rho$, with the property
\[\liminf_{r\to\infty}\frac{T(r,f)}{T(r,f')}>1.\]\\

\noindent\textbf{Problem 2.70} Let $H$ be an entire function, let $f_1, f_2$ be linearly independent solutions
of the differential equation $w'' + Hw = 0$, and let $E = f_1 f_2$. Clearly $f_1$, $f_2$ and $E$ are entire. Also, it is
well-known that if $H$ is a polynomial of degree $n$, then the orders of  $f_1$ and $f_2$ are
 $\rho(f_1)=\rho(f_2)=\frac{1}{2}(n+2)$. Furthermore, the exponent of convergence of the zeros of $E$  is $\lambda(E)=\frac{1}{2}(n+2)$, provided that $n > 1$.
If $H$ is transcendental then
\[ \rho(f_1)=\rho(f_2)=+\infty\]
and one might hope that, by analogy with the previous remarks, $\lambda(E)=+\infty$.\\
However, examples by Bank and Laine \cite{MR664047} showed that this is not necessarily the case if $\rho(H)$ is a positive
integer or $+\infty$. They also asked whether $\lambda(E)=+\infty$ if $\rho(H)$
is non-integral and finite and they showed that this is indeed the case if
$\rho(H)<\frac{1}{2}$, a result improved to  $\rho(H)\leq\frac{1}{2}$ by Rossi \cite{MR831388}. What happens in general?\\
\textit{(S. Hellerstein and J. Rossi)}\\

\noindent\textbf{Update 2.70}
This question became known as the Bank-Laine conjecture. A negative answer was given by Bergweiler and Eremenko in \cite{Bergweiler2014OnTB}
 and \cite{2015arXiv151005731B}: for every $\rho > \frac12$ there exists an entire $H$ for which $\rho(H) = \rho$ but $\lambda(E) < + \infty$.\\

\noindent\textbf{Problem 2.71} It is shown by Hellerstein and Rossi \cite{MR847009}, and Gundersen \cite{MR853961} that if $f_1$ and $f_2$ are two linearly independent solutions to the differential equation $w'' + Hw = 0$, where $H$ is a polynomial and $f_1$ and $f_2$ have only finitely many non-real zeros, then $H$ is a non-negative constant. It is also shown
that, if $H(z) = az + b$, for $a,b\in\mathbb{R}$, then the differential equation admits a solution with
only real zeros (and infinitely many of them). Furthermore, as pointed out by Gundersen \cite{MR853961},
Titchmarsh (\cite[pp. 172-173]{MR0176151}) showed that, if $H(z) = z^4-\beta$ for special choices of $\beta$,
then the differential equation also admits solutions with only real zeros (and infinitely
many of them).\\
Characterise all non-constant polynomials $H$ such that the differential equation
admits a solution with only real zeros (and infinitely many of them).\\
\textit{(S. Hellerstein and J. Rossi)}\\

\noindent\textbf{Update 2.71}
Gundersen \cite{MR1833253} showed that for given $a > 0$ and $b\geq≥ 0$, there exists an
infinite sequence of real numbers when $\lambda_k \to +\infty$, so that the differential equation
\[f'' + (az^4 + bz^2 − \lambda_k)f = 0\]
has solutions with infinitely many zeros, and all these zeros, except finitely
many of them, are real. Rossi and Wang \cite{MR1417690} proved that if an equation
$f'' +P(z)f = 0$, where $P$ is a polynomial, has a solution with infinitely many
zeros, all of them real, then the number of real zeros of $P$ must be less than
than $1+\frac{1}{2}\deg P$, counted with multiplicities. \\
Eremenko and Merenkov \cite{MR2210353}
proved that for every $d$ there exist polynomials $P$ of degree $d$ such that some
solution of the equation $f'' + Pf = 0$ has only real zeros. The zero set of
such $f$ can be infinite if and only if $d\not\equiv 2\pmod 4$.\\
Eremenko and Merenkov \cite{MR2210353} show that for every non-negative integer $d$, there exist differential equations $w'' + Pw = 0$, where $P$ is a polynomial of degree $d$, such that some non-trivial solution $w$ has only real roots. For polynomials of degree $3$, Eremenko and Gabrielov \cite{MR2894426} give a curve $\Gamma_0$, which parametrises all such equations; and a curve $\Gamma_n$ which parametrises all such equations having solutions with exactly $2n$ non-real zeros.\\

\noindent\textbf{Problem 2.72} Let $\{f_1,\ldots,f_n\}$ be a fundamental system for the differential equation
\begin{equation}
\label{2.72}
L_n(w)\equiv w^{(n)}+a_{n-1}(z)w^{(n-1)}+\ldots+a_0(z)=0,
\end{equation}
where $a_o,...,a_{n-1}$ are polynomials. Frank \cite{MR0264160} has proved that each function
$f_1,\ldots,f_n$ has finitely many zeros, if and only if (\ref{2.72}) can be transformed into a
differential equation with constant coefficients by a transformation of the form
$w(z) = \exp(q(z) u(z))$, where $q$ is a suitable polynomial.
Does the same result hold if each function $f_1,\ldots,f_n$ is assumed to have only finitely
many non-real zeros? In view of Hellerstein and Rossi \cite{MR847009} and Gundersen \cite{MR853961}, we may assume that $n\geq3$.\\
\textit{(S. Hellerstein and J. Rossi)}\\

\noindent\textbf{Update 2.72}
Br\"uggeman \cite{MR1057941} and Steinmetz \cite{MR1081702} independently gave a positive answer.
In fact, each of them proved a stronger result than conjectured.\\

\noindent\textbf{Problem 2.73} Let $F(z, a, b)$ be an entire function of three complex variables, and suppose
that $F$ is not of the form
\begin{equation}
\label{2.73}
 F(z,a,b) = G(z,H(a,b))
\end{equation}
for any entire functions $G$ and $H$ of two complex variables. Can 
\[\{F(z, a,b):a,b\in\mathbb{C}\}\]
constitute a normal family of entire functions of $z$? (Put rather loosely, does there
exist a two-parameter normal family of entire functions?)
Notice that if $F$ does have the form (\ref{2.73}) then $F_bF_{a,z} = F_aF_{b,z }$, where subscripts
denote partial differentiation. The purpose of ruling out the form (\ref{2.73}) is to ensure that
there are two honest parameters of the family. Otherwise we could have, say,
\[F(z, a, b) = z + 5a^2 + \sin b,\]
which is really a one-parameter family in disguise. \\
\textit{(L. A. Rubel)}\\

\noindent\textbf{Update 2.73}
No progress on this problem has been reported to us.\\

\noindent The next two problems concern the class of real entire functions of finite order
with only real zeros. We may partition this class as follows: for each integer non-negative $p$,
denote by $V_{2p}$ the class of entire functions $f$ of the form
\[f(z)=g(z)\exp(-az^{2p+2})\]
where $a\geq0$, and $g(z)$ is a constant multiple of a real entire function of genus at most
$2p+ 1$ with only real zeros. Then write $U_0 = V_0$, and $U_{2p} = V_{2p}-V_{2p-2}$ for $p\in\mathbb{N}$.\\

\noindent\textbf{Problem 2.74} Suppose that $f(z) = 1 + a_1z + a_2 z^2 +\ldots \in U_{2p}$. If $p = 0$ (so that $f\in U_0)$ and if $f$ is not a polynomial, it is well-known that $f$ cannot have two consecutive Taylor coefficients
equal to zero.\\
If $p > 0$, can an analogous assertion be made? Is it true, for example, that the
Taylor series of $f$ in $U_{2p}$ cannot have $2p + 2$ consecutive coefficients equal to zero?\\
\textit{(J. Williamson)}\\

\noindent\textbf{Update 2.74}
No progress on this problem has been reported to us.\\

\noindent\textbf{Problem 2.75} Suppose that $f$ is entire of proximate order $\rho(r)$, and that $f$ has a
representation as a Dirichlet series
\[f( z ) = \sum^\infty_{n=1}a_ne^{\lambda_n z},\hspace{1cm} 0\leq\lambda_1<\lambda_2<\ldots\to\infty, \hspace{1cm}a_n>0.\]
Can one give a complete characterisation of the indicator
\[h(\theta,f)=\limsup_{r\to\infty}r^{-\rho(r)}\log|f(re^{i\theta})|\]
of such functions? If $f$ has a representation
\[f(z)=\int^\infty_0 e^{iz}\,dF(t)\]
where $F$ is positive and increasing, Gol'dberg and Ostrovskii \cite{MR883581} gave such a
characterisation. \\
\textit{(A. A. Gol'dberg and I. V. Ostrovskii)}\\

\noindent\textbf{Update 2.75} 
Gol'dberg and Ostrovskii \cite{MR1054303} solved the problem under the following
additional assumption on the sequence of the exponents: there exists an
entire function $L(\lambda)$ of exponential type, such that $L(\lambda_n) = 0$ for $n = 1, 2,\ldots,$
and
\[\lim_{n\to\infty} \lambda_n^{-1} \log(1/|L'(\lambda_n)|) < \infty.\]
If these conditions are satisfied, the sequence $(\lambda_n)$ is said to have \textit{finite index
of concentration}. If the sequence of exponents of a Dirichlet series, $f$ has finite
index of concentration, the only possible indicators are $h(f, \theta) = a(\cos^+ \theta)^\rho$,
where $a > 0$. Gol'dberg and Ostrovskii \cite{MR1054303} also obtained other results, with weaker
conditions on the sequence $(\lambda_n)$, and studied the lower indicators.\\

%\section* *{Iteration theory}

\noindent
We draw to the attention of the reader many other problems in iteration raised
in Douady (\cite{MR728980}, \cite{MR936853}), Herman \cite{MR915567}, Lyubich \cite{MYuLyubich1986Tdor}, Beardon \cite{zbMATH01542658}, Bergweiler \cite{zbMATH00487133} and Milnor \cite{zbMATH06257274}. In this section, $f^n$ denotes the \textit{$n$th iterate} of a function $f$.\\

\noindent\textbf{Problem 2.76} Let $\Omega$ be a component of the normal set of an entire function (under
iteration). Is $\dim(\partial\Omega) > 1$? Or is $\partial\Omega$ a circle/line? \\
\textit{(D. Hamilton)}\\

\noindent\textbf{Update 2.76}
Bishop \cite{Bishop2018} has constructed transcendental entire functions with Julia sets of dimension $1$. So, in particular, the boundary of the components of the Fatou set have dimension $1$, but they are not circles or lines.\\

\noindent\textbf{Problem 2.77} Let $\Omega$ be a component of the normal set of an entire function $f$ (under
iteration). Do there exist such an $f$ and such an $\Omega$ with the following properties:
\begin{enumerate}[(a)]
\item $f^n(\Omega)$ is uniformly bounded, for $n = 0,1,2,\ldots$;
\item $f^n(\Omega)\cap f^m(\Omega)=\emptyset$ for $n\neq m$.
\end{enumerate}
\textit{(I. N. Baker, R. Herman and I. Kra; communicated by D. Hamilton)}\\

\noindent\textbf{Update 2.77}
We note that `uniformly' has been added to the original statement of $(a)$, as otherwise there are examples for such $f$.\\ 
See Update 2.87 for details on the equivalence of Problem 2.77 and Problem 2.87, and the special case that Problem 2.67 is of these problems. There has been no progress reported to us on any of these problems.\\

\noindent In the following problems we denote by $R_d$ the class of rational functions of degree
$d$ on the Riemann sphere $\hat{\mathbb{C}}$, and we suppose that $d\geq2$.\\

\noindent\textbf{Problem 2.78} \textbf{(Fatou's conjecture)} Show that the subset $U$ of functions $g$ in $R_d$, such that all the critical points of $g$ are in the basins of attraction of periodic sinks,
is dense in $R_d$.\\
The property `$g\in U$' is also sometimes called Axiom A. See also Fatou \cite{MR1504820}.\\
\textit{(P. Fatou; communicated by M. R. Herman)}\\

\noindent\textbf{Update 2.78}
An analogue of the Fatou conjecture for the real quadratic family 
\[\{z\mapsto z^2 + c : c \in\mathbb{R}\}\]
(instead of the family $R_d$ of all rational functions of degree
$d$ at least $2$) has been established by Graczyk and \'Swiatek \cite{MR1469316}, and by Lyubich \cite{MR1459261}.\\

\noindent The next three problems were suggested in order to give a positive answer to the
Fatou conjecture. See, for example, Ma\~n\'e \cite{MR807068}.\\

\noindent\textbf{Problem 2.79} 
\begin{enumerate}[(a)]
\item Show that, if a function $g$ in $R_d$ has the property that its Julia set
$J(g) \neq \hat{\mathbb{C}}$, then $g$ does not leave invariant a non-trivial Beltrami form on $J(g)$. Here a
\textit{Beltrami form} $\mu$ means that $\mu\in L^{\infty}( \hat{\mathbb{C}})$ and $\|\mu\|_{L^\infty}<1$; and \textit{trivial} on $J(g)$ means that $\mu(x)=0$ for almost all $x$ in $J(g)$.
\item More generally, is the Lebesgue measure of $J(g)$ in $\hat{\mathbb{C}}$ equal to zero? (This is the
analogous conjecture to the Ahlfors conjecture for finitely-generated Kleinian
groups.)
\end{enumerate}
Negative answers to both $(a)$ and $(b)$ have been proved (by McMullen (uncited) and Eremenko and Lyubich \cite{MR918638}) for the class of transcendental entire functions $g$. Douady has conjectured that the answer is negative, and that a counter-example is the function $P_\lambda(z) = \lambda(z + z^2)$, for some $\lambda$ of modulus one. See also Douady \cite{MR728980}, Lyubich \cite{MYuLyubich1986Tdor} and Ma\~n\'e \cite{MR807068}.\\
\textit{(D. Sullivan; communicated by M. R. Herman)}\\

\noindent\textbf{Update 2.79}
\begin{enumerate}[(a)]
\item has been established for the real quadratic family by McMullen \cite{MR1312365},
and extended to the family $\{z\mapsto z^d +c, c \in\mathbb{R}\}$ by Levin and van Strien \cite{MR1637647}.
\item Buff and Ch\'eritat \cite{zbMATH06093940} showed that the answer is `no'  in general, and that a counter-example is given by
 $\lambda(z+z^2)$ for suitable $\lambda$, as conjectured by Douady.
\end{enumerate}

\noindent\textbf{Problem 2.80} Let the function $g$ in $R_d$ have the property that its Julia set $J(g) = \hat{\mathbb{C}}$. Is the dimension $k$ of the space of Beltrami forms on $\hat{\mathbb{C}}$, invariant under $g$, at most one? Also, are certain of the Latt\`es examples the only rational functions such that
$k\neq0$?\\
For further information, see \cite{MR807068}.\\
\textit{(D. Sullivan; slightly modified and communicated by M. R. Herman)}\\

\noindent\textbf{Update 2.80}
No progress on this problem has been reported to us.\\

\noindent\textbf{Problem 2.81} Let the function $g$ in $R_d$ have the property that its Julia set $J(g) = \hat{\mathbb{C}}$. Is $g$ ergodic for Lebesgue measure? In other words, if $B\subset\hat{\mathbb{C}}$ is a Borel-invariant set under $g$ (that is, $g^{-1}(B) = B$), does it follow that either $B$ or $\hat{\mathbb{C}} \setminus B$ has Lebesgue measure zero?\\
For further information, see \cite{MR807068}.\\
\textit{(D. Sullivan; communicated by M. R. Herman)}\\

\noindent\textbf{Update 2.81}
No progress on this problem has been reported to us.\\

\noindent\textbf{Problem 2.82} Let $L_d$ denote the class of those functions $g\in R_d$ such that every critical
point of $g$ is preperiodic but not periodic. Show that, if the function $g$ in $R_d$ has the
property that $J(g) = \hat{\mathbb{C}}$, then $g$ belongs to the closure of $L_d$ in $R_d$. \\
\textit{(M. R. Herman)}\\

\noindent\textbf{Update 2.82}
No progress on this problem has been reported to us.\\

\noindent\textbf{Problem 2.83} Let a function $f$ in $R_d$ have the property that
\[f(z)=\lambda_\alpha z+O(z^2)\hspace{1cm}\text{ as }z\to0,\]
where $\lambda_\alpha=e^{2\pi i\alpha}$ and $\alpha\in\mathbb{R}\setminus\mathbb{Q}$. Assume also that $f$ is linearisable at $0$, and denote by $S$ its Siegel singular disc.
\begin{enumerate}[(a)]
\item Is $\alpha$ necessarily a Brjuno number? In other words, is it true that
\[\sum^\infty_{n=0}(\log q_{n+1})/q_n<+\infty,\]
where $\{p_n/q_n\}^\infty_{n=0}$ are the convergents of the continued fraction expansion of $\alpha$? Also,
what is the situation here when $f$ is a non-linear entire function?
\item Is $f$ necessarily injective on the boundary $\partial S$ of $S$ in $\hat{\mathbb{C}}$?
\item Is it true that $f$ has no periodic points on $\partial S$?
\end{enumerate}
Both $(b)$ and $(c)$ are open even under the additional hypothesis that $f$ has no
critical point on $\partial S$. A positive answer to $(b)$, under this additional hypothesis, would
imply that when $\alpha$ satisfies a diophantine condition (that is, there exists $\beta$, $\gamma$, $\gamma > 0$ and
$\beta\geq2$ such that, for every number $p/q$ in $\mathbb{Q}$, we have $|\alpha-(p/q)|\geq\gamma q^{-\beta})$ then $f$ has a critical point on $\partial S$.\\
For further information, see Herman \cite{MR796014} and Przytycki \cite{Przytycki}. \\
\textit{(J.-C. Yoccoz and M. R. Herman)}\\

\noindent\textbf{Update 2.83}
For part $(a)$ Yoccoz \cite{MR1367353} has proved that the answer is ‘yes’, if $f$ is a polynomial of degree $2$.\\

\noindent\textbf{Problem 2.84} Does there exist a number $\lambda$ of modulus one that is not a root of unity, such that the positive orbit of $-\frac{1}{2}$ under $P_ \lambda(z) =  \lambda(z + z^2)$ is dense in $J(P_ \lambda)$?\\
\textit{(M. R. Herman)}\\

\noindent\textbf{Update 2.84}
No progress on this problem has been reported to us.\\

\noindent\textbf{Problem 2.85} Suppose that $\lambda$ is of modulus one and not a root of unity, let
$P_ \lambda(z) =  \lambda(z + z^2)$ and
\[h_\lambda(z)=z+O(z^2)\]
is the unique formal power series such that $P_ \lambda(h_\lambda(z)) =  h_\lambda(\lambda z)$. Denote by $R(\lambda)$ the radius of convergence of $h_\lambda$.
\begin{enumerate}[(a)]
\item Calculate (or, at least, estimate up to $\pm10^{-10}$) the value of $m = \sup_\lambda R(\lambda)$.
\item Prove that $m$ is realised by a function $h_\lambda$, where $\lambda = e^{2\pi i\alpha}$ and $\alpha$ is a real algebraic number of degree $2$.
\item If $R(\lambda) = 0$, is the following true: for every positive $\varepsilon$, the function $P_\lambda$ has a
repelling periodic cycle included in $\{|z| < \varepsilon\}$? This property is known to hold for a
dense $G_\delta$-set of numbers $\lambda$ of modulus one, see Cremer \cite{Cremer}. 
\end{enumerate}
\textit{(M. R. Herman and J. -C. Yoccoz)}\\

\noindent\textbf{Update 2.85}
No progress on this problem has been reported to us.\\

\noindent\textbf{Problem 2.86} Let the function $f(z)$, $f(z) = \lambda(e^z-1)$ with $|\lambda| = 1$, have a Siegel singular disc $S_\lambda$ that contains zero.
\begin{enumerate}[(a)]
\item Prove that there exists some number $\lambda$, where $|\lambda| = 1$, such that $S_\lambda$ is bounded
in $\mathbb{C}$.
\item If $S_\lambda$ is unbounded in $\mathbb{C}$, does $-\lambda$ belong to $\partial S_\lambda$?
\end{enumerate}
\textit{(M. R. Herman, I. N. Baker and P. J. Rippon)}\\

\noindent\textbf{Update 2.86}
Rempe-Gillen writes that $(a)$ can be solved by a method of Ghys which is mentioned by Douady \cite{MR936853} and found in detail in \mbox{Rempe-Gillen} \cite{RempePHD}. A positive answer to $(b)$ is given by Rempe-Gillen \cite{MR2069014} with a relatively simple proof.\\

\noindent\textbf{Problem 2.87} Does there exist a non-linear entire function $g$ with wandering domain $W$
such that $\bigcup_{n\geq0}g^n(W)$ is bounded in $\mathbb{C}$?\\
It has been conjectured by Lyubich \cite{MYuLyubich1986Tdor} that if $g$ and $W$ exist, then $g^n(W)$ cannot converge as $n\to\infty$ to a fixed point of $g$.\\
\textit{(M. R. Herman; A. Eremenko and M. Yu. Lyubich)}\\

\noindent\textbf{Update 2.87}
The first question is the same as Problem 2.77, and Problem 2.67 is a special case of
this question. All these problems remain unsolved. Bergweiler writes that in principle, it is conceivable that there exists a function satisfying Problem 2.67, but not Problem 2.77 or Problem 2.87. In other words, the limit functions in the wandering domain are all bounded, but the domains are not.\\
Rempe-Gillen writes that the first part of Problem 2.87 remains an open problem, although results in \cite{MR0419759} imply the absence of all wandering domains for some functions. However, the absence of bounded wandering domains can be deduced from a result of Clark, Trejo and van Strien (unpublished). This problem is about whether there exists a transcendental entire function $g$ with a wandering component $W$ of the Fatou set, such that $\bigcup_n g^n(W)$ is a bounded set. Nicks \cite{MR3008885} shows that there does exist a quasi-analytic mapping of the plane, with an essential singularity at infinity, that has this ‘bounded wandering’ property. This problem asks for a non-linear entire function $g$ with wandering domain $W$ such that $\bigcup_n g^n(W)$ is bounded in $\mathbb{C}$.  Under certain additional assumptions, a negative answer is given by Zheng \cite{MR1772425}.\\
The conjecture that iterates in a
wandering domain cannot converge to a fixed point arises from the results
of P\'erez Marco (\cite{MR1360570}, \cite{MR1607557}).\\

\noindent\textbf{Problem 2.88} Let $B$ denote the boundary of the Mandelbrot set (or, equivalently, the
topological bifurcation set of the family $z\mapsto z^2+c$, $c\in\mathbb{C}$.
\begin{enumerate}[(a)]
\item Is $B$ locally-connected?
\item Prove that $B$ has Hausdorff dimension $2$.
\item Does $B$ have Lebesgue measure zero?
\end{enumerate}
For further information, see Douady \cite{MR728980} and Lyubich \cite{MYuLyubich1986Tdor}.\\
\textit{(A. Douady and J. H. Hubbard; N. Sibony; M. Rees; M. R. Herman)}\\

\noindent\textbf{Update 2.88}
\begin{enumerate}[(a)]
\item It is shown by Hubbard \cite{MR1215974} and Lyubich \cite{MR1755441} that $B$ is locally connected in the neighbourhoods of certain points.
\item Shishikura (\cite{MR1285397}, \cite{MR1626737}) proved that the boundary of the Mandelbrot set $B$
has Hausdorff dimension $2$.\\
\end{enumerate}

\noindent\textbf{Problem 2.89} Let the function $f_0$ in $R_d$ have an invariant Herman singular ring $A_f$ of
rotation number $\alpha$, where $\alpha$ satisfies a diophantine condition. Denote by $H_{d,\alpha}$ the class
of all functions $f_1$ in $R_d$ such that $f_1$ can be joined to $f_0$ by a continuous path $f_t$,
$0\leq t\leq1$, in $R_d$, where each $f_t$ has a Herman singular ring $A_{f_t}$ of rotation $\alpha$, and the
annuli $A_f$ vary continuously with $f$ (in the sense of Carath\'eodory).
\begin{enumerate}[(a)]
\item Is $H_{d,\alpha}$ locally closed in $R_d$?
\item Is the boundary of $H_{d,\alpha}$ in its closure in $R_d$ a topological manifold?
\end{enumerate}
Both $(a)$ and $(b)$ are related to the investigation of rational functions with an
invariant Herman singular ring, when the moduli of their invariant rings tend to zero. See Herman
\cite{MR771920}. \\
\textit{(M. R. Herman)}\\

\noindent\textbf{Update 2.89}
No progress on this problem has been reported to us.\\

\noindent\textbf{Problem 2.90} Does there exist a number $\alpha$ in $\mathbb{R}\setminus\mathbb{Q}$ that does not satisfy a diophantine condition, such that every $\mathbb{R}$-analytic orientation-preserving diffeomorphism of the
circle with rotation number $\alpha$ is $\mathbb{R}$-analytically conjugated to a rotation?\\
For related questions, see Douady \cite{MR936853}. If $\alpha$ satisfies a diophantine condition, the global
analytical conjugacy theorem has been proved. See Herman \cite{MR538680} and Yoccoz \cite{MR777374}. \\
\textit{(M. R. Herman)}\\

\noindent\textbf{Update 2.90} 
No progress on this problem has been reported to us.\\

\section*{New Problems}

\noindent To appear.

% ==================== CHAPTER 3 ==========================================

\chapter{Subharmonic and Harmonic Functions}

\underline{Preface by S. Gardiner (to appear)}\\

\noindent\textbf{Notation} A function $u(z)$ in a domain $D$ of the plane, is said to be \textit{subharmonic} in $D$ if:
\begin{enumerate}[(a)]
\item $u(z)$ is upper semi-continuous in $D$,
\item $-\infty\leq u(z_0)<+\infty$, and $u(z)\not\equiv-\infty$ in $D$,
\item for every $z_0$ in $D$, and all sufficiently small $r$ (depending on $z_0$), we have
\[u(z_0)\leq\frac{1}{2\pi}\int^{2\pi}_0 u(z_0+re^{i\theta})\,d\theta.\]
\end{enumerate}
In space of higher dimensions, subharmonic functions are defined analogously. If $u(z)$ and $-u(z)$ are subharmonic, then $u(z)$ is \textit{harmonic}.
If $f(z)$ is analytic in a domain $D$, and $f(z)\not\equiv0$, then 
\[u(z)=\log |f(z)|\]
is subharmonic in $D$, and many properties of the modulus of analytic functions suggest corresponding properties of subharmonic functions in the plane.\\ 

\noindent Sometimes, analogous properties of subharmonic functions in space also hold.
We recall, from the introduction to Chapter 2, Iversen's theorem \cite{Iversen}, that if $f(z)$ is a non-constant entire function, then there exists a path $\Gamma$, such that $f(z)\to\infty$ as $z\to\infty$ along $\Gamma$. 
The analogue, that a subharmonic function $u(z)$ in the plane, which is not constant, tends to $+\infty$ along a path $\Gamma$ was proved by Talpur \cite{zbMATH03493031}. \\

\noindent Problems 2.6 and 2.7 have immediate analogues for subharmonic functions, and an analogue of Problem 2.8 can be formulated without too much difficulty. If $f(z)$ is analytic and non-zero, then $u(z)=\log|f(z)|$ is harmonic. Thus, Problem 2.9 in its first formulation reduces to a problem on harmonic functions $u(z)$. It should be said that Huber \cite{MR0094452} proved his result (\ref{2.2}) for subharmonic functions, (see Problem 2.10), that for every positive $\lambda$ there exists a path $C_\lambda$ such that
\[\int_{C_\lambda}|f(z)|^{-\lambda}\,|dz|<\infty.\]
That is, with $|f(z)|$ replaced by $\exp(u(z))$ in Problem 2.10, where $u(z)$ is a subharmonic function such that
\[\frac{B(r,u)}{\log r}\to+\infty,\]
where, throughout this section,
\[A(r,u)=\inf_{|z|=r}u(z),\hspace{1cm}B(r,u)=\sup_{|z|=r}u(z).\]
Unless the contrary is explicitly stated, we shall assume that the functions $u(z)$ are defined in a plane domain.\\

\noindent\textbf{Problem 3.1} If $u(z)$ is harmonic in the plane, and not a polynomial, does there exist a path $\Gamma_n$ for every positive integer $n$, such that 
\begin{equation}
\label{3.1}
\frac{u(z)}{|z|^n}\to+\infty
\end{equation}
as $z\to\infty$ along $\Gamma_n$. Does there exist a path $\Gamma_\infty$, such that (\ref{3.1}) holds for every fixed $n$, as $z\to\infty$ along $\Gamma_\infty$?\\
We note that we can apply the result (see Problem 2.6) of Boas, that for every transcendental entire function there exists a path $\Gamma_\infty$ such that for every $n$, $\big|\frac{f(z)}{z^n}\big|\to\infty$ as $z\to\infty$ along $\Gamma$;  to $f(z)=e^{u+iv}$, where $v$ is the harmonic conjugate of $u$, but this only yields
\[\frac{u(z)}{\log|z|}\to+\infty.\]

\noindent\textbf{Update 3.1}
We can ask the analogue of Iversen's theorem (see \cite{Iversen}) for a harmonic function
$u(z)$ (see Problem 2.6). Boas's theorem gives a path on which
\[\frac{u(z)}{\log|z|}\to+\infty,\]
where $u(z)$ is harmonic and not constant in the plane. Barth, Brannan and Hayman \cite{MR507611}
have shown that there exists a path for every $\alpha$, with $\alpha <\frac{1}{2}$ such that
\[\frac{u(z)}{|z|^\alpha}\to+\infty,\]
and a little more. On the other hand, there need not exist such a path if $\alpha\geq\frac{1}{2}$.\\
 
\noindent\textbf{Problem 3.2} If $u(x)$ is harmonic and not constant in space of $3$ or more dimensions, is it true that there exists a path $\Gamma$ such that $u(x)\to+\infty$ as $x\to\infty$ along $\Gamma$?\\
The corresponding result for subharmonic functions is certainly false, since if 
\[r=\Big(\sum^n_{\nu=1}x_\nu^2\Big)^\frac{1}{2}\]
is the distance of $x= (x_1, x_2,\ldots,x_n)$ from the origin, then $u(x)=\max(-1,-r^{2-n})$ is subharmonic and bounded in space of $n$ dimensions, when $n>2$. On the other hand, a bounded harmonic function in space is constant.\\

\noindent\textbf{Update 3.2}
This problem has been largely settled. If $u(x)$ is subharmonic
and has finite least upper bound $M$, then Hayman (\cite[Chp. 4]{MR0460672}) showed that $u(x)\to M$
as $x\to\infty$ along almost all rays through the origin. If $u(x)$ is not bounded
above, and in particular, if $u(x)$ is harmonic, there always exists a path
$\Gamma$ such that
\begin{equation}
\label{B4}
u(x)\to+\infty\hspace{1cm}\text{ as }x\to\infty \text{ along }\Gamma.
\end{equation}
This is a result of Fuglede \cite{MR0379870} which completes earlier results of
Talpur and Hayman \cite{Talpur}. Fuglede uses a deep theorem about Brownian motion
by Nguyen-Xuan-Loc and Watanabe \cite{Watanabe}. The path $\Gamma$ is locally a Brownian
motion and so highly irregular. It is an open question whether a smooth
path or a polygonal path exists satisfying (\ref{B4}). For continuous subharmonic
functions, and in particular for harmonic functions, such paths certainly
exist. Carleson \cite{MR0470233} has proved the existence of a
polygonal path for general subharmonic functions.\\

\noindent\textbf{Problem 3.3} Suppose that $u(z)$ is subharmonic and $u(z)<0$ in the half-plane $|\theta|<\pi/2$, where $z=re^{i\theta}$. Suppose also that 
\[A(r)=\inf_{|\theta|<\pi/2}u(re^{i\theta})\leq-K, \hspace{1cm}0<r<\infty.\]
Is it true then that
\[u(r)\leq-\frac{1}{2}K,\hspace{1cm} 0<r<\infty\,?\]
The result $u(r)\leq-K/3$ is true, and is due to Hall \cite{Hall}.\\

\noindent\textbf{Update 3.3}
The answer is `no'. Hayman (\cite[Chp. 7]{MR0460672}) has constructed an example of a function
satisfying Hall's conditions, and such that $u(r) > -\frac{1}{2}K$ on the whole
positive axis. What is the best constant $\alpha$ such that $u(r)<-\alpha K$ on the whole positive axis?\\

\noindent\textbf{Problem 3.4} Consider the class of functions subharmonic in the unit disc $\mathbb{D}$, and satisfying $u(z)\leq0$ there. Suppose also that $A(r,u)\leq-1$, for $r$ lying on a set $E$ consisting of a finite number of straight line segments. Then it is known (see e.g. Nevanlinna \cite[p.100]{MR0344426}) that $B(r,u)$ is maximal, when $u(z)=u_E(z)$, where $u_E(z)$ is harmonic in $\mathbb{D}$, except on a set $E$ of the positive real axis, and $u(z)$ assumes boundary values $0$ on $|z|=1$, and $-1$ on $E$. 
This is the solution to the so-called Carleman-Milloux problem \cite{MR0344426}. For $0<r<1$, $0<K<1$, let $C(r,K)$ be the set of all $\theta$ such that $u(re^{i\theta})<-K$. \\
Is it true that $C(r,K)$ has minimal length only if $u(z)=u_E(z)$? The special case where $E$ consists of the whole interval $[0,1]$ has particular interest.\\
\textit{(T. K\"ovari)}\\

\noindent\textbf{Update 3.4}
No progress on this problem has been reported to us.\\
  
\noindent\textbf{Problem 3.5} Suppose that $u(z)$ is positive and subharmonic in $\mathbb{D}$, and that there exists a series of arcs $\gamma_n$ tending to the arc $\alpha\leq\theta\leq\beta$ of $|z|=1$, such that
\begin{equation}
\label{3.2}
u(z)\leq M,\hspace{1cm} z\text{ on }\gamma_n,\hspace{1cm} n=1,2.
\end{equation}
If in addition, 
\begin{equation}
\label{3.3}
\int^1_0(1-r)u(re^{i\theta})\,d\theta<+\infty,
\end{equation}
for a set $E$ on $\theta$, which is dense in the interval $(\alpha, \beta)$, then Maclane \cite{MR0148923} proved that $u(re^{i\theta})$ is uniformly bounded for 
\[\alpha+\delta\leq\theta\leq\beta-\delta, \hspace{1cm}0\leq r<1,\]
and any fixed positive $\delta$. These conclusions thus hold in particular, if
\begin{equation}
\label{3.4}
\int^1_0(1-r)B(r,u)\,dr<+\infty.
\end{equation}
Can the growth conditions (\ref{3.3}) and (\ref{3.4}) be weakened without weakening the conclusions?\\

\noindent\textbf{Update 3.5}
Hornblower \cite{MR0296305} and \cite{MR0364642} has shown that (\ref{3.4}) can indeed be replaced
by the much weaker condition
\[\int^1_0\log^+B(r)\,dr<\infty,\]
and that this is more or less best possible. However, it is still an open
question as to whether condition (\ref{3.3}) is sharp.\\

\noindent\textbf{Problem 3.6} It follows from a result of Wolf \cite{MR0000288}, that if 
\[u(re^{i\theta})\leq f(\theta),\hspace{1cm} 0<r<+\infty,\]
where
\[\int^{2\pi}_0\log^+f(\theta)\,d\theta<+\infty,\]
then $u(z)$ is bounded above, and so is constant. What is the $3$-dimensional analogue of this result?\\

\noindent\textbf{Update 3.6} We note that the original statement of Problem 3.6 had $\log^+\log^+f(\theta)$ instead of $\log^+f(\theta)$. Wolf proved that if $u$ is subharmonic in the plane, and satisfies
\[u(re^{i\theta})\leq f(\theta),\hspace{1cm} 0<r<\infty,\]
where
\[\int^{2\pi}_0\log^+f(\theta)\,d\theta<\infty,\] 
then $u$ is bounded above, and so constant. Yoshida has pointed out that the higher-dimensional analogue of this result had been obtained by Carleman \cite{MR1555232}.\\

\noindent\textbf{Problem 3.7} Problem 1.17 can be reformulated for subharmonic functions, if we replace $\log M(r,f)$ by a general subharmonic function $u(z)$. The same positive theorems hold, and the same conclusions are conjectured in the general case.\\

\noindent\textbf{Update 3.7}
The result, generalizing Govorov's theorem to subharmonic
functions in space, has been obtained by Dahlberg \cite{MR0322184}.\\

\noindent\textbf{Problem 3.8} If $\omega(\zeta)$ is continuous and bounded in $\mathbb{D}$, and each point of $\zeta$ of $\mathbb{D}$ is the centre of at least one circle $C_\zeta$ lying in $\mathbb{D}$, such that $\omega(\zeta)$ is equal to the average of the values of $\omega(z)$ on $C_\zeta$, is it true that $\omega(z)$ is harmonic in $\mathbb{D}$? The result is true if $\omega(z)$ is continuous in $|z|\leq1$, and is certainly false if this condition is removed. The corresponding problem with $C_\zeta$ replaced by its interior $D_\zeta$ is also open.\\
\textit{(J. E. Littlewood)}\\

\noindent\textbf{Update 3.8}
This has been completely solved by Hansen and Nadirashvili \cite{doi:10.1112/jlms/50.2.349}. Eremenko writes that Problem 3.8 has two interpretations, depending on the meaning of the word `circle'. If a circle means a disc, then the answer is positive; but if the circle means circumference, then the answer is negative.\\

\noindent\textbf{Problem 3.9} If $D$ is a convex domain in space of $3$ or more dimensions, can we assert any inequalities for the Green's function $g(P,Q)$ of $D$ which generalise the results of $2$ dimensions, that follow from the classical inequalities for convex univalent functions?\\
Gabriel \cite{MR0072959} proved that the level surfaces $G(P,Q)=\lambda>0$ are convex, but the proof is long. It would be interesting to find a simpler proof, and also to obtain definite inequalities for the curvatures. It may be conjectured that half-space gives the extreme case.\\
\textit{(G. E. H. Reuter)}\\

\noindent\textbf{Update 3.9}
No progress on this problem has been reported to us.\\

\noindent In general, there are many problems for harmonic functions in $3$ dimensions, where the proofs of the corresponding two-dimensional results depend on conjugacy arguments, and analytic functions, and so do not readily extend. An example is the following:\\

\noindent\textbf{Problem 3.10} Suppose that $u(X)$ is harmonic on the unit ball $|X|<1$, and remains continuous with partial derivatives of all orders on $|X|=1$, where $X$ is a point $(x_1, x_2, x_3)$ in space, and 
\[|X|^2=x^2_1+x^2_2+x^2_3.\]
Suppose further that there is a set $E$ of positive area on $|X|=1$, such that both $u$ and its normal derivative vanish on $E$. Is it true that $u\equiv0$?\\
Here the two-dimensional analogue is almost trivial, since if $u$ is harmonic in $\mathbb{D}$, and $u$ and its partial derivatives remain continuous on the unit cirlce $\mathbb{T}$, we may consider 
\[f(z)=z\Big(\frac{\partial u}{\partial x}-i\frac{\partial u}{\partial y}\Big), \hspace{1cm}\text{ where }z=x+iy.\]
If $u(z)$ and its normal derivatives both vanish on a set $E$, then $f(z)$ vanishes at all limit points of $E$. Now the Poisson-Jensen formula shows at once that $\log|f(z)|=-\infty,f(z)=0,$ identically in $\mathbb{D}$, provided that the closure of $E$ has positive $1$-dimensional measure.\\
\textit{(L. Bers)}\\

\noindent\textbf{Update 3.10}
The problem was solved by Wolff \cite{zbMATH00799747} who constructed a counterexample.\\

\noindent\textbf{Problem 3.11} If $u(x)$ is a homogeneous harmonic polynomial of degree $n$ in
$\mathbb{R}^m$, what are the upper and lower bounds of
\[-\frac{A(r,u)}{B(r,u)},\]
where $A(r, u) = \inf_{|x|=r} u(x), B(r, u) = \sup_{|x|=r} u(x)$?
If $n$ is odd, it is evident that $A(r, u) = -B(r, u)$, but if
$u(x) = x_1^2 + x^2_2 - 2x^2_3$ in $\mathbb{R}^3$, then $B(r, u) = r^2$, $A(r, u) = -2r^2$. For
transcendental harmonic functions, such that $u(0) = 0$, we can prove that
\[-A(r, u) \leq\frac{(R+r)R^{m-2}}{(R-r)^{m-1}}B(r,u), \hspace{1cm}0 < r < R,\]
by Poisson's formula and this leads to
\begin{equation}
\label{B3.5}
-A(r, u) <B(r) (\log B(r) )^{m - 1 + \varepsilon}
\end{equation}
outside a set of $r$ of finite logarithmic measure. However, (\ref{B3.5}) is
unlikely to be sharp. We note that for $m = 2$, it follows from a classical
result of Wiman that, for any harmonic function $u$,
\[A(r)\sim-B(r)\]
as $r\to\infty$ outside a set of finite logarithmic measure.\\

\noindent\textbf{Update 3.11}
Armitage \cite{MR540060} has solved this problem by obtaining the least upper bound of 
\[\frac{\sup_{|x|=1}u(x)}{\inf_{|x|=1}u(x)},\]
where $u(x)$ is a homogeneous harmonic polynomial of degree $n$ in $\mathbb{R}^m$, $m\geq3$. He shows that the ratio is bounded by a constant depending only on $m$.\\

\noindent\textbf{Problem 3.12}  Consider a domain of infinite connectivity in $\mathbb{R}^3$ whose complement
$E$ lies in the plane $P : x_3 =0$. Suppose further than any disc of positive
radius $R$ in $P$ contains a subset of $E$ having area at least $\varepsilon$,
where $\varepsilon, R$ are fixed positive constants. If $u$ is positive and harmonic
in $D$, continuous in $\mathbb{R}^3$ and zero on $E$, is it true that
\[u = cx_3 + \phi(x), \hspace{1cm}x_3 > 0,\]
where $c$ is a constant and $\phi(x)$ is uniformly bounded? One can also ask
the analogue of this result for $\mathbb{R}^m$ when $m>3$. It is true in $\mathbb{R}^2$ (but
Kjellberg cannot remember who proved it).\\
\textit{(B. Kjellberg)}\\

\noindent\textbf{Update 3.12}
Kjellberg asked for conditions on a closed set in $\mathbb{R}^k$ so that the complement in  $\mathbb{R}^{k+1}$ has one of two Martin boundary points. A fairly complete answer has been given by Benedicks \cite{MR608327}.\\

\noindent\textbf{Problem 3.13} Let $u(x)$ be subharmonic in $\mathbb{R}^m$. One can define the quantities
$n(r, 0), N(r, 0), T(r)$ as in Nevanlinna theory in the plane, taking the
analogue of the case $u(z) = \log|f(z) |$, where $f(z)$ is an entire function (see, for example Hayman \cite{MR0111839}). Define
\[\delta(u) = 1 - \limsup_{r\to\infty}\frac{N(r,u)}{T(r)}.\]
If the order $\rho$ of $u$ is less than $1$, it is possible to obtain the sharp upper
bound for $\delta(u)$ in terms of $\rho$ and $m$. The bound is attained when $u(x)$
has all its mass on a ray (see Hayman and Kennedy \cite[Chp. IV]{MR0460672}). One can ask the corresponding
question for $\rho > 1$.\\
One can ask whether a lower bound $A(\rho)$ can be obtained for
$\delta(u)$ if $\rho > 1$ and all the mass of $u(x)$ lies on a ray, or more generally
on some suitable lower dimensional subspace $S$ of $\mathbb{R}^m$ and $\rho>\rho_0(S)$.
For fixed $S$ we may conjecture by analogy with the case $m = 2$, that
$A(\rho)\to1$ as $\rho\to\infty$. This is proved by Hellerstein and Shea \cite{MR0393486} in the
case $m = 2$.\\
\textit{(D. Shea)}\\

\noindent\textbf{Update 3.13}
Rao and Shea \cite{MR0444974} have extended from $2$ to $m$ dimensions the work of
Hellerstein and Shea \cite{MR0393486} on the deficiency of a subharmonic function $u(x)$ with all its
mass concentrated on a ray. More precisely write
\[N(r,u)= \int_{|x|=r} u(x)\,d\sigma(x),\hspace{1cm} m_2(r) = \int_{|x|=r} u(x)^2\,d\sigma(x),\]
where $\sigma$ denotes surface measure on $|x|=r$, normalised to have total value one. Then
for $m = 2, 3, 4$ they obtain the sharp lower bound for
\[\limsup_{r\to\infty}\frac{m_2(r)}{N(r)}\]
as a function of the order of $u$. It turns out again that the limit is attained when all the
Riesz mass lies on a ray. This leads to nice bounds for the deficiencies of subharmonic
functions in space in terms of the order.\\

\noindent\textbf{Problem 3.14}  Let there be given an integrable function $F$ on $\mathbb{T}$
and a point $z_0$ in $\mathbb{D}$. The problem is to maximise $u(z_0)$,
where $u$ runs through all functions which are subharmonic in $\mathbb{D}$,
equal to $F$ on $\mathbb{T}$, and which satisfy
\[\inf u(re^{i\theta})\leq 0, \hspace{1cm}0 < r < 1.\]
\textit{(A. Baernstein)}\\

\noindent\textbf{Update 3.14}
No progress on this problem has been reported to us.\\

\noindent\textbf{Problem 3.15} Let $D$ be a doubly-connected domain with boundary curves $\alpha$ and
$\beta$ and let $z_0, z_1$ be points of $D$. Let $A, B$ be given real numbers.
The problem is to maximise $u(z_0)$, where $u$ runs through all functions
which are subharmonic in $D$, take the values $A$ and $B$ on $\alpha$ and $\beta$
respectively and are non-positive on some curve connecting $z_1$ to $\alpha$.\\
\textit{(A. Baernstein)}\\

\noindent\textbf{Update 3.15}
No progress on this problem has been reported to us.\\

\noindent\textbf{Problem 3.16} A compact set $E$ in $\mathbb{R}^n$, $n\geq3$ is said to be \textit{thin} at $P_0$ if
\begin{equation}
\label{B3.1}
\int^1_0\frac{c(P_0,r)}{r^{n-1}}\,dr<\infty,
\end{equation}
where $c(P_0,r)= \text{cap}\big[E\cap\{P:|P-P_0| \leq r\}\big]$; this is the integrated form of the Wiener
criterion. It follows from Kellogg's theorem \cite{MR1501602} that the points of $E$ where (\ref{B3.1}) holds,
form a polar set. Can one give a direct proof of this fact, which shows perhaps that
(\ref{B3.1}) is best possible? \\
\textit{(P. J. Rippon)}\\

\noindent\textbf{Update 3.16}
No progress on this problem has been reported to us.\\

\noindent\textbf{Problem 3.17}  Let $D, D'$ be Lipschitz domains in $\mathbb{R}^n$, $n\geq 3$ with $D' \subset D$ and $\partial D'\cap\partial D$ lying compactly in the interior of a set $\Gamma$ in $\partial D$; for any fixed $P_0$ in $D'$, let $H(P_0)$
denote the family of positive harmonic functions $h$ on $D$ that vanish continuously
on $\Gamma$ and satisfy $h(P_0) = 1$. Is there a constant $C$ such that, for all $h_1, h_2$ in $H(P_0)$, we
have $h_1(P)\leq Ch_2(P)$ for all $P$ in $D'$?\\
The problem was first considered by Kemper \cite{MR0293114}, but the proof he gives
contains an error. \\
\textit{(P. J. Rippon)}\\

\noindent\textbf{Update 3.17}
Suppose $D$ is a Lipschitz domain, $P_0$ is a point in $D$, $E$ is a relatively open set on $\partial D$, and $S$ is a subdomain of $D$ satisfying $\partial S\cap\partial D \subseteq E$. Dahlberg \cite{MR0466593} and Wu \cite{MR513884} proved independently that there is a constant $C$ such that whenever $u_1$ and $u_2$ are two positive harmonic functions in $D$ vanishing on $E$, and $u_1(P_0)=u_2(P_0)$, then $u_1(P)\leq CU_2(P)$ for all $P\in S$. Ancona \cite{MR513885} proved the above result for solutions of certain elliptic operators instead of harmonic functions.\\

\noindent\textbf{Problem 3.18} It is known that the set $E$ of least capacity $C$ and given volume is a ball.
If $E$ displays some measure of asymmetry (for instance, if every ball with the same
volume as $E$ in space contains a minimum proportion $\delta$ in the complement of $E$), can
one obtain a lower bound for the capacity of $E$ which exceeds $C$ by some positive
function of $\delta$? \\
\textit{(E. Fraenkel; communicated by W. K. Hayman)}\\

\noindent\textbf{Update 3.18}
Hall, Hayman and Weitsman \cite{zbMATH00020033} obtain for general sets if $n=2$, and for convex sets if $n\geq3$, 
\[\text{cap } E\geq C(1+k\delta^{n+1}),\] 
where $k$ is a constant depending on the dimension $n$ .\\

\noindent\textbf{Problem 3.19} Let $C_0$ be a tangential path in $\mathbb{D}$ which ends at $z=1$, and let $C_\theta$ be any rotation of $C_0$. Littlewood \cite{MR1574417} showed that there exists a function $u(z)$, harmonic and satisfying $0<u(z)<1$ in $\mathbb{D}$, such that
\[\lim_{|z|\to1,\, z\in C_\theta}u(z)\]
does not exist for almost all $\theta$, $0\leq\theta\leq2\pi$. Surprisingly, it seems to be unknown whether there exists a $v(z)$, positive and harmonic in $\mathbb{D}$, such that
\[\lim_{|z|\to1,\, z\in C_\theta}v(z)\]
does not exist for all $\theta$, $0\leq\theta\leq2\pi$.
The corresponding result is known for bounded analytic functions in $\mathbb{D}$. See, for example, Collingwood and Lohwater, \cite[Chp. 2]{MR0231999}. \\
\textit{(K. H. Barth)}\\

\noindent\textbf{Update 3.19}
Aikawa \cite{MR990410} proves a theorem which automatically answer’s Barth's question in the affirmative. Further, he notes that the question can be solved more easily, and does in fact find a positive unbounded harmonic function with the required property.\\

\noindent\textbf{Problem 3.20} Suppose that you have a continuous real function $u(x)$ on $\mathbb{R}^n$, and you want to know whether a homeomorphism $\phi:\mathbb{R}^n\to\mathbb{R}^n$ and a harmonic function $v$ on $\mathbb{R}^n$ exist, such that
\[v(x)=u(\phi(x)).\]
Is it necessary and sufficient that there should exist mappings $\mu_1, \mu_2, \ldots, \mu_n$ so that
\[U=(u, \mu_1, \mu_2, \ldots, \mu_n)\]
is a light open mapping of $\mathbb{R}^n$ into $\mathbb{R}^n$? The case $n=2$ is a result of Stoilow which is in \cite{MR0165476}, for example.\\
\textit{(L. A. Rubel, communicated by D. A. Brannan)}\\

\noindent\textbf{Update 3.20} No progress on this problem has been reported to us.\\

\noindent\textbf{Problem 3.21}  Let $\alpha$ be a continuum in the closure of the unit disc $\mathbb{D}$, and let
$\omega(z) = \omega(z; \mathbb{D}; \alpha)$ be the harmonic measure at $z$ of $\alpha$ with respect to $\mathbb{D}$. Is it true that
$\omega(0)\geq\frac{1}{\pi}\arcsin\frac{1}{2}d$, where $d$ is the diameter of $\alpha$?\\
\textit{(B. Rodin)}\\

\noindent\textbf{Update 3.21}
 FitzGerald, Rodin and Warschawski \cite{MR768733} show that the answer is positive. \\

\noindent\textbf{Problem 3.22}  Let $D$ be a domain containing the origin whose `outer boundary' is
$\mathbb{T}$ and whose `inner boundary' is a closed set $E$ in  $\mathbb{D}$. If every radius of
the unit disc meets $E$, determine the supremum of the harmonic measure at $0$ of
$\mathbb{T}$ with respect to $D$.\\
\textit{(W. H. J. Fuchs)}\\

\noindent\textbf{Update 3.22}
A special case of this is listed as an open problem by Betsakos in \cite{MR1877323}. This problem is to find a continuum in the closed unit disc which meets every radius and whose harmonic measure at the origin is minimal. Jenkins \cite{MR922301} characterises the class of continua in which such an extremal continuum must occur; and in \cite{MR948195}, the essentially unique extremal continuum is determined. A simplified derivation of this solution using the method of the extremal metric is given by Jenkins \cite{MR1162090}. Fuchs asked what the optimal constant $c$ is. Marshall and Sundberg \cite{MR948195} show how close $c$ is to $1$ in a special case.\\

\noindent\textbf{Problem 3.23} Determine whether or not there exists a function $g(r)$, defined for $r\geq0$,
with $g(r)\to0$ as $r\to\infty$, such that the following holds: if $u$ is any Green potential in
$\mathbb{D}$ satisfying $u(0) = 1$, then for every non-negative $r$ the set
\[E_r = \{z : z \in \mathbb{D}, u(z) > r\}\]
can be covered by a family of discs $\{D(a_k; r_k)\}$ (with centres $a_k$ and radii $r_k$),
depending on $r$, such that $\sum_k r_k\leq g(r)$. One can ask the same question with `Green
potential' replaced by `positive harmonic function'.\\
Results of this type are known for ordinary logarithmic potentials (c.f. Cartan's
lemma, see for example \cite{MR0180696}) and the Riesz potentials (in higher dimensions).\\
\textit{(R. Zeinstra)}\\

\noindent\textbf{Update 3.23}
Eiderman points out that an old result of Govorov \cite{MR0269865} shows this to be the case
with
\[g(r)= 
\begin{cases}
1 &\text{if }0 \leq r \leq 28,\\
\frac{27}{r-1}&\text{if }r>28.\\
\end{cases}\]

\noindent\textbf{Problem 3.24} For which positive $p$ does there exist a function $u$, $u\not\equiv0$ harmonic on $\mathbb{R}^3$ and
vanishing on the cone $x^2_1+x^2_2=px^2_3$?\\
\textit{(H. S. Shapiro)}\\

\noindent\textbf{Update 3.24}
A complete characterisation of all quadratic harmonic divisors in $\mathbb{R}^n$, for arbitrary $n$, is given by Agranovsky and Krasnov \cite{MR1799672}.\\

\noindent\textbf{Problem 3.25} Is there a harmonic polynomial $P(x_1, x_2, x_3)$, $P\not\equiv 0$ that is divisible by $x^4_1+x^4_2+x^4_3$?\\
\textit{(H. S. Shapiro)}\\

\noindent\textbf{Update 3.25}
The answer is `no' since by a theorem of Brelot and
Choquet \cite{MR0069968}, any divisor of a harmonic polynomial in $\mathbb{R}^n$ assumes positive and negative
values.\\

\noindent\textbf{Problem 3.26} Given $n, n\geq4$, find a continuous function $f$ on $(0,1)$ such that the
following statement is true: if $u$ is a subharmonic function in the unit ball $B$ of $\mathbb{R}^n$
with $u(0) > 0$ and $0\leq u < 1$ in $B$, then there exists a path $\gamma$ from the origin to $\partial B$
with $u > 0$ on $\gamma$ and
\[\text{length of } \gamma\leq f(u(0)).\]
Such an $f$ exists when $n = 2$ and when $n = 3$, see David and Lewis \cite{MR735222}.
In the particular case $n = 2$, it has been shown by Lewis, Weitsman and Rossi \cite{MR735882} that one can take
\[f(t) = c_1t^{-c_2}, \hspace{1cm} 0 < t < 1,\]
where $c_1, c_2$ are absolute constants. What is the smallest exponent $c_2$ for which
such an $f$ exists (for the case $n = 2$)?\\
\textit{(J. Lewis)}\\

\noindent\textbf{Update 3.26}
For $n = 2$, Wu \cite{MR825926} has obtained such a result with
\[f(t)=C\Big(1+\log\frac{1}{t}\Big).\]
Lewis has reported that Ancona pointed out to him that the same $f$ will do for general
$n$. This follows from a result of Brelot and Choquet \cite{MR0062883}.\\

\noindent\textbf{Problem 3.27} Let $D$ be an unbounded domain in $\mathbb{R}^n$, $n\geq2$. Is there a positive
continuous function $\varepsilon(|x|)$ such that, if $u$ is harmonic in $D$ and $|u(x)| < \varepsilon(|x|)$, then
$u \equiv 0$?\\
For $n = 2$, the answer is `yes'. The answer is also `yes' if we restrict our attention to
positive harmonic functions. For fine domains and finely harmonic functions, it
follows from an example of Lyons \cite{MR576642}, \cite{MR0576643} that the answer is `no'.\\
\textit{(P. M. Gauthier and W. Hengartner)}\\

\noindent\textbf{Update 3.27}
A positive answer has been given by Armitage, Bagby and Gauthier \cite{MR813738}. See also Armitage and Goldstein \cite{MR1274778} and Armitage \cite{MR1277091}.\\

\noindent\textbf{Problem 3.28}  Determine all domains $\Omega$ in $\mathbb{R}^n$, $n\geq2$, satisfying the identity
$\int_\Omega h(x)\,dx = 0$ for every function $h$ harmonic and integrable on $\Omega$.\\
In the case $n = 2$, the answer is given by Sakai in \cite{MR659788}.\\
\textit{(M. Sakai)}\\

\noindent\textbf{Update 3.28}
No progress on this problem has been reported to us.\\

\noindent\textbf{Problem 3.29}  It is known that the Newtonian potential of a uniform mass distribution
spread over an ellipsoid $K$ in $\mathbb{R}^n$, $n\geq2$ is a quadratic function of the coordinates of
$x = (x_1,\ldots,x_n)$ for $x\in K$.\\
Nikliborc (\cite{MR1545320}, \cite{Nikliborc}) and Dive (\cite{Dive1}, \cite{Dive2}) independently proved that for $n = 2$ and
$n = 3$, the ellipsoid is the only body with this property. Prove this converse assertion
for $n > 3$ (preferably by a new method, since Nikliborc and Dive both use methods
involving highly non-trivial calculations).\\
\textit{(H. S. Shapiro)}\\

\noindent\textbf{Update 3.29}
No progress on this problem has been reported to us.\\

\noindent\textbf{Problem 3.30} Let $K(z, z')$ denote the kernel of the double layer potential occurring in
Fredholm's theory where $z, z' \in \Gamma$, $\Gamma$ being a smooth Jordan curve. (Recall that $K(z,z')=\frac{\cos \phi}{|z-z'|}$, where $\phi$ is the angle between the inward normal to $\Gamma$ at $z$ and the
line $(z,z')$.)\\
When $\Gamma$ is a circle, the function $z \mapsto K(z, z')$ is a constant (i.e. the same for each
choice of $z'$); and consequently the integral operator
\[T_\Gamma:f \mapsto \int_\Gamma f(z)K(z,z')\,ds_z\]
is of rank one (as an operator from $C(\Gamma)$ to $C(\Gamma)$). Are there any other $\Gamma$ for which
the rank of $T_\Gamma$ is finite?\\
\textit{(H. S. Shapiro)}\\

\noindent\textbf{Update 3.30}
No progress on this problem has been reported to us.\\

\noindent\textbf{Problem 3.31} Let $D$ be an unbounded domain in $\mathbb{R}^n$, $n\geq 2$. Points in $\mathbb{R}^n$ will be denoted by $x = (x_1,x_2,\ldots,x_n)$, and $|x|$ will denote the Euclidean norm of $x$. The following
result was given by Ess\'en \cite{MR1050784}.
\begin{theorem}
Assume that the least harmonic majorant $\Psi$ of $|x_1|$ in $D$ is such that
$\Psi(x)=O(|x|)$ as $x\to\infty$ in $D$. If $|x|$ has a harmonic majorant in $D$, then $|x_1|\log^+|x_1|$ has
a harmonic majorant in $D$.
\end{theorem}
In the plane, $|x|^p$ (if $p > 0$) has a harmonic majorant in $D$, if and only if
$F\in H^p$, where $F:\{|x| < 1\} \to D$ is a universal covering map with $F(0) = 0$. There are
similar statements for $|x_1|\log^+|x_1|$ and $\text{Re}\, F\in L\log L$. Thus, when $n = 2$, Theorem 1
is closely related to the following.
\begin{theorem}
Suppose that $n = 2$, and that $F\in H^1(\mathbb{D})$ (where $\mathbb{D}$ is the unit disc). Then
$\text{Re}\, F\in L\log L$ if and only if
\begin{equation}
\label{Jstar}
\int^\infty_{-\infty}N(1,iv,F)\log^+|v|\,dv<\infty,
\end{equation}
where $N(1,\omega,F)$ is the Nevanlinna counting function (see \cite[Theorem 1]{MR812321}).
\end{theorem}
Theorem 2 is an extension of a well-known result of Zygmund. Note also that, in
the case $n = 2$, there are functions $F$ in $H^1(\mathbb{D})$ such that $\text{Re}\, F\notin L\log L$, see \cite[Section
6]{MR812321}.
\begin{enumerate}[(a)]
\item In the case $n=2$, what is the relation between the condition on $\Psi$ in
Theorem 1 and condition (\ref{Jstar}) in Theorem 2?
\item Suppose now that $n\geq2$. The assumption on $\Psi$ in Theorem 1 was introduced
in the proof for purely technical reasons. Is this the correct condition needed in
Theorem 1 (we assume always that $|x|$ has a harmonic majorant in $D$)? Does there
exist a domain $D$ such that $|x|$ has a harmonic majorant in $D$ while $|x_1|\log^+|x_1|$ does
not have a harmonic majorant in $D$? 
\end{enumerate}
\textit{(M. Ess\'en)}\\

\noindent\textbf{Update 3.31}
No progress on this problem has been reported to us.\\

\noindent\textbf{Problem 3.32} Let $\Omega$ be an open ball in $\mathbb{R}^n$, $n\geq 2$. It is shown by Armitage \cite{Armitage} that $V\in L^p(\Omega)$ for any positive superharmonic function $V$ on $\Omega$ and any $p$ in $(0,n/(n-1))$. Now suppose that
$\Omega$ is a bounded Lipschitz domain in $\mathbb{R}^n$ for which the interior cones have half-angle
at least $\alpha$. For what values of $p$ do we have $V\in L^p(\Omega)$ for every positive superharmonic
function $V$ on $\Omega$? \\
\textit{(S. J. Gardiner)}\\

\noindent\textbf{Update 3.32}
Aikawa \cite{MR1169019} shows that this is the case, for 
\[0<p<\min\{n/(n+\alpha-2),1/(\alpha-1)\}.\]\\

\noindent\textbf{Problem 3.33}  Let $\Omega$ be a bounded open subset of $\mathbb{R}^n$, $n\geq 2$, and suppose that $y\in\partial\Omega$. Denote the open ball of centre $y$ and radius $r$ by $B(r)$. A point $y$ is said to be \textit{$B$-regular} for $\Omega$ if, for each resolutive function $f$ on $\partial\Omega$ that is bounded in $B(R)\cap\partial\Omega$ for some positive $R$,
the Perron-Wiener-Brelot solution $H^\Omega$ of the Dirichlet problem is bounded in
$B(r)\cap\Omega$ for some positive $r$. Also, a point $y$ is said to be \textit{$lB$-regular} for $\Omega$ if there exists
a positive null sequence $\{r_m\}^\infty_1$, such that $y$ is $B$-regular for $B(r_m)\cap\Omega$ for all $m$.
If $y$ is $lB$-regular for $\Omega$, then it is $B$-regular for $\Omega$, see Sadi \cite{Sadi}. Is the converse true?\\
\textit{(D. H. Armitage and A. Sadi)}\\

\noindent\textbf{Update 3.33}
No progress on this problem has been reported to us.\\

\noindent\textbf{Problem 3.34}  Let $\Omega$ be a bounded domain in $\mathbb{R}^n$, $n\geq 2$, with the property that there exists $\alpha$ in  $(0,\pi]$ such that for every point $y$ in $\partial\Omega$ there is an open truncated cone with vertex
$y$ and angle $\alpha$ contained in $\Omega$. Is there some positive number $p=p(\alpha)$ such that every
positive superharmonic function in $\Omega$ belongs to $L^p(\Omega)$? If so, can such a $p$ be
characterised in terms of $\alpha$? \\
\textit{(D. H. Armitage)}\\

\noindent\textbf{Update 3.34}
This problem is also answered by Aikawa \cite{MR1169019}. \\

\noindent\textbf{Problem 3.35} For $r_1<r_2$, we will call the set $\{x\in\mathbb{R}^n:n\geq3,r_1<\|x\|<r_2\}$ an annulus and its closure a closed annulus. Let $\Omega$ be a non-empty subset of $\mathbb{R}^n$, $n\geq3$, such that
$\lambda(\overline{\Omega})<+\infty$, where $\lambda$ denotes $n$-dimensional Lebesgue measure. Then if, for each
point $x\in\mathbb{R}^n\setminus\overline{\Omega}$ we have that
\[\frac{1}{\lambda(\overline{\Omega})}\int_{\overline{\Omega}}\|x-y\|^{2-n}\,d\lambda(y)\]
equals the mean-value of the function $y\mapsto\|x—y\|^{2-n}$ over the unit sphere of $\mathbb{R}^m$, it can
be shown that $\overline{\Omega}$ is a closed annulus.\\
If, throughout the hypotheses, we replace $\overline{\Omega}$ by $\Omega$, can we conclude that $\Omega$ is an
annulus?\\
For details of closely-related work, see \cite{10.2307/2000876} and \cite{MR820533}.\\
\textit{(D. H. Armitage and M. Goldstein)}\\

\noindent\textbf{Update 3.35}
No progress on this problem has been reported to us.\\

\section*{New Problems}

\noindent To appear.

% ==================== CHAPTER 4 ==========================================

\chapter{Polynomials}

\underline{Preface by E. Crane (to appear)}\\

\noindent\textbf{Problem 4.1} Let $\{z_n\}, 1\leq n<\infty$ be an infinite sequence such that $|z_n|=1$. Define
\[A_n=\max_{|z|=1}\prod^n_{i=1}|z-z_i|.\]
Is it true that $\limsup_{n\to\infty}A_n=\infty$, and if so, how quickly must $A_n$ tend to infinity?\\
We may define $z_n$ inductively as follows, $z_1=1, z_2=-1$ and if $z_\nu$ has already been defined for $1\leq\nu\leq2^k$, then we define for $1\leq p\leq2^k$,
\[z_{p+2^k}=z_p\exp\Big(\frac{\pi i}{2^k}\Big).\]
With this definition, we easily see that $A_n\leq n+1$, with equality if, and only if, $n=2^k-1$ for some integer $k$. Is this example extreme?\\
\textit{(P. Erd\"os)}\\

\noindent\textbf{Update 4.1}
Wagner \cite{MR571725} has proved that   
\[\limsup_{n\to\infty}\frac{\log A_n}{\log\log n}>0\]
thus answering the first question  positively. Linden \cite{MR0437725} proved that $A_n$ can be less than $n^{1-\delta}$ for all large $n$, where $\delta$ is a positive constant thus giving a negative answer to the second question. A follow-up question by Erd\"os asks about the behaviour of 
\[\big(\prod^N_{n=1}A_n\big)^{1/N}.\]\\

\noindent\textbf{Problem 4.2} Let $p(z)=a_0+a_1z+\ldots+a_nz^n$ be a polynomial, all of whose zeros are on $|z|=1$. If 
\[A=\max_{0\leq k\leq n}|a_k|,\hspace{1cm}M=\max_{|z|=1}|p(z)|,\]
is $M\geq2A$?\\

\noindent\textbf{Update 4.2}
The original problem was mistated as $M\leq2A$. The inequality \mbox{$|a_\nu|\leq\frac{1}{2}M$} has been proved by Saff and Sheil-Small \cite{MR0382609}, except when $n$ is even and $\nu=\frac{1}{2}n$. Sheil-Small points out that Kristiansen \cite{MR551458} has proved that $M\geq2A$ in general.\\

\noindent\textbf{Problem 4.3} Let $P_N(z)$ be a polynomial with $N$ terms, satisfying $|P_N(z)|\leq1$ on $|z|=1$.
How large can $P_n(z)$ be if $P_n(z)$ is a partial sum of $P_N(z)$?\\
\textit{(P. Erd\"os)}\\

\noindent\textbf{Update 4.3}
No progress on this problem has been reported to us.\\

\noindent\textbf{Problem 4.4} Is there a function $f(k)$ of the positive integer $k$, so that the square of every polynomial having at least $f(k)$ terms has a least $k$ terms? Erd\"os proved that at any rate, $f(k)>k^{1+c}$ for a positive constant $c$.\\
\textit{(P. Erd\"os)}\\

\noindent\textbf{Update 4.4}
No progress on this problem has been reported to us.\\
  
\noindent\textbf{Problem 4.5} Let $P(z)$ be a polynomial whose zeros $z_1, z_2, \ldots, z_n$ lie in $|z|\leq1$. Is it true that $P'(z)$ always has a zero in $|z-z_1|\leq1$?\\
\textit{(Bl. Sendov)}\\

\noindent\textbf{Update 4.5} 
This problem is commonly known as the Sendov conjecture, although it had originially been erroneously attributed to Ilieff. With over eighty papers published on this problem, it would be a difficult task to list them all here, but we note that it was proved by Meir and Sharma \cite{MR0249587} for $n\leq5$, and obtained for general $n$ if $|z_1|=1$ by Rubinstein \cite{MR0237753}.\\
Gundersen writes that \cite{MR3051172} contains a useful survey, which draws attention to some highlights. Brown \cite{MR1081693} proves the conjecture for polynomials of degree $n=6$ with respect to a zero of the polynomial of modulus not exceeding $\frac{63}{64}$; and an asymptotic proof of a version of the conjecture is obtained by D\'egot \cite{MR3162254}. Kasmalkar \cite{MR3225918} strengthens a theorem by Chijiwa on this conjecture, and offers results which sharpen some of the estimates. Sendov \cite{MR2722582} also announces a stronger conjecture; proves it for polynomials of degree $n=3$; and  also announces a number of other conjectures, including a variation of Smale's mean value conjecture.  \\

\noindent\textbf{Problem 4.6}  If $H_\nu(z)$ is the $\nu$-th Hermite polynomial, so that 
\[H_\nu(z)e^{-z^2}=(-1)^\nu\big(\frac{d}{dz}\big)^\nu e^{-z^2},\]
is it true that the equation 
\[1+H_1(z)+aH_n(z)+bH_m(z)=0,\]
where $2\leq n<m$, and $a, b$ are complex, has at least one zero in the strip $|\text{Im}\, z|\leq c$, where $c$ is an absolute constant? (This is true if $b=0$, see Makai and Tur\'an \cite{MR0166418}.)\\
\textit{(P. Tur\'an)}\\

\noindent\textbf{Update 4.6}
No progress on this problem has been reported to us, although progress on Problem 4.5 was erroneously listed under this number in \cite{C}.\\

\noindent In the following five problems we use the notation
\[E_f^{(n)}=\{z\in\mathbb{C}:|f(z)|\leq1\}\]
where $n$ is the degree of the polynomial $f$.\\

\noindent\textbf{Problem 4.7} Let $f(z)=z^n+a_1z^{n-1}+\ldots+a_n$ be a polynomial of degree $n$. Cartan \cite{MR1509288} proved that the set $|f(z)|\leq1$, which we call $E^{(n)}_f$ can always be covered by discs, the sum of whose radii is at most $2e$. It seems likely that $2e$ can be replaced by $2$. If $E^{(n)}_f$ is connected, this was proved by Pommerenke \cite{MR0120481}, who also proved the general result, with $2.59$ instead of $2$.\\

\noindent\textbf{Update 4.7}
Mo \cite{Moguodan} attempts to prove this conjecture with a sharp
constant $2^{1-1/n}$ for each fixed $n$, instead of $2$. Unfortunately there are gaps in the
argument which so far it has proved impossible to fill. Eremenko and Hayman \cite{eremhayman1999ae} proved that the boundary of $E_f^{(n)}$ has length less than $9.173n$. See also Borwein \cite{zbMATH00741404}.\\

\noindent\textbf{Problem 4.8} Assume that $E^{(n)}_f$ is connected. Is it true that
\begin{equation}
\label{4.8eqn}
\max_{z\in E^{(n)}_f}|f'(z)|\leq\frac{1}{2}n^2\,?
\end{equation}
Pommerenke \cite{MR0109208} proved this with $\frac{1}{2}en^2$ instead of $\frac{1}{2}n^2$.\\

\noindent\textbf{Update 4.8}
Eremenko writes that the inequality (\ref{4.8eqn}) is not true as stated, as it is violated by Chebyshev's polynomials. The correct inequality is
\[\max_{z\in E_f^{(n)}}|f'(z)|\leq 2^{1/n-1}n^2\]
which is best possible. This has been proved by Eremenko and Lempert \cite{eremlempert}. There is also a generalisation of this inequality by Eremenko \cite{erem2006ae}.\\

\noindent\textbf{Problem 4.9} Is it true that to every positive $c$, there exists an $A(c)$ independent of $n$, such that $E_f^{(n)}$ can have at most $A(c)$ components of diameter greater than $1+c^2$?\\
\textit{(P. Erd\"os)}\\

\noindent\textbf{Update 4.9}
Erd\"os wrote that Pommerenke \cite{MR0151580} showed that to every positive $\varepsilon$
and integer $k$ there is an integer $n_0$ so that for $n > n_0$ there is a polynomial $f_n(z)$
of degree $n$ for which $E^{(n)}_f$ has at least $k$ components of diameter
greater than $4 - \varepsilon$, thus showing that the original conjecture is false.\\
One could try to estimate the number of
components of $E^{(n)}_f$ having diameter greater than $1 + c$. Is it
$o(n)$ as $n\to\infty$? This seems certain, but it could be $o(n^\varepsilon)$ as $n\to\infty$.\\

\noindent\textbf{Problem 4.10} Is it true that the length of the curve $|f_n(z)|=1$ is maximal for $f_n(z)=z^n-1$?\\
\textit{(P. Erd\"os)}\\

\noindent\textbf{Update 4.10}
 Eremenko and Hayman \cite{eremhayman1999ae} proved this for $n=2$. They also proved
that the set $L(f_n)=\{ z:|f_n(z)|=1\}$ is connected for
extremal polynomials. Further progress was made by Fryntov and Nazarov \cite{zbMATH05580784}
who proved that the length of $L(f_n)$ is
at most $(2+o(1))n$ which is the optimal asymptotics.
They also proved that the polynomial $z^n-1$ gives a local maximum
for this problem.\\

\noindent\textbf{Problem 4.11} If $|z_i|\leq1$, estimate from below, the area of $E^{(n)}_f$. Erd\"os, Herzog and Piranian \cite{MR0101311} prove that, given positive $\varepsilon$, the area of $E^{(n)}_f$ can be made less that $\varepsilon$, if $n>n_0(\varepsilon)$.\\

\noindent\textbf{Update 4.11}
No progress on this problem has been reported to us.\\

\noindent\textbf{Problem 4.12} If $-1\leq z_1\leq z_2\leq \ldots\leq z_n\leq1$, where the $z_i$ are the zeros of $f(z)$, it is known that the diameter of the set on the real line, for which $|f(x)|\leq1$, is at most $3$. Is it true that the measure of this set is at most $2\sqrt{2}$ (see also Problem 4.7)?\\

\noindent\textbf{Update 4.12}
This conjecture has been proved by Elbert \cite{MR0209422} and \cite{MR0204964}.\\

\noindent\textbf{Problem 4.13} It is known that there exists a polynomial $P(z)$
\[P(z)=\sum^n_{k=1}\varepsilon_k z^k, \hspace{1cm}\varepsilon_k=\mp1\]
for which
\begin{equation}
\label{4.1}
\max_{|z|=1}|P(z)|<C_1\sqrt{n}.
\end{equation}
(See Clunie \cite{MR0107715}). Is it necessarily true that $C_1>1+A$ if (\ref{4.1}) holds, where $A$ is a positive absolute constant?\\

\noindent\textbf{Update 4.13}
No progress on this problem has been reported to us.\\

\noindent\textbf{Problem 4.14} Does there exist a polynomial of the type in Problem 4.13, for which
\begin{equation}
\label{4.2}
\min_{|z|=1}|P(z)|>C_2\sqrt{n}
\end{equation}
for every $n$? More generally, does there exist such a polynomial satisfying both (\ref{4.1}) and (\ref{4.2})?\\

\noindent\textbf{Update 4.14}
An affirmative answer has been given by Beller and Newman \cite{MR0355015} for the case when $|\varepsilon_k|\leq1$ and by K\"orner \cite{MR572106} with $|\varepsilon_k|=1$. The case $\varepsilon_k=\pm1$ remains open.\\

\noindent\textbf{Problem 4.15} If again $\varepsilon_k=\mp1$, is it true that, for large $n$, all but $o(2^n)$ polynomials $P(z)=\sum^n_{k=1}\varepsilon_k z^k$ have just $n/2+o(n)$ roots in $\mathbb{D}$?\\

\noindent\textbf{Update 4.15}
No progress on this problem has been reported to us.\\

\noindent\textbf{Problem 4.16} Is it true that for all but $o(2^n)$ polynomials $P(z)$
\[\min_{|z|=1}|P(z)|<1,\]
or, if not, what is the corresponding correct result?\\

\noindent\textbf{Update 4.16}
No progress on this problem has been reported to us.\\

\noindent\textbf{Problem 4.17} It is shown by Salem and Zygmund \cite{MR0065679} that there exist positive constants $C_3, C_4$ such that for every positive $\varepsilon$, we have
\[(C_3-\varepsilon)\big(n\log n\big)^\frac{1}{2}<\max_{|z|=1}|P(z)|<(C_4+\varepsilon)\big(n\log n\big)^\frac{1}{2}\]
apart from $o(2^n)$ polynomials $P(z)$. Is this result true with $C_3=C_4$ and if so, what is the common value?\\

\noindent\textbf{Update 4.17}
Hal\'asz \cite{MR0367545} has proved this conjecture with $C_3 = C_4 = 1$ together
with the analogous result for trigonometric polynomials.\\

\noindent\textbf{Problem 4.18} If $f$ is any polynomial or rational function of degree $N$, find the least upper bound $\phi(N)$ of
\[\frac{1}{r}\int^r_0dt\int^\pi_{-\pi}\frac{|f'(re^{i\theta})|}{1+|f|^2}\,d\theta,\]
for varying $r$ and $f$. It is known only that $\phi(N)=O(N^{\frac{1}{2}})$, $\phi(N)\neq O(\log N)^{\frac{1}{2}-\varepsilon}$, as $N\to\infty$. The upper bound for varying $f$ of 
\[\int^\infty_0dt\int^\pi_{-\pi}\frac{|f'|}{1+|f|^2}\,dr\, d\theta\]
should not be much larger than $\phi(N)$, but there is no proof that it is not $+\infty$. The special cases $N=2, 3, \ldots$ would be of interest. See Littlewood \cite{MR0049315}.\\

\noindent\textbf{Update 4.18}
If $f$ is a polynomial or rational function, Littlewood asked for bounds depending on $N$ only, where $N$ is the degree or order of $f$, of the functionals
\[J_1(f)=\int^1_0r\,dr\int^\pi_{-\pi}\frac{|f'(re^{i\theta})|}{1+|f|^2}\,d\theta,\hspace{1cm} J_2(f)=\int^\infty_0\,dr\int^\pi_{-\pi}\frac{|f'(te^{i\theta})|}{1+|f|^2}\,d\theta.\]
Hayman \cite{MR581802} has shown that the sharp bounds $\psi_1(N), \psi_2(N)$ for $J_1(f)$ and $J_2(f)$ and rational functions $f$ of degree $N$ have orders $N^{1/2}$ and $N$ respectively.\\
The corresponding questions for polynomials are more difficult and more interesting. If $\phi_1(N)$, $\phi_2(N)$ are the bounds, then Hayman \cite{MR581802} has proved
\[\phi_1(N)>A\log N, \hspace{1cm}\phi_2(N)<A(N\log N)^{1/2},\]
thus slightly improving previous results by Littlewood \cite{MR0049315} and Chen and Liu \cite{MR0276445} respectively.
If $f$ is a polynomial of degree $N$, let $\phi(N)$ be the supremum of
\[\frac{1}{r}\int^r_0dt\int^\pi_{-\pi}\frac{|f'(re^{i\theta})}{1+|f|^2}\,d\theta,\]
for varying $r$ and $f$. Littlewood conjectured that
\[\phi(N)=O(N^{\frac{1}{2}-\eta})\]
for some positive $\eta$, and this has been proved by Lewis and Wu \cite{MR942832} by the same
method which they used for their attack on Problem 1.6. It depends on some
harmonic measure estimates of Bourgain \cite{MR880292}. Earlier Eremenko and Sodin \cite{MR897006} had shown
that
\[\phi(N)=o(N^\frac{1}{2}).\]
A consequence of these results is that an entire function of finite order has `almost
all roots of almost all equations $f(z) = a$' in a set of zero density in the plane.\\
D. Beliaev writes that Problems 4.18, 6.5, 6.7, and 6.8 are essentially equivalent, and the current best estimates for the optimal exponent in all of these problems are: lower bound $0.23$ (\cite{MR2680427}, \cite{MR2425327}), upper bound $0.46$ (\cite{MR2130416}).  See Problem 6.5 for details of the equivalence of Problems 6.5, 6.7 and part of 6.8. To see how Problem 4.18 is also essentially equivalent we need a proof of Binder and Jones that was never published.\\

\subsection*{Trigonometric polynomials}

\noindent\textbf{Problem 4.19} 
Littlewood conjectured that if $n_1, n_2, \ldots, n_k$ are distinct positive integers then
\begin{equation}
\label{4.3}
\int^{2\pi}_0\Big|\sum^k_{i=1}\cos (n_1 x)\Big|\,dx>c\log k.
\end{equation}
The best result known in this direction is due to Davenport \cite{MR0124681} who proved (\ref{4.3}) with
$\big(\log k/\log \log k\big)^{\frac{1}{4}}$ instead of $\log k$, thus sharpening an earlier result of Cohen \cite{MR0133397}.\\

\noindent\textbf{Update 4.19}
This was solved independently by Konyagin \cite{zbMATH03806164} and by McGehee, Pigno and Smith \cite{zbMATH03741921}.\\

\noindent\textbf{Problem 4.20} Let $f_n(\theta)$ be a trigonometric polynomial of degree $n$, all of whose roots are real. Is it true that
\[\int^\pi_0\big|f_n(\theta)\big|\,d\theta\leq 4?\]
See Erd\"os \cite{MR0003595}.\\
\textit{(P. Erd\"os)}\\

\noindent\textbf{Update 4.20}
This problem was originally misattributed to F. R. Keogh. It has been settled affirmatively
by Saff and Sheil-Small \cite{MR0382609}. They have also obtained the sharp bounds
for
\[\int^\pi_{-\pi}\big|f_n(re^{i\theta})\big|^p\,d\theta\]
whenever $f(z) = z^n + \ldots$ is a polynomial of degree $n$ all of whose zeros
lie on $|z|=1$ and $p > 0$.\\

\noindent\textbf{Problem 4.21} If $a_k=\mp1, k=0,\ldots,n$ and
\[b_k=a_na_{n-k}+a_{n-1}a_{n-k-1}+\ldots+a_ka_0,\]
is it true that
\[\sum^n_1|b_k|^2>An^2,\]
where $A$ is an absolute constant?\\
If $p(z)=a_0+a_1z+\ldots+a_nz^n$, then
\[|p(z)|^2=n+2\sum^n_1b_k\cos(k\theta),\]
so that the truth of the above inequality would imply
\[\frac{1}{2\pi}\int^{2\pi}_0\big|p(e^{i\theta})\big|^4\,d\theta\geq n^2(1+A).\]
\textit{(F. R. Keogh)} \\

\noindent\textbf{Update 4.21}
No progress on this problem has been reported to us.\\

\noindent\textbf{Problem 4.22} Using the notation of Problem 4.7, if $|z_i|\leq1$, Clunie
and Netanyahu (personal communication) showed that a path exists joining the origin to $| z| = 1$ in
$E^{(n)}_f$. What is the shortest length $L^{(n)}_f$ of such a path? Presumably
$L^{(n)}_f$ tends to infinity with $n$, but not too fast.\\
\textit{(P. Erd\"os)}\\

\noindent\textbf{Update 4.22}
No progress on this problem has been reported to us.\\

\noindent\textbf{Problem 4.23} Some of the Problems 4.7 to 4.12 extend naturally to the space of
higher dimensions. Let $x_i$ be a set of $n$ points in $\mathbb{R}^m$ and let $E^{(m)}_n$
be the set of points for which
\[\prod^n_{i=1}|x-x_i|\leq1.\]
When is the maximum volume of $E^{(m)}_n$ attained and how large
can it be? Piranian observed that the ball is not extreme for $m = 3$,
$n = 2$. If $E^{(m)}_n$ is connected, can it be covered by a ball of radius $2$?
For $m = 2$, this was proved by Pommerenke \cite{MR0109875}.\\
\textit{(P. Erd\"os)}\\

\noindent\textbf{Update 4.23}
No progress on this problem has been reported to us.\\

\noindent\textbf{Problem 4.24}  Let
\[P(z)=\sum^n_0a_kz^k\]
be a self-inversive polynomial, i.e. if $\zeta$ is a zero of $P(\zeta)$ with multiplicity
$m$, then $1/\zeta$ is also a zero with multiplicity $m$. Is it true that
$w = P(z)$ maps $\mathbb{D}$ onto a domain containing a disc of radius
$A= \max_{0\leq k\leq n} |a_k|$?\\
\textit{(T. Sheil-Small)}\\

\noindent\textbf{Update 4.24}
No progress on this problem has been reported to us.\\

\noindent\textbf{Problem 4.25} Determine
\[\inf\int^\pi_{-\pi}\big|1-e^{i\theta}\big|^{2\lambda}\big|P(e^{i\theta})\big|^2\,d\theta,\hspace{1cm}\lambda>0,\]
where $P(z)$ ranges over all polynomials with integer coefficients and leading coefficient
unity. (The solution would have number-theoretic applications.)\\ 
\textit{(W. H. J. Fuchs)}\\

\noindent\textbf{Update 4.25}
No progress on this problem has been reported to us.\\

\noindent\textbf{Problem 4.26} Let $P_n$ denote the class of polynomials $p(z)$, $p(0) = 1$, of degree at most $n$ 
and of positive real part in $\mathbb{D}$. Find
\[\max_{p\in P_n}\int^{2\pi}_0|p(e^{i\theta})|^2\,d\theta.\]
\textit{(F. Holland)}\\

\noindent\textbf{Update 4.26}
Let $P_n$ denote the class of polynomials $p(z)$ of degree $n$, having positive
real part and such that $p(0) = 1$. It was asked to find
\[\Lambda_n=\max_{p\in P_n}\sum^n_0|a_\nu|^2.\]
Goldstein and McDonald \cite{MR734993} have proved that $\Lambda_n\leq n +1$ and
calculated $\Lambda_2$ and $\Lambda_5$. Their research suggests that
\[\frac{\Lambda_n}{n}\to\Lambda\hspace{1cm}\text{ as }n\to\infty,\text{ where }\frac{2}{3}\leq\Lambda\leq1.\]\\

\noindent\textbf{Problem 4.27} Let $p(x)$ be a real polynomial of degree $n$ in the real variable $x$ such that
$p(x) = 0$ has $n$ distinct (real) rational roots. Does there necessarily exist a (real)
non-zero number $t$ such that $p(x)- t = 0$ has $n$ distinct (real) rational roots? (I can
prove this for $n = 1,2, 3$.)\\ 
\textit{(L. A. Rubel)}\\

\noindent\textbf{Update 4.27}
No progress on this problem has been reported to us.\\

\noindent\textbf{Problem 4.28} Suppose that $P$ is a non-linear polynomial with real coefficients. Show that $P^2(z)+P'(z)$ has non-real zeros. We conjecture that the lower bound for the number of non-real zeros is $\deg(P)-1$. If $P$ itself has only real zeros, this is proved by P\'olya and Szeg\"o \cite{MR1492448}. For the origin of this problem, see Problem 2.64.\\
\textit{(S. Hellerstein)}\\

\noindent\textbf{Update 4.28}
This has been solved by Sheil-Small \cite{zbMATH04097742}. An alternative solution was proposed by Eremenko (unpublished) which is reproduced in a paper by Bergweiler, Eremenko and Langley \cite{zbMATH02190101}.\\

\noindent\textbf{Problem 4.29} Yang \cite{MR0457415} claims to prove the following: let $P(z), Q(z)$ be monic polynomials such that $(i)$ $P(z)=0 \iff Q(z)=0$, and $(ii)$ $P'(z)=0 \iff Q'(z)=0$. Then there exist positive integers $m, n$ such that $P(z)^m\equiv Q(z)^n$. As pointed out by Rubinstein (personal communication), Yang's proof is incorrect since the inequalities \cite[p. 597]{MR0457415} are wrong. The problem then is to settle the Yang conjecture. If we let $\{z_1, z_2, \ldots, z_\nu\}$ be the distinct points at which $P$ (and therefore $Q$) has zeros, then the conjecture is easily established if $\nu\leq5$, and also in the case that $\nu$ is arbitrary and all the points $z_i$ are collinear.\\
\textit{(E. B. Saff)} \\

\noindent\textbf{Update 4.29}
A counter-example has been given by Roitman \cite{MR692530}. \\

\noindent\textbf{Problem 4.30} Let $\mathcal{P}$ denote the set of all polynomials of the form
\[p(z)=\prod^n_{\nu=1}(z-\zeta_\nu),\]
where $n\geq2$ and $|\zeta_\nu|\leq1$, $\nu=1, 2, \ldots, n$. The Sendov conjecture (see Problem 4.5) states: if $p(z)\in\mathcal{P}$ then each disc 
\[\{z:|z-\zeta_\nu|\leq1\},\hspace{1cm} \nu=1, 2, \ldots, n,\]
 contains at least one zero of $p'(z)$. Schmeisser \cite{MR0486436} proved this conjecture for certain subsets of $\mathcal{P}$. In all but two of these special cases, the proof show that the following stronger result is true:\\
\textit{If $\zeta$ is an arbitrary point of the convex hull of the zeros of $p(z)$, then the disc \mbox{$\{z:|z-\zeta|\leq1\text{ contains at least one zero of }p'(z)\}$}}.\\
We ask:
\begin{enumerate}[(a)]
\item is this stronger result true for all $p(z)$ in $\mathcal{P}_1$, where $\mathcal{P}_1$ is the subset of all polynomials in $\mathcal{P}$ which vanish at $0$?
\item is this stronger result true for all $p(z)$ in $\mathcal{P}_2$, where $\mathcal{P}_2$ is the subset of all polynomials in $\mathcal{P}$ which are of the form $p(z)=z^n+a_{n-1}z^{n-1}+\ldots+a_0$, $a_\nu\leq0$, $\nu=0, 1, \ldots, n-1$?
\end{enumerate}
\textit{(G. Schmeisser)}\\

\noindent\textbf{Update 4.30}
No progress on this problem has been reported to us.\\

\noindent\textbf{Problem 4.31} Erd\"os and Newman conjectured that if 
\begin{equation}
\label{D1}
f(z)=\sum^n_{k=0}a_kz^k,\hspace{1cm} |a_k|=1, \hspace{1cm}0\leq k\leq n,
\end{equation}
then there is an absolute constant $c$ such that
\begin{equation}
\label{D2}
\max_{|z|=1}|f(z)|>(1+c)n^{1/2}.
\end{equation}
The weaker form of our conjecture stated that (\ref{D2}) holds if we assume \linebreak \mbox{$a_k=\pm1$} (Problem 4.13). A stronger form would be that (\ref{D2}) holds even if \mbox{$f(z)=\sum^n_{k=1}a_{n_k}z^{n_k}$}, $n_k$ natural numbers, $|a_{n_k}|=1$. However, (\ref{D2}) was disproved by Kahane (no citation). In fact, he showed that, given that $\varepsilon>0$, there are polynomials of the form (\ref{D1}) for which
\begin{equation}
\label{D3}
\max_{|z|=1}|f(z)|<n^{1/2}+O(n^{(3/10)+\varepsilon})\hspace{1cm}\text{ as }n\to\infty.
\end{equation}
Show that $n^{(3/10)+\varepsilon}$ cannot be replaced by $n^\varepsilon$ in (\ref{D3}). Is there any $n$-th degree polynomial of the form (\ref{D1}) for which
\begin{equation}
\label{D4}
\min_{|z|=1}|f(z)|>(1-\varepsilon)n^{1/2}
\end{equation}
for every positive $\varepsilon$ if $n>n_0(\varepsilon)$? Perhaps there is an $n$-th degree polynomial of the form (\ref{D1}), for which, for all $z$, $|z|=1$,
\begin{equation}
\label{D5}
(1-\varepsilon)n^{1/2}<|f(z)|<(1+\varepsilon)n^{1/2}.
\end{equation}
In fact, (\ref{D5}) could possibly hold with $n^{1/2}+O(1)$ on the left and right. It would be worthwhile to determine if (\ref{D2}) holds for $a_k=\pm1$ and for $a_{n_k}=\pm1$.\\
\textit{(P. Erd\"os and D. Newman)} \\

\noindent\textbf{Update 4.31} No progress on this problem has been reported to us.

\section*{New Problems}

\noindent To appear.

%======================= CHAPTER 5 ====================================

\chapter{Functions in the Unit Disc}

\underline{Preface by L. R. Sons (to appear)}\\

\noindent\textbf{Notation} In this chapter, $\mathbb{D}$ denotes the open unit disc $|z|<1$, $\mathbb{T}$ denotes the circumference $|z|=1$ and $\overline{\mathbb{D}}=\mathbb{D}\cup \mathbb{T}$.
There are a number of problems concerned with the growth of the coefficients, means and maximum modulus of functions in $\mathbb{D}$ omitting certain values. Let
\[f(z)=\sum^\infty_{n=0}a_nz^n\]
be an analytic function in $\mathbb{D}$ with coefficients $a_n$, and suppose that
\[f(z)\neq w_k,\]
where $w_k$ is complex, $r_k=|w_k|$ is monotonic increasing and $r_k\to\infty$, as $k\to\infty$.\\
Write 
\[M(r,f)=\sup_{|z|=r}|f(z)|,\]
\[I_\lambda(r,f)=\Big\{\frac{1}{2\pi}\int^{2\pi}_0\big|f(re^{i\theta})\big|^\lambda\,d\theta\Big\}^{1/\lambda}.\]
While the theory for $M(r,f)$ is well-worked out, (see e.g. Cartwright \cite{MR1512983} and  Hayman \cite{hayman_1948}, \cite{PLMS:PLMS0450}) very little that is not trivial is known about $I_\lambda(r,f)$ and the coefficients $a_n$. If 
\[r_{k+1}\leq Cr_k,\]
it is known (see Littlewood \cite{MR1574417}) that
\[M(r,f)=O(1-r)^{-A(C)}\]
where $A(C)$ depends on $C$ only. If
\begin{equation}
\label{5.1}
\frac{r_{n+1}}{r_n}\to1,
\end{equation}
Cartwright \cite{MR1512983} proved that
\begin{equation}
\label{5.2}
M(r,f)=O(1-r)^{-2-\varepsilon},\hspace{1cm}\text{ as }r\to1
\end{equation}
for every positive $\varepsilon$. If
\begin{equation}
\label{5.3}
\sum\Big(\log\frac{r_{n+1}}{r_n}\Big)^2<+\infty,
\end{equation}
then Hayman \cite{PLMS:PLMS0450} proved that
\begin{equation}
\label{5.4}
M(r,f)=O(1-r)^{-2}.
\end{equation}
These two results are essentially best possible.\\

\noindent\textbf{Problem 5.1} Is it true that (\ref{5.1}) implies
\[I_1(r,f)=O(1-r)^{-1-\varepsilon}\]
and
\[|a_n|=O(n^{1+\varepsilon})\,?\]\\

\noindent\textbf{Update 5.1} For Problem 5.1 to 5.3, substantial progress has been made by Pommerenke \cite{MR0313496}. He
has shown that if $f(z)$ satisfies
\[\Big|\frac{f'(z)}{f(z)}\Big|\leq\frac{\alpha}{1-|z|},\hspace{1cm} r_0<|z|<1,\]
then
\[I_1(r,f)=O((1-r)^{-\lambda(\alpha)}),\hspace{1cm} r\to1\]
and hence
\[a_n=O(n^{\lambda(\alpha)}),\]
where $\lambda(\alpha)=\frac{1}{2}(\sqrt{1+4\alpha^2}-1)$. In particular, if $f$ is weakly univalent we
may take $\alpha = 2 + \varepsilon$, $\lambda(\alpha) = \frac{1}{2}(\sqrt{17}-1)+\varepsilon=1. 562\ldots + \varepsilon$.\\
Suppose that
\[f(z)=\sum^\infty_{n=0} a_nz^n\]
is analytic in $\mathbb{D}$ and $f(z)\neq w_n$ there, where $w_n$ is a sequence of complex values,
such that
\[w_n\to\infty,\hspace{1cm} \Big|\frac{w_{n+1}}{w_n}\Big|\to1,\hspace{1cm}\text{ as }n\to\infty.\]
Then Cartwright \cite{MR1512983} had shown that
\[M(r,f) = O(1 -r)^{-2-\varepsilon},\hspace{1cm}\text{ as }r\to\infty\]
for every positive $\varepsilon$. Baernstein and Rochberg \cite{MR0437765} have obtained the analogous
results for the means
\[I_\lambda(r,f)=\frac{1}{2\pi}\int^{2\pi}_0\big|f(re^{i\theta})\big|^\lambda \,d\theta,\]
namely
\[I_\lambda(r,f)= O(1 -r)^{1-2\lambda-\varepsilon}, \hspace{1cm}\lambda\geq\frac{1}{2},\]
\[I_\lambda(r,f)= O(1 -r)^{-\varepsilon},\hspace{1cm} \lambda<\frac{1}{2},\]
as $r\to1$. From the case $\lambda=1$, they deduce that
\[a_n=O(n^{1+\varepsilon}).\]
Baernstein and Rochberg \cite{MR0437765} also obtain the analogous result when the equations \mbox{$f(z)=w_n$} have
at most $p$ roots for some positive integer $p$.\\

\noindent It is natural to ask whether, under stronger hypotheses on the $w_n$, it is
possible to get rid of the $\varepsilon$ in the above results. Hansen \cite{MR0262512}, Baernstein \cite{MR0417406}, and Hayman and Weitsman \cite{MR0361049} each did so on the hypothesis that $f(z)$ fails to take some value on every circle
$|w|=R$ (or every circle apart from a set of $R$ of finite logarithmic measure). However,
it ought to be possible to obtain the same conclusion under weaker hypotheses.\\

\noindent\textbf{Problem 5.2} Is it true that (\ref{5.3}) implies that
\begin{equation}
\label{5.5}
I_1(r,f)=O(1-r)^{-1}
\end{equation}
and 
\begin{equation}
\label{5.6}
|a_n|=O(n)?
\end{equation}\\

\noindent\textbf{Update 5.2}
See Update 5.1.\\

\noindent\textbf{Problem 5.3} An even stronger hypothesis than (\ref{5.3}) is that $f(z)$ is \textit{weakly univalent} (see Hayman \cite{MR0045210}) i.e. for every $r$ with $0<r<\infty$, either $f(z)$ assumes every value on $|w|=r$ exactly once, or there exists a complex $w=w_r$, such that $|w_r|=r$ and $f(z)\neq w_r$. Even with this assumption, nothing stronger than the results
\[I_\lambda(r,f)=O(1-r)^{-2},\]
\[|a_n|=O(n^2)\]
are known (which are trivial consequences of (\ref{5.4})). It would be interesting to obtain some sharpening of these results even if it is not possible to deduce the full strength of (\ref{5.5}), (\ref{5.6}).\\

\noindent\textbf{Update 5.3}
See Update 5.1.\\

\noindent\textbf{Problem 5.4} If the sequence $w_n$ satisfies
\begin{equation}
\label{5.7}
\arg w_n=O\big(|w_n|^{\frac{1}{2}}\big)
\end{equation}
and 
\begin{equation}
\label{5.8}
|w_{n+1}- w_n|=O(|w_n|^{\frac{1}{2}})
\end{equation}
then it is known (see Hayman \cite{MR0070716}) that (\ref{5.5}) and hence (\ref{5.6}) hold. It is interesting to ask whether the method will yield the same conclusions under somewhat weaker hypotheses, such as for instance, replacing the index $\frac{1}{2}$ in (\ref{5.7}) and (\ref{5.8}) by a smaller positive number.\\

\noindent\textbf{Update 5.4}
If the sequence $w_n$ satisfies (\ref{5.7}) and $|w_{n+1}-w_n|=O\big(|w_n|^\lambda\big)$, where $\lambda<1$, then
Higginson (unpublished) has shown that this is indeed possible. Hayman had previously obtained the
corresponding result when $\lambda=\frac{1}{2}$. Higginson \cite{MR0447549} has also shown that if
\[f(z)=z^p+a_{p+1}z^{p+1}+\ldots\]
is \textit{weakly $p$-valent} in $\mathbb{D}$, that is, if $f(z)$ either assumes, for every $R$, all values on
$|w| = R$ exactly $p$ times, or assumes at least one such value less than $p$ times, then
\begin{equation}
\label{C5.4.4}
\limsup_{n\to\infty}\frac{|a_n|}{n^{2p-1}}<\frac{1}{\Gamma(2p)}
\end{equation}
unless $f(z) \equiv z^p(1-ze^{i\theta})^{-2p}$. The sharp bounds for the means of this class of
functions had previously been obtained by Baernstein \cite{MR0417406}. In contrast to the $p$-valent
or mean $p$-valent case, it is not true that the upper limit in (\ref{C5.4.4}) exists as a limit in general.
It does so, if and only if, there is exactly one radius of greatest growth.\\

\noindent\textbf{Problem 5.5} If $f(z)=u+iv$ assumes only values in the right half-plane, then subordination shows that
\begin{equation}
\label{5.9}
a_n=O(1).
\end{equation}
It is of interest to ask what other hypotheses on the values assumed by $f(z)$ result in (\ref{5.9}). Let $d(r)$ be the radius of the largest disc whose centre lies on $|w|=r$, and every value of whose interior is assumed by $f(z)$. If $d(r)\leq d$, then it is shown by Hayman \cite{MR0070716} that (\ref{5.9}) holds. It is fair to ask whether this conclusion still holds if $d(r)\to\infty$ sufficiently slowly.\\

\noindent\textbf{Update 5.5} For Problem 5.5 to 5.7,  let $D$ be a domain and let $w = f(z) = \sum^\infty_0a_n z^n$ be a function in the unit
disc $\mathbb{D}$ with values in $D$. It was asked under what conditions on $D$ we can conclude
that
\begin{enumerate}[(a)]
\item $a_n = o(1)$,
\end{enumerate}
and
\begin{enumerate}[(b)]
\item  $a_n = O(1)$.
\end{enumerate}
Fernandez \cite{MR769161} has shown that if the complement of $D$ has capacity zero, then $(a)$ is never
true, and $(b)$ holds if and only if $D$ is a \textit{Bloch domain}, that is, $D$ does not contain
arbitrarily large discs. Pommerenke \cite{MR0284574} has shown that if the complement of $D$ has
positive capacity and $D$ is a Bloch domain, then $(a)$ holds.\\

\noindent\textbf{Problem 5.6} It is known that there exist functions which fail to take any of the values $2\pi ik$, $-\infty<k<+\infty$ and which do not satisfy (\ref{5.9}), and in fact $|a_n|\leq \log\log n$ (see Littlewood \cite[p.205]{MR0012121},  and Hayman \cite{MR0070716}). However (\ref{5.9}) holds if $f(z)$ omits all but a finite interval of the imaginary axis (again by subordination). This suggests that (\ref{5.9}) might still hold if the omitted values $w_n$ cluster near $\infty$ sufficiently close to the imaginary axis.\\

\noindent\textbf{Update 5.6} See Update 5.5.\\

\noindent\textbf{Problem 5.7} If $c_k$ is a sequence of positive numbers such that
\[\sum c_k=S<+\infty,\]
and $n_k$ is an arbitrary sequence of positive integers, then
\[f(z)=\sum^\infty_{n=0}c_kz^{n_k}=\sum^\infty_{n=0}a_nz^n\]
is bounded in $\mathbb{D}$, and so takes no values outside a fixed disc. This shows that no conditions on the omitted values $w$ can imply more than
\begin{equation}
\label{5.10}
a_n=o(1).
\end{equation}
Clearly (\ref{5.10}) holds if $f(z)$ is bounded, since then
\[I_2(r,f)=\Big(\sum^\infty_{n=0}|a_n|^2\Big)^{\frac{1}{2}}<+\infty.\]
It would be of interest to obtain a non-trivial condition on the values omitted by $f(z)$, which would imply (\ref{5.10}). Such a condition might be $d(r)\to0$, where $d(r)$ is defined as in Problem 5.5.\\

\noindent\textbf{Update 5.7}
The work of Hansen \cite{MR0262512} and of Hayman and Weitsman \cite{MR0361049} also
leads to conditions on a set of omitted values such that the corresponding function
satisfies
\begin{equation}
\label{C5.7.5}
a_n\to0.
\end{equation}
However all our conditions require the set of omitted values to have positive
capacity. The biggest functions omitting a set $E$ of complex values from the point of
view of the means and maximum modulus, are the functions $F$ which map $\mathbb{D}$
onto the infinite covering surface over the complement of $E$. This led Pommerenke
to ask (see Problem 5.33), whether these functions $F(z)$ satisfy (\ref{C5.7.5}). Patterson, Pommerenke
and Hayman \cite{MR0486496} have been able to give an affirmative answer to this question, by showing
that (\ref{C5.7.5}) holds for $F(z)$ if $E$ is a lattice or if $E$ is thick at $\infty$ in the sense that the distance
$d(w)$ of any point in the $w$-plane from $E$ tends to zero as $w\to\infty$. This answers
Problem 5.7 for the functions $F(z)$. However the case of the subordinate
functions $f(z)$ remains open.
See Update 5.5 also.\\

%\section* *{Bloch's Theorems}

\noindent\textbf{Problem 5.8} Suppose that $f(z)=z+a_2z^2+\ldots$ is analytic in $\mathbb{D}$. Then $f(z)$ maps some sub-domain of $\mathbb{D}$ univalently into a disc of radius at least $B$, where $B$ is Bloch's constant. What is the value of $B$? The best results known are $B\geq\sqrt{3}/4>0.433$, due to Ahlfors \cite{MR1501949}, and $B<0.472$ due to Ahlfors and Grunsky \cite{MR1545698}. The upper bound is conjectured to be the right one. Heins \cite{MR0130974} has shown that $B>\sqrt{3}/4$.\\

\noindent\textbf{Update 5.8}
Tur\'an \cite{MR0219491} has shown that if $f(z)=\sum^\infty_{n=1}a_nz^n$ is univalent in $\mathbb{D}$
with $\sum^{2p}_p|a_n|^2\geq1$ for some $p$, then the corresponding Bloch's constant is at least $(32e^2)^{-1}$. An analogous result holds with $2p$ replaced by $Cp$, where
$C$ is any constant greater than $1$, but not if $C$ is allowed to tend to $\infty$
with $p$, however slowly. This latter result is due to Petruska \cite{MR0262473}. Bonk \cite{zbMATH04173340} was the first to prove that $B>\sqrt{3}/4$. The best current estimate from below is 
$B>\sqrt{3}/4+2.10^{-4}$ by Chen and Gauthier \cite{zbMATH00971867}.\\

\noindent\textbf{Problem 5.9} With the hypotheses of Problem 5.8, it follows that $f(z)$ assumes all values in some disc of radius $L$, $L\geq B$. What is the value of $L$? The best lower bound for the Landau constant $L$ is $L\geq\frac{1}{2}$, due to Ahlfors \cite{MR1501949}.\\

\noindent\textbf{Update 5.9}
An upper bound for $L$ is $L<0.54326$, which has been been obtained by Rademacher \cite{MR0008257}. Rademacher mentions that this bound has also been found, but not published, by Robinson.\\

\noindent\textbf{Problem 5.10} If, in addition, $f(z)$ is univalent in $\mathbb{D}$, the conclusions of Problem 5.8 and Problem 5.9 follow with a constant $S$, $S\geq L$, known as the schlicht Bloch's constant. What is the value of $S$? We may also ask the same question when, in addition, $f(z)$ is star-like, thus obtaining a still larger constant $S_1$. If $f(z)$ is convex, the correct value of the constant is $\pi/4$, attained for $f(z)=\frac{1}{2}\log\frac{1+z}{1-z}$, which maps $\mathbb{D}$ onto the strip $|\text{Im}\, z|<\pi/4$.\\

\noindent\textbf{Update 5.10}
No progress on this problem has been reported to us.\\

%\section* *{Normal families}

\noindent
Various conditions on functions $f(z)$ in the plane result in these functions being constant. This is the case, for instance, if $f(z)$ is analytic and bounded (Liouville's Theorem) or more generally, if $f(z)$ is meromorphic and fails to take three fixed values (Picard's Theorem), or if $f(z)$ is analytic, $f(z)\neq0, f'(z)\neq1$ (see Hayman \cite{MR0110807}). Classes of functions satisfying one of the above conditions in $\mathbb{D}$ form normal families, the results being due to Montel \cite{Montel} and Miranda \cite{MR1505035} respectively.\\

\noindent This suggests investigating other conditions  for functions in the unit disc, where the corresponding conditions for functions in the plane lead to constants.\\

\noindent The interested reader is directed to a more recent result in this area, namely Zalcman's Lemma \cite{MR0379852}, as stated by Schwick \cite{MR1762263}: Let $\mathcal{F}$ be a family of meromorphic functions on $\mathbb{D}$ which are not normal at $0$. Then there exist sequences $f_n$ in $\mathcal{F}, z_n, \rho_n$, and a non-constant function $f$ meromorphic in the plane such that
\[ f_n(z_n+\rho_nz)\to f(z),\] 
locally and uniformly (in the spherical sense) in $\mathbb{C}$, where $z_n\to0$ and $\rho_n\to0$.\\
By means of this lemma, a positive answer to several of the questions in this subsection can be given.\\

\noindent Do the following classes form normal families, possibly after suitable normalisations?\\

\noindent\textbf{Problem 5.11} $f(z)$ meromorphic in $\mathbb{D}$, $f(z)\neq0, f^{(l)}(z)\neq1$, where $l\geq1$.\\

\noindent\textbf{Update 5.11}
This has been shown to be the case by Ku Yung Hsing \cite{MR511294}.\\

\noindent\textbf{Problem 5.12} $f(z)$ analytic in $\mathbb{D}$, $f'(z)f(z)^n\neq1$, where $n\geq1$.\\

\noindent\textbf{Update 5.12}
This has been answered affirmatively for $n\geq2$ by Yang and Zhang \cite{MR576357}. The
case $n = 1$ has also been answered affirmatively by Oshkin \cite{MR650781}.\\

\noindent\textbf{Problem 5.13} $f(z)$ meromorphic in $\mathbb{D}$, $f'(z)f(z)^n\neq1$, for $n\geq3$.\\

\noindent\textbf{Update 5.13}
This is solved for $n\geq 2$ by Pang \cite{MR1057999}; and for $n\geq 1$ this follows from
combining the results of Pang \cite{zbMATH04159040} and Bergweiler and Eremenko \cite{MR1344897}. See also Wang and Fang \cite{zbMATH01159053}.\\

\noindent\textbf{Problem 5.14} $f'-f^n\neq a$, where $a$ is some complex number, and $n\geq 5$ if $f$ is meromorphic, $n\geq3$ if $f$ is entire.\\
The corresponding results for functions in the place are proved by Hayman \cite{MR0110807}, except for the case $n=1$ of Problem 5.12, which is a result of Clunie \cite{MR0214769}.\\

\noindent\textbf{Update 5.14}
Eremenko writes that this problem with $n\geq3$ is equivalent to Problem 5.13 with $n\geq1$.\\

%\section* *{Symmetrisation}

\noindent\textbf{Notation} Let $D$ be a domain in the plane. The \textit{circularly symmetrised domain} $D^*$ of $D$ is defined as follows. For every $r$, with $0<r<\infty$, the intersection of $D^*$ and the circle $|z|=r$ is
\begin{enumerate}[(i)]
\item The whole circle $|z|=r$, if $D$ contains the whole circle $|z|=r$;
\item null if $D$ does not meet $|z|=r$;
\item the arc $|\arg z|<\frac{1}{2}\ell(r)$ if neither $(a)$ nor $(b)$ holds, but $D$ meets $|z|=r$ in a set of arcs of total length $r\ell(r)$.
\item In addition $D^*$ contains the origin, if and only if $D$ does.
\end{enumerate}
Suppose now that $a_0$ is a point of $D$ on the positive real axis. Let
\[f(z)=a_0+\sum^\infty_{n=1}a_nz^n\]
be analytic in $\mathbb{D}$, and assume only values lying in $D$, and let
\[f^*(z)=a_0+\sum^\infty_1a_n^*z^n\]
map $\mathbb{D}$ onto the infinite covering surface over $D^*$. The following facts are known:
\begin{enumerate}[(a)]
\item If $D^*$ is simply connected, so that $f^*(z)$ is univalent, then $|a_1|\leq|a_1^*|$.\\ 
From this, it is not difficult to deduce that:
\item If $m_0(r,f)=\inf_{|z|=r}|f(z)|$, and $M(r,f)=\sup_{|z|<r}|f(z)|$, then
\[\mu(r,f^*)\leq m_0(r,f)\leq M(r,f)\leq M(r,f^*),\]
and, in special cases, rather more (see Hayman \cite{MR0108586}, Ch. 4).
\end{enumerate}
This leads to the following questions:\\

\noindent\textbf{Problem 5.15} Is it possible to remove the restriction that $D^*$ is simply connected in $(a)$ and $(b)$ above? It might be possible to start with the case when $D$ and $D^*$ are both doubly-connected.\\

\noindent\textbf{Update 5.15}
Another way to phrase this question, is to ask whether the Poincar\'e metric is decreased by symmetrisation. This result was known to be true if the symmetrised domain is simply connected. Lai \cite{MR532327} and Hempel \cite{MR623485} have independently proved that this is indeed the case, when the complement of the domain $D$ consists of the origin and one other point. This leads to sharp forms of the theorems of Schottky and Landau. The result has also been proved by Weitsman \cite{MR554396} when the domain is already symmetrical about some ray.\\
Weitsman \cite{MR847954} has proved  in general that the Poincar\'e metric of a domain $D$ is decreased by symmetrisation. It follows that the maximum modulus of a function $f$ mapping $\mathbb{D}$ into $D$ is dominated by
that of a function $F$ mapping $\mathbb{D}$ onto the universal cover surface over $D^*$. The
corresponding conclusion for the means of $f$ and $F$ is still open, unless $D^*$ is simply
connected, when it had been proved by Baernstein \cite{MR0417406}.\\

\noindent\textbf{Problem 5.16} Do corresponding results to Problem 5.15(a) apply to the means
\[I_\lambda(r,f)=\Big\{\frac{1}{2\pi}\int^{2\pi}_0\big|f(re^{i\theta})\big|^\lambda \,d\theta\Big\}^{1/\lambda}, \hspace{1cm}0<\lambda<\infty\]
or the Nevanlinna characteristic
\[T(r,f)=\frac{1}{2\pi}\int^{2\pi}_0\log^+\big|f(re^{i\theta})\big|\,d\theta?\]
This is known to be the case when $D$ is already symmetrical (but possibly multiply-connected) so that $D=D^*$, as a consequence of the theory of subordination (see e.g Littlewood \cite[Theorem 210, p. 164]{MR0012121}).\\

\noindent\textbf{Update 5.16}
These results have been proved by Baernstein \cite{MR0374429} if $D$ (and so
$D^*$) is simply-connected. See also \cite{MR0387590}.\\

\noindent\textbf{Problem 5.17} Let $D=D_0$ be a domain, $g(z,a_0)$ be the Green's function of $D$ with respect to a point $a_0$ on the positive real axis, and let $D_\lambda$ be the part of $D$ where $g>\lambda$ for $0<\lambda<\infty$. Is it true, at least in some simple cases, that $(D_\lambda)^*\subset(D^*_\lambda)$? The cases where $D^*$ consists of the plane or the unit disc cut along the negative real axis are of particular interest. A positive answer to this problem for simply-connected domains $D$ would lead to a positive answer of Problem 5.16 for the same class of domains, using a formula of Hardy-Stein-Spencer (see Hayman \cite[Chapter 3, p. 42]{MR0108586}).\\
With the general notation of the introduction above, when can we assert that $|a_n|\leq|a_n^*|$ for general $n$? This is true, for instance, when $D^*$ is the plane cut along the negative real axis, so that we obtain the Littlewood conjecture that $|a_n|\leq4|a_0|n$ for non-zero univalent $f$. In fact, $|a_1|\leq4|a_0|$ in this case, by symmetrisation and subordination, and $|a_n|\leq n|a_1|$ by De Branges' proof of the Bieberbach conjecture (see Update 6.1).\\
If $D$ is convex and $D=D^*$, the exact bound for $|a_n|$ is $|a_1|$, but $|a_n^*|\leq|a_1|$ in general (see Hayman \cite[Theorem 1.7, p. 12]{MR0108586}).\\

\noindent\textbf{Update 5.17}
In the original formulation of this question, Littlewood's conjecture was left open. There has been no further progress on this problem reported to us.\\

\noindent\textbf{Problem 5.18} Let $f(z)=\lambda+a_1z+\ldots$ be analytic in $\mathbb{D}$, where $0<\lambda<1$. Find the best constant $B(\lambda)$ such that if
\[F(r)=\lambda+|a_1|r+|a_2|r^2+\ldots,\]
then 
\[F[B(\lambda)]\leq1\]
for all $f$. It is known that  (see Bombieri \cite{MR0162918}):
\[B(\lambda)= 
\begin{cases}
(1+2\lambda)^{-1} &\text{if }\frac{1}{2}\leq\lambda<1,\\
1/\sqrt{2}&\text{if }\lambda=0,\\
\end{cases}\]
and that
\[B(\lambda)>[(1-\lambda)/2]^{\frac{1}{2}},\hspace{1cm}\text{if }0<\lambda<\frac{1}{2}.\]

\noindent\textbf{Update 5.18}
No progress on this problem has been reported to us.\\

\noindent\textbf{Problem 5.19} A function meromorphic in $\mathbb{D}$ which has no asymptotic value, assumes every value infinitely often in the disc. Every point of the circumference $\mathbb{T}$ is a Picard point for such a function, i.e. a point such that all except perhaps two values are taken in every neighbourhood. Functions with no exceptional values in the global sense are known. Can locally exceptional values occur? (See Cartwright and Collingwood  \cite{MR0050010}, \cite{MR0130379}).\\
\textit{(E. F. Collingwood)}\\

\noindent\textbf{Update 5.19}
Collingwood noted that a function meromorphic in $\mathbb{D}$ and
without asymptotic values, must assume every value infinitely often in the disc, and
that near any point of $\mathbb{T}$, all values with at most $2$ exceptions must be
taken infinitely often. He asked whether such local exceptional values can occur. \\
A strongly positive answer has been given by Eremenko (unpublished), who noted that one can
construct such a function, for which the only limit point of poles and zeros is the
point $z = 1$, by slightly modifying an example of Barth and Schneider \cite{MR0293089}.\\

\noindent\textbf{Problem 5.20} Plessner \cite{MR1581142} proves after Privaloff that if $f$ is analytic in $\mathbb{D}$, almost all points $P$ of the boundary are of two kinds. Either
\begin{enumerate}[(a)]
\item $f$ tends to a finite limit, as $z\to P$ in any Stolz angle $S$, lying in $\mathbb{D}$, or
\item as $z\to P$ in any $S$, $f$ takes (infinitely) often all values of a dense set.
\end{enumerate}
Can `dense set' be replaced by anything bigger here; e.g. the complement of a set of measure zero?\\
\textit{(E. F. Collingwood)}\\

\noindent\textbf{Update 5.20}
No progress on this problem has been reported to us.\\

\noindent\textbf{Problem 5.21} Corresponding to each function $f$ analytic in $\mathbb{D}$, and each value $w$, with $(|w|<1)$, write
\[f_w(z)=f\Big(\frac{z-w}{1-wz}\Big)=\sum^\infty_{n=0}a_n^{(w)}z^n,\]
\[\|f_w\|=\sum^\infty_{n=0}|a^{(w)}_n|\]
and let $W_f$ be the set of all values $w$, with $(|w|<1)$ for which
\[\|f_w\|<\infty,\]
Since $a^{(w)}_n$ is a continuous function of $w$, it follows that $W_f$ is a set of type $F_\sigma$. What more can be said? For example, if $W_f$ is everywhere dense in $D$ (or uncountable, or of positive measure), is $W_f$ the unit disc? It is known that $W_f$ may be a proper non-empty subset of $D$. Is $W_f$ either empty, or all of $D$ if $f$ is univalent?\\

\noindent\textbf{Update 5.21}
No progress on this problem has been reported to us.\\

%\section* *{Spaces of analytic functions}

\noindent\textbf{Problem 5.22} Let $H^p$ be the space of functions $f(z)=\sum^\infty_{n=0}a_nz^n$ analytic in
$\mathbb{D}$, and such that
\[\int^{2\pi}_0\big|f(re^{i\theta})\big|^p\,d\theta\]
remains bounded as $r\to1$. We define $H^\infty$ to be the class of bounded
functions in $\mathbb{D}$. For $0 < p < 1$ describe the coefficient multipliers
from $H^p$ to $H^p$. That is, for each such $p$ describe the sequences
$\lambda_n$ such that
\[\sum\lambda_n a_n z^n\in H^p\hspace{1cm}\text{whenever }\sum a_nz^n\in H^p.\]
\textit{(P. L. Duren)}\\

\noindent\textbf{Update 5.22}
No progress on this problem has been reported to us.\\

\noindent\textbf{Problem 5.23} Describe similarly the coefficient multipliers from $S$ to $S$,
where $S$ is the class of functions $\sum^\infty_{n=1} a_nz^n$ univalent in $\mathbb{D}$ either
\begin{enumerate}[(a)]
\item with the normalisation $a_1 = 1$, or\
\item generally.
\item What are the multipliers of the space of close-to-convex functions
into itself?
\item What are the multipliers of $S$ into the class $C$ of convex
functions?
\item What are the multipliers from the class $N$ of functions of bounded
characteristic into itself? The analogous problem for the class $N^+$ may
be more tractable. (The definition of $N^+$ is too lengthy for this work, but the reader is directed to Duren \cite[p. 25]{MR0268655} for more details).
\end{enumerate}
Ruscheweyh and Sheil-Small in solving Problem 6.9 have shown
that $(\lambda_n)$ is a multiplier sequence from $C$ into itself if and only if
$\sum\lambda_nz^n\in C$. In general, one can obtain only some sufficient conditions.
Thus in most cases $f(z) = \sum a_n z^n$ belongs to a class $A$ if $a_n$ is sufficiently
small, and conversely if $f\in A$, then $a_n$ cannot be too big. For
example, $\sum^\infty_2n|a_n|\leq1$ is a sufficient condition for $f(z)$ in $S$, and
$|a_n|\leq n\sqrt{7/6}$ is a necessary condition (see Fitzgerald \cite{MR0335777}). Similarly, if
$\sum^\infty_1|a_n|<\infty$, then $f(z)$ is continuous in $\overline{\mathbb{D}}$, and so belongs to $H^p$ for every positive $p$ and to $N$, whilst if $f\in N$, then $|a_n|\leq\exp(cn^\frac{1}{2})$ for
some constant $c$. Again, if $f(z)$ belongs to one of the above classes, then
so does $\frac{1}{t}f(tz)$ for $0 < t < 1$, so that the sequence $(t^{n-1})$ is a multiplier
sequence. In other cases, negative results are known. Thus
Frostman \cite{MR0012127} showed that $(n)$ is not a multiplier sequence from $N$ to $N$,
and Duren \cite{MR0268655} showed that $\big(\frac{1}{n+1}\big)$ is not such a sequence either.\\
\textit{(P. L. Duren, except for $(e)$ which is due to A. Shields)}\\

\noindent\textbf{Update 5.23}
For part $(a)$, Sheil-Small observes that $\sum^\infty_{n=1}b_na_nz^n$ is close-to-convex for every function $\sum^\infty_{n=1}a_nz^n$ which is close-to-convex, if and only if, $\sum^\infty_{n=1}b_nz^n$ is convex in $\mathbb{D}$. The sufficiency was shown by Ruscheweyh and Sheil-Small (\cite{MR0328051} and \cite{MR0447552}), while Sheil-Small notes that the necessity is a simple application of the duality principle.\\

\noindent\textbf{Problem 5.24} Is the intersection of two finitely generated ideals in $H^\infty$ finitely
generated?\\
\textit{(L. A. Rubel)}\\

\noindent\textbf{Update 5.24}
For $H^\infty(\mathbb{D})$, an affirmative answer was given by McVoy and Rubel \cite{MR0410393}. This was extended to any finitely connected domain by Barnard (no citation). Amar (no citation) has given a negative answer in $H^\infty(D^n)$ or $H^\infty(B^n)$ where $D^n$ is the polydisc and $B^n$ is the unit ball, for $n\geq3$. The case $n=2$ is still open.\\

\noindent\textbf{Problem 5.25} Let $W^+$ be the Banach algebra of power series $f(z)=\sum^\infty_{n=0}a_nz^n$  absolutely convergent in $|z|\leq1$, with $\|f\|=\sum^\infty_{n=0}|a_n|$.\\
Which functions generate $W^+$? More precisely, for which functions
$f$ is it true that the polynomials in $f$ are dense in $W^+$? It is clear
that a necessary condition is that $f$ be univalent in $\overline{\mathbb{D}}$. Newman \cite{MR0170206} has shown that if in addition, $f'\in H^1$ then $f$ generates
$W^+$. Hedberg \cite{MR0219729} and Lisin \cite{MR0229836} have shown (independently) that if $f$ is univalent
and $\sum n|a_n|^2(\log n)^{1+\varepsilon}<\infty$ for some positive $\varepsilon$, then $f$ generates $W^+$.
Neither Newman's condition nor Hedberg-Lisin's condition implies the
other. Is univalentness enough?\\
\textit{(L. Zalcman)}\\

\noindent\textbf{Update 5.25}
No progress on this problem has been reported to us.\\

\noindent\textbf{Problem 5.26} Let $B$ be the \textit{Bergman space of square integrable functions} in
$\mathbb{D}$, that is, those functions $f(z) = \sum^\infty_{n=0} a_nz^n$ for which $\sum^\infty_{n=1} n^{-1}|a_n|^2<\infty$.
A subspace $S$ is said to be \textit{invariant} if $zf \in S$ whenever $f \in S$. What
are the invariant subspaces of $B$? The corresponding problem for $H^2$
was solved by Beurling \cite{MR0027954} and uses the inner-outer factorisation
of $H^p$ functions, a tool unavailable in the present context.\\
\textit{(L. Zalcman)}\\

\noindent\textbf{Update 5.26}
No progress on this problem has been reported to us.\\

\noindent\textbf{Problem 5.27} \textbf{(The corona conjecture)} Let $D$ be an arbitrary domain in the
plane that supports non-constant bounded analytic functions. Suppose
that $f_1(z),\ldots, f_n(z)$ are bounded and analytic in $D$ and satisfy
\[\sum^n_{\nu=1}|f_\nu(z)|\geq\delta>0\]
in $D$. Can one find bounded analytic functions $g_\nu(z)$ in $D$ such that
\[\sum^n_{\nu=1}f_\nu(z)g_\nu(z)\equiv1\]
in $D$?\\
When $D$ is a disc, Carleson \cite{MR0141789} proved that the answer is `yes', and the
result extends to finitely connected domains. The result is also known to
be true for certain infinitely connected domains (see Behrens \cite{MR0256166}, Gamelin \cite{MR0276742}), but
false for general Riemann surfaces of infinite genus (Cole, unpublished). Presumably
the answer for the general plane domain is negative. Proofs of all
positive results depend on Carleson's theorem.\\
\textit{(L. Zalcman)}\\

\noindent\textbf{Update 5.27}
A positive answer for certain classes of domains has been given by
Carleson \cite{MR730079} and by Jones \cite{MR597809}.\\
It has also proved for domains $D$ whose complement lies on the real axis,
by Garnett and Jones \cite{MR793236}.\\

\noindent\textbf{Problem 5.28} Let $f$ be continuous in $\overline{\mathbb{D}}$ and analytic in $\mathbb{D}$. Let
\[\omega(f,\delta)=\sup|f(z)-f(w)|,\hspace{1cm}\text{for }|z-w|\leq\delta\text{ and }z, w\in \mathbb{D},\]
\[\tilde{\omega}(f,\delta)=\sup|f(z)-f(w)|,\hspace{1cm}\text{for }|z-w|\leq\delta\text{ and }|z|=|w|=1.\]
Is it true that
\[\lim_{\delta\to0}\frac{\omega(f,\delta)}{\tilde{\omega}(f,\delta)}=1?\]
The Ruben, Shields and Taylor (unpublished) have shown that there is an
absolute constant $C$ such that
\[\omega(f,\delta)\leq C\tilde{\omega}(f,\delta),\]
but that one may not take $C = 1$.\\
\textit{(L. A. Rubel, A. Shields)}\\

\noindent\textbf{Update 5.28}
A negative answer was given by Rubel, Shields and Taylor \cite{MR0412395}.\\

\noindent\textbf{Problem 5.29} A \textit{$G_\delta$ set} is a subset of a topological space that is a countable intersection of open sets. Let $E$ be a $G_\delta$ set of measure zero on $|z|=1$. Then does there 
exist an $f$ in $H^{\infty}$, $f\neq0$ such that $f = 0$ on $E$, and every point of the unit
circle is a Fatou point of $f$?\\
\textit{(L. A. Rubel)}\\

\noindent\textbf{Update 5.29}
An affirmative solution has been given by Danielyan \cite{MR3525402}.\\

%\section* *{Bloch functions}

\noindent\textbf{Problem 5.30} Let $\mathcal{B}$ be the space of Bloch functions, that is the space of
functions analytic in $\mathbb{D}$ with
\[\|f\|_\mathcal{B}=|f(0)|+\sup_\mathbb{D}(1-|z|^2)|f'(z)|<\infty.\]
Let $\mathcal{B}_S$ be the space of functions of the form
\begin{equation}
\label{B5.11}
f(z) = \log g'(z),\hspace{1cm} g \in S
\end{equation}
where $S$ is the class of functions as in Problem 5.23(a). Let $\mathcal{B}_Q$ be the space of
functions $g$ in $S$ that have a quasi-conformal extension to the closed plan, see Anderson, Clunie and Pommerenke \cite{MR0361090}.
\begin{enumerate}[(a)]
\item Is $\mathcal{B}_S$ connected in the norm topology?
\item Is $\mathcal{B}_Q$ dense in $\mathcal{B}_S$ in the norm topology?
\end{enumerate}
\textit{(L. Bers)}\\

\noindent\textbf{Update 5.30}
For part $(b)$ Becker points out that the answer is `no'. If $f=\log g'\in\mathcal{B}_S$, then define a new norm by $\|f\|=\sup\big\{(1-|z|^2)^2|S_g(z)|:|z|<1\big\}$, where $S_g$ is the Schwarzian of $g$. It is easy to show that
\begin{equation}
\label{D5.28.1}
\|f_n-f\|_B\to0\hspace{0.5cm}\implies\hspace{0.5cm}\|f_n-f\|\to0
\end{equation}
(compare with Becker \cite{MR623464}, Lemma 6.1). Gehring \cite{MR0499134} disproved the Bers conjecture, which means that $\mathcal{B}_Q$ is not dense in $\mathcal{B}_S$ with respect to $\|\cdot\|$. By (\ref{D5.28.1}) it follows immediately that $\mathcal{B}_Q$ is not dense in $\mathcal{B}_S$ with respect to $\|\cdot\|_\mathcal{B}$. Baernstein notes that it is an interesting open problem to determine geometrically the closure of $\mathcal{B}_Q$.\\

\noindent\textbf{Problem 5.31} It was shown by Becker \cite{MR0299780} that
\[\big\{f:\|f\|_\mathcal{B}<1\big\}\subset \mathcal{B}_Q.\]
Is the radius 1 best possible? Is it true that for $f \in \mathcal{B}_S$,
\[ \limsup_{|z|\to1}\,(1 - |z|^2)|f'(z)| < 1\hspace{0.5cm}\implies\hspace{0.5cm} f \in\mathcal{B}.\]
\textit{(Ch. Pommerenke)}\\

\noindent\textbf{Update 5.31}
Campbell \cite{Campbell} showed that the answer to the second part is `no'. The reason is that one of the assumptions has been omitted, namely that the image domain of the function $g(f=\log g')$ must be a Jordan domain. However, if the image domain of the functions $g(f=\log g')$ is a Jordan domain, then even the condition
\[\limsup_{|z|\to1}\Big|(1-|z|^2)zf'(z)-c\Big|<1,\]
where $|c|<1$, is enough to imply that $f$ is in $\mathcal{B}_Q$. Becker's proof of this latter statement consists in noting that \[g(z,t)=g(e^{-t}z)+(1-c)^{-1}(e^t-e^{-t})zg'(e^{-t}z)\] 
is a subordination chain for small $t$, hence univalent for small $t$ with $g(z,t)$ giving a quasiconformal extension of $g(z,0)$ onto a larger disc.\\

\noindent\textbf{Problem 5.32} Suppose that $f_n\in\mathcal{B}_S$. What does $\|f_n-f\|_\mathcal{B}\to0$ as $n\to\infty$ mean
geometrically for the functions $g_n$ related to $f_n$ by (\ref{B5.11})?\\
\textit{(Ch. Pommerenke)}\\

\noindent\textbf{Update 5.32}
No progress on this problem has been reported to us.\\

\noindent\textbf{Problem 5.33} Let $L$ be a regular triangular lattice in the plane. Let $f(z)$ map
$\mathbb{D}$ onto the universal covering surface over the complement of $L$. Is it true that the coefficients $a_n$ of $f$ tend to $0$ as $n\to\infty$?\\
\textit{(Ch. Pommerenke)}\\

\noindent\textbf{Update 5.33}
See Update 5.7.\\
An affirmative answer was given by Hayman, Pommerenke and Patterson \cite{MR0486496}. They proved that $a_n=O\big((\log n)^{-1/2}\big)$ in this case.\\

\noindent\textbf{Problem 5.34} It was proved by Hall that every Bloch function has (possibly
infinite) angular limits on an uncountably dense subset of $|z | = 1$. Do
there always exist angular limits on a set of positive measure relative
to some fixed Hausdorff measure, such as logarithmic measure for
example?\\
\textit{(J. E. McMillan, Ch. Pommerenke)}\\

\noindent\textbf{Update 5.34}
No progress on this problem has been reported to us.\\

\noindent\textbf{Problem 5.35} Let $F$ be any discontinuous group of M\"obius transformations of
$\mathbb{D}$. Does there always exist a meromorphic function automorphic
with respect to $\Gamma$ and normal, i. e. such that
\[(1-|z|^2)\frac{|f'(z)|}{1+|f(z)|^2}\]
is bounded in $\mathbb{D}$?\\
\textit{(Ch. Pommerenke)}\\

\noindent\textbf{Update 5.35}
Pommerenke \cite{MR0357779} has shown that there does exist a character-automorphic normal
function $f$, i.e. one for which $|f|$ is invariant under $\Gamma$. The original question for a
normal $f$, which is invariant under $\Gamma$, remains open.\\

%\section* *{Power series with gaps}

\noindent\textbf{Problem 5.36} Let $(n_k)$ be a sequence of positive integers such that
\begin{equation}
\label{B5.12}
n_{k+1}>\lambda n_k,\text{ where }\lambda>1,
\end{equation}
and suppose that
\begin{equation}
\label{B5.13}
f(z)=\sum^\infty_{k=0}a_kz^{n_k}
\end{equation}
is analytic in $\mathbb{D}$. Is it true that if $\sum^\infty_{k=0}|a_k|=\infty$, then $f(z)$ assumes
every finite value
\begin{enumerate}[(a)]
\item at least once;
\item infinitely often;
\item in every angle $\alpha<\arg z<\beta$ of $| z | < 1$?
\end{enumerate}
See Weiss and Weiss \cite{MR0145052}.\\
\textit{(J. P. Kahane)}\\

\noindent\textbf{Update 5.36}
Murai \cite{zbMATH03682094} has proved that $f(\mathbb{D})\neq\mathbb{C}$ implies $\sum_{k=0}^\infty|a_k|<\infty$. Hence (a) is true. Sons points out that the best known result is due to I-Lok Chang \cite{MR0315097} who showed that $\sum^\infty_{k=0}|a_k|^{2+\varepsilon}=\infty$ for some positive $\varepsilon>0$ implies infinitely many zeros in any sector.\\

\noindent\textbf{Problem 5.37} Suppose that $f(z)$ is a function as in (\ref{B5.13}) and define
\[\mu=\limsup_{r\to1}\frac{\log\log M(r,f)}{-\log(1-r)},\]
where $M(r,f)$ is the maximum modulus of $f(z)$ on $|z| = r$. We do not
now assume (\ref{B5.12}), but let $N^0(t)$ be the number of $n_k$ not greater than $t$.
If $N(r, a)$ is the function of Nevanlinna theory (see Chapter 1) it is known
that
\[\limsup_{r\to1}\frac{N(r,0)}{\log M(r,f)}=1\]
provided that either
\begin{enumerate}[(a)]
\item $\mu > 0$ and
\[\liminf_{k\to\infty}\frac{\log(n_{k+1}-n_k)}{\log n_k}>\frac{1}{2}\Big(\frac{2+\mu}{1+\mu}\Big),\]
(this is implicit in Wiman, see Sunyer and Balaguer \cite{MR0048559}), or
\item $\mu>\frac{1-\beta}{\beta}$ and $N^0(t)=O(t^{1-\beta})$ as $t\to\infty$, where $0<\beta<1$, see Sons \cite{MR0203030}.
\end{enumerate}
If $0<\mu<\frac{1-\beta}{\beta}$ with $N^0(t)=O(t^{1-\beta})$ as $t\to\infty$ we ask $(a), (b)$ and
$(c)$ of the preceding problem, at least for those cases not covered above.\\
In particular we may consider the cases $n_k = [k^\alpha]$, where $1 < \alpha < \frac{3}{2}$.\\
\textit{(L. R. Sons)}\\

\noindent\textbf{Update 5.37}
The best result on this problem so far is by Nicholls and Sons \cite{MR620424}.

\subsection*{Functions in the unit disc}

For the rest of this chapter, and for Chapter 6, the notation $f$ in $S$ means that $f$ is analytic and univalent in $\mathbb{D}$ with
\[f(z)=z+\sum^\infty_{n=2}a_nz^n.\]

\noindent In the next three questions $\prec$ means `is subordinate to'. If $f$ and $g$ are analytic functions in $\mathbb{D}$, then $g$ is \textit{subordinate to} $f$ if $g=f\circ \phi$ where $\phi$ is analytic in $\mathbb{D}$ and satisfies $\phi(\mathbb{D})\subset\mathbb{D}$, $\phi(0)=0$. See also Pommerenke \cite[Chp. 2]{MR0507768}.\\

\noindent\textbf{Problem 5.38} Tao-Shing Shah \cite{MR0100094} has shown that if $g(z)\prec f(z)$ in $\mathbb{D}$, $g'(0)/f'(0)$
is real, and $f$ in $S$ then
\begin{equation}
\label{C5.1}
|g(z)|\leq|f(z)|\text{ for }|z|\leq\frac{1}{2}(3-\sqrt{5}),
\end{equation}
\begin{equation}
\label{C5.2}
|g'(z)|\leq|f'(z)|\text{ for }|z|\leq3-\sqrt{8}.
\end{equation}
Both constants are `best-possible'. Shah's proofs are technically very involved and
it would be nice to have simpler proofs. Goluzin \cite{MR0043209} gave simpler proofs but with worse
constants in each case. His methods appear to be incapable of yielding (\ref{C5.1}) and
(\ref{C5.2}).\\ 
\textit{(P. L. Duren)}\\

\noindent\textbf{Update 5.38}
Shah's \cite{MR0100094} majorisation results were generalised by Campbell (\cite{MR0299769}, \cite{MR0315120} and \cite{MR0349987}) as well as similar questions by Biernacki, MacGregor and Lewandowski. The proper setting for such questions is not univalent function theory but locally univalent functions of finite order. In particular, Shah's sharp results of $|g(z)|\leq|f(z)|$ in $(3-\sqrt{5})/2$ and $|g'(z)|\leq |f'(z)|$ in $3-\sqrt{8}$ hold for all functions in $U_2$, the \textit{universal linear invariant family} of order $2$, which contains $S$ as a proper subclass, and which contains functions of infinite valence.\\

\noindent\textbf{Problem 5.39} Goluzin \cite{MR0043209} has shown that, if $g(z)\prec f(z)$ in $\mathbb{D}$, then
\[M_2(r,g')\leq M_2(r,f'),\hspace{1cm}0\leq r\leq\frac{1}{2}.\]
Here, for $p > 0$,
\[M_p(r,h)=\Big(\frac{1}{2\pi}\int^{2\pi}_0|h(re^{i\theta})|^p\,d\theta\Big)^{1/p}.\]
The result is not necessarily true if $r > \frac{1}{2}$, as $f(z) = z, g(z) = z^2$ shows, though it
follows from a theorem of Littlewood that, for all $p$,
\[M_p(r,g)\leq M_p(r,f),\hspace{1cm} 0 < r < 1.\]
Find the largest number $r_p, 0 < r_p < 1$, independent of $f$ and $g$ so that
\[M_p(r,g')\leq M_p(r,f'),\hspace{1cm} 0<r<r_p,\]
if $g\prec f$.\\
Note: If $g(z) \prec f(z)$ in $\mathbb{D}$, then $g(z) = f(\phi(z))$ so that
\[|g'(z)|\leq |f'(\phi(z))|,\hspace{1cm}|z|\leq\sqrt{2}-1,\]
(see Carath\'eodory \cite[p. 19]{MR0064861}). Thus if $h(z) = f'(\phi(z))$ then $h\prec f'$ so that, by Littlewood's
theorem, $M_p(r,g') \leq M_p(r,h) \leq M_p(r,f')$ for $p > 0, r \leq\sqrt{2}-1$. This improves
what one gets if one applied (\ref{C5.2}) to the $p$-th means since $\sqrt{2}-1>3-\sqrt{8}$; but it
may not be best possible.\\ 
\textit{(P. L. Duren)}\\

\noindent\textbf{Update 5.39}
No progress on this problem has been reported to us.\\

\noindent\textbf{Problem 5.40} Suppose that $f(z) = \sum^\infty_0 a_nz^n$ and that $F(z)$ is analytic in $\mathbb{D}$, with $f\prec F$. What non-trivial conditions on $F$ imply that
\begin{equation}
\label{C5.3}
a_n\to0\hspace{1cm}\text{ as }n\to\infty?
\end{equation}
In particular, is (\ref{C5.3}) implied by
\begin{equation}
\label{C5.4}
F\text{ is a Bloch function, }\textit{and}
\end{equation}
\begin{equation}
\label{C5.5}
\int^{2\pi}_0|f'(re^{i\theta})|^2\,d\theta=o\Big(\frac{1}{1-r}\Big)^2,\hspace{1cm}\text{ as }r\to1^-\,?
\end{equation}
Is it true that (\ref{C5.5}) by itself is preserved under subordination?\\
It is known that (\ref{C5.3}) holds if
\[(1-r)|F'(re^{i\theta})|\to0\]
uniformly in $\theta$ as $r\to1$ and that (\ref{C5.3}) holds if both (\ref{C5.4}) is satisfied, and given any positive $\varepsilon$,
$F$ can be written in the form
\[F(z)=F_1(z)+F_2(z)\]
where $(1-|z|^2)|F_1'(z)|\leq\varepsilon$ and $F_2$ has bounded characteristic in $\mathbb{D}$.\\
\textit{(W. K. Hayman)}\\

\noindent\textbf{Update 5.40}
The problem asked whether various conditions on a function $F$ might imply that $(\ref{C5.3})$
above holds when $f$ is subordinate to $F$. The examples of Fernandez \cite{MR769161} show that the
answer is `no'. In particular the condition
\begin{equation}
\label{K(c)}
\int^{2\pi}_0|F'(re^{i\theta})|^2\,d\theta=\frac{o(1)}{(1-r)^2}
\end{equation}
is not preserved under subordination, even if $F$ is a Bloch function. For instance, if
$D$ is a domain such that the distance $d(w)$ of $w$ from the complement of $D$ tends to zero
as $w\to\infty$, then $(\ref{C5.3})$ and $(\ref{K(c)})$ hold for the functions $F$ mapping $\mathbb{D}$ into the universal
covering surface over $D$, by a theorem of Hayman, Patterson and Pommerenke \cite{MR0486496},
but by Fernandez's examples neither $(\ref{C5.3})$ nor $(\ref{K(c)})$ hold for the functions $f$ analytic in
$\mathbb{D}$ and with values in $D$, if the complement of $D$ is countable or more generally has
capacity zero.\\

\noindent\textbf{Problem 5.41} Let $k$ be a positive integer, and define
\[\Phi_k=\{\phi:\phi\in C[0,2\pi],\phi\uparrow,\phi(0)=0,\phi(2\pi)=2k\pi\}.\]
Let $e^{i\phi(t)}$ with $\phi\in\Phi_k$ have the Fourier expansion
\[\sum^\infty_{-\infty}C_ne^{int}.\]
\begin{enumerate}[(a)]
\item Shapiro conjectures that 
\begin{equation}
\label{C5.6}
\sum^k_{n=0}|C_n|^2\geq\delta_k,
\end{equation}
where $\delta_k$ is a positive constant depending only on $k$. (For $k=1$ this is known,
and it can be proved by combining results of Heinz (uncited) and Rado (uncited).
Inequalities of this type arise in studying the distortion of harmonic
mappings.)
\item Conceivably (\ref{C5.6}) holds even with $\delta_k\geq\delta>0$ for some absolute constant $\delta$.
In any case Shapiro conjectures also that, given only that $\phi\in C(0,2\pi)$ and $\phi\uparrow$ then
\[\sum^\infty_{n=0}|C_n|^2\geq\delta>0.\]
\item Is it true that, for $\phi\in\Phi_k$,
\[\sum^k_{n=1}|C_n|^2\geq\sum^{-1}_{n=-k}|C_n|^2\,?\]
(This is known to be true for $k = 1$). 
\end{enumerate}
\textit{(H. S. Shapiro)}\\

\noindent\textbf{Update 5.41}
Sheil-Small \cite{MR803610} has proved this conjecture. The case  $k=2$ had been previously settled by Hall \cite{MR657497}. Hall \cite{MR848860} also proved that $\delta_k\to0$ as $k\to\infty$, and gave an example for which
\begin{equation}
\label{(d)}
\sum^k_{n=1}|C_n|^2\leq\frac{1}{k^2}\sum^{-1}_{n=-k}|C_n|^2,
\end{equation}
thus answering in the negative two other questions of Shapiro. He conjectured that
the factor $k^{-2}$ in (\ref{(d)}) is extreme.\\

\noindent\textbf{Problem 5.42} Suppose that
\[f(z)=\sum^\infty_{n=0}a_nz^n\]
is analytic in $\mathbb{D}$, with
\[\sum^\infty_{n=0}|a_n|=1,\hspace{1cm} |f(z)|\geq\delta>0\text{ in }\mathbb{D},\hspace{1cm}\text{ and}\hspace{1cm}\frac{1}{f(z)}=\sum^\infty_{n=0}b_nz^n.\]
The following facts are known about $M=\sum^\infty_{n=0}|b_n|$:
\begin{enumerate}[(a)]
\item $M < +\infty$;
\item if $\delta<\frac{1}{2}, M$ cannot be bounded in terms of $\delta$ (Katznelson (no citation)).
\item if $\delta>2^{-\frac{1}{2}}, M$ can be bounded in terms of $\delta$ (Katznelson (no citation), Newman (no citation)).
\end{enumerate}
What is the infimum of those $\delta$ such that $M$ can be bounded in terms of $\delta$? (A likely
guess is $\frac{1}{2}$).\\ 
\textit{(H. S. Shapiro)}\\

\noindent\textbf{Update 5.42}
No progress on this problem has been reported to us.\\

\noindent\textbf{Problem 5.43} Determine the Laurent coefficient bodies for analytic functions taking values
of modulus at most unity in a given annulus
\[A_r = \{z : r < |z| < 1\}.\]
Determine the extremal functions.\\ 
\textit{(M. Heins)}\\

\noindent\textbf{Update 5.43}
No progress on this problem has been reported to us.\\

\noindent\textbf{Problem 5.44} Suppose that $0 < \alpha < 1$ and
\[\frac{(1+xz)^\alpha}{1-z}=\sum^\infty_{n=0}A_n(x)z^n,\hspace{1cm} A_0(x)=1,\hspace{1cm} |x|=1.\]
Is it true that
\[|A_{2n+1}(x)|\leq|A_{2n+1}(1)|,\hspace{1cm}n\geq2,\hspace{1cm}|x|=1?\]
The above is true for $n = 1$; the corresponding result is false for all \mbox{$A_{2m}(x)$} with $m \geq 1$, see Brannan \cite{MR0412411}.
More generally, one can ask the same question for the coefficients of
\[(1+xz)^\alpha(1-z)^{-\beta}, \hspace{1cm}|x|=1, \hspace{1cm}\alpha>0,\hspace{1cm}\beta>0.\]
Here it is not even known if $|A_3(x)| \leq A_3(1)$.\\ 
\textit{(D. A. Brannan)}\\

\noindent\textbf{Update 5.44}
No progress on this problem has been reported to us.\\

\noindent\textbf{Problem 5.45} If $A$ is any analytic subset of the Riemann sphere it was shown by Kierst \cite{Kierst}
that $A$ is (exactly) the set of asymptotic values of a function meromorphic in $\mathbb{D}$. If $\infty\in A$, Kierst also proved that $A$ is the set of asymptotic values of a
function $f$ analytic in $\mathbb{D}$. However, there exist analytic sets which are not the set
of asymptotic values of a function analytic in $\mathbb{D}$. Ryan \cite{MR0203034} has characterised
those subsets of the Riemann sphere which are the set of asymptotic values of a
function analytic in $\mathbb{D}$. Can one find a simpler characterisation?\\ 
\textit{(K. F. Barth)}\\

\noindent\textbf{Update 5.45}
No progress on this problem has been reported to us.\\

\noindent\textbf{Problem 5.46} A non-constant function $f$ analytic in $\mathbb{D}$, is said to be in
the \textit{MacLane class} $\mathcal{A}$ if the set of points of the unit circle $\mathbb{T}$ at which $f$ has an
asymptotic value is a dense subset of $\mathbb{T}$ (see MacLane \cite{MR0148923}). A function $f$ is said to have an
\textit{arc tract} if there exists a sequence $\{\gamma_n\}$ of arcs, $\gamma_n\subset \mathbb{D}$, and a non-degenerate subarc
$\gamma$ of $\mathbb{T}$ such that $\gamma_n\to\gamma$ (in the obvious fashion) and $\min\{|f(z)| :z\in\gamma_n\}\to\infty$ as
$n\to\infty$. Does there exist a function $f\in\mathcal{A}$ with an arc tract and with non-zero
derivative (see MacLane \cite[p. 281]{MR0274765})?\\ 
\textit{(K. F. Barth)}\\

\noindent\textbf{Update 5.46}
No progress on this problem has been reported to us.\\

\noindent\textbf{Problem 5.47} Let
\[f(z)=\sum^\infty_{n=0}a_nz^{\lambda_n}\]
be analytic in $\mathbb{D}$ and have Hadamard gaps, i.e. $\lambda_{n+1}/\lambda_n\geq q>1$. Need $f$ have
any radial limits (finite or infinite)? If not, need $f$ have any asymptotic value on a
path ending at a single point? For the case $q\geq3$, see MacLane \cite{MR0148923}.\\
\textit{(J. M. Anderson and R. Hornblower)}\\

\noindent\textbf{Update 5.47}
Murai \cite{MR692343} has shown that, for $q>1$, a Hadamard function has the
asymptotic value $\infty$ if the coefficients are unbounded. It is still not known whether
every such function has any radial limits whatsoever.\\
There are several partial results: Gnuschke and Pommerenke (see \cite[Theorem 1]{zbMATH03908769}) have
shown that, if $q > 1$ and $f$ is unbounded and if furthermore
\[|a_k|/(|a_1|+\ldots+|a_k|)\to0\hspace{1cm}\text{ as }k\to\infty,\]
then $\text{Re}\, f$ has the angular limit $+\infty$ on a set of positive Hausdorff dimension: this
improves a (proof and) result of Csordas, Lohwater and Ramsey \cite{MR674281}.
Gnuschke and Pommerenke (see \cite[Theorem 3]{zbMATH03908769}) have also proved that there exists a
Hadamard function $f$ (with $q=\frac{33}{32}$) such that $\text{Re}\, f$ and $\text{Im}\, f$ oscillate between
$-\infty$ and $+\infty$ on every radius. This does not exclude the possibility that $f$ spirals to
$\infty$ on some radii. Furthermore, there is a result of Hawkes \cite{MR576895} that, for instance,
$f(z) = \sum z^{2^k}$ has the radial limit $\infty$ on a set of Hausdorff dimension one.\\

\noindent\textbf{Problem 5.48} Let
\[f(z)=\sum^\infty_{k=1}a_{n_k}z^{n_k}\]
be analytic in $\mathbb{D}$ and have Hadamard gaps. Characterise sets $S$ in $\mathbb{D}$ with the
property that if $f$ is bounded on $S$ then $f(z)$ is bounded for $z\in\mathbb{D}$. \\
\textit{(K. G. Binmore)}\\

\noindent\textbf{Update 5.48}
No progress on this problem has been reported to us.\\

\noindent\textbf{Problem 5.49} What kind of gaps can the Taylor expansion of a non-constant automorphic
function have? For example, can it have Hadamard gaps? (Presumably the sharp 
answer would depend on the group concerned.) This is closely related to a theorem of
R\'enyi (no citation) that a non-constant periodic entire function cannot have more than half of
its coefficients zero.\\ 
\textit{(L. A. Rubel)}\\

\noindent\textbf{Update 5.49}
Nicholls and Sons \cite{MR620424} have shown that $f(z)=\sum^\infty_{k=0}c_kz^{n_k}$ cannot have 
\[n^{-1}_k(n_{k+1}-n_k)\log n_k\to\infty.\] 
Furthermore, if $f$ is automorphic and of the above form, and 
\[N_0(t)\equiv\max\{k:n_k<t\},\] 
then $N_0(t)\neq o(\log\log t)$ as $t\to\infty$.\\
Hwang \cite{MR726986} has proved that, if $f$ is automorphic with respect to
any group containing a parabolic element, then $f$ does not have Hadamard gaps.
The proof that this holds for any Fuchsian group appears to be false.\\

\noindent\textbf{Problem 5.50} A function $f$ analytic in $\mathbb{D}$, is said to be \textit{annular} if there
exists a sequence $\{J_n\}$ of Jordan curves in $\mathbb{D}$ such that 
\begin{enumerate}[(a)]
\item $J_n$ lies in the inside of $J_{n+1}$, 
\item for each positive $\varepsilon$ there exists a number $N(\varepsilon)$ such that for if $n > N(\varepsilon)$, $J_n$
lies in the domain $\{z: 1 - \varepsilon < |z| < 1\}$, 
\item $\min\{|f(z)| : z \in J_n\}\ \to\infty$ as $n\to\infty$.
\end{enumerate}
Let $\mathbb{T}$ denote the unit circle and let $Z'(f)$ denote the set of limit points of the zeros of an
annular function $f$. Write 
\[S(f) = \{a:a\in\mathbb{C}\text{ and }Z'(f-a) \neq T\},\] 
and let $|S(f)|$ denote the cardinality of this set. Can $|S(f)| = \mathfrak{N}_0$? The cases $|S(f)| = 1$ (Barth
and Schneider \cite{MR0269834}) and $|S(f)| = 2$ (Osada (no citation)) are known, however neither construction
can be easily adapted to the general case.\\
\textit{(K. F. Barth and D. D. Bonar)}\\

\noindent\textbf{Update 5.50}
Barth observes that Carroll \cite{MR555223} has constructed a strongly annular function $f$ such that the set of singular values of $f$ is countably infinite. An analytic function $f$ in $\mathbb{D}$ is said to be \textit{strongly annular} if there exists a sequence $\{r_n\}$, $0<r_n<1$ such that the minimum modulus of $f$ on the circle $|z|=r_n$ tends to infinity as $n\to\infty$.\\

\noindent\textbf{Problem 5.51} It can be shown that there exists a Blaschke product $B(z)$ with \mbox{$B(0) = 0$}
such that
\[\frac{1+B(z)}{1-B(z)}\]
is a Bloch function. Give an explicit construction for such a product.\\ 
\textit{(F. Holland)}\\

\noindent\textbf{Update 5.51}
No progress on this problem has been reported to us.\\

\noindent\textbf{Problem 5.52} Let F be a Fuchsian group in $\mathbb{D}$, and let $B$ be a set on $\mathbb{T}$ such that
$B\cap \gamma(B)=\emptyset$ for $\gamma\in\Gamma,\gamma\neq I$. (This is the case, for instance, if
\[B=(\partial F\cap\mathbb{T})\setminus C\hspace{1cm}(C\text{ countable})\]
where $F$ is a normal fundamental domain.) If cap $B > 0$, does it follow that $F$ is of
convergence type? (Conversely, it is known that, for every group $F$ of convergence
type, there exists such a set $B$ with cap $B > 0$).\\
\textit{(Ch. Pommerenke)}\\

\noindent\textbf{Update 5.52}
No progress on this problem has been reported to us.\\

\noindent\textbf{Problem 5.53} The ratio $R$ of two Blaschke products $B(z,a_n), B(z,b_n)$ is of bounded
characteristic, but need not be a normal meromorphic function. That is, if $a_n$ is `close' to
$b_n$ for infinitely many $n$ we can arrange that the spherical derivative of $R$ is too large
for $R$ to be normal. When is $R$ a normal function? Cima and Colwell \cite{MR0227423}  have shown
that if $\{a_n\}$ and $\{b_n\}$ are both interpolating sequences, then $R$ is normal if and
only if $\{a_n\} \cup \{b_n\}$ is also an interpolating sequence. What happens if the sequences
are not interpolating sequences?\\ 
\textit{(J. M. Anderson)}\\

\noindent\textbf{Update 5.53}
No progress on this problem has been reported to us.\\

\noindent\textbf{Problem 5.54} Suppose that
\[f(z)=\sum^\infty_{n=0}a_nz^n\]
is convergent in $\mathbb{D}$, that $|z_0| = 1$, and that neither of the series
\[\sum^\infty_{n=0}(\text{Re}\, a_n)z_0^n,\hspace{1cm}\sum^\infty_{n=0}(\text{Im}\, a_n)z_0^n,\]
is absolutely convergent. Let $S$ be the set of complex values assumed by the series
\[\sum^\infty_{n=0}\varepsilon_n a_n z_0^n,\]
where $\varepsilon=\pm1$, and is allowed to vary over all possible choices. Is it true that $S = \mathbb{C}$?\\
\textit{(A. C. Offord; communicated by J. G. Clunie)}\\

\noindent\textbf{Update 5.54}
Partial results have been obtained by Jakob and Offord \cite{MR709720}. In particular, if $R_n=\sum^\infty_n|a_n|^2$ is such that
\[\sum n^{-1}(\log n)^{-1/2}R_n^{1/2}\]
is convergent, then almost all the functions $\sum^\infty_0\pm a_nz^n$ will be bounded. Related results are due to Salem and Zygmund \cite{MR0065679}. Notice that $\sum|a_n|^2<\infty$ is in itself not sufficient to secure boundedness.\\

\noindent\textbf{Problem 5.55} Let the function $h(z)$ be analytic in $\mathbb{D}$, $G_h = \{h(\mathbb{D})\}$, $K_h$ a compact
subset of $G_h$. The function $h(z)$ will be said to have the \textit{$L$-property} (`$h(z_n)$ leaves
the range of $h$') on the sequence $\{z_n\}^\infty_1$ with $\lim |z_n|=1$, if only finitely many points
of $\{h(z_n)\}^\infty_1$ lie in $K_h$ for every $K_h$. We will define two functions $f,g$ analytic in $\mathbb{D}$
to be an \textit{ordered $L$-pair} if, on any sequence for which $f$ has the $L$-property, $g$ also has
the $L$-property. A non-constant function $\alpha(z)$ analytic in $\mathbb{D}$ will be called an
\textit{$L$-atom} if, corresponding to any function $f$ such that $f, \alpha$ is an ordered $L$-pair, there
exists a function $\phi_f$ analytic on $G_\alpha$ such that $f=\phi_f\circ\alpha$. Prove or disprove the conjecture
that $\alpha$ is an $L$-atom if and only if it is univalent. What happens if $\mathbb{D}$ is
replaced by more general domains?\\ 
\textit{(L. A. Rubel)}\\

\noindent\textbf{Update 5.56}
No progress on this problem has been reported to us.\\

\noindent\textbf{Problem 5.56} By a \textit{first-order property} we shall mean a ring-theoretic property in the
first-order predicate calculus. For functions $f$ analytic in $\mathbb{D}$, are the following
first-order properties:
\begin{enumerate}[(a)]
\item $f$ is constant?
\item $f$ is bounded ?
\item $f$ is admissible, i.e. $T(r,f) \neq O(\log(1-r)^{-1})?$
\end{enumerate}
What is the situation in more general domains?\\ 
\textit{(L. A. Rubel)}\\

\noindent\textbf{Update 5.57}
No progress on this problem has been reported to us.\\

\noindent\textbf{Problem 5.57} Suppose that $f$ is analytic in $\mathbb{D}$. Plessner's theorem \cite{MR1581142} asserts that
at almost all points $e^{i\theta}$ on the unit circle, either $f$ has an angular limit, or else the
image $f(S)$ of every Stolz angle $S$ with vertex at $e^{i\theta}$ is dense in $\mathbb{C}$. How much can
`dense' be improved? In particular, is it true that at almost all $e^{i\theta}$ either $f$ has an
angular limit, or else $f(S)$ is all of $\mathbb{C}$, except perhaps for a set of zero logarithmic
capacity? This result would be best possible, since if $E$ is any closed set of capacity
zero then the universal covering map of the disc onto $\mathbb{C} \setminus E$ has angular limits almost
nowhere.\\ 
\textit{(A. Baernstein II)}\\

\noindent\textbf{Update 5.57}
No progress on this problem has been reported to us.\\

\noindent\textbf{Problem 5.58} Suppose that $f$ is univalent and zero-free in $\mathbb{D}$. It has
been shown by Baernstein \cite{baernstein1976} that, for each $p\in(0,\frac{1}{2})$, $f$ admits a factorisation $f=B_p(F_p)^{1/p}$, where
$B_p\in H^\infty, 1/B_p\in H^\infty$, and $\text{Re}\, F_p > 0$. Is it possible to pass to the limit $p = \frac{1}{2}$, and thus
factor $f$ into a bounded function times a function subordinate to a map onto a slit
plane?\\ 
\textit{(A. Baernstein II)}\\

\noindent\textbf{Update 5.58}
While this conjecture is still open,
Wolff \cite{MR1979771} has shown that a certain natural stronger conjecture is false; on the other
hand, the factorisation given in Problem 5.58 is possible when $f$ is a monotone slit mapping. A \textit{monotone slit mapping} is a function $f(z)$ which is analytic and univalent in $\mathbb{D}$ and for which $f(0)=0$ and $f'(0)=1$, whose image domain is the complement of a path $\Gamma(t)$ on $[0,\infty)$ for which $|\Gamma(t_1)|<\Gamma(t_2)|$ if $t_1<t_2$.\\

\noindent\textbf{Problem 5.59} \textbf{(Subordination and extreme point problem)} Let $g(z)=\sum^\infty_{n=0}B_nz^n$ be analytic in $\mathbb{D}$. Denote by $S_g$ the family of functions $f(z)$ subordinate to $g$. Find general conditions on $g$ so that the only extreme points of the closed convex hull of $S_g$ are the functions $g(ze^{it})(0\leq t<2\pi)$. This is known for certain functions $g$, see Clunie \cite{MR623468}, for example:
\begin{enumerate}[(a)]
\item $g(z)=[(1+cz)/(1-z)]^\alpha$ where $|c|\leq1,\alpha\geq1$;
\item $g(z)=\exp[(1+z)/(1-z)]$.
\end{enumerate}
Sheil-Small can prove it for a functions of the form \linebreak\mbox{$g(z)=h((1+z)/(1-z))$}, where $h(w)$ is a univalent quadratic polynomial in $\text{Re}\, w>0$. One might expect that the conclusion would hold  for a functions $g$ with positive coefficients increasing in a suitably regular manner.\\
\textit{(T. Sheil-Small)}\\

\noindent\textbf{Update 5.59}
No progress on this problem has been reported to us.\\

\noindent\textbf{Problem 5.60} \textbf{(Hadamard convolutions)} Suppose that $\alpha\geq1, \beta\geq1$ and that $\phi$ is analytic in $\mathbb{D}$, and satisfies
\[\phi(z)\ast\frac{(1+xz)^\alpha}{(1-z)^\beta}\neq0,\hspace{1cm}|x|=1, |z|<1.\]
Is it true that
\[\phi(z)\ast\frac{(1+xz)^{\alpha-1}}{(1-z)^\beta}\neq0,\hspace{1cm}|x|=1, |z|<1\,?\]
This is true when $\alpha$ is a natural number. For the proof and further background, see Sheil-Small \cite{MR623496}.\\
\textit{(T. Sheil-Small)}\\

\noindent\textbf{Update 5.60}
No progress on this problem has been reported to us.\\

\noindent\textbf{Problem 5.61} Let $w(z)$ be analytic in $\mathbb{D}$ with $w(0)=0$. If \mbox{$|w(z)+zw'(z)|<1$}, for $|z|<1$, then a simple application of Schwarz's lemma shows that $|w(z)|<1$, for $|z|<1$. Miller and Mocanu \cite{MR506307} showed that
\[|w(z)+zw'(z)+z^2w''(z)|<1,\hspace{1cm}|z|<1\]
implies that $|w(z)|<1$, for $|z|<1$. Is it true that
\[|w+zw'+z^2w''+\ldots+z^nw^{(n)}|<1\hspace{0.5cm}\implies\hspace{0.5cm}|w(z)|<1,\]
for $n=1, 2, 3, \ldots$?\\
\textit{(S. Miller)} \\

\noindent\textbf{Update 5.61}
This has been proved to be the case by Goldstein, Hall, Sheil-Small and Smith \cite{MR679930}.\\

\noindent\textbf{Problem 5.62} Let $u$ be a continuous real-valued function on the unit circle $\mathbb{T}$. Give a necessary and sufficient condition on $u$ such that $u$ is the real part of a function $f$ in the disc algebra $A(\overline{\mathbb{D}})$.\\
Remarks: 
\begin{enumerate}[(a)]
\item A solution would have applications in the algebraic ideal theory of $A(\overline{\mathbb{D}})$.
\item An answer to the analogous problem for $L^p(T), H^p(\mathbb{D})$ is the Burkholder-Gundy-Silverstein Theorem (see Peterson \cite[p.13]{Peterson}).
\end{enumerate}
\textit{(M. von Renteln)}\\

\noindent\textbf{Update 5.62}
No progress on this problem has been reported to us.\\

\noindent\textbf{Problem 5.63} One of the many equivalent norms on BMOA on $\mathbb{D}$ is defined by
\[\|f\|_h=\inf_q\sup_{z\in\mathbb{D}}|f(z)+\overline{q(z)}|,\]
the infimum being taken over all functions analytic in $\mathbb{D}$.\\
Given $f$ in BMOA, does there exist a $q$ such that
\[|f(z)+\overline{q(z)}|\equiv\|f\|_h\text{ a.e.  on }|z|=1?\]
The answer is `yes' when $f$ is a rational function.\\
\textit{(J. A. Hempel)}\\

\noindent\textbf{Update 5.63}
No progress on this problem has been reported to us.\\

\noindent\textbf{Problem 5.64} Let $f$ be analytic in $\mathbb{D}$ with
\begin{equation}
\label{D5.64.1}
|f(z)|=O((1-|z|)^{-k}),\hspace{1cm}k\geq0.
\end{equation}
Then $f$ induces a distribution on $C^\infty(T)$, as follows. For $\phi\in C^\infty(T)$
\[\lim_{r\to1}\Big(\frac{1}{2\pi}\int f(re^{i\theta})\phi(e^{-i\theta})\,d\theta\Big)=\Lambda_f (\phi).\]
What can be said about the order of the distributions satisfying (\ref{D5.64.1}) and
\[\lim_{|z|\to1}|f(z)|(1-|z|)^k=0\]
\textit{(J. A. Cima)}\\

\noindent\textbf{Update 5.64}
No progress on this problem has been reported to us.\\

\noindent\textbf{Problem 5.65} Does there exist a non-constant function $f$ in the disc algebra such that $f(e^{i\theta})\in f(\mathbb{D})$ for almost all $\theta$? Caution: the Rudin-Carleson theorem (see Bishop \cite{MR0133462}) allows construction of a good candidate, but it is not immediately clear whether it does work.\\
\textit{(K. Stephenson)}\\

\noindent\textbf{Update 5.65}
Stegenga and Stephenson (unpublished) have shown this to be the case for almost all $f$ in a category sense. A solution has
also been obtained by Gol'dberg \cite{MR760854} using a simple construction and without the use of
the Rudin-Carleson Theorem.\\

\noindent\textbf{Problem 5.66} Let $B$ be an infinite Blaschke product in $\mathbb{D}$. Does there exist a positive $\delta$, depending on $B$, such that, for every $w$, $|w|<\delta$, the set $B^{-1}(\{w\})$ is infinite? Stephenson has obtained some related results.\\
\textit{(K. Stephenson)}\\

\noindent\textbf{Update 5.66}
Stephenson (unpublished) has answered this
question negatively by constructing an inner function $f$ such that $f^{-1}(w)$ is finite on a
dense subset of $|w|< 1$. Then
\[f_\alpha=\frac{f-\alpha}{1-\overline{\alpha}f}\]
is a Blaschke product for almost all $\alpha$.\\

\noindent\textbf{Problem 5.67} Let the function $f$ in $\mathbb{D}$ be given by
\[f(z)=\sum^\infty_{k=0}a_kz^{n_k},\hspace{1cm}\frac{n_{k+1}}{n_k}\geq\lambda>1,\hspace{1cm} k\geq0,\]
with
\[m_0(r,f) \equiv \max_{k\geq0} |a_k|r^{n_k}\to\infty\hspace{1cm}\text{ as }r\to 1^-.\]
Define
\[E=\Big\{\theta:\liminf_{r\to1^-}\frac{|f(re^{i\theta})|}{m_0(r,f)}>0\Big\}.\]
Is it true that $E$ has measure $0$?\\
This has been proved for the case that $\mu(\frac{1}{2}(1+r))/m_0(r,f)\leq$ constant.\\
\textit{(D. Gnuschke and Ch. Pommerenke)}\\

\noindent\textbf{Update 5.67}
No progress on this problem has been reported to us.\\

\noindent\textbf{Problem 5.68} Let the function $f$ where
\[f(z)=1+\sum^\infty_{n=1}a_nz^n,\hspace{1cm}|z|\leq1,\]
be a Bloch function with positive real part in $\mathbb{D}$. Determine the rate of growth
(as $N\to\infty$) of the sequence $\big\{\sum^N_{n=1}|a_n|^2\big\}$.\\
\textit{(F. Holland)}\\

\noindent\textbf{Update 5.68}
No progress on this problem has been reported to us.\\

\noindent\textbf{Problem 5.69} Let the function $f$ where
\[f(z)=1+\sum^\infty_{n=1}a_nz^n,\hspace{1cm}|z|\leq1,\]
be a Bloch function with positive real part in $\mathbb{D}$ and such that each $a_n\geq0$. Does it
follow that $\sum^\infty_{n=1}a^2_n<\infty$?\\
An equivalent formulation of the problem is the following: let $\mu$ be a probability
measure in Zygmund's class $\Lambda_*$ on the circle, and let
\[\hat{\mu}(n)=\int^{2\pi}_0e^{-inx}\,d\mu(x)\geq0,\hspace{1cm}n\in\mathbb{Z}.\]
Is it true that $\hat{\mu}\in l_2$? A counter-example, if one exists, cannot be constructed using
Riesz products (see Duren \cite{MR0185057}, and Holland and Twomey \cite{MR0430651}).\\
An affirmative answer would mean that, if $f(z) = 1 + \sum^\infty_{n=1}a_nz^n$ is a Bloch
function in $\mathbb{D}$ with positive real part, then $\sum^\infty_{n=1}|a_n|^4<\infty$. Even if this last
inequality is false, perhaps it is still true in the general case that there exists $p > 4$,
such that $\sum^\infty_{n=1}|a_n|^p<\infty$.\\
\textit{(F. Holland)}\\

\noindent\textbf{Update 5.69}
No progress on this problem has been reported to us.\\

\noindent\textbf{Problem 5.70} Barth and Clunie \cite{MR660605} have constructed a bounded analytic function in $\mathbb{D}$ with a level set component of infinite length; this component is highly
branched. Can one construct a bounded analytic function with an unbranched level
set component of infinite length?\\
\textit{(K. F. Barth and J. G. Clunie)}\\

\noindent\textbf{Update 5.70}
No progress on this problem has been reported to us.\\

\noindent\textbf{Problem 5.71} Suppose that 
\[f(z) = \sum^\infty_{k=1}a_kz^{n_k}, \hspace{1cm}n_{k+1}/n_k\geq q>1,\] 
is an analytic function in $\mathbb{D}$ with Hadamard gaps, such that $T(r,f)\to\infty$ as $r\to1$. Does $\delta(w,f)=0$ hold for
every (finite) complex number $w$?\\ 
\textit{(T. Murai)}\\

\noindent\textbf{Update 5.71}
No progress on this problem has been reported to us.\\

\noindent\textbf{Problem 5.72} Let the function $f$ have the power series $f(z)=\sum^\infty_{n=0}a_nz^n$ of radius of
convergence $1$; let $E$ be the singular set on $\mathbb{T}$, and suppose that
\[\sup_{\xi\in E}\, \sup_{N\geq0}\Big|\sum^N_{n=0}a_n\xi\Big|<\infty.\]
It has been shown by Allan, O'Farrell and Ransford \cite{MR915430} that, if $E$ has measure zero, then $\sum^\infty_{n=0}a_nz^n$ converges to $f(z)$
at each point $z$ in $\mathbb{T}\setminus E$.\\
Does the conclusion remain true if $E$ has positive measure? (Nothing appears to
be known either way.)\\ 
\textit{(T. J. Ransford)}\\

\noindent\textbf{Update 5.72}
No progress on this problem has been reported to us.\\

\noindent\textbf{Problem 5.73} Let $0<\alpha<1$ and let $R_\alpha$ denote the set of all Riesz potentials $p(x)$ of finite
positive Borel measures $\mu$ on $\mathbb{R}$:
\[p(x)=\int^\infty_{-\infty}|x-t|^{-\alpha}\,d\mu(t).\]
Characterise those non-negative measurable functions $f(x)$ on $\mathbb{R}$ that are dominated by some
function $p$ in $R_\alpha$.\\ 
\textit{(B. Korenblum)}\\

\noindent\textbf{Update 5.73}
No progress on this problem has been reported to us.\\

\noindent\textbf{Problem 5.74} Characterise those non-negative measurable functions $f$ on the unit circle that
are dominated almost everywhere by moduli of the boundary values of an analytic function in
the unit disc with positive real part. (Problem 5.73 might be considered a step
towards solving this problem.)\\ 
\textit{(B. Korenblum)}\\

\noindent\textbf{Update 5.74}
No progress on this problem has been reported to us.\\

\noindent\textbf{Problem 5.75} Does there exist a bounded analytic function in $\mathbb{D}$ such that the
image of every radius has infinite length? See, for example, Anderson \cite{MR0287295}, and Rudin \cite{MR0079093}.\\
\textit{(W. Rudin; communicated by K. F. Barth)}\\

\noindent\textbf{Update 5.75}
No progress on this problem has been reported to us.\\

\noindent\textbf{Problem 5.76} Let $\gamma$ be a non-tangential arc that lies in $\mathbb{D}$ except for one endpoint
at $z = 1$, and define
\[\gamma_\theta=e^{i\theta\gamma},\hspace{1cm}\theta\in[0,2\pi).\]
Does there exist a function $g\in H^\infty$ such that
\[\lim_{z\to e^{i\theta},\,z\in\gamma_\theta}g(z)\]
exists for no value of $\theta$? (See Rudin \cite{MR516096}.)\\
If $\gamma$ is tangential, the answer is `yes', see Collingwood and Lohwater \cite[p. 43]{MR0231999}.\\
\textit{(W. Rudin; communicated by K. F. Barth)}\\

\noindent\textbf{Update 5.76}
No progress on this problem has been reported to us. We note that this problem is not clearly formulated, as the general form of Fatou's theorem says that the angular limit exists for almost all point on $\mathbb{T}$, see Collingwood and Lohwater \cite[p. 21]{MR0231999}.\\

\noindent\textbf{Problem 5.77} In general, the radial behaviour of the derivative of a bounded analytic
function in $\mathbb{D}$ can be pretty arbitrary; in fact, even under much stronger
restrictions than just bounded, it can still be quite arbitrary.\\
Is it true that, given any measurable function $m(\theta)$ on $[0,2\pi)$, there exists a
function $f(z)$, continuous on $\overline{\mathbb{D}}$  and univalent on $\mathbb{D}$, such that
\[\lim_{r\to1}f'(re^{i\theta})=m(\theta)\]
for almost all $\theta$? (In terms of known results, the conjecture seems quite plausible. Ortel and Schneider \cite{MR813642} showed that the conjecture is true under slightly strengthened hypotheses;
and Lohwater, Piranian and Rudin \cite{MR0072218} showed that the conjecture is true with a slightly weakened
conclusion.)\\ 
\textit{(W. J. Schneider)}\\

\noindent\textbf{Update 5.77}
No progress on this problem has been reported to us.\\

\noindent\textbf{Problem 5.78} Let $H^1$ denote Hausdorff one-dimensional measure on $\mathbb{C}$, and $\mathbb{T}$;
let $g:\mathbb{T}\to[-\infty,\infty]$ denote an arbitrary Borel function. Does there exist a
corresponding function $f$, analytic in $\mathbb{D}$  and with bounded
Taylor coefficients, such that
\[\lim_{r\to1^-}f(rz)=g(z)\]
for $H^1$-almost all $z$ in $\mathbb{T}$?\\
For work on related questions, see Ortel and Schneider \cite{MR813642}.\\ 
\textit{(M. Ortel and W. J. Schneider)}\\

\noindent\textbf{Update 5.78}
No progress on this problem has been reported to us.\\

\noindent\textbf{Problem 5.79}
\begin{enumerate}[(a)]
\item Let $f$ be a non-constant analytic function in $\mathbb{D}$, $m$ be a positive integer,
and define $\psi = (f)^mf'$. Then it is shown by Sons (\cite[Theorem 2]{MR889119}) that, when $f$ and $\psi$ both belong
to MacLane's class $A$ \cite{MR0148923}, either
\begin{enumerate}[(i)]
\item $f$ has finite asymptotic values on a dense subset of $\mathbb{T}$, or
\item $\psi$ assumes every finite value infinitely often.
\end{enumerate}
What replacements can be found for $(i)$ (to give another correct theorem)? Can $\psi$ be
taken to be of the form
\[\psi=(f)^{m_0}(f')^{m_1}\ldots(f^{(k)})^{m_k},\]
whenever $k,m_0,m_1,\ldots,m_k\in\mathbb{N}\cup\{0\}$?\\
\item It is known that the conclusion in part $(a)$ is true when $\psi$ is replaced by
\[\psi=f^{(l)}+\sum^{i-1}_{\nu=0}a_\nu f^{(\nu)},\]
where $l\in\mathbb{N}$ and the $a_\nu$ are analytic functions in $\{|z| < 1 +r\}$, for some positive $r$. What
replacements can be found in this case for $(i)$ to give another correct theorem? Can
the $a_\nu$ be taken as functions analytic in $\mathbb{D}$ with $T(r,a_\nu) = o(T(r,f))$ as $r\to1$?\\
(For both parts, compare Problem 1.38.)
\end{enumerate}
\textit{(L. R. Sons)}\\

\noindent\textbf{Update 5.79}
No progress on this problem has been reported to us.\\

\section*{New Problems}

\noindent To appear.

%======================= CHAPTER 6 ====================================

\chapter{Univalent and Multivalent Functions}

\underline{Preface by Ch. Pommerenke (to appear)}\\

\noindent\textbf{Notation} A function $w=f(z)$ analytic or meromorphic in a domain $D$, is said to be \textit{schlicht} or \textit{univalent} in $D$, if $f(z)$ assumes different values at different points of $D$. Such a function maps $D$ univalently onto a domain $\Delta$ in the $w$ plane. Among classes of univalent functions, the following play a particularly important role. We denote by $S$ the class of functions of the form
\begin{equation}
\label{6.1}
f(z)=z+\sum^\infty_{n=2}a_nz^n
\end{equation}
analytic and univalent in the unit disc $\mathbb{D}$ and by $\Sigma$ the class of functions
\begin{equation}
\label{6.2}
F(z)=z+\sum^\infty_{n=1}b_nz^{-n}
\end{equation}
analytic and univalent for $1<|z|<\infty$.\\
An account of various elementary extremal problems for the class $S$ is given by Hayman \cite[Chps. 1, 6]{MR0108586}. 
It is known that for many problems in $S$, the extremal functions are the Koebe functions
\[w=f_\theta(z)=\frac{z}{(1-ze^{i\theta})^2}=z+\sum^\infty_{n=2}nz^ne^{i(n-1)\theta},\]
which map $\mathbb{D}$ onto the plane cut along a ray going from $w=-\frac{1}{4}e^{-i\theta}$ to $\infty$ in a straight line away from the origin. The corresponding functions of $\Sigma$ are
\[w=F_\theta(z)=f_\theta(z^{-2})^{-\frac{1}{2}}=z-\frac{e^{i\theta}}{z},\]
which map $|z|>1$ onto the plane cut along the segment from $-2ie^{i(\theta/2)}$ to $2ie^{i(\theta/2)}$.\\
Many of the most interesting problems for $S$ and $\Sigma$ relate to the size of the coefficients $a_n$ and $b_n$. Pre-eminent among these is the famous Bieberbach conjecture.\\

\noindent\textbf{Problem 6.1} \textbf{The Bieberbach conjecture} Is it true that $|a_n|\leq n$ for $f$ in $S$ with equality only for $f(z)\equiv f_\theta(z)$? The result is known to be true for $n=2$ (see Bieberbach \cite{bieberbach1916koeffizienten}), $n=3$ (see L\"owner \cite{Lowner1923}), and $n=4$ (see Garabedian and Schiffer \cite{MR0070712}, Charzynski and Schiffer \cite{MR0114930}).\\

\noindent\textbf{Update 6.1} De Branges' \cite{MR772434} proof of the Bieberbach conjecture has
transformed this field. Suppose that $S$ is the class of functions
\[f(z)=z+\sum^\infty_{n=2}a_nz^n\]
univalent in $\mathbb{D}$. We also define
\[\log\frac{f(z)}{z}=2\sum^\infty_{k=1}\gamma_kz^k\]
and
\[\Big(\frac{f(z)}{z}\Big)^p=\sum^\infty_{n=0}a_{n,p}z^n,\hspace{1cm}p>0.\]
Then de Branges proved Milin's conjecture (see \cite{MR772434} and \cite{MR0158056}) that
\[\sum^N_{n=1}\sum^n_{k=1}k|\gamma_k|^2\leq\sum^N_{n=1}\sum^n_{k=1}\frac{1}{k},\]
(Problem 6.42), and this in turn implies Robertson's conjecture (see Sheil-Small \cite{MR0320761})
\[\sum^n_{k=0}|a_{k,\frac{1}{2}}|^2\leq n+1\]
(Problem 6.39), which implies the Bieberbach conjecture $|a_n|\leq n$ (Problem 6.1), as well
as the corresponding result for functions subordinate to $f$, and a positive answer to
Problems 6.2 and 6.85.\\ 
Among earlier proofs, the result for $n = 6$ was proved independently by Pederson \cite{MR0239072}
and Ozawa \cite{MR0239073}, and for $n = 5$ by Pederson and Schiffer \cite{MR0311891}. Among later proofs, that of Fitzgerald and Pommerenke \cite{MR792819} may
be mentioned.\\

\noindent\textbf{Problem 6.2} Define $A_n=\sup_{f\in S}|a_n|$. It is shown by Hayman \cite{MR0096810} that
\[\frac{A_n}{n}\to K_0,\hspace{1cm}\text{ as }n\to\infty.\]
Is it true that $K_0=1$? The best known result so far is $K_0<1.243$ due to Milin \cite{MR0172991}.\\

\noindent\textbf{Update 6.2} See Update 6.1.\\
With regard to earlier weaker results, Fitzgerald \cite{MR0335777} had proved that
\[|a_n|\leq Kn,\]
where $K=(7/6)^{\frac{1}{2}}=1.080\ldots$. Horowitz \cite{MR0480979} sharpened this to $K<1.0657$, and Fitzgerald has remarked that the best his method might be expected to yield is $K<(7/6)^{\frac{1}{4}}=1.039\ldots$. Further, Hamilton \cite{MR684512} proved that $A_{n+1}-A_n\to K_0$, and  also \cite{MR732121} that $K_0=1$ is implied by, and so is equivalent to, Littlewood's conjecture $|a_n|\leq4dn$ where $d$ is the distance from the origin to the complement of the image of $\mathbb{D}$ by $f(z)$. By de Branges' Theorem (see Update 6.2) $A_n=n$ for all $n$, so that $K_0=1$, and Littlewood's conjecture holds.\\

\noindent\textbf{Problem 6.3} If $f(z)$ in $S$ it is shown by Bombieri \cite{bombieri2}, that there exist constants $c_n$ such that for $f(z)$ in $S$
\[|\text{Re}\,(n-a_n)|\leq c_n\text{Re}\,(2-a_2).\]
What is it the exact size of the constants $c_n$? Is it true that there exists $d_n$ such that 
\[\big|n-|a_n|\big|\leq d_n\big(2-|a_2|\big)\,?\]
\textit{(E. Bombieri)}\\

\noindent\textbf{Update 6.3} No progress on this problem has been reported to us.\\

\noindent\textbf{Problem 6.4} Suppose that $F(z)=z+\sum^\infty_{n=1}b_nz^{-n}\in\Sigma$ and that 
\[G(z)=z+\sum^\infty_{n=1}c_nz^{-n}\]
is the inverse function of $F(z)$. What are the exact bounds for $|b_n|$ and $|c_n|$? The following results are known\\
\indent$|b_1|=|c_1|$ (Bieberbach \cite{bieberbach1916koeffizienten}).\\
\indent$|b_2|=|c_2|\leq \frac{2}{3}$ (Schiffer \cite{MR1505083}, Goluzin \cite{Goluzin}).\\
\indent$|b_3|\leq \frac{1}{2}+e^{-6}$ (Garabedian and Schiffer \cite{MR0066457}).\\
\indent$|c_3|\leq 1$ (Springer \cite{MR0041935}). \\
These are sharp.\\

\noindent\textbf{Update 6.4} Springer \cite{MR0041935} had conjectured that 
\[|c_{2n-1}|\leq\frac{(2n-2)!}{n!(n-1)!},\hspace{1cm}n=1,2,\ldots.\]
Kubota \cite{MR0437735} has proved this conjecture for $n=3, 4, 5$, and has also \cite{MR0369683} obtained the sharp upper bounds of $|b_4|$ and $|b_5|$ among all $F$ with real coefficients.\\

\noindent\textbf{Problem 6.5} If it proves too difficult to obtain sharp bounds for all of the coefficients in Problem 6.4, we ask for the orders of magnitude. An area principle shows that
\[\sum^\infty_{n=1}n|b_n|^2\leq1\]
and hence
\[b_n=o(n^{-\frac{1}{2}}).\]
Clunie and Pommerenke \cite{MR0207974} have shown that 
\[|b_n|=O(n^{-\frac{1}{2}-\frac{1}{300}}).\]
In the opposite direction, examples due to Clunie \cite{MR0107715} show that
\[|b_n|>n^{-1+\delta}\]
is possible for indefinitely many $n$ and a fixed $F(z)$ in $\Sigma$, where $\delta$ is an absolute constant.\\

\noindent\textbf{Update 6.5} Pommerenke \cite{MR0220920} has proved that
\[|b_n|\geq n^{0.17-1}\]
can hold for infinitely many $n$.\\
As stated above, Clunie and Pommerenke \cite{MR0207974} proved that 
\[|b_n|=O\big(n^{-\frac{1}{2}-\delta}\big)\]
where $\delta=1/300$. A corresponding improvement to the coefficient estimates for bounded functions in $S$ follows from this method.\\
See also Update 4.18.\\

\noindent\textbf{Problem 6.6} What are the orders of magnitude of the $c_n$ in Problem 6.4? Springer \cite{MR0041935} obtained the estimate 
\[|c_n|\leq\frac{2^n}{n}\]
and also showed that, given $\varepsilon>0$, 
\[|c_{2n-1}|>(1-\varepsilon)2^{2n-2}e/(\pi n^3)^{\frac{1}{2}}\] 
is possible for all sufficiently large $n$.\\

\noindent\textbf{Update 6.6} No progress on this problem has been reported to us.\\

\noindent\textbf{Problem 6.7} If $f(z)$ in $S$ and is bounded, i.e. satisfies $|f(z)|<M$ for $z\in\mathbb{D}$, we again ask for the order of magnitude of the coefficients $a_n$. Since the area of the image of $\mathbb{D}$ by $f(z)$ is at most $\pi M^2$, we deduce that $\sum n|a_n|^2\leq M^2$, so that again
\[|a_n|=o(n^{-\frac{1}{2}}),\hspace{1cm}\text{ as }n\to\infty.\]
Here also Clunie and Pommerenke \cite{MR0207974} have shown that 
\[|a_n|=O(n^{-\frac{1}{2}-\frac{1}{300}}).\]
Examples in the opposite direction due to Littlewood \cite{Littlewood2}, show again that for a sufficiently small positive $\delta$, we can have
\[|a_n|>n^{-1+\delta}\]
for infinitely many $n$ and a fixed $f(z)$. The problems for this class of functions seem very analogous to the corresponding problems for $\Sigma$, see Problem 6.5.\\

\noindent\textbf{Update 6.7} Carleson and Jones \cite{MR1162188} established that this problem is equivalent to Problem 6.5; they also have shown that, as conjectured above, the bounds for the coefficient $|a_n|$ in $\Sigma$ and those for bounded functions in $S$ are of the form $n^{-\gamma+o(1)}$ where $\gamma$ is an absolute constant, which they conjecture to be $\frac{3}{4}$. See also Update 4.18 and Update 6.8.\\

\noindent\textbf{Problem 6.8} We write
\[I_\lambda(r,f)=\Big\{\frac{1}{2\pi}\int^{2\pi}_0|f(re^{i\theta}|^\lambda\,d\theta\Big\}^{1/\lambda}.\]
What are the exact bounds for $I_\lambda(r,f)$ and $I_\lambda(r,f')$ when $f$ in $S$ or $f$ in $\Sigma$? If $f(z)$ in $S$, it is known that for fixed $\lambda$, the orders of magnitude of $I_\lambda(r,f)$ and $I_\lambda(r,f')$ are maximal when $f(z)$ is the Koebe function. For the best results in this direction, see Bazilevi\v c \cite{MR0039804}.\\
If $f$ in $S$ and $|f|<M$, or if $f$ in $\Sigma$, it is almost trivial from the area principle that
\[I_1(r,f')=o(1-r)^{-\frac{1}{2}},\hspace{1cm}\text{ as }r\to1.\]
Clunie and Pommerenke \cite{MR0207974} have improved this to
\[I_1(r,f')=O(1-r)^{-\frac{1}{2}+\frac{1}{300}}.\]

\noindent A function $f(z)$ in $S$ is said to be \textit{convex} if the image of $\mathbb{D}$ by $w=f(z)$ is a convex domain $D$ in the $w$-plane, i.e. for any two points $w_1, w_2$ in $D$, the straight line segment $w_1, w_2$ also lies in $D$. A function $F(z)$ in $\Sigma$ is said to be \textit{convex} if the complement of the image of $|z|>1$ by $F(z)$ is convex.\\

\noindent\textbf{Update 6.8} Baernstein \cite{MR0387590} has shown that if $\phi(R)$ is any convex function of
$\log R$ then for each $r$, with $0 < r < 1$,
\[\int^{2\pi}_0\phi\big(|f(re^{i\theta})|\big)\,d\theta\]
is maximised in the class $S$ by the Koebe function. Thus he completely settled this classical problem for class $S$.
The problem for class $\Sigma$ remains open.
Carleson and Jones \cite{MR1162188} proved that problems 6.5, 6.7 and 6.8 for class $\Sigma$
are equivalent.\\
See also Update 4.18.\\

\noindent\textbf{Problem 6.9} \textbf{(Schoenberg's conjecture)}  If $f(z)=\sum^\infty_{n=1}a_nz^n$ and $g(z)=\sum^\infty_{n=1}b_nz^n$ are convex, and $f$, $g$ belong to $S$, is it true that
\[f\ast g=\sum^\infty_{n=1}a_nb_nz^n\]
is also a convex function in $S$? See P\'olya and Schoenberg \cite{MR0100753}.\\

\noindent\textbf{Update 6.9} Ruscheweyh and Sheil-Small (see \cite{MR0328051} and \cite{MR0447552}) have proved this result.\\

\noindent\textbf{Problem 6.10} If $F(z)$, $G(z)$ are convex functions in $\Sigma$, it is known that for $0<\lambda<1$,
\[H(z)=\lambda F(z)+(1-\lambda)G(z)\in\Sigma,\]
see Pommerenke \cite{MR0158058}. Is it true that $H(z)$ is also convex?\\
\textit{(Ch. Pommerenke)}\\

\noindent\textbf{Update 6.10} Contrary to previous updates, no progress on this problem has been reported to us.\\

\noindent\textbf{Problem 6.11} If $f(z)$, $g(z)$ are convex functions in $S$, is it true that for $0<\lambda<1$, $\lambda f+(1-\lambda)g$ is star-like and univalent? A function $w=f(z)$ in $S$ is said to be \textit{star-like} if the image domain $D$ is star-like with respect to the origin $O$, that is, if for any point $P$ in $D$ the straight line segment $OP$ lies in $D$. It is known that $f(z)$ is convex if and only if $zf'(z)$ is star-like, see e.g. Hayman \cite[Chp. 1]{MR0108586}.\\

\noindent\textbf{Update 6.11} A counter-example has been given by Macgregor \cite{MR0236369} who has also
found the largest disc in which $\lambda f+(1-\lambda)g$ is star-like.\\

\noindent\textbf{Problem 6.12} If $f(z)=z+\sum^\infty_{k=2}a_{n_k}z^{n_k}\in S$, and
\[\liminf_{k\to\infty}\frac{n_{k+1}}{n_k}>1,\]
then Pommerenke \cite{MR0167611} has proved that 
\begin{equation}
\label{6.3}
a_n=o\Big(\frac{1}{n}\Big)
\end{equation} 
and this is sharp. If
\[\liminf_{k\to\infty}(n_{k+1}-n_k)>4,\]
then Hayman \cite{MR0213535} has shown that
\begin{equation}
\label{6.4}
a_n=o(n^{-\frac{1}{2}}).
\end{equation} 
Are there intermediate gap conditions which allow us to interpolate between (\ref{6.3}) and (\ref{6.4})?\\

\noindent\textbf{Update 6.12} No progress on this problem has been reported to us.\\

\noindent\textbf{Problem 6.13} Suppose that $f(z)$ in $S$, and that positive integers $k, m, n,$ are given. It is known that there exist complex numbers $c_0, c_1,\ldots,c_m,$ depending on $k, m, n$ and $f$, such that $c_0=1$, $|c_m|\geq4^{-m}$ and
\begin{equation}
\label{6.5}
|c_0a_{n+j}+c_1a_{n+j+1}+\ldots+c_ma_{n+j+m}|\leq Kn^{\alpha_m},\hspace{1cm}0\leq j\leq k,
\end{equation} 
where $K$ is a constant depending on $k, m$ and $f$, and $\alpha_m=8/\sqrt{m-\frac{1}{2}}$. What are the best possible values for the $\alpha_m$? See Pommerenke \cite{MR0185105}. It is known that $\alpha_1=0$, and it may be conjectured that $\alpha_2=-\frac{1}{3}$, $\alpha_3=-\frac{1}{2}+\varepsilon$, and $\alpha_m<-\frac{1}{2}$ for $m>4$.\\
If, in addition, $f(z)$ is star-like, then Pommerenke has shown that (\ref{6.5}) holds with $\alpha_m=-1+\frac{2}{m+1}$. The result is sharp.\\
 
\noindent\textbf{Update 6.13} Pommerenke \cite{MR0220920} has made some progress with this, showing that for large $m$,
\[\alpha_m<-\frac{1}{2}+\frac{2}{m^2}.\]\\

\noindent\textbf{Problem 6.14} If $f(z)$ in $S$, set 
\[A^{(k)}_n=
\begin{vmatrix}
a_n,&a_{n+1},&\ldots,&a_{n+k-1}\\
\hdotsfor{4}\\
a_{n+k-1},&a_{n+k},&\ldots,&a_{n+2k-2}
\end{vmatrix}\]
Using (\ref{6.5}), Pommerenke \cite{MR0185105} proved that if $k$ is fixed
\begin{equation}
\label{6.6}
|A^{(k)}_n|^{1/k}=O(n^{j_k}),\hspace{1cm}\text{ as }n\to\infty,
\end{equation} 
with $j_k=-\frac{1}{2}+16/\sqrt{k}$. Here the conjecture is that $j_k<-\frac{1}{2}$ for large $k$. It can be shown that for $k=2$, (\ref{6.6}) is false with $j_2=0$, and holds with $j_2=\frac{1}{4}$, see Hayman \cite{MR0219715}.\\
If $f(z)$ is in addition star-like, then Pommerenke \cite{MR0185105} showed that the best possible values of $j_k$ are $j_k=-1+2/k$.\\

\noindent\textbf{Update 6.14} Pommerenke \cite{MR0215976} has proved that for $k=2, 3, \ldots$,
\[A_n^{(k)}=O\big(n^{\frac{3}{2}-\frac{1}{2}k-k\beta}\big),\hspace{1cm}\text{ as }n\to\infty\]
for some constant $\beta$, where $\beta>\frac{1}{4000}$.\\

\noindent\textbf{Problem 6.15} If $f(z)$ in $S$, write
\[f_\alpha(z)=\int^z_0f'(\zeta)^\alpha\, d\zeta.\]
For what values of $\alpha$, is it true that $f_\alpha(z)\in S$? The result is known to hold for $\alpha\leq(\sqrt{5}-2)/3$ (see Duren, Shapiro and Shields \cite{MR0199359}), but not for $\alpha>\frac{1}{3}$, see Royster \cite{MR0183866}.\\
If for $z$ in $\mathbb{D}$, $f(z)$ is analytic and $|f''(z)/f'(z)|\leq c/(1-|z|^2)$, then $f(z)$ is univalent if 
\begin{equation}
\label{6.15eqn}
c\leq2(\sqrt{5}-2),
\end{equation}
 see Duren, Shapiro and Shields \cite{MR0199359}; but not if $c>2$, see Hille \cite{MR0030000}. What is the best value of $c$?\\
\textit{(P. L. Duren)}\\ 

\noindent\textbf{Update 6.15} Becker \cite{MR0299780} has shown that the right hand side of (\ref{6.15eqn}) could be improved to $1$. Becker and Pommerenke \cite{zbMATH03857435} showed that this is best possible.\\

\noindent\textbf{Problem 6.16} Let $S^*$ be the class of all star-like functions $f(z)$ in $S$. Marx \cite{MR1512785} conjectured that for each fixed $z_0$, $|z_0|<1$, the set of all numbers $f'(z_0)$ for $f\in S^*$ coincides with the set of all number $k'(z)$, $|z|\leq |z_0|$, where
\[k(z)=\frac{z}{(1-z)^2}\]
is the Koebe function. This is known to be true for $|z|\leq 0.736$, see Duren \cite{MR0176044}.\\
\textit{(P. L. Duren)}\\ 

\noindent\textbf{Update 6.16} This is a version of the Marx conjecture \cite{MR1512785}. A counter-example has been given by Hummel \cite{MR0304629}.\\

\noindent\textbf{Problem 6.17} If $f(z)=z+\sum^\infty_{n=2}a_nz^n$ in $S$, then
\[ A=\pi\sum^\infty_{n=1}n|a_n|^2\]
is the area of the image domain. What is the minimum value of $A$ when $a_2$ is given? Clearly $A\geq\pi(1+2|a_2|^2)$ always, but this bound is not sharp if $|a_2|>\frac{1}{2}$, since in this case $f(z)=z+a_2z^2$ is not univalent in $\mathbb{D}$.\\
\textit{(H. S. Shapiro)}\\ 

\noindent\textbf{Update 6.17} No progress on this problem has been reported to us.\\

\noindent\textbf{Problem 6.18} If $F(z)=z+\sum^\infty_{n=1}b_nz^{-n}$ in $\Sigma$, then
\[A(F)=\pi-\pi\sum^\infty_{n=1}n|b_n|^2\]
is the area of the set of values not assumed by $F(z)$ in $\mathbb{D}$. If $F_n(z)\in\Sigma$, and
\[F_n(z)\to F(z),\hspace{1cm}\text{ as }n\to\infty,\]
for $|z|>1$, under what additional hypotheses is it true that
\begin{equation}
\label{6.7}
A(F_n)\to A(F),\hspace{1cm}\text{ as }n\to\infty\,?
\end{equation}
It is suggested that (\ref{6.7}) might be true under some hypotheses on $(1-|z|^2)^2\{F(z),z\}$ where
\[\{F(z),z\}=\Big(\frac{F''}{F'}\Big)^\prime-\frac{1}{2}\Big(\frac{F''}{F'}\Big)^2\]
is the Schwarzian derivative of $F(z)$.\\
\textit{(L. Bers)}\\ 

\noindent\textbf{Update 6.18} We note that there were errors in the original statement of this problem. No progress on this problem has been reported to us.\\

\noindent\textbf{Problem 6.19} If $f(z)=\sum^\infty_{n=1}a_nz^n$ is analytic in $\mathbb{D}$ and $\sum^\infty_{n=1}|a_n|<+\infty$, can $f(z)$ map the unit circle $\mathbb{T}$ onto a curve of positive two-dimensional measure if
\begin{enumerate}[(a)]
\item $f(z)$ in $S$,
\item more generally, $f'(z)\neq0$ in $\mathbb{D}$?
\end{enumerate}

\noindent\textbf{Update 6.19} No progress on this problem has been reported to us.\\

\noindent\textbf{Problem 6.20} Let $C$ be a closed curve inside the unit circle $\mathbb{T}$. Under what conditions on $C$ does there exist a univalent function $f$ in $\mathbb{D}$ such that $f(C)$ and $f(\mathbb{T})$ are both convex?\\

\noindent\textbf{Update 6.20} No progress on this problem has been reported to us.\\

\noindent\textbf{Problem 6.21} A function $f(z)$ analytic in $\mathbb{D}$ is said to be \textit{typically real} if $f(z)$ is real, when and only when $z$ is real, see Rogosinski \cite{MR1545292}.\\
If $f(z)=z+\sum^\infty_{n=2}a_nz^n$ is typically real in $\mathbb{D}$, then
\[f(z)=\frac{z}{1-z^2}P(z),\]
where $P(0)=0$, $\text{Re}\, P(z)>0$ in $\mathbb{D}$. What other conditions must $P(z)$ satisfy to make $f(z)$ univalent in $\mathbb{D}$?\\

\noindent\textbf{Update 6.21} No progress on this problem has been reported to us.\\

\noindent\textbf{Problem 6.22} If $f(z)=z+\sum^\infty_{n=2}a_nz^n$ is univalent and star-like of order $\frac{1}{2}$ in $\mathbb{D}$, i.e. 
\[\text{Re}\,\frac{zf'(z)}{f(z)}\geq\frac{1}{2},\]
find the radius of the largest disc $|z|<r$ in which $f(z)$ is convex. In other words, when is 
\[\min_{|z|=r}\text{Re}\,\Big\{\frac{P+1}{2}+z\frac{P'(z)}{P+1}\Big\}>0,\]
where $\text{Re}\, P(z)>0$ in $\mathbb{D}$ and $P(0)=1$?\\

\noindent\textbf{Update 6.22} No progress on this problem has been reported to us.\\

%\section* *{Multivalent Functions}

\noindent\textbf{Notation} Let $f(z)$ be analytic in a domain $D$ and let $n(w)$ be the number of roots of the equation $f(z)=w$ in $D$. Then $f(w)$ is said to be \textit{mean $p$-valent} in $D$ if
\begin{equation}
\label{6.8}
\int^{2\pi}_0\int^R_0n(\rho e^{i\phi})\rho \,d\rho \,d \phi\leq \pi pR^2,\hspace{1cm}0<R<\infty.
\end{equation} 
Frequently, it is possible to prove analogous results to those for univalent functions for the wider class of mean $p$-valent functions. Counter-examples for mean $p$-valent functions are in general much easier to find. We discuss briefly how some of the preceding problems could be modified for mean $p$-valent functions. In this connection, the following counter-example is valuable, see Pommerenke \cite{MR0133448}.\\
If 
\[f(z)=1+\sum^\infty_{m=1}a_mz^m,\]
where
\begin{enumerate}[(a)]
\item $\sum^\infty_{m=1}|a_m|<\varepsilon$,
\item $\sum^\infty_{m=1}m|a_m|^2<\varepsilon$, 
\end{enumerate}
then $f(z)$ is bounded and mean $p$-valent in $\mathbb{D}$ with $p<\varepsilon/(1-\varepsilon)^2$. In fact
\[n(\rho e^{i\phi})=0,\hspace{1cm}\rho<1-\varepsilon,\]
and
\[\int^{2\pi}_0\int^R_0 n(\rho e^{i\phi})\,d\rho \,d\phi\leq\pi\sum^\infty_{n=1}n|a_n|^2<\pi\varepsilon<\pi\frac{\varepsilon}{(1-\varepsilon)^2}R^2,\]
if $R>1-\varepsilon$. \\
On the other hand, $(a)$ and $(b)$ permit isolated relatively large coefficients so that no boundedness or gap-condition for mean $p$-valent functions can imply more than $|a_n|=o(n^{-\frac{1}{2}})$ for this class. In particular, the results of Pommerenke and Clunie quoted in Problems 6.5, 6.7, 6.8 and 6.12 have no analogue for mean $p$-valent functions, though (\ref{6.4}) holds for the wider class, provided that 
\[\limsup (n_{k+1}-n_k)>4p,\]
and this condition is sharp, see Hayman \cite{MR0213535}.\\

We discuss modifications of some of the other problems for $p$-valent functions.\\
Since a polynomial of degree $p$ is automatically mean $p$-valent for any positive integer $p$, it is clear that we cannot in general expect to limit the growth of mean $p$-valent functions
\[f(z)=\sum^\infty_{n=0}a_nz^n\]
unless we have information about the coefficients $a_0$ to $a_p$. We write
\[\mu_p=\max_{0\leq\nu\leq p}|a_\nu|,\]
and denote by $S(p)$ the class of functions $f(z)$ analytic and mean $p$-valent in $\mathbb{D}$, and such that $\mu_p=1$. It is known (see e.g. Hayman \cite[Theorem 3.5, p. 50]{MR0108586}) that, if $f(z)$ in $S(p)$ with $p>\frac{1}{4}$,
\[|a_n|<A(p)n^{2p-1}.\]
Our counter-example shows that this is false for $p<\frac{1}{4}$. It was also shown by Spencer (in a letter to Hayman) that $|a_3|>3$ is possible for $f(z)=z+\sum^\infty_{n=2}a_nz^n$ in $S(1)$, so that the  analogue of Bieberbach's conjecture fails for $n=3$. However, Spencer \cite{MR0004903} proved that $|a_2|\leq2$ in this case.\\ 

\noindent\textbf{Problem 6.2$^\prime$} If $A^{(p)}_n=\sup_{f\in S_p}|a_n|$ is it true that
\[\frac{A^{(p)}_n}{n^{2p-1}}\to K_p,\hspace{1cm}\text{ as }n\to\infty,\]
and if so, what is $K_p$? It is known that for a fixed $f$ in $S(p)$, the limit
\[\alpha_f=\lim_{n\to\infty}\frac{|a_n|}{n^{2p-1}}\]
exists if $p>\frac{1}{4}$, see Hayman \cite{MR0071536} and Eke \cite{MR0004903}.\\

\noindent\textbf{Update 6.2$^\prime$} No progress on this problem has been reported to us.\\

\noindent\textbf{Problem 6.7$^\prime$} Here our counter-example shows that $|a_n|=o(n^{-\frac{1}{2}})$ is best possible for bounded $f(z)$ in $S(p)$.\\

\noindent\textbf{Problem 6.8$^\prime$} If we ask the analogous problems to those of Problem 6.8 for the class $S(p)$, the correct orders of magnitude are again known in many cases, but not the exact bounds. If $f$ in $S(p)$ and $f$ is bounded, then
\[I_1(r,f')=o(1-r)^{-\frac{1}{2}}\]
and this is sharp for the class, by our introductory example.\\

\noindent\textbf{Problem 6.13$^\prime$} The results of Pommerenke \cite{MR0185105} were proved in fact for mean $p$-valent functions, and if $f(z)=\sum^\infty_{n=0}a_nz^n$ is mean $p$-valent with $p>\frac{1}{4}$, then (\ref{6.5}) holds with $\alpha_m=-\frac{1}{2}+8p^\frac{3}{2}/\sqrt{m}$. This result is probably far from the best possible, though clearly $\alpha_m\geq-\frac{1}{2}$ in all cases. Lucas \cite{lucasnew} has shown that $\alpha_1=2p-2$ if $p\geq1$, and $\alpha_1\leq 2p-2\sqrt{p}$ for $\frac{1}{4}<p<1$. It is fair to conjecture that the correct value of $\alpha_1$ is $p-1$ for $\frac{1}{2}<p<1$ and $-\frac{1}{2}$ if $p<\frac{1}{2}$. It might also be conjectured that $\alpha_m=-\frac{1}{2}$ for $m+1>4p$.\\
The functions
\[f(z)=\Big(\frac{1+z^{m+1}}{1-z^{m+1}}\Big)^{2p/(m+1)}=1+\sum^\infty_{n=2}a_nz^n\]
provide counter-examples. These functions are mean $p$-valent, and $a_n=0$ except when $(m+1)$ divides $n$, and for the remaining values of $n$ we have, as $n\to\infty$,
\[|a_n|\sim Cn^{2p/(m+1)-1},\hspace{1cm}\text{ where }C\text{ is a constant}.\]
Thus (\ref{6.5}) cannot hold with $\alpha_m<2p/(m+1)-1$.\\
 
\noindent\textbf{Update 6.13$^\prime$} Some work on this probem has been done by Noonan and Thomas \cite{MR0306479}.\\

\noindent\textbf{Problem 6.14$^\prime$} Here again the main conclusions extend to mean $p$-valent functions. In this case (\ref{6.6}) holds with $j_k=-\frac{1}{2}+16(p^3/k)^{\frac{1}{2}}$. This is still unlikely to be best possible.\\

\noindent\textbf{Update 6.14$^\prime$} Some work on this probem has been done by Noonan and Thomas \cite{MR0306479}.\\

\noindent\textbf{Problem 6.23} A related problem concerns upper bounds for $|a_{n+1}|-|a_n|$ when $f(z)$ is mean $p$-valent. Lucas \cite{lucasnew} has proved that
\[\big||a_{n+1}|-|a_n|\big|=O(n^{j_d}),\]
where $j_d=2p-2$, if $p\geq1$;  $j_d\leq2p-2\sqrt{p}$ if $\frac{1}{4}<p<1$; and $j_d=-\frac{1}{2}$, if $p<\frac{1}{4}$. The result for $\frac{1}{4}<p<1$ is probably not sharp. A similar question may be asked for symmetric mean $p$-valent functions of the type
\[f(z)=\sum^\infty_{n=0}a_nz^{an+b}.\]
The coefficients of such functions behave rather like those of functions $\sum a_nz^n$ which are mean $(p/a)$-valent. In particular, if $f(z)=z+\sum a_nz^{2n+1}$ is mean univalent, then Lucas \cite{lucasnew} proved that
\[\big||a_{n+1}|-|a_n|\big|=O(n^{1-\sqrt{2}}).\]

\noindent\textbf{Update 6.23} No progress on this problem has been reported to us.\\

\noindent\textbf{Problem 6.24} If $f(z)=z+\sum^\infty_{n=2}a_nz^n\in S(1)$, prove that on $|z|=r$,
\[|f(z)|\leq\frac{r}{(1-r)^2}.\]
It is shown by Garabedian and Royden \cite{MR0060020} that $f(z)$ assumes in $\mathbb{D}$ each value $w$ such that $|w|<\frac{1}{4}$ and hence that
\[|f(z)|\geq\frac{r}{(1+r)^2},\hspace{1cm}|z|=r.\]
The question of sharp bounds for $|f'(z)|$ and $|f'(z)|/|f(z)|$ is also open. The corresponding results for $S$ are elementary, see e.g. Hayman \cite[Chp. 1]{MR0108586}.\\

\noindent\textbf{Update 6.24} No progress on this problem has been reported to us.\\

%\section* *{$p$-valent and circumferentially mean $p$-valent functions}

\noindent\textbf{Notation} If $f(z)$ is analytic in a domain $D$, and $n(w)$ is the number of roots of the equation $f(z)=w$ in $D$, then $f(z)$ is called \textit{$p$-valent} in $D$ if $p$ is an integer and $n(w)\leq p$ always. If $p$ is any positive number and
\begin{equation}
\label{6.9}
\frac{1}{2\pi}\int^{2\pi}_0n(Re^\phi)\,d\phi\leq pR,\hspace{1cm}0<R<\infty,
\end{equation}
then $f(z)$ is called \textit{circumferentially mean $p$-valent}. Clearly (\ref{6.9}) implies ({\ref{6.8}), that is, circumferentially mean $p$-valence implies mean $p$-valence.\\

\noindent\textbf{Problem 6.25} Suppose that $p$ is an integer and $f(z)=\sum^\infty_{n=0}a_nz^n$ is $p$-valent in $\mathbb{D}$. It is conjectured by Goodman \cite{MR0023910} that
\[|a_n|\leq\sum^p_{k=1}|a_k|D(p,k,n),\]
where
\[D(p,k,n)=\frac{2kn\prod^p_{\alpha=1}(n^2-\alpha^2)}{(p+k)!(p-k)!(n^2-k^2)},\hspace{1cm}1\leq k\leq p<n.\]
This result, containing the Bieberbach conjecture as a special case, is likely to be extremely difficult. The inequality if true would be sharp in all cases. No counter-examples are known and the conjecture has been proved only if $a_k=0$ for $k=1, 2, \ldots, (p-1)$ and $n=p+1$, see Spencer \cite{MR0004903}; or $n=p+2$, see Jenkins \cite[p. 160]{MR0096806}. The conjecture is true and sharp for certain classes of $p$-valent functions, namely those which are typically real of order $p$, see Goodman and Robertson \cite{MR0040430}.\\
We recall that the conjecture is definitely false for areally mean $p$-valent functions, if $p=1$ and $n=3$ by the example of Spencer \cite{MR0004903}, though for circumferentially mean $p$-valent functions it remains true in this case, see Jenkins \cite[p. 160]{MR0096806}.\\

\noindent\textbf{Update 6.25} No progress on this problem has been reported to us.\\

\noindent\textbf{Problem 6.26} Suppose that $f(z)=\sum^\infty_{n=0}a_nz^n$ is circumferentially mean $p$-valent and $f(z)\neq0$ in $\mathbb{D}$. (This latter condition is a consequence of mean $p$-valency if $p<1$). It is conjectured that in this case, at least if $p\geq 1$, we have
\begin{equation}
\label{6.10}
|a_n|\leq A_{n,p}
\end{equation}
where
\begin{equation}
\label{6.11}
F(z)=a_0\Big(\frac{1+z}{1-z}\Big)^{2p}=\sum^\infty_{n=0}A_{n,p}z^{2n}.
\end{equation}
The conjecture (\ref{6.10}) is known to be true for $n=1$ and all $p$; and for $n=2, 3$ if $p=1$, by the results quoted in Problem 6.25. It is certainly false for small positive $p$, and large $n$, since it would imply $a_n=O(n^{\varepsilon-1})$ for every positive $\varepsilon$ as $n\to\infty$, for a bounded univalent function. For if $g(z)$ is bounded and univalent, and $\varepsilon>0$, then if $K$ is a sufficiently large positive constant $g(z)+K$ is circumferentially mean $\varepsilon$-valent.\\
On the other hand, if $p>\frac{1}{4}$ and $f(z)$ is fixed, (\ref{6.10}) holds for all sufficiently large $n$, see Hayman \cite[Theorem 5.10, p. 113]{MR0108586}. In the special case when $p=1$ and $f(z)$ is univalent, (\ref{6.10}) reduces to the Littlewood conjecture $|a_n|\leq4|a_0|n$. This conjecture is somewhat weaker than the Bieberbach conjecture. It was shown by Nehari \cite{MR0083569} (see also Bombieri \cite{MR0185108}) that we have at any rate $|a_n|\leq 4K_0|a_0|n$ if $f(z)$ is univalent, where $K_0$ is the constant of Problem 6.2. Thus, Littlewood's conjecture holds since $K_0=1$, as was pointed out in Update 6.2.\\
It also seems likely that for $\lambda p>1$
\[I_\lambda(r,f)\leq I_\lambda(r,F)\]
if $f(z)$ satisfies the above hypotheses, and $F(z)$ is given by (\ref{6.11}).\\

\noindent\textbf{Update 6.26} See Update 6.2. Nothing beyond this is known.\\

\noindent\textbf{Problem 6.27} Suppose that
\[g(z) = z + b_0 + b_1z^{-1} + \ldots\]
is univalent in $|z|>1$. Is it true that for each positive $\varepsilon$ we have
\[n|b_n|=O(n^\varepsilon)\big\{\max_{0<|\nu-n|<\frac{n}{2}}(\nu(|b_\nu|+1)\big\},\hspace{1cm}\text{ as }n\to\infty\,?\]
This is suggested by some results of Clunie and Pommerenke.\\
\textit{(J. Clunie, Ch. Pommerenke)}\\

\noindent\textbf{Update 6.27} There were errors in the original statement of this problem. No progress on this problem has been reported to us.\\

\noindent\textbf{Problem 6.28} Suppose that $f(z) = z+\sum^\infty_{n=2}a_nz^n$ in $S$ and that \mbox{$P(z) = \sum^n_{k=0}b_kz^k$} is a polynomial of degree at most $n$. Is it true that
\[\max_{|z|=1} |P(z)\ast f(z)| \leq n \max_{|z|=1} |P(z )|?\]
Here $P \ast f = \sum^n_{k=0}a_kb_kz^k$. The above result would imply Rogosinski's
generalised Bieberbach conjecture but is weaker than Robertson's conjecture, see Sheil-Small \cite{MR0320761}.\\
\textit{(T. Sheil-Small)}\\

\noindent\textbf{Update 6.28} No progress on this problem has been reported to us.\\

\noindent\textbf{Problem 6.29} With the above notation $f (z)$ in $S$ if and only if for each pair of
numbers $\xi_1, \xi_2$ satisfying $|\xi_1|\leq1$, $|\xi_2|\leq1$, we have
\[f(z)\ast\frac{z}{(1-\xi_1z)(1-\xi_2z)}\neq0,\hspace{1cm}0<|z|<1.\]
On the other hand it is true that if $F(z) \ast f(z) \neq 0$, $0 <|z| < 1$ whenever
$f$ in $S$, then $F(z)$ is star-like. What is the complete class of star-like
functions having this property? $F(z) = z + z^n/n$ has the property since the Bieberbach conjecture holds.\\
\textit{(T. Sheil-Small)}\\

\noindent\textbf{Update 6.29} No progress on this problem has been reported to us.\\

\noindent\textbf{Problem 6.30} If $f$ in $S$, Baernstein has shown that
\[\int^{2\pi}_0|f(re^{i\theta})|^p\,d\theta\leq\int^{2\pi}_0|k(re^{i\theta})|^p\,d\theta,\hspace{1cm}0<r<1,\hspace{1cm}0<p<\infty,\]
where $k(z)$ is the Koebe function.\\
Does the corresponding inequality hold for integral means of the
derivatives at least for certain values of $p$? The best we can hope for is
that it holds for $p\geq\frac{1}{3}$ because $k'(z)\in H^p$ for $p<\frac{1}{3}$, and there exist
functions $f$ in $S$ for which $f'(z)$ does not belong to any $H^p$. (An
example is due to Lohwater, Piranian and Rudin \cite{MR0072218}). For close-to-convex
functions it was proved by MacGregor \cite{MR0222272} that the result holds for $p\geq1$,
and in fact the corresponding inequality holds for derivatives of all orders.\\
\textit{(A. Baernstein)}\\

\noindent\textbf{Update 6.30} Leung \cite{MR554395} proved that 
\[\int^{2\pi}_0|f'(re^{i\theta})|^p\,d\theta\leq\int^{2\pi}_0|k'(re^{i\theta})|^p\,d\theta,\hspace{1cm}0<r<1,\hspace{1cm}0<p<\infty,\]
for all functions in a certain subclass of the Bazilevi\v{c} functions which includes the close-to-convex functions. For $p<1$ this was previously unknown even for star-like functions. Feng and MacGregor \cite{Feng1976} have shown that in the full class $S$ the inequality is correct in order of magnitude as $r\to1$ if $p>2/5$. The asymptotic problem therefore remained open for $1/3\leq p\leq 2/5$, but Makarov \cite{MR818450} has given a counter-example if $p-\frac{1}{3}$ is sufficiently small. This is related to the Brennan conjecture. See Problem 6.96.\\

\noindent\textbf{Problem 6.31} Duren \cite{MR0364623} has shown that if $f(z) = \sum^\infty_{n=0}a_nz^n$ in $S$ and if
\[(1 - r )^2f( r ) = \lambda + O\big(( 1 - r )^\delta\big),\hspace{1cm}\text{ as }r\to1-,\]
for some $\lambda$, $\delta$, where $\lambda\neq0$, $\delta>0$, then
\[\frac{a_n}{n}=\lambda+O\Big(\frac{1}{\log n}\Big), \hspace{1cm}\text{ as }n\to\infty.\]
To what extent can this estimate be improved?\\
\textit{(P. L. Duren)}\\

\noindent\textbf{Update 6.31} No progress on this problem has been reported to us.\\

\noindent\textbf{Problem 6.32} Let $S_\alpha$, $0 < \alpha \leq 1$ be the subclass of $S$ of functions $f$ such that $\mathbb{C}\setminus f(\mathbb{D})$ is a single piecewise analytic slit from some finite point $\omega_0$ to $\infty$ that makes an angle at most $\alpha \pi/2$ with the radii vectors. What can be said about the Taylor coefficients of
functions in $S_\alpha$? If $f(-1) = \infty$ and $f(e^{i\phi_f}) = \omega_0$, $-\pi< \phi < \pi$, find $\sup_{f\in S_\alpha} |\phi_f|$.\\
\textit{(K. W. Lucas)}\\

\noindent\textbf{Update 6.32} No progress on this problem has been reported to us.\\

\noindent\textbf{Problem 6.33} The same questions as in Problem 6.32 can be asked under the alternative hypothesis
that $\mathbb{C}\setminus\{f(\mathbb{D})\}$ is a single piecewise analytic slit from some finite point ($\omega_0$, say)
to $\infty$ lying in an infinite sector with opening $\alpha\pi$ $(0<\alpha\leq2)$ and vertex $\omega_0$.\\
\textit{(K. W. Lucas and D. A. Brannan)}\\

\noindent\textbf{Update 6.33} No progress on this problem has been reported to us.\\

\noindent\textbf{Problem 6.34} A function $f(z) = z + a_2z^2 +\ldots$ analytic in $\mathbb{D}$ is said to belong to
\textit{Ruscheweyh's class} $M$ if the $*$ (i.e. Hadamard) convolution of $f$ with every (normalised)
convex function is univalent. All close-to-convex functions lie in $M$. Suppose that
$g(z) = z + b_2z^2 +\ldots$ is analytic in $\mathbb{D}$ and satisfies the condition
\begin{equation}
\label{C6.1}
\text{Re}\,\Big\{\frac{\phi*(gF)}{\phi*g}\Big\}>0\hspace{1cm}|z|<1,
\end{equation}
for all normalised convex functions $\phi$, and all normalised functions $F$ of positive real
part in $\mathbb{D}$. If $g$ is star-like, then (\ref{C6.1}) certainly holds; is (\ref{C6.1}) true for any other $g$,
or (maybe) for some significant larger family of $g$?\\
If $g$ satisfies (\ref{6.2}) and the condition 
\[\text{Re}\, \Big(\frac{zf'}{g}\Big) > 0,\hspace{1cm} |z| < 1,\]
then $f\in M$; does this classify $M$?\\
\textit{(T. Sheil-Small)}\\

\noindent\textbf{Update 6.34} Sheil-Small observes that the suggested classification for $M$ is invalid. The duality method shows that if $g(z)=z+\ldots$ is analytic and satisfies
\[\phi(z)*\frac{1+xz}{1-yz}g(z)\neq0,\hspace{1cm}0<|z|<1,\]
for $|x|=|y|=1$ and every convex function $\phi$, then $g$ is star-like. Duality also shows that $M$ is a linear invariant family. Its geometrical classification remains unknown.\\

\noindent\textbf{Problem 6.35} Let $\mathbb{O}$ be a subset of $\mathbb{D}=\{|\omega|< 1\}$. Find a characterisation of those $\mathbb{O}$ that are of the form $(\mathbb{C}\setminus f(\mathbb{D}))\cap\mathbb{D}$ for some $f$ in $S$. What is the maximum area of $\mathbb{O}$ ?
Given that $\omega_1, \omega_2\in\{|\omega|<1\}$, how does one tell whether there exists such an $f$ with
$\omega_1, \omega_2 \in\mathbb{O}_f$? The same question could be asked for $\omega_1, \omega_2, \omega_3$ etc.\\ 
\textit{(A. W. Goodman)}\\

\noindent\textbf{Update 6.35} This problem was originally wrongly attributed to D. A. Brannan. Barnard and Suffridge have observed that an extremal domain $g(\mathbb{D})$, $g$ in $S$, for this problem must be circularly symmetric and have as its boundary (up to rotation) the negative reals up to $-1$, an arc $\lambda$ of the unit circle which is symmetric about $-1$, and a Jordan curve $\gamma$ symmetric about the real axis, lying in the unit disc $\mathbb{D}$ and connecting the endpoints of $\lambda$. If $\gamma$ is assumed to be piecewise twice continuously differentiable, then it can be proved that $|zg'(z)|$ is constant on $g^{-1}(\gamma)$. See also Barnard, Pearce, and Solynin \cite{MR3002600} for a summary of progress.\\

\noindent\textbf{Problem 6.36} Suppose that $f$ in $S$ and define
\[f_p(z)=[f(z)]^p=z^p+\sum^\infty_{n=p+1}a_{n, p}z^n.\]
What can be said about bounds for $a_{n,p}$? If $|a_{n,1}|\leq Kn$ for all $n$ and fixed $K$, then,
for integral $p$,
\[|a_{n,\,p}|\leq K^p\frac{2p(2p+1)\ldots(n+p-1)}{(n-p)!},\]
but it might be easier to obtain bounds for $f_p$ than for $f$.\\
\textit{(W. K. Hayman)}\\

\noindent\textbf{Update 6.36} It was observed independently by Aharonov (unpublished), Grinshpan \cite{MR697766} and Hayman and Hummel \cite{MR877651}, that de Branges' theorem yields the sharp bounds
\[|a_{n,\,p}|\leq\frac{\Gamma(2p+n)}{\Gamma(n+1)\Gamma(2p)}=b_{n,\,p},\hspace{1cm}n=0,1,2,\ldots,\hspace{1cm}p\geq1.\]
The result holds trivially for all $p$ if $n=0$ or $1$, but, as Grinshpan \cite{MR697766} has shown, it fails for $p<1$ whenever $n\geq2$. However, Hayman and Hummel \cite{MR877651} showed that at least for $\frac{1}{4}\leq p<1$, if
\[A_{n,\,p}=\sup|a_{n,\,p}|,\]
then
\[\frac{A_{n,\,p}}{b_{n,\,p}}\to K_p,\hspace{1cm}\text{ as }n\to\infty.\]
Also $K_p>1$ if $p<0.499$. We conjecture that $K_p=1$ if and only if $p\geq\frac{1}{2}$.\\

\noindent\textbf{Problem 6.37} Suppose that $f(z)=z + c_3z^3 + c_5z^5 +\ldots$ is an odd univalent function in $\mathbb{D}$, and
let $d_n = |c_{2n+1}|-|c_{2n-1}|$. It is known that $d_n\to0$, the best known estimate being
$d_n =O(n^{1-\sqrt{2}})$; can this be improved to 
\begin{equation}
\label{J6.37}
d_n=O(n^{-\frac{1}{2}})\,?
\end{equation}
(Nothing better is possible, as
is shown by the fourth-root transform of the Koebe-function).\\
Milin (see \cite[p. 107]{MR0427620}; \cite{MR0450541}; and \cite{MR0369684}) proves that $d_n\leq K(\alpha)n^{-\frac{1}{2}}$ for functions $f$ such that
$g(z) = [f(z^{\frac{1}{2}})]^2 = z + a_2z^2 +\ldots$ is univalent in $\mathbb{D}$ and has positive Hayman
number $\alpha = \lim_{n\to\infty} n^{-1}|a_n|$; but $K(\alpha)\to\infty$ as $\alpha\to0$. Levin \cite{MR1545451} showed that
$d_n=O(n^{-\frac{1}{2}}\log n)$ if $c_n$ vanishes for $n\not\equiv1\pmod 4$.\\ 
\textit{(P. L. Duren)}\\

\noindent\textbf{Update 6.37} For general univalent functions, there has been no improvement on Lucas' upper bound $d_n =O(n^{1-\sqrt{2}})$, however (\ref{J6.37}) has been proved by Elhosh \cite{MR741643} for close-to-convex functions.\\

\noindent\textbf{Problem 6.38} With the notation of Problem 6.37, is it true that
\[\sum^\infty_{n=1}n^{-\beta}d_n^2<\infty\]
where $\beta=(\sqrt{2}-1)^2$?\\ 
\textit{(K. W. Lucas)}\\

\noindent\textbf{Update 6.38} Milin \cite{MR0231992} has proved that for every positive $\beta$,
\[\sum^\infty_{n=1}n^{-2\beta}d_n^2<50\sum^\infty_{n=1}\frac{1}{n^{1+2\beta}}.\]\\

\noindent\textbf{Problem 6.39} Suppose $f$ in $S$ and define $h(z) = \{f(z^2)\}^{\frac{1}{2}} = z + c_3z^3 + c_5z^5 +\ldots$ Robertson's conjecture (see Sheil-Small \cite{MR0320761}) asserts that
\[1+|c_3|^2+|c_5|^2+\ldots+|c_{2n-1}|^2\leq n.\]
This is known to be true if $f$ is star-like; is it true if $f$ has real coefficients, or if $f$ is close-to-convex?\\
\textit{(P. L. Duren)}\\

\noindent\textbf{Update 6.39} See Update 6.1.\\

\noindent\textbf{Problem 6.40} If $f(z)$ in $S$ and if the $a_n$ are real, then
\begin{equation}
\label{C6.2}
1+a_3+\ldots+a_{2n-1}\geq a_n^2,\hspace{1cm}n\geq1.
\end{equation}
The Bieberbach conjecture for such functions follows easily, see e.g. Fitzgerald \cite{MR0335777}.
It was pointed out by Clunie and Robertson that (\ref{C6.2}) holds for
normalised typically-real functions; the inequality is clear from the representation
formula for these functions.\\
The Bieberbach conjecture for $S$ would follow analogously if we could prove
that if $f$ in $S$,
\begin{equation}
\label{C6.3}
1+|a_3|+\ldots+|a_{2n-1}|\geq|a_n|^2,\hspace{1cm}n\geq1.
\end{equation}
Bshouty (unpublished) has shown that if $f$ in $S$, then there exists an $N(f)$ such that
(\ref{C6.3}) holds for $n > N(f)$.\\ 
\textit{(C. Fitzgerald)}\\

\noindent\textbf{Update 6.40} No progress on this problem has been reported to us.\\

\noindent\textbf{Problem 6.41} Let $K(\alpha)$ and $S^*(\alpha)$ be those subsets of $S$ consisting of the class of functions
convex in $\mathbb{D}$ of order $\alpha$ i.e.
\[\text{Re}\,\Big[1+\frac{zf''(z)}{f'(z)}\Big]\geq\alpha,\hspace{1cm}|z|<1,\]
and star-like of order $\alpha$ in $\mathbb{D}$ i.e.
\[\text{Re}\,\Big[z\frac{f'z)}{f(z)}\Big]\geq\alpha,\hspace{1cm}|z|<1,\]
respectively.
\begin{enumerate}[(a)]
\item Prove that (see Goel \cite{MR0367172})
\[\min_{f\in K(\alpha)}\,\min_{|z|=r}\Big|\frac{zf'(z)}{f(z)}\Big|=\min_{f\in K(\alpha)}\,\min_{|z|=r}\Big[\text{Re}\,\frac{zf'(z)}{f(z)}\Big].\]
\item Show that the functions
\[(2\alpha-1)^{-1}[1-(1-z)^{2\alpha-1}],\hspace{1cm}\alpha\neq\frac{1}{2}\hspace{1cm}\text{ and }\hspace{1cm}-\log(1-z)\]
are star-like of order
\[4^\alpha(2\alpha-1)[4-4.2^\alpha]^{-1}\hspace{1cm}\text{ and }\hspace{1cm}(\log 4)^{-1}\]
respectively, see MacGregor \cite{MR0367175}.
\end{enumerate}
Either $(a)$ or $(b)$, combined with the work of I. S. Jack \cite{MR0281897} would solve the
following problem of F. R. Keogh: find
\[\max_{f\in K(\alpha)}\{\beta: f\in S^*(\beta)\}.\]
\textit{(D. Benjamin)}\\

\noindent\textbf{Update 6.41} Sheil-Small observes that Feng and Wilken \cite{MR575386} have established Jack's assertion \cite{MR0281897} which therefore completes the solution of the order of star-likeness of a function convex of order $\alpha$. Fournier \cite{zbMATH06729355} gives a new proof of Jack's assertion, and discusses an extension for polynomials.\\

\noindent\textbf{Problem 6.42} If $f$ in $S$, write
\[\log[f(z)/z]=2\sum^\infty_{k=1}\gamma_kz^k.\]
If $f$ is star-like then $|\gamma_k| \leq1/k$; this is false in general, even in order of magnitude.
Milin (see \cite[p. 151]{MR708494}) has shown that
\[\sum^n_{k=1}k|\gamma_k|^2\leq\sum^n_{k=1}\frac{1}{k}+\delta\]
where $\delta<0.312$, and that $\delta$ cannot be reduced to $0$; Milin conjectured that
\begin{equation}
\label{C6.4}
\sum^N_{n=1}\sum^n_{k=1}k|\gamma_k|^2\leq\sum^N_{n=1}\sum^n_{k=1}\frac{1}{k},
\end{equation}
which would imply Robertson's conjecture (see Problem 6.39). Inequality (\ref{C6.4}) is
known to be true for $N = 1, 2, 3$ (see Grin\v span \cite{MR0328061}). Is it true in general if $f$ has real
coefficients, or if $f$ is close-to-convex?\\ 
\textit{(P. L. Duren)}\\

\noindent\textbf{Update 6.42} See Update 6.1.\\

\noindent\textbf{Problem 6.43} Using the notation of Problem 6.42, it is well-known that
\[\Big|\sum^\infty_{k=1}k\gamma_kz^k\Big|=O\Big(\frac{1}{1-r}\Big),\hspace{1cm}r\to1-,\]
for $|z| = r < 1$; is it true that
\[\sum^\infty_{k=1}k|\gamma_k|r^k=O\Big(\frac{1}{1-r}\Big),\hspace{1cm}r\to1-\,?\]
\textit{(D. Aharonov)}\\

\noindent\textbf{Update 6.43} Hayman \cite{MR568639} has constructed a function $f$ in $S$ for which
\[\sum^\infty_{k=1}k^2|\gamma_k|^2r^k\neq o\Big(\frac{1}{1-r}\log\frac{1}{1-r}\Big),\]
but the general problem remains open.\\

\noindent\textbf{Problem 6.44} Let $f$, $g$ be formal power series
\[\sum^\infty_{n=0}a_nz^n,\hspace{1cm} \sum^\infty_{n=0}b_nz^n\]
respectively, and define 
\[(f\otimes g)(z)=\sum^\infty_{n=1}a_nb_nn^{-1}z^n.\]
Let $S_R$ denote the class of functions in $S$ with real coefficients. Prove (or disprove)
that $f,g\in S_R$ implies that $f\otimes g\in S_R$. (Robertson (uncited) has proved the corresponding
result for typically-real functions.\\ 
\textit{(J. G. Krzyz)}\\

\noindent\textbf{Update 6.44} Bshouty \cite{MR577757} has shown that the result is false by producing a counter-example using a deep result of Jenkins (see \cite[Corollary 4.8 and Example 4.5, p. 120]{MR0507768}).\\

\noindent\textbf{Problem 6.45} Let $S^*(\alpha)$ be the \textit{class of $\alpha$-strongly-star-like functions} $f$, that is, those $f$ in $S$ for
which
\[\Big|\arg\Big(\frac{zf'(z)}{f(z)}\Big)\Big|<\frac{\alpha\pi}{2}\hspace{1cm}\text{ for } |z|<1,\]
where $0 < \alpha < 1$.
\begin{enumerate}[(a)]
\item Prove (or disprove) that $S^*(\alpha)$ is closed under $\otimes$ (see Problem 6.44).
\item Prove (or disprove) that, if $f\in S^*(\alpha)$ and $g \in S^*(\beta)$, then $f\otimes g\in S^*(\gamma)$ where
$\gamma=\gamma(\alpha, \beta)<1$.
\end{enumerate}
One could ask the same question for different convolutions in place of $\otimes$.\\
\textit{(J. G. Krzyz)}\\

\noindent\textbf{Update 6.45} Barnard and Kellogg \cite{MR555840} have observed that $S^*(\alpha)$ is closed under $*$ and that $\gamma(\alpha,\beta)=\min[\alpha, \beta]$ and involves the Ruscheweyh--Sheil-Small theorem \cite{MR0328051}. Suppose that $f, g \in S^*(\alpha)$ and define
\[\phi=g*\log\frac{1}{1-z}.\]
 Then $\phi$ is convex and $f$ star-like, so that
\[\Big|\arg\Big(\frac{\phi*(zf')}{\phi*f}\Big)\Big|<\frac{\alpha\pi}{2}\]
as required.\\

\noindent\textbf{Problem 6.46} Suppose that $f$ in $S$ and is star-like. Is it true that
\begin{equation}
\label{C6.5}
\big||a_{n+1}|-|a_n|\big|\leq1?
\end{equation}
This is certainly true if $\lim_{r\to1} (1-r)M(r,f) > 0$, (D. A. Brannan, unpublished). T. Sheil-Small (uncited) has obtained an upper bound $2$ in (\ref{C6.5}). Notice that in addition to the `obvious' extremal functions $z(1-z^2)^{-1}$ and $z(1-z)^{-2}$ we have $z(1+z+z^2)^{-1}$.\\
\textit{(J. G. Clunie)}\\

\noindent\textbf{Update 6.46} Leung \cite{MR0480980} has solved this problem by proving that
\[\big||a_{n+1}|-|a_n|\big|\leq1,\hspace{1cm}n=1. 2, \ldots,\]
for all star-like functions $f$ in $S$. By a similar method, Leung \cite{MR534395} has also proved a conjecture of Robertson that
\[\big|n|a_n|-m|a_m|\big|\leq|n^2-m^2|,\hspace{1cm}n,m=1,2,\ldots,\]
for all close-to-convex functions $f$ in $S$. For the full class $S$, Grinspan \cite{grinspan1976} gives
\[-2.97\leq|a_{n+1}|-|a_n|\leq3.61,\hspace{1cm}n=1,2,\ldots.\]
Duren \cite{MR540059} has shown that
\[\big||a_{n+1}|-|a_n|\big|\leq e^\delta\alpha^{-1/2}<1.37\alpha^{-1/2},\hspace{1cm}n=1,2,\ldots,\] 
for all $f$ in $S$ with positive Hayman index $\alpha$, where $\delta$ is Milin's constant, $\delta<0.312$. Hamilton \cite{DHHamilton} has proved Leung's results independently.
 \\

\noindent\textbf{Problem 6.47} If $f$ in $S$ and $f'$ is also univalent in
$\mathbb{D}$, what can be said about $\max|a_n|$, $n \geq 2$? The function $z(1-z)^{-1}$ shows that
$\max|a_n|\geq1$ whenever $n\geq2$. \\
\textit{(J. G. Clunie)}\\

\noindent\textbf{Update 6.47} Define $S'=\{f:f\in S\text{ and }f'\text{ is univalent}\}$. It was noted by Barnard \cite{Barnardproblem} that if $f(z)=\sum a_nz^n\in S'$ and $A_2=\max_{f\in S'}|a_2|$, then \mbox{$|a_n|\leq 2A_2(n-1)/n$}, whenever the Bieberbach conjecture holds for $(f'-1)/2A_2$. It was conjectured that $A_2=2/(\pi-2)\approx1.75$. However, it has been shown by Barnard and Suffridge \cite{MR694924} that the function $F$ defined by
\[F'(z)=\frac{1+z}{(1-z)^2}\exp\Big\{\frac{1}{i\pi}\int^z_0\log\Big(\frac{1+iw}{1-iw}\Big)\,\frac{dw}{w}\Big\}=1+2B_2z+\ldots\]
is in $S'$ with $B_2=\frac{3}{2}+\frac{1}{\pi}\approx1.82$. $F$ has the property that it takes the left half of the unit circle onto a slit along the negative reals such that \mbox{$-\infty<F(i)\leq F(e^{i\theta})\leq F(-1)$} for $\pi/2\leq\theta\leq3\pi/2$; while $F'$ takes the right half of the unit circle onto a slit along the negative reals such that $-\infty\leq F'(e^{i\theta})\leq F'(i)<0$ for $-\pi/2\leq\theta\leq\pi/2$. If $C'=\{g\in S':g\text{ and }g'\text{ are convex}\}$, then Barnard and Suffridge \cite{MR694924} have shown that if $g(z)=\sum b_nz^n \in C'$ then $|b_n|\leq4n/3$. This is sharp since the function defined by $G'(z)=1+4z/3(1-z)$ is in $C'$.\\

\noindent\textbf{Problem 6.48} Suppose that $f$ in $S$. The coefficient problem, except in certain cases, remains
open for each of the following subclasses of univalent functions. (We limit ourselves to
one-parameter families).
\begin{enumerate} [(a)]
\item $B(\alpha)$, the class of Basilevi\^c functions comprising functions $f$ such that
\[f(z)=\Big[\int^z_0 p(t)s^\alpha(t)t^{-1}\,dt\Big]^{1/\alpha},\]
where $\alpha> 0$, $p(t) = 1+p_1t+\ldots$ is analytic and of positive real part in
$|t|<1$ and $s(t) = t + s_2t^2 +\ldots$ is star-like in $|t| < 1$.
\item $M(\alpha)$, the class of Mocanu-Reade functions comprising functions $f$ such that
\[\text{Re}\,\Big[\alpha\Big(1+\frac{zf''}{f'}\Big)+(1-\alpha)\frac{zf'}{f}\Big]>0,\]
where $0 <\alpha < 1$.
\item $S^*(\alpha)$, the class of strongly-star-like functions, comprising functions $f$ such that
\[\Big|\arg\Big(\frac{zf'}{f}\Big)\Big|<\frac{\alpha\pi}{2},\]
where $0 < \alpha < 1$.
\end{enumerate} 
\textit{(D. A. Brannan)}\\

\noindent\textbf{Update 6.48} No progress on this problem has been reported to us.\\

\noindent\textbf{Problem 6.49} What are the extreme points of the following classes of functions?
\begin{enumerate} [(a)]
\item Basilevi\^c functions (see Problem 6.48).
\item $S^*(\alpha)$ (see Problem 6.48).
\item Close-to-convex functions of order $\alpha$, $0 < \alpha < 1$.
\item Functions of boundary rotation $k\pi$, $2 < k < 4$.
\end{enumerate}
\textit{(D. A. Brannan)}\\

\noindent\textbf{Update 6.49} No progress on this problem has been reported to us.\\

\noindent\textbf{Problem 6.50} If $0 \le \alpha \le 1$, and $f(z)$, $g(z)\in \Sigma$, and if we define $F(z)$ by
\begin{eqnarray}
F(z)&=&f(z)^{1-\alpha}g(z)^\alpha, \hspace{1cm}|z|>1;\\
&=&z+\sum^\infty_{n=0}A_nz^{-n}, \hspace{1cm} |z|\text{ sufficiently large};
\end{eqnarray}
is it true that
\[\sum^\infty_{n=1}n|A_n|^2\le 1?\]
See Thomas \cite{MR0357765}. \\
\textit{(D. K. Thomas)}\\

\noindent\textbf{Update 6.50} The original statement of this problem contained errors which now have been corrected. No progress on this problem has been reported to us.\\

\noindent\textbf{Problem 6.51} Let $D$ be a domain in $\mathbb{C}$ (containing the origin) of connectivity $n$, and let
$S(D)$ be the class of analytic univalent functions in $D$ with $f(0) = 0$, $f'(0) = 1$. Find
the functions $f$ in $S(D)$ that minimise
\[\int_D\int|f'(z)|^2\,d\sigma_z.\]
\textit{(D. Aharonov)}\\

\noindent\textbf{Update 6.51} No progress on this problem has been reported to us.\\

\noindent\textbf{Problem 6.52} Suppose that $f(z)$ is analytic in $\mathbb{D}$, and has the whole complex plane as
its range. Does there necessarily exist a bounded univalent function
$g(z)$ in $\mathbb{D}$ such that $f(z)+g(z)$ has the whole complex plane as its range?\\ 
\textit{(L. A. Rubel)}\\

\noindent\textbf{Update 6.52} No progress on this problem has been reported to us.\\

\noindent\textbf{Problem 6.53}  Hentgartner and Schobe \cite{zbMATH03406294} and Goodman and Saff \cite{zbMATH03634586} have shown that
if $f(z) = z + a_2z^2 +\ldots$ maps $\mathbb{D}$ univalently onto a domain $G_1$ that is convex in the
direction of the imaginary axis (CIA), then $G_r = \{f(|z| < r)\}$ is not necessarily CIA
for all $r$, $r < 1$, or even for $r$ bigger than some constant.
\begin{enumerate} [(a)]
\item Find $\sup\{r : G_r \text{ is necessarily CIA}\}$. (A result of Goodman and Saff \cite{zbMATH03634586} suggests
that this is $\sqrt{2}-1$.)
\item Find reasonable sufficient conditions on $G_1$ (or, equivalently, on $f$) that
imply that $G_r$ is in CIA for all $r$ in $(0,1)$.
\item Suppose that we slit $G_1$ along the real axis and let $G'$, $G''$ be those two components
of the resulting family of domains that have $0$ on their boundaries. If
$G_1 = G' \cup G''$, does it necessarily follow that $G_r$ is CIA for all $r$ in $(0,1)$?
\item If $G_1$ is CIA but $G_{r_0}$ is not CIA, is it true that $G_r$ is not CIA for $r_0 < r < 1$?
\end{enumerate}
\textit{(D. A. Brannan)}\\

\noindent\textbf{Update 6.53}  Brown \cite{MR891142} proves a weaker form of the Goodman-Saff conjecture and answers Problem 6.43(d) by Brannan negatively. Ruscheweyh and Salinas \cite{MR997971} completely verify the Goodman-Saff conjecture.\\

\noindent\textbf{Problem 6.54} Let $D$ be a Jordan domain with boundary $C$, $\{F_n(z)\}^\infty_1$ the sequence of
Faber polynomials for $D$, and $S(D)$ the class of univalent functions
\[f(z)=F_1(z)+\sum^\infty_{n=2}a_nF_n(z)\]
in $D$. Does the coefficient region of at least one $a_n$ have the same shape as $C$, or is it
at least in the subclass of star-like functions? Royster \cite{MR0107008} has shown that, if $f$ is
star-like in an ellipse, then there exists a direction $\theta_f$ and sequences
\[\{\lambda_n\}^\infty_2,\hspace{1cm}\{\mu_n\}^\infty_2,\hspace{1cm}\lambda_n>\mu_n>0,\]
such that each coefficient $a_n$ lies in an ellipse of centre $0$, major axis $\lambda_n$, minor axis
$\mu_n$, inclined at an angle $\theta_f$ to the real axis. (The $\lambda_n, \mu_n$ are known.)\\ 
\textit{(K. W. Lucas)}\\

\noindent\textbf{Update 6.54} No progress on this problem has been reported to us.\\

\noindent\textbf{Problem 6.55} Let $f(z)$ be a normalised bounded star-like function in $\mathbb{D}$, and set
\[f(\xi)=\lim_{r\to1-}f(r\xi),\]
where $|\xi|=1$, $\xi\in E$, $E\subset\{|z|=1\}$. Is it true that $f(E) = \{f(\xi):\xi\in E\}$ has zero linear
(= one-dimensional Hausdorff) measure if $\text{cap } E = 0$, or at least if $E$ has zero logarithmic
measure? (It is known that $\text{meas} \{\arg f(\xi):\xi\in E\}=0$ if $\text{cap } E=0$.)\\
\textit{(Ch. Pommerenke)}\\

\noindent\textbf{Update 6.55} No progress on this problem has been reported to us.\\

\noindent\textbf{Problem 6.56} Let $S_R(q)$ be the class of normalised univalent functions in $\mathbb{D}$ with real
coefficients that admit a quasi-conformal extension to the whole plane with complex
dilatation bounded (in modulus) almost everywhere by $q$, $q < 1$. Following the notation in Problem 6.44, prove (or disprove)
that, if $f\in S_R(p)$ and $g\in S_R(q)$, then $f\otimes g\in S_R(r)$.\\ 
\textit{(J. G. Krzyz)}\\

\noindent\textbf{Update 6.56} No progress on this problem has been reported to us.\\

\noindent\textbf{Problem 6.57} It is known that functions $f$ in $\Sigma(q)$, i.e. those $f$ analytic and univalent in
$\mathbb{D}$, normalised by $f(z) = z + b_0 + b_1z^{-1}+\ldots$ and having a quasi-conformal
extension to the whole plane with $|f_{\overline{z}}/f_z|\leq q<1$ almost everywhere, satisfy
\begin{equation}
\label{C6.5a}
{\int\int}_{|\zeta|>1}|u(z,\zeta)|^2\,d\sigma_\zeta\leq\frac{\pi q^2}{(|z|^2-1)^2},
\end{equation}
where
\[u(z,\zeta)=\frac{f'(z)f'(\zeta)}{f(z)-f(\zeta)}-\frac{1}{(z-\zeta)^2}.\]
Prove (or disprove) that (\ref{C6.5a}) is also sufficient for $f$ to have a quasiconformal extension
to the whole plane, possibly for $q$ sufficiently small.\\ 
\textit{(J. G. Krzyz)}\\

\noindent\textbf{Update 6.57} This has been answered affirmatively by Zuravlev \cite{MR517191}.\\

\noindent\textbf{Problem 6.58} Following the notation in Problem 6.57, the well-known Golusin inequality for functions $f$ in $\Sigma(q)$ (defined in
Problem 6.57) is:
\begin{equation}
\label{C6.6}
\Big|\log\frac{f'(z)f'(\zeta)(z-\zeta)^2}{[f(z)-f(\zeta)]^2}\Big|\leq q\log\frac{|z\overline{\xi}-1|^2}{(|z|^2-1)(|\zeta|^2-1)}.
\end{equation}
Prove (or disprove) that (\ref{C6.6}) is also sufficient for $f$ to have a quasi-conformal
extension to the whole plane, possibly for $q$ sufficiently small.\\ 
\textit{(J. G. Krzyz)}\\

\noindent\textbf{Update 6.58} No progress on this problem has been reported to us.\\

\noindent\textbf{Problem 6.59} Let $D$ be a plane domain containing $\infty$. Let there be given a continuous
assignment of numbers (thought of as angles) to the components of $\mathbb{C}\setminus D$. Consider
conformal mappings of $D$ onto the complement in the extended plane of straight line
segments. Show that one of these mappings is such that the straight line segments
make angles with the positive real axis equal to the corresponding preassigned
angles. In other words, under the mapping, each component of $\mathbb{C}\setminus D$ is associated
with a slit in its preassigned direction.\\
\textit{(B. Rodin; communicated by C. FitzGerald)}\\

\noindent\textbf{Update 6.59} Fitzgerald and Weening \cite{MR1694289} show that Rodin's continuity assumption is sufficient to imply the existence of a rectilinear slit map achieving any angle assignment with finite range.\\

\noindent\textbf{Problem 6.60} Let $C$ be a closed Jordan curve. Then if $f(z) = z + a_2z^2 +\ldots$,
$g(z) =z^{-1}+b_0 +b_1z + \ldots$ map $\mathbb{D}$ onto the inside and outside of $C$ respectively,
the area principle shows that $C$ is the unit circumference $\mathbb{T}$. If we remove the normalisation on $g$ and replace $g$ by 
$g_1(z) = b_{-1}z^{-1} + b_0 + b_1z + \ldots$, what can be said about
the connection between $f$ and $g_1$? For example, what about the asymptotic behaviour
of their coefficients, or their Lipschitz behaviour on $\mathbb{T}$, or about prime ends
(if we make $C$ slightly non-Jordan)?\\ 
\textit{(D. A. Brannan)}\\

\noindent\textbf{Update 6.60} Lesley \cite{MR545203} has obtained best possible results on the Lipschitz behaviour.\\

\noindent\textbf{Problem 6.61} Let $D_1, D_2$ be Jordan domains bounded by rectifiable curves $C_1, C_2$ of equal
length. Suppose that an isometric sewing of $C_1$ and $C_2$ is everywhere conformally
admissible, thus generating a Riemann surface equivalent to a sphere $S$. Is the curve
$C$ on $S$ which corresponds to $C_1$ (and $C_2$) necessarily rectifiable?\\ 
\textit{(A. Huber)}\\

\noindent\textbf{Update 6.61}  No progress on this problem has been reported to us.\\

\noindent\textbf{Problem 6.62} Let $D_1$ and $D_2$ be bounded Jordan domains, bounded by curves $C_1$ and $C_2$
of bounded boundary rotation (in the sense of Paatero, see e.g. Noonan \cite{NOONAN1972721}); then $C_1$ and $C_2$ are rectifiable,
and we shall assume that they have the same length. It is known that in this
case, every isometric sewing is conformally admissible and generates a Riemann surface
which is equivalent to a sphere $S$. It follows from results of Aleksandrov (no citation) and
Reshetnjak (no citation) that the curve $C$ on $S$ which corresponds to $C_1$ (and $C_2$) is of
bounded boundary rotation. Can one find a function-theoretic proof of this?\\
\textit{(A. Huber)}\\

\noindent\textbf{Update 6.62} No progress on this problem has been reported to us.\\

\noindent\textbf{Problem 6.63} Let $\alpha$ be a homeomorphic mapping of $(0, \infty)$ onto $(\alpha(0), \infty)$, $\alpha(0)\geq 0$,
such that $x\to\alpha(x)+i$ defines a conformal sewing of the half-strip
\[H = \big\{z:\text{Re}\, z > 0, 0 \leq \text{Im}\,(z) \leq 1\big\}.\]
If $\alpha$ is hyperbolic (that is, if the generated Riemann surface has a hyperbolic end at $\infty$),
does there exist a positive continuous function $\varepsilon: (0,\infty)\to\mathbb{R}$ with the
property that each conformal sewing $x\to\beta(x)+i$ of $H$ satisfying the inequality
\[|\beta(x)-\alpha(x)|<\varepsilon(x),\hspace{1cm}0<x<\infty,\]
is also hyperbolic?\\ 
\textit{(C. Constantinescu; communicated by A. Huber)}\\

\noindent\textbf{Update 6.63} No progress on this problem has been reported to us.\\

\noindent\textbf{Problem 6.64} Let $\alpha$ be real and suppose that $f(z)=z+\sum^\infty_{n=2}a_nz^n$ is analytic in $\mathbb{D}$ with $f(z)f'(z)/z\neq0$. We say $f$ is in $M_\alpha$, the \textit{class of alpha-convex functions}, if 
\[\text{Re}\,\Big[(1-\alpha)\frac{zf'}{f}+\alpha\Big(1+\frac{zf''}{f'}\Big)\Big]>0\]
for $|z|<1$. Note that $M_0=S^*$, the class of star-like functions and $M_1=K$, the class of convex functions. It is known that $M_\alpha\subset S^*$ for all $\alpha$.\\
Clunie and Keogh \cite{MR0110814} have shown the following:
\begin{enumerate}[(a)]
\item If $\sum^\infty_{n=2}n|a_n|<1$ then $f\in M_0$
\item If $\sum^\infty_{n=2}n^2|a_n|<1$ then $f\in M_1$.
\end{enumerate}
What generalisation of these conditions implies that $f$ belongs to $M_\alpha$?\\
\textit{(S. Miller)}\\

\noindent\textbf{Update 6.64} No progress on this problem has been reported to us.\\

\noindent\textbf{Problem 6.65} Given $M$, $1<M<\infty$, let $S^*(M)$ be the class of star-like univalent functions $f$ in $\mathbb{D}$ with $f(0)=0$, $f'(0)=1$, and $|f(z)|\leq M$ for $|z|<1$. If $f(z)=z+\sum^\infty_{n=2}a_nz^n$, find $\sup\{|a_3|:f\in S^*(M)\}$ for $e<M<5$. Barnard and Lewis \cite{MR0379816} have found the supremum for other values of $M$.\\
\textit{(J. Lewis)}\\

\noindent\textbf{Update 6.65} No progress on this problem has been reported to us.\\

\noindent\textbf{Problem 6.66} Describe the extreme points of the class $\Sigma_0$ consisting of all functions $g$ in $\Sigma$ with constant term $b_0=0$. Springer (see Pommerenke \cite[p. 183]{MR0507768}) showed that every $g$ in $\Sigma_0$ whose omitted set has measure zero is an extreme point. Is this condition also necessary?\\
\textit{(P. L. Duren)}\\

\noindent\textbf{Update 6.66} Hamilton \cite{MR656110} answers this question positively.\\

\noindent\textbf{Problem 6.67} Let $f$ be univalent in $\mathbb{D}$ and let $f(\mathbb{D})$ be a Jordan domain. Does the condition
\begin{equation}
\label{6.67}
\limsup_{|z|\to1}(1-|z|^2)|f''(z)/f'(z)|<2
\end{equation}
imply that $f$ has a quasi-conformal extension over the unit circle?\\
It is known that, for $|c|<1$, the condition 
\[\limsup_{|z|\to1}\big|(1-|z|^2)zf''(z)/f'(z)-c\big|<1\]
is sufficient for $f$ to have a quasiconformal extension. It is also known that (\ref{6.67}) is sufficient if $f$ satisfies in addition
\[\limsup_{|z|\to1}(1-|z|^2)^2|S_f(z)|<2,\]
where $S_f$ denotes the Schwarzian derivative. It might also be asked whether the last condition alone is already sufficient. In both cases, the constant $2$ on the right-hand side would be best possible.\\
\textit{(J. Becker)}\\

\noindent\textbf{Update 6.67} The first question has been answered positively by Gehring and Pommerenke \cite{MR749106}.\\

\noindent\textbf{Problem 6.68} Let $\Sigma$ be the class of univalent functions in $\{|z|>1\}$ with the usual normalisation $f(z)=z+\sum^\infty_{n=0}b_nz^{-n}$. Let $S_f$ denote the Schwarzian derivative of $f$. Does
\[\sup_{|z|>1}(|z|^2-1)^2\big|S_{f_n}(z)-S_f(z)\big|\to0\hspace{1cm}\text{ as }n\to\infty\]
imply that 
\[\sup_{|z|>1}(|z|^2-1)\Big|\frac{f_n''(z)}{f_n'(z)}-\frac{f''(z)}{f'(z)}\Big|\to0\hspace{1cm}\text{ as }n\to\infty,\]
if $f, f_n\in\Sigma$ $(n=1, 2, \ldots)$?\\
This is known to be true if $f$, $f_n$ have quasiconformal extensions onto the plane.\\
\textit{(J. Becker)}\\

\noindent\textbf{Update 6.68} No progress on this problem has been reported to us. \\ 

\noindent\textbf{Problem 6.69} Let $B$ be the Banach space of analytic functions $\phi$ in $\{|z|>1\}$ with finite norm
\[\|\phi\|:=\sup_{|z|>1}(|z|^2-1)|z\phi(z)|.\]
Let $S$ and $T$ be the subsets defined by $S:=\{f''/f':f\in \sum\}$ and $T:=\{f''/f':f\in\sum, f\text{ has a quasiconformal extension to }\mathbb{C}\}$.\\
It is known that $T$ is topologically equivalent to the universal Teichm\"uller space. From results of Ahlfors and Gehring (no citations), it follows that $T$ is a subdomain in $B$ and that $S\setminus\overline{T}\neq\emptyset$, where $\overline{T}$ denotes the closure of $T$.\\
Is it also true (analogous to another result of Gehring (no citation)) that $T=S^0$, where $S^0$ denotes the interior of $S$?\\
This would follow if the answer to Problem 6.68 were affirmative. The problem is closely related to a characterisation of quasicircles given by Gehring.\\
\textit{(J. Becker)}\\

\noindent\textbf{Update 6.69} This has been proved by Astala and Gehring \cite{MR861687}.\\

\noindent\textbf{Problem 6.70} Is every extreme point of $S$ a support point? Is every support point an extreme point?\\
\textit{(P. L. Duren)}\\

\noindent\textbf{Update 6.70} Hamilton \cite{MR915532} has constructed examples of extreme points which are not support points
and a different construction has been given by Duren and Leung \cite{MR861691}. The converse
problem seems still to be open. \\ 

\noindent\textbf{Problem 6.71} For each $f$ in $S$, it can be shown that
\[\int^{2\pi}_0\Big|\frac{f'(Re^{i\theta})}{f(Re^{i\theta})}\Big|^2d\theta=O\Big(\frac{1}{1-R}\log\frac{1}{1-R}\Big)\]
as $R\to1$. Hayman \cite{MR568639} has constructed an example showing that `$O$' cannot be replaced by `$o$'. For certain subclasses of $S$, such as the star-like functions and the functions with positive Hayman index (including all support points of $S$), the estimate can be improved to $O(1/(1-R))$. Can the same improvement be made for the extreme points of $S$?; for close-to-convex functions?\\
\textit{(P. L. Duren)}\\

\noindent\textbf{Update 6.71} No progress on this problem has been reported to us.\\

\noindent\textbf{Problem 6.72} Let $\Gamma$ be the analytic arc omitted by a support point of $S$. Must $\Gamma$ have monotonic argument? Must the angle between the radius and tangent vectors be monotonic on $\Gamma$? (The second property implies the first. Brown \cite{MR766220} has shown that the support points associated with point-evaluation functions have both properties.)\\
\textit{(P. L. Duren)}\\

\noindent\textbf{Update 6.72} No progress on this problem has been reported to us.\\

\noindent\textbf{Problem 6.73} Let $f(z)=z+\sum^\infty_{n=2}a_nz^n$ be in $S$. Is it true that
\[\limsup_{n\to\infty}\big||a_{n+1}|-|a_n|\big|\leq1?\]
Hamilton \cite{DHHamilton} has proved that this is true for odd functions in $S$, functions of maximal growth in $S$, and spiral-like functions.\\
\textit{(D. H. Hamilton)}\\

\noindent\textbf{Update 6.73} Hamilton \cite{MR732121} has proved more precisely that \mbox{$\big||a_{n+1}|-|a_n|\big|\leq 1$}, $n>n_0(f)$, with equality only for 
\[f(z)=\frac{z}{(1-ze^{-\theta_1})(1-ze^{i\theta_2})},\] 
when $(\theta_2-\theta_1)$ is a rational multiple of $\pi$.\\

\noindent\textbf{Problem 6.74} Suppose $f(z)=z+\sum^\infty_{n=2}a_nz^n$ is univalent and bounded by $M$ in $\mathbb{D}$. Find 
\[\sup_t \max_{0\leq t\leq 2\pi}|s_n(e^{it})|,\]
where
\[s_n(z)=z+\sum^n_{k=2}a_kz^k.\]
\textit{(F. Holland)}\\

\noindent\textbf{Update 6.74} No progress on this problem has been reported to us.\\

\noindent\textbf{Problem 6.75} Let $\mathcal{P}_n$ be the class of polynomials
\[P_n(z)=z+a_2z^2+\ldots+a_nz^n\]
univalent in $\mathbb{D}$, and let 
\[A_m(n)=\max_{\mathcal{P}_n}|a_m|.\]
If $n$ is fixed, is it true that as $m$ increases, the quantity $A_m(n)$ increases strictly up to some index $n_0$, and then decreases strictly? Again, for fixed (but arbitrary) $n$ determine the least $n_1$ such that
\[A_m(n)\leq1,\hspace{1cm}n_1\leq m\leq n.\]
\textit{(D. A. Brannan)}\\

\noindent\textbf{Update 6.75} No progress on this problem has been reported to us.\\

\noindent\textbf{Problem 6.76} Let $\mathcal{V}_n$ denote the class of polynomials
\[P_n(z)=z+a_2z^2+\ldots+a_nz^n\]
analytic and bi-univalent in $\mathbb{D}$ (that is, $P_n$ and $P^{-1}_n$ are both univalent in the unit disc). Determine $\max_{\mathcal{V}_n}|a_2|$ and $\max_{\mathcal{V}_n}|a_n|$.\\
\textit{(D. A. Brannan)}\\

\noindent\textbf{Update 6.76} No progress on this problem has been reported to us.\\

\noindent\textbf{Problem 6.77} Let $\mathcal{P}_n$ be the class of polynomials
\[p_n(z)=z+a_2z^2+\ldots+a_nz^n\]
univalent in $\mathbb{D}$. Determine 
\[\max_{p\in\mathcal{P}_n}\int^{2\pi}_0|p(e^{it})|^q\,dt,\hspace{1cm}0<q<\infty.\]
\textit{(F. Holland)}\\

\noindent\textbf{Update 6.77} No progress on this problem has been reported to us.\\

\noindent\textbf{Problem 6.78} Suppose that $f$ in $S$. Consider the region $\mathbb{D}(f)$ on the Riemann sphere which is the stereographic projection of the image of the unit disc under $f$. We can associate with each $f$ in $S$ the spherical area of $\mathbb{D}(f)$. What is
\[\min_{f\in S}\{\text{area of }\mathbb{D}(f)\},\]
and what is the extremal function?\\
\textit{(Y. Avci)}\\

\noindent\textbf{Update 6.78} No progress on this problem has been reported to us.\\

\noindent\textbf{Problem 6.79} Let $S_k(\infty)$ denote the class of all analytic and univalent functions $f(z)=z+a_2z^2+\ldots$ defined in $\mathbb{D}$ which admit a $k$-quasiconformal extension $(0<k<1)$ to the whole plane, with $f(\infty)=\infty$. Prove or disprove:
\[f(z)\in S_k(\infty)\hspace{0.5cm}\implies\hspace{0.5cm}\frac{f(rz)}{r}\in S_k(\infty),\hspace{1cm}0<r<1.\]
\textit{(D. Bshouty)}\\

\noindent\textbf{Update 6.79} This has been proved by Krushkal \cite{MR924982}.\\

\noindent\textbf{Problem 6.80} If $f$ is univalent analytic in $\mathbb{D}$, then it is well known (see Pommerenke \cite[p. 262]{MR0507768}) that both $f$ and its first derivative $f'$ must be normal, while the higher derivatives $f^{(n)}$ need not be normal if $n\geq2$. Setting $f^{(-1)}(z)=\int^z_0f(t)\,dt$ and $f^{(-n-1)}(z)=\int^z_0f^{(-n)}(t)\,dt$, it is easy to verify that $f^{(-2)}$ is Bloch; while if $n\geq3$, $f^{(-n)}$ is bounded (hence Bloch, and therefore normal). Thus, if $f$ is univalent analytic in $\mathbb{D}$, then the functions $f^{(n)}$ must be normal for $n=1, 0, -2, -3,\ldots,$ and need not be normal for $n=2, 3, \ldots$. Must $f^{(-1)}(z)=\int^z_0f(t)\,dt$ be normal if $f(z)$ is univalent analytic in $\mathbb{D}$?\\
\textit{(D. Campbell)}\\

\noindent\textbf{Update 6.80} No progress on this problem has been reported to us.\\

\noindent\textbf{Problem 6.81} Let $G$ be the set of functions analytic and \textit{not} univalent in $\mathbb{D}$. Set, for $f\in G$,
\[M_f=\sup\{|f'(z)|:|z|<1\},\hspace{1cm}m_f=\inf\{|f'(z)|:|z|<1\}\]
and put 
\[\gamma=\inf\Big\{\frac{M_f}{m_f}:f\in G\Big\}.\]
Find $\gamma$.\\
John \cite{MR0422606} has proved that $\gamma\geq e^{\pi/2}\approx4.7$ and Yamashita \cite{MR0507738} has shown that $\gamma\leq e^\pi\approx 23.1$.\\
\textit{(G. M., J. M. and T. M. Rassias)}\\

\noindent\textbf{Update 6.81} No progress on this problem has been reported to us.\\

\noindent\textbf{Problem 6.82} The above definition of a bi-univalent function is difficult to understand. What is also meant is that the inverse function $f^{-1}$ has an analytic continuation to the unit disc where it is univalent. Let $\sigma$ denote the class of bi-univalent functions, namely the class of functions $f(z)=z+a_2z^2+\ldots$ analytic and univalent in $\mathbb{D}$, such that their inverses $f^{-1}$ are also analytic univalent in the $\mathbb{D}$. Brannan conjectures that
\[a^*_2\equiv\sup_\sigma|a_2|=\sqrt{2}.\]
It is known that $a^*_2<1.51$ (see Lewin \cite{MR0206255}) and that $a^*_2>\frac{4}{3}+0.02$ (see Styer and Wright \cite{MR609659}).\\
Prove or disprove the statement that, if $f\in\sigma$, then
\[d_n\equiv\big||a_{n+1}|-|a_n|\big|\leq A\]
for some absolute constant $A$ with $0<A<1$.\\
Alternatively, what is $\sup_\sigma d_n$ for $n\geq2$?\\
\textit{(D. A. Brannan)}\\

\noindent\textbf{Update 6.82} Brannan's conjecture is verified for some classes by Sivasubramanian, Sivakumar, Kanas, and Kim \cite{MR3320475}.\\

\noindent\textbf{Problem 6.83} Let $S$ be the usual class of normalised univalent functions in the unit
disc $\mathbb{D}$. Characterise those sequences $\{z_n\}$ of points in $\mathbb{D}$ for which $f(z_n) = g(z_n)$ for
two different functions $f, g$ in  $S$.\\
Notice that a necessary condition is that $\sum^\infty_{n=1}(1-|z_n|) < \infty$, because $(f-g)\in H^p$
for all $p <\frac{1}{2}$.\\
\textit{(P. L. Duren)}\\

\noindent\textbf{Update 6.83} Overholt \cite{MR1749956} has provided a partial solution to this problem.\\

\noindent\textbf{Problem 6.84} If $f(z)$ in $S$, write
\[\log\frac{f(z)}{z}=2\sum^\infty_{n=1}\gamma_nz^n\]
and
\[f(z^p)^{1/p}=z+\sum^\infty_{n=1}c^{(p)}_nz^{pn+1}\hspace{1cm}(p=1,2,\ldots).\]
Szeg\"o's conjecture asserts that $c^{(p)}_n=O(n^{2/p-1})$ as $n\to\infty$. This has been proved for
$p = 1,2$ and $3$, but Pommerenke \cite{MR0220920} has shown that it is false for $p\geq 12$; his
example also has $\gamma_n\neq O(1/n)$. Milin \cite{MR804182} has shown that, if a function $f$ in
$S$ has the property that $\gamma_n=O(1/n)$, then $c^{(p)}_n=O(n^{2/p-1})$ for every $p$. Is the converse
true?\\
\textit{(P. L. Duren)}\\

\noindent\textbf{Update 6.84} Szeg\"o's conjecture has been proved by Baernstein \cite{MR846491} for $p=4$. The cases $5\leq p\leq11$ remain open.\\

\noindent\textbf{Problem 6.85} Each function $f$ in S that maximises $\text{Re}\, \{L(g) :g \in S\}$ for some continuous
linear functional $L$ must map the unit disc onto the complement of an arc $\Gamma$ that is
asymptotic at infinity to the half-line
\[w=\frac{L(f^3)}{3L(f^2)}-L(f^2)t,\hspace{1cm}t\geq0.\]
For the coefficient functional $\Lambda_n(f)=a_n$, show that the asymptotic half-line is
radial; that is, that $\Lambda_n(f^3)/[\Lambda_n(f^2)]^2$ is real. (This is true for $n = 2, 3, 4, 5, 6$.)\\
\textit{(P. L. Duren)}\\

\noindent\textbf{Update 6.85} See Update 6.1.\\

\noindent\textbf{Problem 6.86} Sundberg notes that it is well known fact (see Hayman \cite[p. 5]{MR0108586}) that, for each fixed $z_0$ in $\mathbb{D}$, 
\[\Big|z_0\frac{f''(z_0)}{f'(z_0)}-\frac{2\rho^2}{1-\rho^2}\Big|\leq\frac{4\rho}{1-\rho^2},\]
 and asks if this can be improved, if we have that $f$ is real on the real axis, or equivalently, if $f$ has real coefficients?\\
\textit{(C. Sundberg; communicated by P. L. Duren)}\\

\noindent\textbf{Update 6.86} This problem has been restated. No progress on this problem has been reported to us.\\

\noindent\textbf{Problem 6.87} Let $L_1$, $L_2$ be two complex-valued continuous linear functionals on $H(\mathbb{D})$,
the space of all analytic functions on the unit disc $\mathbb{D}$, that are not constant on the set
$S$ of all normalised univalent functions on $\mathbb{D}$. Assume, in addition, that $L_1 \neq  tL_2$ for
any positive $t$. If a function $f$ in $S$ maximises both $\text{Re}\, \{L_1\}$ and $\text{Re}\, \{L_2\}$ in $S$, must $f$ be a
rotation of the Koebe function?\\
\textit{(P. L. Duren)}\\

\noindent\textbf{Update 6.87} Bakhtin \cite{MR865514} has answered this positively for certain classes of coefficient functionals.\\

\noindent\textbf{Problem 6.88} Let the function $f = z + a_2z^2 + \ldots$ in $S$ map $\mathbb{D}$ onto a domain
with finite area $A$. Then Bieberbach's inequality $|a_2|\leq2$ can be sharpened to the
following:
\begin{equation}
\label{J6.88}
|a_2|\leq2-cA^{-1/2}
\end{equation}
where $c$ is an absolute constant. What is the best value of $c$?\\
Aharonov and Shapiro \cite{MR0412400} have shown that (\ref{J6.88}) holds for some $c$, and have a
conjecture concerning the sharp constant $c$ and the extremal function for (\ref{J6.88}). They
also conjecture that $|a_2|\leq 2-c_1l^{-1}$, where $l$ is the length of $\partial f(\mathbb{D})$. For some
background information see Aharonov and Shapiro (\cite{MR0412400} and \cite{AharonovShapiro}) and Abarov, Shapiro and Solynin  (\cite{zbMATH01369959} and \cite{zbMATH01640257}).\\
\textit{(H. S. Shapiro)}\\

\noindent\textbf{Update 6.88} No progress on this problem has been reported to us.\\

\noindent\textbf{Problem 6.89} Let $S^*(\frac{1}{2})$ denote the class of functions $g$ analytic in $\mathbb{D}$ and such that
$\text{Re}\, (zg'/g) > \frac{1}{2}$ in $\mathbb{D}$.\\
Is it true that, if $f\in S^*(\frac{1}{2})$, $r\in (0,1)$ and $\theta\in[0,2\pi)$, then
\begin{equation}
\label{J6.89}
\frac{1}{|f(re^{-i\theta})|}\int^r_0|f'(te^{i\theta})|\,dt\leq\frac{\arcsin r}{r}?
\end{equation}
(The left-hand-side of (\ref{J6.89}) is the ratio of the length of the image of a radius and the
distance between the endpoints of that image (as the crow flies).) Inequality (\ref{J6.89}) is
true when $f(z) = z/(1-z)$; in addition, the left-hand-side of (\ref{J6.89}) never exceeds $\pi/2$ for
any $r$ or $f$.\\
Let $z_1,z_2$ lie in $\mathbb{D}$. Then, in the smaller class of convex functions $f$, what
can be said about
\[\frac{1}{|f(z_1)-f(z_2)|}\int^{z_2}_{z_1}|f'(t)|\,|dt|?\]
\textit{(R. R. Hall)}\\

\noindent\textbf{Update 6.89} No progress on this problem has been reported to us.\\

\noindent\textbf{Problem 6.90} Let $E$ be a set of positive logarithmic capacity on the unit circle $\mathbb{T}$. Is $E$
necessarily a set of uniqueness for functions univalent in the unit disc $\mathbb{D}$?\\
Carleson \cite{MR0050011} has shown that this is false for functions $f$ analytic in $\mathbb{D}$ for
which ${\int\int}_{|z|<1}|f'(z)|^2\,dx\,dy<\infty$. Beurling \cite{MR0001370} has shown that univalent functions
cannot have constant boundary values on a set of positive capacity.\\
\textit{(D. H. Hamilton)}\\

\noindent\textbf{Update 6.90} Hamilton \cite{}[EMAILLED 11/4/18] has shown that this is false: there exists a set of positive logarithmic capacity on $\mathbb{T}$ which is not a set of uniqueness for functions univalent in $\mathbb{D}$.\\

\noindent\textbf{Problem 6.91} Let $\Omega$ be an arbitrary domain in $\mathbb{C}$. Does there necessarily exist a set
$E$ in $\partial\Omega$, of full harmonic measure, with the following property: for each $z$ in $E$, there
exist circular arcs $C_r$ in $\Omega$, of radius $r$ (where $r$ is small) and centred at $z$, for which
\[\lim_{R\to0}\Big(\frac{1}{\pi R^2}\int^R_0\theta(C_r)r\,dr\Big)=\frac{1}{2},\]
where $\theta(\,\cdot\,)$ denotes angular measure?\\
For simply-connected domains $\Omega$, this is a theorem of McMillan \cite{MR0294614}. For general
domains $\Omega$, at least there is a sequence $\{r_n\}$ decreasing to $0$ for which $\theta(C_{r_n})$
approaches $\pi/4$.\\
\textit{(D. Stegenga and K. Stephenson)}\\

\noindent\textbf{Update 6.91} No progress on this problem has been reported to us.\\

\noindent\textbf{Problem 6.92} If $\mathbb{R}^2_+=\{(x,y)\in\mathbb{R}^2:y>0\}$, suppose that $E\subset\mathbb{R}^2_+$, and let $f:\mathbb{R}^2_+\to B^2$ be analytic and conformal where $B^2 = \{(x, y)\in\mathbb{R}^2_+: x^2 + y^2 < 1\}$; assume also that
\[\lim_{x\to0,\,x\in E}f(x)=\alpha.\]
Lindel\"of's theorem shows that $\alpha$ is an angular limit of $f$ at $0$ if $E$ is a curve
terminating at $0$. Vuorinen \cite{MR608334} has shown that the same conclusion holds under
much weaker hypotheses on $E$; for instance, the condition that
\[\liminf_{r\to0}m(A\cap(0,r))/r>0,\]
where $A = \{|x|: x \in E\}$, is sufficient (see Vuorinen \cite[pp. 169, 176]{MR608334}). This problem concerns a
converse result.\\
Denote by $\mathcal{F}$ the class of all analytic and conformal maps of $\mathbb{R}^2_+$ into $B^2$. Let $K$
be a subset of $\mathbb{R}^2_+$ with the property that, whenever
\[\lim_{x\to0,\,x\in E}f(x)=\alpha\hspace{1cm}\text{ and }\hspace{1cm}f\in\mathcal{F},\]
then necessarily $f$ has an angular limit $\alpha$ at $0$. Since $f$ is conformal it follows
(see Vuorinen \cite[p. 176]{MR608334}) that
\[\lim_{x\to0,\,x\in K_1}f(x)=\alpha,\]
where $K_1 = \{(x,y)\in\mathbb{R}^2_+ : \rho((x,y), K) < l\}$ and $\rho$ is the Poincar\'e metric of $\mathbb{R}^2_+$.
What can be said about the thickness of $K_1$ at $0$? Is it true that
\begin{equation}
\label{J6.92}
\limsup_{r\to0}m\big(B^\prime\cap(0,r)\big)/r>0,
\end{equation}
where $B^\prime = \{(x^2 + y^2)^{1/2} : (x,y)\in K_1\}$? It is easy to see that (\ref{J6.92}) must hold if $K$ is
contained in an angle with vertex at $0$ whose closure lies inside $\mathbb{R}^2_+$.\\
\textit{(M. Vuorinen)}\\

\noindent\textbf{Update 6.92} No progress on this problem has been reported to us.\\

\noindent\textbf{Problem 6.93} Let the function $f(z) = z + a_2z^2 + \ldots$ map $\mathbb{D}$ univalently onto a
domain $\Omega$, and let $F: \Omega\to\mathbb{D}$ denote the inverse of $f$. Is it true that
\[\int_{\Omega\cap \mathbb{R}}|F'(x)|^p\,dx<\infty\hspace{1cm}\text{ for }1\leq p<2?\]
Hayman and Wu \cite{MR639358} and Garnett, Gehring and Jones \cite{MR721565} have shown that the
answer is `yes' for $p = 1$. For $p = 2$, it is possible that $\int_{\Omega\cap\mathbb{R}}|F'(x)|^2\,dx = \infty$; for
example, when $\Omega = \{w: |w| < R, R > 1\}\setminus L$ where, for a suitably chosen $R_1$,
$L = \{(u,0): -R\leq u \leq-R_1 < 0\}$.\\
\textit{(A. Baernstein II)}\\

\noindent\textbf{Update 6.93} No progress on this problem has been reported to us.\\

\noindent\textbf{Problem 6.94} Let $\Omega$ be a simply-connected domain in $\mathbb{C}$ with at least two boundary
points, and let the function $\phi$ map $\Omega$ analytically and conformally onto $\mathbb{D}$.
For which values of $p$ is it true that
\begin{equation}
\label{J6.94}
{\int\int}_\Omega|\phi'|^p\,dx\,dy<\infty?
\end{equation}
If $\Omega$ is star-like or close-to-convex, (\ref{J6.94}) holds for $\frac{4}{3} < p < 4$ and this is sharp. More
generally, it is known that there is a universal constant $\tau$, independent of $\Omega$, with
$0 < \tau < 1$, such that (\ref{J6.94}) holds whenever $\frac{4}{3} < p < 3 + \tau$. Is $\frac{4}{3}< p < 4$ the correct
range for general types of $\Omega$? (For background material and additional information,
see \cite{MR509942}.)\\
\textit{(J. E. Brennan)}\\

\noindent\textbf{Update 6.94} Pommerenke \cite{MR813740} proved that (\ref{J6.94}) is true for $\frac{4}{3}<p<4$. \\

\noindent\textbf{Problem 6.95} Determine an intrinsic characterisation for the class $\mathcal{H}$ of functions $h$
analytic in $\mathbb{D}$ that admit a decomposition of the form $2h = f + f^{-1}$ for some
function $f$ in the class $S$ of normalised univalent functions in $\mathbb{D}$. (Here $f^{-1}$ denotes the function inverse to $f$ and it is assumed that $f^{-1}$ has an analytic continuation to $\mathbb{D}$.) Notice that the functions $z\mapsto z$ and $z\mapsto z/(1-z^2)$ both belong to $\mathcal{H}$. In particular, does the function 
$z \mapsto z + a \sin (2\pi z)$ belong to $\mathcal{H}$ for any non-zero constant $a$?\\
\textit{(F. Holland)}\\

\noindent\textbf{Update 6.95} No progress on this problem has been reported to us.\\

%\section* *{Univalent and multivalent functions}

\noindent\textbf{Problem 6.96} For $-\infty<p<+\infty$ let
\[B(p):=\sup\{\beta_f(p):f\text{ conformal map of }\mathbb{D}\text{ into }\mathbb{D}\}\]
where
\[\beta_f(p)=\limsup_{r\to1}\Big(\int_{|\zeta|=r}|f'(r\zeta)|^p\,|d\zeta|\Big)/\log\frac{1}{1-r}.\]
The  BCJK-conjecture states that
\[
B(p)=
\begin{cases}
-p-1  & \text{for } p\leq-2,\\
p^2/4 & \text{for } -2\leq p\leq2,\\
p-1 & \text{for } p\geq2.\\
\end{cases}
\]
The claim that $B(p)=|p|-1$ is the famous Brennan conjecture. Carleson and Jones \cite{MR1162188} proved that $B(p)=p-1+O((p-2)^2)$ as $p\to2$, $p>2$. Based on extensive computer experiments Kraetzer conjectured that $B(p)=p^2/4$ for $|p|\leq2$. See Garnett and Marshall \cite[p. 305]{zbMATH05296346}.
\textit{(Ch. Pommerenke)}\\

\noindent\textbf{Update 6.96} No progress on this problem has been reported to us.\\

\noindent\textbf{Problem 6.97} Goodman \cite{MR2937804} conjectured that if $f(z) = \sum^\infty_{n=1}a_nz^n$ is $p$-valent in
$\mathbb{D}$, then for each $n> p$, we have
\[|a_n|\leq\sum^p_{k=1}\frac{2k(n+p)!}{(p+k)!(p-k)!(n-p-1)!(n^2-k^2)}|a_k|.\]
If the bound is true, then it is sharp in all of the variables $p$, $n$, $a_1,\ldots, a_p$.\\
The conjecture has been proved for large subclasses (\cite{MR0040430}, \cite{MR0243054}, \cite{MR0051308}), but is still open
in general. The simplest case, namely
for all $2$-valent functions, is still open; but Watson \cite{MR0215980} has made important
contributions towards proving this inequality.\\ 
\textit{(A. W. Goodman)}\\

\noindent\textbf{Update 6.97} No progress on this problem has been reported to us.\\

\noindent\textbf{Problem 6.98} The coefficients of a $p$-valent function are bounded by some function of its
zeros. In particular, let the function
\[f(z)=z^q+\sum^\infty_{n=q+1}a_nz^n\]
be $p$-valent in $\mathbb{D}$ and have $s$ zeros $\beta_k$, $k= 1,2, \ldots, s$, where $0 < |\beta_k|<1$. Goodman \cite{MR0041920} conjectured that, under these hypotheses, $|a_n|\leq |A_n|$ where $A_n$ is defined
by the identity
\[F(z)=\frac{z^q}{(1-z)^{2q+2s}}\Big(\frac{1+z}{1-z}\Big)^{2t}\prod^s_{k=1}\Big(1+\frac{z}{|\beta_k|}\Big)\big(1+|\beta_k|z\big)=z+\sum^\infty_{n=2}A_nz^n,\]
where $t = p-q-s\geq0$.\\
The conjecture has been proved if $t = 0$ and $f(z)$ is $p$-valent star-like with respect
to the origin \cite{MR0033886}. However, it is still open in general.\\ 
\textit{(A. W. Goodman)}\\

\noindent\textbf{Update 6.98} No progress on this problem has been reported to us.\\

\noindent\textbf{Problem 6.99} A function $f(z) = z + a_2 z^2 +\ldots$ is said to belong to the \textit{class $CV(R_1,R_2)$} if
it is univalent and convex in $\mathbb{D}$, and if on $f(\{|z|=1\})$ the curvature $\rho$ satisfies the
inequalities $R_1\leq\rho\leq R_2$. So far, little progress has been made on the study of this
class of functions (see Goodman \cite{MR760942} and \cite{MR835886}); we do not even know the sharp bound for $|a_2|$
in $CV(R_1,R_2)$.\\ 
\textit{(A. W. Goodman)}\\

\noindent\textbf{Update 6.99} No progress on this problem has been reported to us.\\

\noindent\textbf{Problem 6.100} Given two functions $f$, $g$ in the (usual) class $S$, we can form the new
functions (arithmetic and geometric mean functions)
\[F(z)=\alpha f(z)+\beta g(z)\hspace{1cm}\text{ and }\hspace{1cm} G(z)=z\Big(\frac{f(z)}{z}\Big)^\alpha\Big(\frac{g(z)}{z}\Big)^\beta,\]
where $\alpha$, $\beta\in(0,1)$ and $\alpha+\beta=1$. It is known \cite{MR0247057} that, if
\[0.042\simeq\frac{1}{1+e^\pi}<\alpha,\hspace{1cm}\beta<\frac{e^\pi}{1+e^\pi}\simeq0.988,\]
then there are functions $f$, $g$ in $S$ such that $F$ and $G$ have valence infinity in $\mathbb{D}$.
What can be said about the `fringes' of the interval $(0,1)$? Is there some bound on the
valence of $F$ and $G$ that is a function of $\alpha$ for $0 < \alpha\leq 1/(1+e^\pi)$?\\
\textit{(A. W. Goodman)}\\

\noindent\textbf{Update 6.100} No progress on this problem has been reported to us.\\

\noindent\textbf{Problem 6.101} Let $K$ be a closed set of points in $\mathbb{C}$, and let $F(K)$ denote the family of
functions $f$ of the form
\[f(z)=\sum^n_{k=1}\frac{A_k}{z-a_k},\]
where $A_k > 0$ and $a_k\in K$, $k=1, 2,\ldots, n$. Find a maximal domain of $p$-valence for the
class $F(K)$.\\
For $p=1$, this problem was completely solved by Distler \cite{MR0167610}; however, for $p> 1$,
we do not even have a good conjecture. (See also Goodman \cite{MR1005663}.)\\ 
\textit{(A. W. Goodman)}\\

\noindent\textbf{Update 6.101} No progress on this problem has been reported to us.\\

\noindent\textbf{Problem 6.102} Let $\{v_n\}^\infty_1$ be a sequence of positive integers (which may include $\infty$);
the sequence is called a \textit{valence sequence} if there is a function $f(z)$, analytic
in $\mathbb{D}$, such that $f^{(n)}(z)$ has valence $v_n$ in $\mathbb{D}$, for $n = 1,2,\ldots$.\\
Find interesting necessary conditions for $\{v_n\}$ to be a valence sequence. Also, find
sufficient conditions for $\{v_n\}$ to be a valence sequence. (For some results of this type,
see Goodman \cite{MR567984}.)\\ 
\textit{(A. W. Goodman)}\\

\noindent\textbf{Update 6.102} No progress on this problem has been reported to us.\\

\noindent The next five problems concern harmonic univalent functions, that is, functions that can be written in the form
\[f(z)=\sum^\infty_{n=0}a_nz^n+\sum^\infty_{n=1}a_{-n}\overline{z}^n\]
and are univalent. For a good
introduction to these, see Clunie and Sheil-Small\cite{MR752388}. We will denote by $S$ the (usual) class of normalised
analytic univalent functions in $\mathbb{D}$, by $S_H$ the class of normalised harmonic
univalent functions in $\mathbb{D}$ (with $a_0 = 0$ and $a_1 = 1$), and by $S^0_H$ the subclass of $S_H$ for
which $a_{-1} = 0$.\\

\noindent\textbf{Problem 6.103} The function
\[k(z)=2\text{Re }\Big(\frac{z+\frac{1}{3}z^3}{(1-z)^3}\Big)=\sum^\infty_{n=1}\frac{1}{3}(2n^2+1)r^n(e^{in\theta}+e^{-in\theta}),\]
where $z = re^{i\theta}$, lies in the closure of $S_H$. Prove that $k$ is extremal for the coefficient
bounds in $S_H$.\\
It is known by Clunie and Sheil-Small \cite{MR752388} that $|a_n| <\frac{1}{3}(2n^2+1)$ for functions in $S_H$ with real coefficients, and
that $|a_n| \leq\frac{1}{3}(2n^2+1)$ for functions in $S_H$ for which the domain $f(\mathbb{D})$ is close-to-convex.\\
\textit{(T. Sheil-Small)}\\

\noindent\textbf{Update 6.103} No progress on this problem has been reported to us.\\

\noindent\textbf{Problem 6.104} It is known \cite{MR752388} that, for functions $f$ in $S^0_H$, $\{|w| < \frac{1}{16}\}\subset f(\mathbb{D})$. Prove that the correct value $d$, of the Koebe constant for the class $S^0_H$ is $\frac{1}{6}$.\\
Note that the function
\[k_0(z)=\text{Re }\Big(\frac{z+\frac{1}{3}z^3}{(1-z)^3}\Big)+i\text{Im }\Big(\frac{z}{(1-z)^2}\Big)\]
belongs to $S^0_H$ and maps $\mathbb{D}$ onto the plane cut along the real axis from $-\frac{1}{6}$ to $-\infty$,
and that $\frac{1}{6}$ is the correct constant for those functions $f$ in $S^0_H$ for which $f(\mathbb{D})$ is close-to-convex.\\
Also, determine the number
\[\alpha=\sup\{|a_2|:f\in S_H\}.\]
The best known estimate for $\alpha$ is $\alpha < 57.05$, see Sheil-Small \cite{MR1083443}. It is also known that $d\geq 1/(2\alpha)$.\\
\textit{(T. Sheil-Small)}\\

\noindent\textbf{Update 6.104} No progress on this problem has been reported to us.\\

\noindent\textbf{Problem 6.105} What are the convolution multipliers $\phi^*:K_H\to K_H$, where $K_H$ is the
subclass of functions $f$ in $S_H$ with convex images $f(\mathbb{D})$?\\
A particularly interesting case is the radius of convexity problem: for which
values of $r$ in $(0,1)$ is the function $z\mapsto f(rz)$ convex in $\mathbb{D}$, when $f$ is convex in $\mathbb{D}$? (It is
known that $r\leq\sqrt(2)-1$, see Clunie and Sheil-Small \cite{MR752388}.)\\ 
\textit{(T. Sheil-Small)}\\

\noindent\textbf{Update 6.105} No progress on this problem has been reported to us.\\

\noindent\textbf{Problem 6.106} Let $J$ be a Jordan curve in $\mathbb{C}$ bounding a domain $D$. Suppose that
$f:e^{it}\mapsto f(e^{it})$ is a sense-preserving homeomorphism of the unit cirlce $\mathbb{T}$ onto $J$, and that the
harmonic extension of $f$ to $\mathbb{D}$ satisfies the relation $f(\mathbb{D})\subset D$. Prove that $f$ is a
homeomorphism of $\mathbb{D}$ onto $D$.\\
This is known to be true if $J$ is convex, when the hypothesis $f(\mathbb{D})\subset D$ is
automatically satisfied because of the positivity of the Poisson kernel by the
Kneser-Rado-Choquet theorem, see Duren \cite{zbMATH02117216}. The result is also known to be true when $\partial f/\partial t$ is
continuous and non-zero (on $\mathbb{T}$), and when the co-analytic and analytic parts of
$f$ have continuous derivatives on $\overline{\mathbb{D}}$.\\ 
\textit{(T. Sheil-Small)}\\

\noindent\textbf{Update 6.106} No progress on this problem has been reported to us.\\

\noindent\textbf{Problem 6.107} Prove that, for $f\in S^0_H$,
\begin{equation}
\label{J6.107}
\big||a_n|-|a_{-n}|\big|\leq n,\hspace{1cm} n=2,3,4,\ldots.
\end{equation}
(This is a generalisation of the Bieberbach conjecture for $S$.)\\
It is known that (\ref{J6.107}) holds in the following cases:
\begin{enumerate}[(a)]
\item when $f$ has real coefficients, see Clunie and Sheil-Small \cite{MR752388};
\item when $f(\mathbb{D})$ is star-like with respect to the origin, see Sheil-Small \cite{MR1083443};
\item when $f(\mathbb{D})$ is convex in one direction, see Sheil-Small \cite{MR1083443}.
\end{enumerate}
\textit{(T. Sheil-Small)}\\

\noindent\textbf{Update 6.107} No progress on this problem has been reported to us.\\

\noindent\textbf{Problem 6.108} Let $f$ be analytic univalent in $\mathbb{D}$, and consider
\[I_\lambda(r,f')=\Big(\frac{1}{2\pi}\int^{2\pi}_0\big|f'(re^{i\theta})\big|^\lambda\,d\theta\Big)^{1/\lambda}\]
where $\lambda > 0$.
\begin{enumerate}
\item What is the maximal order of magnitude of $I_\lambda(r,f')$ as $r\to1$, where
$0\leq \lambda\leq\frac{2}{5}$?\\
If $\lambda>\frac{2}{5}$, it is known that $I_\lambda(r,f') = O\big(1/(1-r)^{3-1/\lambda}\big)$ with equality when $f$ is the
Koebe function. The case $\lambda >\frac{1}{2}$ follows easily from classical facts, while the case
$\frac{2}{5}\leq \lambda \leq\frac{1}{2}$ is due to Feng and MacGregor \cite{Feng1976}. (It was formerly conjectured that the
Koebe function would still be extremal for $\lambda>\frac{1}{3}$; but an example of Makarov \cite{MR818450} has
shown that this is not the case for $\lambda\leq\frac{1}{3}+\varepsilon$, for some positive $\varepsilon$.)
\item Now normalise the function $f$ to belong to the usual class $S$. For $\lambda>\frac{2}{5}$, we
know that
\[I_\lambda(r,f')\leq C_\lambda I_\lambda(r,k'),\]
where $k$ is the Koebe function. For which $\lambda$ is the best constant $C_\lambda$ equal to $1$? (It
follows from de Branges' theorem \cite{MR772434} that $C_\lambda = 1$ for \mbox{$\lambda = 2,4,6,\ldots$}. Presumably
$C_\lambda=1$ for $\lambda\geq2$, but the proposer knows of no proof of this.)
\end{enumerate} 
\textit{(A. Baernstein)}\\

\noindent\textbf{Update 6.108} Related questions have become of great interest, see Problem 6.96. No progress on this problem has been reported to us.\\

\noindent\textbf{Problem 6.109} Let $f$ be analytic univalent in $\mathbb{D}$, and consider
\[I_{-\lambda}(r,f')=\Big(\frac{1}{2\pi}\int^{2\pi}_0\big|f'(re^{i\theta})\big|^{-\lambda}\,d\theta\Big)^{1/\lambda},\]
where $\lambda>0$. Except for the elementary case $\lambda=\infty$, the maximal order of magnitude
as $r\to1$ is not known for any positive $\lambda$. A particularly interesting case is when $\lambda=2$. It
seems possible that
\begin{equation}
\label{J6.109}
I_{-2}(r,f')=O\big(I_{-2}(r,k')\big)=O\big((1-r)^{-1/2}\big)
\end{equation}
might be true, where $k$ is the Koebe function. This equation (\ref{J6.109}) is slightly stronger than an earlier
conjecture of Brennan (see  Problem 6.94). The best known bounds are due to
Pommerenke \cite{MR875238}.\\
Related information appears in \cite{MR951982}.\\ 
\textit{(A. Baernstein)}\\

\noindent\textbf{Update 6.109} No progress on this problem has been reported to us.\\

\noindent\textbf{Problem 6.110} Let $\Omega$ be a simply-connected domain in the finite plane whose complement 
contains $n$ disjoint closed balls with centres on the interval $[0,1]$ and common radius $\varepsilon$, 
$\varepsilon<1/n$, (such domains are the so-called `ball and chain domains'). Let $z_0$ in $\Omega$ be a
point at a distance at least $1$ from each ball, and let $w(z_0,\Omega)$ denote the harmonic
measure at $z_0$ of the union of the balls, relative to $\Omega$. Is it true that, for every positive $\delta$,
there is an estimate
\begin{equation}
\label{J1.110}
w(z_0,\Omega)\leq C_\delta(n\varepsilon)^{\frac{1}{2}-\delta},
\end{equation}
where $C_\delta$ depends only on $\delta$.\\
An affirmative answer would imply the $L^p$ extension, for $1 < p < 2$, of the
Hayman-Wu theorem \cite{MR639358} mentioned in Problem 6.93. In the case $n = 1$, the
Beurling-Nevanlinna projection theorem shows that (\ref{J1.110}) is true with exponent $\frac{1}{2}$;
however, for large $n$ there are examples that show that this is then false.\\
\textit{(A. Baernstein)}\\

\noindent\textbf{Update 6.110} No progress on this problem has been reported to us.\\

\noindent\textbf{Problem 6.111} Let $A$ denote the class of functions $f(z) = z + a_2z^2 +\ldots$ analytic in
$\mathbb{D}$. For $\delta\geq0$ and $T= \{T_k\}^\infty_2$ a sequence of non-negative real numbers,
define a \textit{$T$-$\delta$-neighbourhood} of $f\in A$ by
\[TN_\delta(f)=\Big\{g:g(z)=z+b_2z^2+\ldots\in A,\sum^\infty_{k=2}T_k|a_k-b_k|\leq\delta\Big\}.\]
When $T = \{k\}^\infty_2$, we call $TN_\delta(f) = N_\delta(f)$ a \textit{$\delta$-neighbourhood} of $f$ ($\delta$-neighbourhoods were introduced by Ruscheweyh \cite{MR601721}, who used them to generalise the result \cite{MR1503516} that
$N_1(z)\subset St$, the class of normalised star-like functions in $\mathbb{D}$.)\\
Now let $K[A,B]$ denote the class of univalent functions
\[\{f:f\in S,1+zf''(z)/f'(z)\prec(1+Az)/(1+Bz),z\in \mathbb{D}\},\]
where $-1\leq B<A\leq1$, introduced by Janowski \cite{MR0315110}. See also Pommerenke \cite[Chp. 2]{MR0507768} for discussion on $\prec$ and \textit{subordination}. The proposers have shown \cite{MR981503}
that if $f\in K[A,B]$, and either
\[(a)\hspace{1cm}-\frac{1}{4}(\sqrt{3}+2)\leq B<A\leq1\hspace{1cm}\text{ or }\hspace{1cm} (b)\hspace{1cm} -1\leq-A\leq B<A\leq1,\]
then $N_\delta(f)\subset St$, where
\[\delta=
\begin{cases}
(1-B)^{(A-B)/B},&\hspace{1cm} B\neq0,\\
e^{-A},&\hspace{1cm} B=0,
\end{cases}\]
and that this value of $\delta$ is best possible. Is this conclusion still valid without the
hypotheses $(a)$ and $(b)$?\\ 
\textit{(T. Sheil-Small and E. M. Silvia)}\\

\noindent\textbf{Update 6.111} No progress on this problem has been reported to us.\\

\noindent\textbf{Problem 6.112} If  $f$ in $A$ and $\delta > 0$, define a \textit{$\Sigma_\delta(f)$ neighbourhood} of $f$ to be
\[\Big\{g:g\in A,\big|(g'(z)-f'(z))-\frac{1}{z}(g(z)-f(z))\big|+\big|(g'(z)-f'(z))+\frac{1}{z}|g(z)-f(z))\big|<2\delta\Big\}.\]
Clearly, $N_\delta(f)\subset\sum_\delta(f)$. Then it is known by Sheil-Small and Silvia \cite{MR981503} that $\sum_1(z)\subset St$, and \mbox{$\sum_\frac{1}{4}(f)\subset St$} for convex functions $f$, with the notation in Problem 6.111.\\
Given a normal family $\mathcal{F}$ in $A$, the \textit{dual}, $\mathcal{F}^*$, of $\mathcal{F}$ is the set
\[\big\{f:f\in A,f\ast g\neq0\text{ for all }g\in\mathcal{F},0<|z|<1\big\},\]
where $\ast$ denotes the Hadamard product, see Ruscheweyh \cite{MR0382626}. The Bieberbach conjecture (see de Branges \cite{MR772434}) is
equivalent to the statement that $N_1(z)\subset S^*$ (where $S^*$ is the dual of $S$). It seems
reasonable therefore to ask whether $\sum_1(z)\subset S^*$. (Recall that $S^*\subset St$, see Ruscheweyh \cite{MR0318470}.)\\
\textit{(T. Sheil-Small and E. M. Silvia)}\\

\noindent\textbf{Update 6.112} No progress on this problem has been reported to us.\\

\noindent\textbf{Problem 6.113} Following the notation of Problem 6.111 and 6.112, it is known \cite{MR981503} that,
if $|x|\leq\rho\leq1$ and $\gamma=1/(1+\rho)^2$, then
\[N_\gamma\Big(\frac{z}{1-xz}\Big)\subset S^*.\]
Is it true that
\[\Sigma_\gamma\Big(\frac{z}{1-xz}\Big)\subset S^*?\]
Notice that, since not all convex functions belong to $S^*$ (see Sheil-Small and Silvia \cite{MR981503}), it follows that we
cannot replace in the above, $N_\gamma(z/(1-xz))$ by $N_\gamma(g)$, where $g$ is an arbitrary convex
function.\\ 
\textit{(T. Sheil-Small and E. M. Silvia)}\\

\noindent\textbf{Update 6.113} No progress on this problem has been reported to us.\\

\noindent\textbf{Problem 6.114} Let $\Gamma$ be a regular curve and $f$ an analytic and conformal function in the
open unit disc. Does $f^{-1}(\Gamma)$ necessarily have finite length? (A curve $\Gamma$ is said to be
\textit{regular} if the intersection of $\Gamma$ with a disc of radius $r$ has one-dimensional measure
at most $Cr$, where $C$ is a constant independent of $r$.)\\
\textit{(J. L. Fernandez and D. Hamilton)}\\

\noindent\textbf{Update 6.114} No progress on this problem has been reported to us.\\

\noindent\textbf{Problem 6.115} Let $\Gamma$ be a rectifiable curve, and let $E$ be a subset of $\Gamma$  having zero length. If $\Omega$ is any simply-connected domain and $z\in\Omega$, is it true that the harmonic
measure satisfies the equation $\omega(z,E,\Omega)=0$?\\
\textit{(B. \O{}ksendal; R. Kaufman and J.-M. G. Wu; communicated by D. Hamilton)}\\

\noindent\textbf{Update 6.115} No progress on this problem has been reported to us.\\

\noindent\textbf{Problem 6.116} Let $D$ be a domain in $\mathbb{C}$ containing the origin $0$; for $t>0$, let $\Omega_t$ be the
component of $D \cap\{|z|\leq t\}$ containing $0$. In the usual notation for harmonic measure,
define the function $\omega = \omega_D:[0, \infty)\to[0,1]$ by the formula
\[\omega(t)=\omega\big(0,\partial\Omega_t\cap\{|z|=t\},\Omega_t\big).\]
What can be said about $\omega$ and about its relations with $D$? For instance:
\begin{enumerate}[(a)]
\item What are necessary and sufficient conditions on a function to be $\omega_D$ for
some $D$?
\item If $\omega_{D_1 }\equiv \omega_{D_2}$, are $D_1$ and $D_2$ essentially the same? (That is, will $D_1$ be a
rotation or reflection of $D_2$? Or will $D_1$ and $D_2$ differ only on sets of capacity zero?)
What happens if  $\omega_{D_1 }\equiv \omega_{D_2}$ only on some subinterval of $[0,1]$?
\item Given a function $\omega_D$, can one `reconstruct' $D$?
\item Can one infer such properties as connectivity of $D$ fromthe  behaviour of $\omega_D$? One
might start by dealing with domains $D$ that are circularly symmetric.
\end{enumerate}
\textit{(K. Stephenson)}\\

\noindent\textbf{Update 6.116} Snipes and Ward \cite{MR2173369} and Walden and Ward \cite{MR1867264} provide results  which answer the analogues of Stephenson's questions, for their function $h(r)$. Barton and Ward \cite{MR3261731} establish several sets of sufficient conditions on a function $f$, for $f$ to arise as a harmonic measure distribution function.\\

\noindent\textbf{Problem 6.117}
Let $G$ be a domain in $\mathbb{C}$ that contains the origin $0$ and is \textit{axially-symmetric
with respect to the real axis}, that is, if a point $z\in G$ then the line segment with endpoints $z$ and $\bar{z}$ also lies in $G$. (In particular it follows that $G$ is simply-connected.) For
$t>0$, let $G_t$ and $G_t'$ be the components of the domains
\[\{z:z\in G,\text{Im }z< t\}\hspace{1cm}\text{ and }\hspace{1cm}\{z:z\in G,\text{Re }z< t\}\]
containing $0$. Let $\omega_t$ and $\omega_t'$ denote the harmonic measure at $0$ of the sets
\[\partial G_t\cap\{\text{Im }z=t\}\hspace{1cm}\text{ and }\hspace{1cm}\partial G_t'\cap\{\text{Re }z=t\}\]
with respect to the domains $G_t$ and $G_t'$ respectively.\\
Now let $D_1$ and $D_2$ be two domains in $\mathbb{C}$ of this type, and use for each the notation just described. Then:
\begin{enumerate}[(a)]
\item If $\omega_t(D_1)\equiv \omega_t(D_2)$ for each $t>0$, is it true that $D_1 = D_2$?
\item If $\omega_t'(D_1)\equiv \omega_t'(D_2)$ for each $t>0$, is it true that $D_1 = D_2$?
\item If both $(a)$ and $(b)$ fail, is it true that $D_1 = D_2$ if $\omega_t(D_1)\equiv \omega_t(D_2)$ for each $t>0$ and $\omega_t'(D_1)\equiv \omega_t'(D_2)$ for each $t>0$?
\end{enumerate}
If the answer to $(a)$ or $(b)$ is `yes', in that case, how few $t$'s does one need for the
conclusion to hold? (For example, infinitely many $t$'s such that $\{\text{Im }z = t\}$ or $\{\text{Re }z = t\}$
meet the domains? Or do the $t$'s need to be dense?) If the answer to $(a)$, $(b)$ or $(c)$
is `yes', in that case can one replace `axially-symmetric' by `simply-connected', or
perhaps drop this requirement completely?\\ 
\textit{(D. A. Brannan)}\\

\noindent\textbf{Update 6.117} The original statement of this problem contained errors which have now been fixed. No progress on this problem has been reported to us.\\

\section*{New Problems}

\noindent To appear.

%======================= CHAPTER 7 ====================================

\chapter{Miscellaneous}

\underline{Preface by D. Sixsmith (to appear)}\\

\noindent\textbf{Problem 7.1} Let $E$ be the compact plane set of transfinite diameter ($=$capacity) $d(E)=1$ and let
\[d_n(E)^{n(n-1)/2}=\max_{w_\nu\in E}\prod_{1\leq\mu<\nu\leq n}|w_\mu-w_\nu|.\]
It is known that $d_n(E)$ decreases with $n$ and $d_n(E)\to d(E)$ as $n\to\infty$, so that $d_n(E)\geq1$ is trivial. It is known by Pommerenke \cite{MR0168772}, that $d_n(E)\geq n^{2/(n-1)}$ if $E$ is connected. Is this inequality true in general?\\
Further, is it true that 
\[d_n(E)^{(n-1)/2}\leq Kn\]
if $E$ is connected, where $K$ is some absolute constant? It is known that 
\[d_n(E)^{(n-1)/2}\leq\Big(\frac{4}{e}\log n +4\Big)n,\]
in this case.\\
\textit{(Ch. Pommerenke)}\\

\noindent\textbf{Update 7.1}  No progress on this problem has been reported to us.\\

\noindent\textbf{Problem 7.2} Let $f(z)$ be analytic in a simply-connected domain $D$. It is known that $f(z)$ can be expanded in a series of Faber polynomials 
\[f(z)=\sum^\infty_{n=0}a_nP_n(z).\]
Find the domain of variability $V$ of $a_n$ as $f(z)$ runs through all functions analytic in $D$, and having positive real part there. It is known that if $D$ is a circle, $V$ is a circle; and if $D$ is an ellipse, then $V$ is an ellipse, see Royster \cite{MR0203036}. \\
\textit{(W. C. Royster)}\\ 

\noindent\textbf{Update 7.2} No progress on this problem has been reported to us.\\

\noindent\textbf{Problem 7.3} Let $z_i$, $1\leq i\leq n$ be a finite sequence of complex numbers such that $|z_i|\leq1$. Set
\[S_k=\sum^n_{i=1}z^k_i.\]
Can we have
\begin{equation}
\label{7.1}
\max_{2\leq k\leq n+1}|S_k|<A^{-n},
\end{equation}
where $A$ is an absolute constant greater than one? If we assume $z_1=1$, $|z_i|\leq1$, $2\leq i\leq n$, then (\ref{7.1}) can be satisfied. See Tur\'an \cite{MR0060548}.\\
\textit{(P. Erd\"os)}\\

\noindent\textbf{Update 7.3} No progress on this problem has been reported to us.\\

\noindent\textbf{Problem 7.4} If $z_1=1$, and the $z_i$ are arbitrary complex numbers for $2\leq i\leq n$, then Atkinson \cite{MR0126418} proved that
\[\max_{1\leq k\leq n}|S_k|>c\]
with $c=\frac{1}{3}$. What is the best value for the constant $c$?\\

\noindent\textbf{Update 7.4} No progress on this problem has been reported to us.\\
  
\noindent\textbf{Problem 7.5} Let $z_i$, $1\leq i\leq n$ be $n$ complex numbers such that $|z_i|\geq1$. Then there exists an absolute constant $c$ so that the number of sums
\begin{equation}
\label{7.2}
\sum^n_{i=1}\varepsilon_i z_i,\hspace{1cm}\varepsilon_i=\mp1,
\end{equation} 
which fall into the interior of an arbitrary circle of radius $1$, is less than $(c2^n\log n)/n^{\frac{1}{2}}$, see Littlewood and Offord \cite{MR0009656}. If the $z_n$ are real, Erd\"os \cite{MR0014608} proves that the number of the sums (\ref{7.2}) which fall into the interior of any interval of length $2$ is at most $\begin{pmatrix} n\\ [\frac{n}{2}]\end{pmatrix}$ and this is sharp. He asks whether this estimate remains true for complex $z_i$, where his proof gives only $c2^n/\sqrt{n}$, or more generally, vectors of Hilbert space of norm at least one, when he can prove that the number of sums (\ref{7.2}) falling into an arbitrary unit sphere is $o(2^n)$.\\

\noindent\textbf{Update 7.5} The plane case has been settled independently by Katona \cite{MR0205864} and Kleitman \cite{MR0184865}. The general case was also settled by Kleitman \cite{MR0184865}.\\ 

\noindent\textbf{Problem 7.6} We consider the range of the random function
\[F(z)=\sum^\infty_{n=0}\pm a_nz^n\]
($F$ chosen at random in the natural way) defined in $\mathbb{D}$, where $\sum|a_n|^2=\infty$. Is the image of $w=F(z)$ with probability one
\begin{enumerate}[(a)]
\item everywhere dense in the plane?
\item the whole plane?
\item does it contain any given point with probability one?
\end{enumerate}
If $a_n=n^\lambda$, $(b)$ holds if $\lambda>\frac{1}{2}$, and $(a)$ holds if $-\frac{1}{2}<\lambda<+\frac{1}{2}$.\\
\textit{(J. P. Kahane)}\\ 

\noindent\textbf{Update 7.6} Partial results were obtained by Offord \cite{MR0307393}. This problem has been solved by Nazarov, Sodin and Nishry \cite{zbMATH06373401} who proved that (b) always holds.\\

\noindent\textbf{Problem 7.7} Let $f(z)=\sum^\infty_{n=0}a_nz^n$ be analytic in $\mathbb{D}$,
\[f_n(z)=\sum^n_{k=0}a_kz^k,\]
and let $\Gamma_n$ be the modulus of the largest zero of $f_n(z)$. Then Kakeya (no citation) showed that $\liminf_{n\to\infty}\Gamma_n\leq A$, with $A=2$. What is the best possible value of $A$? Clunie and Erd\"os \cite{MR0222262} have shown that $\sqrt{2}<A<2$.\\

\noindent\textbf{Update 7.7} This has been settled by Buckholtz \cite{MR0259076}. He has found a method by which the constant in question can be calculated with any desired accuracy.\\

\noindent\textbf{Problem 7.8} Is it possible to express each $K$-quasiconformal map in $3$-space as the composition of two quasiconformal maps with maximal dilatation less than $K$? The corresponding plane result is true.\\
\textit{(F. W. Gehring)}\\

\noindent\textbf{Update 7.8} No progress on this problem has been reported to us.\\

%\section* *{Quasiconformal mapping}

\noindent\textbf{Problem 7.9} Suppose that $f$ is a plane \textit{$K$-quasiconformal} mapping
of the unit disc $\mathbb{D}$ onto itself. Show that there exists a finite constant
$b = b(K)$ such that
\[m(f(E)) \leq b \{ m(E)\}^{1/K}\]
for each measurable set $E \subset D$. Here $m$ denotes plane Lebesgue
measure. Such an inequality is known (see \cite{MR0201635}) with the
exponent $\frac{1}{K}$ replaced by a constant $a=a(K)$.\\
\textit{(F. W. Gehring)}\\

\noindent\textbf{Update 7.9} This has been proved by Astala \cite{zbMATH00714744}. His proof has been substantially simplified by Eremenko and Hamilton \cite{zbMATH00868218}.\\

\noindent\textbf{Problem 7.10} It was proved by Boyarski\u\i\, \cite{MR0071620} that the partial derivatives of a plane
$K$-quasiconformal mapping are locally $L$-integrable for $2\leq p < 2+c$, where
$c = c(K) > 0$. Show that this is true with $c =\frac{2}{K-1}$.\\
The example $f(z) = |z|^{\frac{1}{K}-1}.z$ shows that such a result would be
sharp.\\
\textit{(F. W. Gehring) }\\

\noindent\textbf{Update 7.10} Gehring points out that the proof published by Okabe \cite{MR501158}, who claims to solve this problem, contains an error. Details of this are found in Gehring's review of Okabe's article for \textit{Mathematical Reviews}. See follow-up papers by Okabe (\cite{MR561991}, \cite{MR581203}) also.\\

\noindent\textbf{Problem 7.11} Show that each quasiconformal mapping of $\mathbb{R}^n$ onto $\mathbb{R}^n$ has a quasiconformal extension
to $\mathbb{R}^{n+1}$ . This has been established by Ahlfors when $n = 2$ and by
Carleson when $n = 3$.
\textit{(F. W. Gehring)}\\

\noindent\textbf{Update 7.11} This has been established in all cases by Tukia and V\"ais\"al\"a \cite{MR647809}.\\

\noindent\textbf{Problem 7.12} Suppose that $f$ is an $n$-dimensional $K$-quasi-analytic function.
Show that the partial derivatives of $f$ are locally $L$-integrable for
$n\leq p\leq n + c$, where $c = c(K, n) > 0$. This was shown by Gehring \cite{MR0402038}
to be true if $f$ is injective.\\
\textit{(F. W. Gehring)}\\

\noindent\textbf{Update 7.12} This has been solved by Elcrat and  Meyers \cite{MR0417568} modulo a slight problem observed by Gehring, which was straightened out by Stredulinsky \cite{MR570689}.\\

%\section* *{Riemann surfaces and conformal mapping}

\noindent\textbf{Problem 7.13} One part of Nevanlinna theory is devoted to the following problem.
How does the geometric structure of a simply connected Riemann covering
surface of the sphere influence the value distribution of the meromorphic
function generating the surface? It is suggested that one also consider halfsheets among
the constituent pieces of such surfaces, which have an
infinite number of branch points on their boundary.\\
A typical example is the covering surface generated by $e^z - z$.
It contains a right half-plane with a second order branch point at $2\pi in$ 
for each integer $n$.\\
\textit{(F. Huckeman)}\\

\noindent\textbf{Update 7.13} No progress on this problem has been reported to us.\\

\noindent\textbf{Problem 7.14} \textbf{(Boundary values of Cauchy integrals).} Let $\gamma$ be a $C^1$ curve in
the plane and let $f$ be a continuous function on $\gamma$. Put
\[F(z)=\int_\gamma\frac{f(t)}{t-z}\,dt\hspace{1cm}(z\notin\gamma).\]
Does $F$ have non-tangential boundary values almost everywhere on $\gamma$? This is a very
old question which has been studied quite extensively, especially by
Russian mathematicians. It is known that the answer is `yes' if slightly
greater smoothness is assumed for $\gamma$ or $f$ (see e.g. Havin \cite{MR0192059} for an
outstanding contribution to this subject).\\
Suppose now that $\gamma$ is a Jordan curve and let $\phi$ be a conformal
map from the unit disc onto the inside of $\gamma$. Does $F\circ\phi$  belong to $H^p$
for some $p < 1$, or perhaps to the class $N$ of functions with bounded
Nevanlinna characteristic?\\
\textit{(A. Baernstein)}\\

\noindent\textbf{Update 7.14} Calder\'on \cite{MR0466568} has proved that Cauchy integrals of measures on $C^1$ curves have boundary values almost everywhere. The second part of the problem, about possibly pulling back to an $H^p$ function, is still open. For non-rectifiable $\gamma$ one needs to consider $F(z)=\int_\gamma (t-z)^{-1}\,d\mu(t)$ for measures $\mu$ since $F(z)=\int_\gamma (t-z)^{-1}f(t)\,d\mu(t)$ does not make sense.\\

\noindent\textbf{Problem 7.15} Let $D$ be a domain in the extended complex plane. A finite
point $z$ on the boundary $\partial D$ of $D$ is called \textit{angular} (relative to $D$)
if there exists $\varepsilon > 0$ such that every component domain of
$D\cap\{|z - z_0| < \varepsilon\}$ which has $z_0$ as a boundary point, is contained in
an angle less than $\pi$ with vertex at $z_0$. Angularity at $\infty$ is similarly
defined.\\
Let $A = A(D)$ be the set of angular points of $\partial D$ relative to $D$.
Obviously $A$ not empty implies that $\partial D$ has positive capacity. The set
$A(D)$ can have positive linear measure. E.g. let $C$ be a Cantor set on
$|z|=1$ and let $D$ consist of the open unit disc from which have been
deleted all points $rz$ with $z \in C$ and $\frac{1}{2}\leq r < 1$. Then $A(D) = C$
which can, of course, have positive linear measure.\\
Yet the following holds for arbitrary domains: $A(D)$ is either
empty or its harmonic measure, relative to any point of $D$, is zero.
This result follows easily from an unpublished theorem on Brownian
paths $\omega(t)$ in the complex plane. This says that almost all such paths
have the property that for every real $t_0$ and every $\varepsilon > 0$ the set of
numbers $\Big|\frac{\omega(t)-\omega(t_0)}{\omega(t)-\omega(t_0)}\Big|$ with
$t_0<t<t_0+\varepsilon$ fills at least an open arc
of length $\pi$ on the unit circle. This theorem is not easy and it would
be desirable to give a direct proof of the above result on $A(D)$.\\
Moreover, the Brownian paths approach will certainly not yield
a similar result for the set $B_\alpha(D)$ of $\partial D$ whose points are defined by
replacing the angles less than $\pi$ with translates of $\{x+iy : 0<x< |y |^\alpha\}$
for a given $\alpha$ with $\frac{1}{2}< \alpha < 1$. (For $\alpha=\frac{1}{2}$ the result is false as can
be seen by taking $D$ to be a disc.)\\
Problem: For which $\alpha\in(\frac{1}{2},1)$ is the harmonic measure of the set $B_\alpha(D)$
always zero? (Of course we may assume that the capacity of $\partial D$ is
positive.)\\
A much more difficult problem would be to characterise the monotone
functions $f(y)$ having the property that the set obtained on replacing
the angles by translates of $\{x+iy : 0 < x < f(|y|)\}$ has necessarily
harmonic measure zero.\\
Similar questions can be asked for Riemann surfaces and for
$n$-dimensional space.\\
\textit{(A. Dvoretzky)}\\

\noindent\textbf{Update 7.15} No progress on this problem has been reported to us.\\

%\section* *{Approximation}

\noindent\textbf{Problem 7.16} Let $\gamma$ be a Jordan arc, $d\mu$ a measure on $\gamma$. Does the Laplace
transform
\begin{equation}
\label{I7.16}
f(z) = \int_\gamma e^{z\zeta}\,d\mu(\zeta)
\end{equation}
always have `asymptotically analytic' growth as $|z|\to\infty$. The answer
might be `no'. However it could be true anyway that the zeros of $f(z)$
have `measurable distribution' in the sense of A. Pfluger \cite{MR0015488}.\\
\textit{(J. Korevaar)}\\

\noindent\textbf{Update 7.16} On regularity of growth and zero distribution of Laplace transforms $f(z)=\int^\gamma e^{z\zeta}\,d\mu(\zeta)$ along arcs, Dixon and Korevaar \cite{MR507458} have obtained some results for arcs of limited slope. A negative answer has been given by Wiegerinck (personal communication) using Borel transforms, to the question whether the Laplace transform (\ref{I7.16}) always has asymptotically regular growth as $|z|\to\infty$.\\

\noindent\textbf{Problem 7.17} Let $f(z)$ be analytic and bounded for $\Re z> 0$. Suppose that
$|\alpha|<\frac{1}{2}\pi$ and that $(r_n)$ is a sequence of positive integers with $\sum\frac{1}{r_n}=\infty$.
Show that the exponential type of $f$ on the sequence $z_n = r_ne^{i\alpha}$ is equal
to the type of $f$ on the ray $z = re^{i\alpha}$. (Proofs by Boas \cite{MR0059365} and by Levinson \cite{MR1501941}
for $\alpha = 0$ do not seem to work for $|\alpha|>\frac{1}{4}\pi$.)\\
\textit{(J. Korevaar)}\\

\noindent\textbf{Update 7.17} This has been proved by Korevaar and Zeinstra \cite{zbMATH03989737}. See also Zeinstra \cite{zbMATH04212730}.\\

\noindent\textbf{Problem 7.18} Let $\Gamma$ be a Jordan curve and suppose that $z = 0$ lies inside it.
Wermer \cite{Wermer} showed that when $\Gamma$ has infinite length, the powers $z^n$, $n\neq0$
span all of $C(\Gamma)$. One can indicate conditions on $\Gamma$ under which the
powers $z^n$, $n\neq n_1,\ldots, n_k$, form a spanning set (see Korevaar and Pfluger, \cite{MR0379810}). Under what conditions on $\Gamma$ can one omit an infinite set
of powers, and still have a spanning set?\\
\textit{(J. Korevaar)}\\

\noindent\textbf{Update 7.18} No progress on this problem has been reported to us.\\

\noindent\textbf{Problem 7.19} For what sets $\Omega$ of lattice points $(m_k, n_k)$ do the monomials
$x^{m_k}y^{n_k}$ span $L^2$ or $C_0$ on the unit square $0\leq x\leq1$, $0\leq y\leq1$? It
is conjectured that the condition $\sum\frac{1}{m_kn_k}=\infty$ is sufficient for sets $\Omega$
in an angle $\varepsilon x \leq y \leq x/\varepsilon$, $\varepsilon>0$. Hellerstein \cite{MR0285724} has shown that the condition is not necessary.\\
\textit{(J. Korevaar)}\\

\noindent\textbf{Update 7.19} No progress on this problem has been reported to us.\\

%\section* *{Other problems}

\noindent\textbf{Problem 7.20} \textbf{(Two constant theorems for the polydisc)} Let $F(z_1, z_2)$ be
defined for $|z_1\leq1, |z_2|\leq2$ except when $z_1 = z_2$ and $|z_1|=|z_2|=1$. Suppose that $F$ is plurisubharmonic 
in $|z_1|<1, |z_2|<1$ and that $F(z_1,z_2)\leq\log\frac{1}{|z_1-z_2|}$, whenever the left hand side is defined.
Further, suppose that $F(z_1, z_2) \leq 0$ for $\{|z_1| = |z_2| = 1, z_1 \neq z_2\}$. Does it follow that 
$F(z_1,z_2)\leq0$  for $|z_1|<1$, $|z_2|<1$?\\
\textit{(L. A. Rubel, A. Shields)}\\

\noindent\textbf{Update 7.20} No progress on this problem has been reported to us.\\

\noindent\textbf{Problem 7.21} Suppose that $|z_k | = 1$ $(1 \leq k < \infty)$. Put
\[A_l=\limsup_{m\to\infty}\big|\sum^m_{k=1}z_k^l\big|.\]
It is easy to see that there is a sequence $z_k,$ for which $A_l<cl$ for
all $l$, and Clunie \cite{MR0213317} proved $A_l>cl^{\frac{1}{2}}$ for infinitely many $l$. Is there a
sequence for which $A_l = o(l)$ as $l\to\infty$?\\
\textit{(P. Erd\"os)}\\

\noindent\textbf{Update 7.21} No progress on this problem has been reported to us.\\

\noindent\textbf{Problem 7.22} Suppose that $A$ and $B$ are disjoint linked Jordan curves in $\mathbb{R}^3$
which lie at a distance $1$ from each other. Show that the length of $A$ is
at least $2\pi$. The corresponding result with a positive absolute constant
instead of $2\pi$ is due to Gehring \cite{MR0361008}.\\
\textit{(F. W. Gehring)}\\

\noindent\textbf{Update 7.22} This problem was solved by several people. Osserman's  survey article \cite[p. 1226]{zbMATH03638660} contains details of its solution and some extension. Cantarella, Fu, Kusner, Sullivan and Wrinkle \cite{MR2284052} state that this was soon solved by Marvin Ortel, but since his elegant solution was never published, they reproduce it with his permission. Also, in refining Gehring's problem to deal with links in a fixed link-homotopy class: they prove that ropelength minimisers exist and introduce a theory of ropelength criticality.\\

\noindent\textbf{Problem 7.23} The expression $u(z, \zeta)$ of Problem 6.57 is closely related to the Schwarzian
derivative $\{f(z), z\}$, e.g. it is invariant under compositions with M\"obius transformations
and
\[\lim_{\zeta\to z}u(z,\zeta)=\frac{1}{6}\{f(z),z\}\]
K\"uhnau (no citation) showed that if $f(z)$ is analytic in $\mathbb{D}$ and has a quasiconformal extension
to the whole plane with
\[|f_{\overline{z}}/f_z|\leq q<1\hspace{1cm}\text{ a.e.}\]
then
\begin{equation}
\label{C7.1}
|u(z,\zeta)|\leq q(1-|z|^2)^{-1}(1-|\zeta|^2)^{-1}.
\end{equation}
Show that (\ref{C7.1}) is also sufficient for $f$ to have a quasi-conformal extension to the
whole plane (possibly with $\frac{1}{3}\leq q < 1$). \\
\textit{(J. G. Krzyz)}\\

\noindent\textbf{Update 7.23} No progress on this problem has been reported to us.\\

\noindent\textbf{Problem 7.24} Let $E = \{|z|< 1\}$, let $0$, $a > 0$, $b = |b|e^{i\beta}(-\pi < \beta\leq\pi)$ be distinct points
of $E$; and let $\mathcal{K}$ be the family of continua $K$ with the properties that
$\{a, b\} \subset K \subset E\setminus\{0\}$ and $E\setminus K$ is connected. On any $K\in\mathcal{K}$ there is a continuous
function $\arg z$; we thus subdivide $\mathcal{K}$ into homotopy classes $\{\mathcal{K}_n\}$ according to the
value of
\[V(K)=\arg b-\arg a\in\{\beta+2n\pi,n\in\mathbb{Z}\}.\]
Let us call $K\in\mathcal{K}$ a \textit{natural continuum} if $K$ is a trajectory of a quadratic differential
$a$ with the following properties: 
\begin{enumerate}[(a)]
\item $0, a, b$ are simple poles of $\sigma$, and there is no other
pole of $\sigma$ in $E$; 
\item $\sigma$ is real on $\partial E$. 
\end{enumerate}
There are many problems that then arise, for
example:
\begin{enumerate}[(i)]
\item Do all homotopy classes contain natural continua? (There are many in
$\mathcal{K}_n$ when $|\beta+2n\pi| < 2\pi$.)
\item Find all natural continua in $\mathcal{K}$.
\item How does the modulus of $E\setminus K$ vary when $K$ runs through the natural
continua in $\mathcal{K}_n$? 
\end{enumerate}
\textit{(F. Huckemann)}\\

\noindent\textbf{Update 7.24(i)} An affirmative answer has been given by D. H. Hamilton \cite{}[EMAILED 16/4/18].\\

\noindent\textbf{Problem 7.25} Let $K$ be a compact set of positive measure in $\mathbb{C}$. Does there necessarily
exist a non-constant analytic function in $\mathbb{C}\setminus K$ with $f(\infty) = 0$ such that
$[f(z)-f(\zeta)]/[z-\zeta]\neq\pm1$ for any $z, \zeta\in\mathbb{C}\setminus K$? Conceivably, the hypotheses even
imply the existence of non-linear analytic functions $f$ with $|f(z) —f(\zeta)| / |z — \zeta|$ \textit{bounded},
but this is a well-known unsolved problem. Of course one can pose more general
problems again such as requiring that the difference quotient omit all values in some
pre-assigned plane set. (These problems arise in the use of variational methods.)\\
\textit{(D. Aharonov and H. S. Shapiro)}\\

\noindent\textbf{Update 7.25} No progress on this problem has been reported to us.\\

\noindent\textbf{Problem 7.26} Is there a homeomorphism of the open unit ball in $\mathbb{R}^3$ onto $\mathbb{R}^3$, whose
coordinate functions are harmonic? In other words, do there exist $u_1,u_2,u_3$ harmonic
in $|\textbf{x}|<1$, $\textbf{x} = (x_1, x_2, x_3)$, such that
\[ (x_1, x_2, x_3)\to (u_1, u_2, u_3)\]
is a homeomorphism of $|\textbf{x}| < 1$ onto all of $\mathbb{R}^3$? The analogous problem in $\mathbb{R}^2$ is
answered negatively; the result is due to T. Rado (no citation), and is an important lemma in the
theory of minimal surfaces.\\ 
\textit{(H. S. Shapiro)}\\

\noindent\textbf{Update 7.26} No progress on this problem has been reported to us.\\

\noindent\textbf{Problem 7.27} For a domain $D$ in $\mathbb{C}$, define
\[\rho(x,y) = \sup\{|f(x)-f(y)| : x,y\in D; f \text{ analytic in }D; |f'|\leq1\text{ in }D\}.\]
If $D$ is convex, then $\rho(x,y) = |x-y|$, but not otherwise. Clearly $\rho(x,y)\leq L(x,y)$,
the infimum of the lengths of paths in $G$ that join $x$ to $y$. What can be said about
$\rho$ for general $D$, in terms of the geometry of $D$?\\ 
\textit{(L. A. Rubel)}\\

\noindent\textbf{Update 7.27} No progress on this problem has been reported to us.\\

\noindent\textbf{Problem 7.28} Suppose that $f(z)$ is continuous in a domain $D$ and that either
\begin{enumerate}[(i)]
\item $\int_{|\zeta-z|=r}f(\zeta)\,d\zeta=0$ for all  $z\in D$ and $0<r\leq r(z)$,
\end{enumerate}
 or (weaker),
\begin{enumerate}[(ii)]
\item $\lim_{r\to0}\Big[r^{-2}\int_{|\zeta-z|=r}f(\zeta)\,d\zeta\Big]=0$ for all $z\in D.$
\end{enumerate}
Does it follow that $f(z)$ is analytic in $D$? (see e.g. \cite{MR0335835}, \cite{MR0348084}).\\
\textit{(D. Gaier and L. Zalcman)}\\

\noindent\textbf{Update 7.28} No progress on this problem has been reported to us.\\

\noindent\textbf{Problem 7.29} Let $f(z)$ be continuous on $\mathbb{D}$, and let $\alpha$ be a fixed number with
$0 < \alpha \leq 1$. If, for each $z$ in $\{|z|< 1\}$,
\[\int_{|\zeta-z|<\alpha(1-|z|)}f(\zeta)\,d\zeta=0,\]
is $f$ necessarily analytic in $\mathbb{D}$? What happens if we are given that $f$ is continuous
only in $\mathbb{D}$?\\ 
\textit{(L. Zalcman)}\\

\noindent\textbf{Update 7.29} No progress on this problem has been reported to us.\\

\noindent\textbf{Problem 7.30} Let $u(z)$ be a real bounded continuous function on $U = \{|z|< 1\}$, and
suppose that to each $z \in U$ there corresponds a real number $r(z)$ with $0<r(z)< 1-|z|$
such that
\begin{equation}
\label{C7.2}
\frac{1}{2\pi}\int^{2\pi}_0u(z+r(z)e^{i\theta})\,d\theta=u(z).
\end{equation}
Must $u(z)$ be harmonic on $U$? Volterra (no citation) showed that this was true in the case that
$u(z)$ is given to be continuous on $U$; the case in which (\ref{C7.2}) is replaced by an areal-mean-value (and the continuity condition on $u(z)$ is relaxed) has been studied by
Veech (see \cite{MR0310269} and \cite{MR0289800}) and others.\\ 
\textit{(L. Zalcman)}\\

\noindent\textbf{Update 7.30} This is the same as Problem 3.8. A complete answer was obtained by Hansen and Nadirashvili \cite{doi:10.1112/jlms/50.2.349}. See Update 3.8. \\

\noindent\textbf{Problem 7.31} Suppose that
\[a_1>0,\hspace{1cm}0\leq a_n\leq n,\hspace{1cm} n\geq1,\hspace{1cm} b_n=\sum^n_{\nu=1}a_\nu,\hspace{1cm} c_n=\sum^n_{\nu=1}b_\nu.\]
Then $\sum(a_n/c_n)^\alpha<\infty$ if $\alpha > \frac{1}{2}$. For what other functions $f(t)$ is it true that
$\sum f(c_n/a_n)<\infty$? Is it true, for instance, that (under some smoothness condition
on $f$) $\sum f(c_n/a_n)$ converges with $\sum f(n^2)$? The analogous result for $c_n/b_n$ was obtained
by Borwein \cite{MR0180776}, that if $xf(x)$ is positive and non-increasing for $x\geq a>0$ and
$\sum f(n)<\infty$, then $\sum f(c_n/b_n)<\infty$.\\
\textit{(W. K. Hayman)}\\

\noindent\textbf{Update 7.31} No progress on this problem has been reported to us.\\

\noindent\textbf{Problem 7.32} Let $\mu(t)$ be a continuous monotonic increasing function of $t$ for $t\in[0,1]$
and let $\omega_1(h,\mu)$, $\omega_2(h,\mu)$ denote its modulus of continuity and modulus of smoothness
respectively. It is known that, if
\[\omega_1(h)=O(h)\hspace{1cm}(h\to0),\]
or
\[\omega_2(h)=O(h(\log 1/h)^{-c})\hspace{1cm}(h\to0),\]
where $c > \frac{1}{2}$, then $\mu(t)$ is absolutely continuous (w.r.t. Lebesgue measure). It is also
known that each of these conditions is essentially best-possible. Are they simultaneously
best possible? More precisely, is it true that given any function
$\phi(t)\uparrow\infty$ as $t\downarrow0$, there is a continuous, monotonic increasing \textit{singular} function $\mu(t)$
such that
\[\omega_1(h)=O(h\phi(h))\hspace{1cm}(h\to0)\]
and
\[\omega_2(h)=O(h(\log 1/h)^{-\frac{1}{2}})\hspace{1cm}(h\to0)?\]
\textit{(J. M. Anderson)}\\

\noindent\textbf{Update 7.32} An affirmative answer has been given by Anderson, Fern\'andez and Shields \cite{MR979966}.\\

\noindent\textbf{Problem 7.33} Let $P(\theta)=\sum^N_{n=1}e^{i\lambda_n\theta}$ be a finite Dirichlet series with exponents $\gamma_m\neq\gamma_n$ for $m\neq n$. What can be said about $\mu\equiv\inf|P(\theta)|$? A trivial argument shows that $\mu\leq (N-1)^{1/2}$. In fact, $|P|^2=N+2\sum_{m\neq n}\cos(\lambda_m-\lambda_n)\theta$. If $w(\theta)=|P|^2-N\geq-c$, then let $h(\theta)=c+w(\theta)=2\sum b_n\cos(\delta_n\theta)$, where the $b_n$ are positive integers. Then 
\[b_n=\lim_{T\to\infty}\frac{1}{2T}\int^T_{-T}h(t)\cos(\delta_nt)\,dt\leq\lim_{T\to\infty}\frac{1}{2T}\int^T_{-T}h(t)\,dt=c.\]
Thus $c\geq1$. The problem arises in prediction theory where $\mu\leq(N-1)^{1/2}$ is adequate. If the $\lambda_n$ are rational, then the problem reduces to a problem on polynomials with coefficients $0$ and $1$.\\
\textit{(S. Rudolfer and W. K. Hayman)}\\

\noindent\textbf{Update 7.33} Goddard \cite{MR1101984} calculates $f(4)$ explicitly and, in the process, discovers some examples of Newman polynomials with few terms, but large minimum modulus.\\

\noindent\textbf{Problem 7.34} Let $\beta_j\in\mathbb{R}^+$ and $\zeta_j\in\mathbb{C}$, and for suitable small $z$ define
\begin{equation}
\label{7.341}
f(z)=\prod^n_{j=1}(1-\zeta_jz)^{\beta_j}=1+a_1z+a_2z^2+\ldots
\end{equation}
Suppose further that the $a_k$ are all real. Then there exists $N=N(\beta_1, \beta_2, \beta_3, \ldots, \beta_n)<\infty$, independent of the $\zeta_j$, such that $\min(a_1, a_2, \ldots, a_N)\leq0$. Find a sharp, or good, upper bound for $N$ as a function of the $\beta$'s. This would be significant for Tur\'an's power sum method. To see that $N<\infty$: if the $\beta_j$'s are all integers, then $f$ is a polynomial and $N=\beta_1+\beta_2+\beta_3+\ldots+\beta_n+1$. (The sort of estimate wanted in the general case.) If not, assume that $|\zeta_1|\leq|\zeta_2|\leq\ldots\leq|\zeta_n|\leq1$, and let $m=\max\{j:\beta_j\in\mathbb{Z}\}$. The radius of convergence of (\ref{7.341}) is $|\zeta_m|^{-1}$ and if the $a_k$ were all positive (or even greater than or equal to $0$). then by Pringshein's theorem (see Titchmarsh \cite{MR3155290}), $|\zeta_m|^{-1}$ would be a singularity of $f$. Since $f$ only has the singularities $\{\zeta^{-1}_j:\beta_j\notin\mathbb{Z}\}$ we deduce that $\zeta_m>0$. We multiply through by $(1-\zeta_m)^{-\beta_m}$ which has positive coefficients, and repeat the arguments. After at most $n$ steps, we obtain a contradiction. This gives a finite $N$ for any particular set of $\zeta_j$ and the uniformity is straightforward.\\
\textit{(R. Hall)}\\

\noindent\textbf{Update 7.34} No progress on this problem has been reported to us.\\

\noindent\textbf{Problem 7.35} According to Fefferman's theorem (no citation), a real function $u$ on the unit circle which has bounded mean oscillation, can be decomposed as $u=b_1+\tilde{b_2}$ where $b_1$ and $b_2$ are bounded functions and $\tilde{b_2}$ denotes the conjugate of $b_2$. Given $u$, what is the smallest possible $\|b_2\|_\infty$? This problem is discussed by Baernstein \cite{MR623463} (Section 10). An affirmative answer to \cite[Section 10, Conjecture 2]{MR623463} would prove the conjecture about factoring non-zero univalent functions that Baernstein made (see Problem 5.58).\\
\textit{(A. Baernstein)}\\

\noindent\textbf{Update 7.35} No progress on this problem has been reported to us.\\

\noindent\textbf{Problem 7.36} This problem is equivalent to Problem 7.9 due to Gehring and Reich about best bounds for area distribution under quasiconformal mapping. Let $E$ denote a measurable subset of the unit disc $\mathbb{D}$, $m$ denote $2$-dimensional measure, and define
\[f_E(z)=-\frac{1}{\pi}\int\int_E\frac{dm(w)}{(w-z)^2},\hspace{1cm}z\in\mathbb{C}\setminus E.\]
Thus, $f_E$ is the $2$-dimensional Hilbert transform of the characteristic function of $E$. It follows from the Calder\'on-Zygmund theory of singular integrals that there are constants $a$ and $b$ such that
\[\int\int_{\mathbb{D}\setminus E}|f_E|\,dm\leq am(E)\log\frac{\pi}{m(E)}+bm(E)\]
for every $E$. The problem is to find the smallest possible $a$ (for which there exists some $b$ such that the inequality holds for every $E\subset\mathbb{D}$). Consideration of $E=\{z:|z|<\delta\}$ for small $\delta$ shows that $a=1$ would be best possible, and this is conjectured by Gehring and Reich.\\
An analogous sharp inequality for sets $E\subset[-1,1]$ and $1$-dimensional Hilbert transforms
\[f_E(x)=-\frac{1}{\pi}\int_E\frac{dt}{x-t}\]
is known, and can be proved either by a subordination argument or by use of a theorem of Stein and Weiss.\\
\textit{(A. Baernstein)}\\

\noindent\textbf{Update 7.36} This has been solved by Eremenko and Hamilton \cite{zbMATH00868218}.\\

\noindent\textbf{Problem 7.37} Let $T$ denote the class of all rational functions $g$ of the form
\[g(z)=\sum^n_{j=1}\frac{\lambda_j}{(z-z_j)^2},\]
where the constants $\lambda_j$ satisfy $\lambda_j>0$ and $\sum^n_{j=1}\lambda_j=1$. Prove (or disprove): There is a constant $C$ with the property that for each $g$ in $T$ we can find a set $S=S(g)$ with $m(S)=\pi$ such that
\[\int\int_{\Delta(R)\setminus S}|g|\,dx\,dy\leq2\pi\log R+C\]
for every $R$ in $(1,\infty)$. Here $\Delta(R)=\{z:|z|<R\}$. This assertion, if true, would imply that the inequality of Problem 7.36 holds with $a=1$, and thus solve the area problem of Gehring and Reich. Problems of this sort have been considered by Fuchs and MacIntyre \cite{MR0002823}.\\
\textit{(A. Baernstein)}\\

\noindent\textbf{Update 7.37} No progress on this problem has been reported to us.\\

\noindent\textbf{Problem 7.38} The Hankel matrices of a function $f$ having a Taylor expansion
\[f(z)=a_0+a_1z+\ldots\]
are defined by
\[H^{(n)}_p=(a_{ij});\hspace{1cm}a_{ij}=a_{n+1+j-2};\hspace{1cm}1\leq i,\, j\leq p+1.\]
If $f$ belongs to the Pick-Nevanlinna class ($\text{Det } H^{(n)}_p\geq0$, all $n, p$), then all the poles of $f$ are simple and they lie on the positive real axis. Denote by 
\[\varepsilon^{(n)}_1\geq\varepsilon^{(n)}_2\geq\ldots\geq\varepsilon^{(n)}_p\geq0\]
the eigenvalues of $H^{(n)}_p$. Then
\[\limsup_{n\to\infty}(\varepsilon^{(n)}_j)^{1/n}=\frac{1}{\lambda_j},\]
where $\lambda_j$ is the $j$-th pole, where the poles are numbered in order of increasing modulus. What can be said about the eigenvalues under less restrictive conditions?\\
\textit{(R. Bouteiller)}\\

\noindent\textbf{Update 7.38} No progress on this problem has been reported to us.\\

\noindent\textbf{Problem 7.39} \textbf{(Subadditivity problem for analytic capacity)} Prove or disprove the existence of a constant $M$ such that
\[\gamma(K_1\cup K_2)\leq M\{\gamma(K_1)+\gamma(K_2)\}\]
for all compact sets $K_1, K_2$ in $\mathbb{C}$. Even the case where $\gamma(K_2)=0$ is open. For background, see Garnett \cite{MR0454006}.\\
\textit{(J. Korevaar)}\\

\noindent\textbf{Update 7.39} This has been proved by Tolsa \cite{zbMATH02115780}, and is the same as Problem 7.75(c).\\

\noindent\textbf{Problem 7.40} Let $D(z, r)=\{w:|w-z|\leq r\}$. Given a sequence $\{D(z_j, r)\}^N_{j=1}$ of disjoint closed balls all contained in $|z|\leq\frac{1}{2}$, put
\[\Omega=\{z:|z|<1\}\setminus\cup^N_{j=1}D(z_j,r).\]
Let $h$ be the function in $|z|\leq1$ satisfying: $h$ is harmonic in $\Omega$, $h=1$ on $|z|=1$, and $h=0$ on $\cup^N_{j=1}D(z_j,r)$. Does there exist $\delta>0$ such that whenever $r\leq\delta$ and $N\geq[1/r]^{2-\delta}$, then $\int^{2\pi}_0h(z_j+2re^{i\theta})\,d\theta\leq r^{2+\delta}$, holds for at least one $D(z_j,r)$? Here $[1/r]$ denotes the greatest integer less than or equal to $1/r$. An affirmative answer should imply a weak version of Arakelian's conjecture for entire functions (Problem 1.6).\\
\textit{(J. Lewis)}\\

\noindent\textbf{Update 7.40} This has been proved by Lewis and Wu \cite{MR942832}, and provides the key to their solution of Problems 1.6 and 4.18.\\

\noindent\textbf{Problem 7.41} Let $\{z_\nu\}$ be a sequence of distinct points in the unit disc $\mathbb{D}$ such that $\sum(1-|z_\nu|)<\infty$. Let $B$ be the Blaschke product corresponding to this sequence. If $0<t<\infty$, define $W_t=\{z\in\mathbb{C}:|B(z)|<t\}$. Denote the space of all bounded analytic functions in $W_t$ by $H^\infty(W_t)$, and put $S=\{z_\nu\}^\infty_{\nu=1}$. Is 
\begin{equation}
\label{I7.41}
H^\infty(\mathbb{D})|_S=H^\infty(W_t)|_S
\end{equation}
 for all $t$? $H^\infty(W_t)|_S$ denotes the restrictions to $S$ of the functions in $H^\infty(W_t)$.\\ 
We note the following:
\begin{enumerate}[(a)]
\item For an interpolating sequence $H^\infty(\mathbb{D})|_S=l^\infty$ this is true, and can easily be deduced from Earl's proof \cite{MR0284588} of Carleson's interpolation theorem.
\item For \textit{any} Blaschke sequence, the result is true when $0<t<1$ for then the result is contained in Carleson's original proof of the corona theorem for $H^\infty$.
\end{enumerate}
\textit{(A. Stray)}\\

\noindent\textbf{Update 7.41} Stray \cite{MR825925} has proved (\ref{I7.41}).\\

\noindent\textbf{Problem 7.42} Define the Harnack function $H_{z_0}(z)$ for a Green domain $D$ relative to $z_0\in D$ to be the supremum of all positive harmonic functions $h$ on $D$ to be the supremum of all positive harmonic functions $h$ on $D$ which satisfy $h(z_0)\leq1$. If $K_\zeta(z)$ is the Martin kernel for $D$ relative to $z_0$, $K(\zeta_0)=1$ for all $\zeta$ in the Martin boundary $\Delta_1$, then
\[H_{z_0}(z)=\sup\{H_\zeta(z):z\in D, \zeta\in\Delta_1\}.\]
In particular for $D$ the unit disc, the boundary $\Delta_1$ may be identified with the unit circle $\mathbb{T}$ and, if $z_0=0$, then $K_\zeta(z)$ is the Poisson kernel with pole at $\zeta\in T$. In this case, $K_\zeta=H_0$ along the radius to $\zeta$. Using the Riemann mapping functions, this property may be described for simply connected $D$ as follows: for each $\zeta\in \Delta_1$, $K_\zeta$ touches $H_{z_0}$ along a Green line for $D$ issuing from $z_0$.\\
Does this property continue to hold for $D$ multiply connected?\\
\textit{(M. Arsove and G. Johnson)}\\

\noindent\textbf{Update 7.42} No progress on this problem has been reported to us.\\

\noindent\textbf{Problem 7.43} A bounded simply-connected domain $D$ is said to be \textit{conformally rigid} if there is some $\varepsilon>0$ such that if $f$ is a conformal self-map of $D$ satisfying $|f(z)-z|<\varepsilon$, then $f(z)\equiv z$. Clearly, if each prime end of $D$ is a singleton, then $D$ is \textit{not} conformally rigid. Show the converse.\\
\textit{(P. M. Gauthier)}\\

\noindent\textbf{Update 7.43} Gaier \cite{zbMATH03834335} has shown that the converse is false.\\

\noindent\textbf{Problem 7.44} Let $U$ be an open set in the plane and $\lambda^U_a$ be the harmonic measure at the point $a$ with respect to $U$. Then \O{}ksendal \cite{MR631668} showed that $\lambda^U_a$ is singular with respect to area measure. Is it also true that  $\lambda^U_a$ is singular with respect to $\beta$-dimensional Hausdorff measure for all $\beta>1$?\\
The same question can be asked for the Keldysh measure $\mu^U_a$ at $a\in K$ with respect to a compact set $K$. \O{}ksendal showed that $\mu^U_a$ is singular with respect to area measure.\\
\textit{(B. \O ksendal)}\\

\noindent\textbf{Update 7.44} This was proved by Makarov \cite{zbMATH03914673} in the case where $U$ is simply connected; and by Jones and Wolff \cite{zbMATH04091817} in the general case.\\

\noindent\textbf{Problem 7.45} Let $D$ be the unit disc cut along $p$ radial slits from the outer boundary, all of the same given length. Let $u$ be the harmonic measure of $\{|z|=1\}$ in $D$. Find the configuration of slits which makes $u(0)$ minimal, when $p$ is fixed.\\
\textit{(A. A. Gonchar, communicated by M. Ess\'en)}\\

\noindent\textbf{Update 7.45} Dubinin \cite{MR746071} has proved that this occurs when the slits are symmetrically distributed. \\

\noindent\textbf{Problem 7.46} \textbf{(A `Universal' Phr\'agmen-Lindel\"of Theorem)} Let $D$ be an arbitrary unbounded plane domain. Suppose that $f(z)$ is analytic on $D$ and continuous on $\overline{D}$. If $|f(z)|\leq1$ on $\partial D$ and $f(z)=o(|z|)$ at $\infty$, show that $|f(z)|\leq1$ throughout $D$. The $o(|z|)$ would then be the `right' condition since that is what is needed for the case $D=\{z:|z|>1\}$.\\
\textit{(D. J. Newman)}\\

\noindent\textbf{Update 7.46} Fuchs \cite{zbMATH03739887} has shown that $|f(z)|\leq1$ throughout $D$. A slightly more general result, where $z^\alpha f(z)$ is one-valued for some real $\alpha$, was obtained by Gehring, Hayman and Hinkkanen \cite{MR663670}, and was applied to problems involving moduli of continuity. \\

\noindent\textbf{Problem 7.47} Let $K\subset\mathbb{C}$ be compact and let $x_0\in K$ be a \textit{non-peak} point for $R(K)$, the uniform limits on $K$ of rational functions with poles outside $K$. Will there always exist a continuous curve $\Gamma$ in $K$ terminating at $x_0$? It is easy to see that if $\sum^\infty_{n=1}2^n M_1(A_n(x_0)\setminus K)<\infty$, where $A_n(x_0)=\{z:2^{-n-1}\leq|z-x_0|\leq2^{-n}\}$, and $M_1$ denotes $1$-dimensional Hausdorff content, then $\Gamma$ can be chosen to be a straight line segment. In this case, \O ksendal \cite{MR568976} showed that the integrated Brownian motion starting at $x_0$ stays inside $K$ for a positive period of time, almost surely. There is an example of a compact set $K$ and a non-peak point $x_0\in K$ such that no straight line segment terminating at $x_0$ is included in $K$.\\
\textit{(B. \O ksendal)}\\

\noindent\textbf{Update 7.47} No progress on this problem has been reported to us.\\

\noindent\textbf{Problem 7.48} A domain $D\subseteq \mathbb{R}^n$ is said to be \textit{linearly accessible}, if each point in the complement of $D$ can be joined to $\infty$ by a ray which does not meet $D$. Let $g(\,\cdot\,,x_0)$ be the Green's function for $D$ with pole at $x_0$ in $D$. Is $\{x:g(x,x_0)>t\}$ linearly accessible for $0<t<\infty$. This conclusion is valid in $\mathbb{R}^2$.\\
\textit{(J. Lewis)}\\

\noindent\textbf{Update 7.48} No progress on this problem has been reported to us.\\

\noindent\textbf{Problem 7.49} Let $x=(x_1, x_2, \ldots, x_n)$ be a point in Euclidean $n$-space, $n\geq3$, with $|x|=\Big(\sum^n_{i=1}x^2_i\Big)^{1/2}$. A function $u$ on $\mathbb{R}^n$ is said to be \textit{homogeneous} of degree $m$ if $u(\lambda x)=\lambda^m u(x)$ for all $\lambda>0$. If $u$ is differentiable, put $\triangledown u=(u_{x_1}, u_{x_2},\ldots, u_{x_n})$. Prove there are no homogeneous polynomials $u$ with real coefficients, $m\geq2$, such that
\begin{equation}
\label{7.49.1}
\triangledown u\cdot(|\triangledown u|^{p-2}\triangledown u)\equiv0,\hspace{1cm}p\text{ fixed},\hspace{1cm}p\neq2,\hspace{1cm}1<p<\infty.
\end{equation}
In $\mathbb{R}^2$ there are no polynomial solutions. A proof of the above would imply that if $f$ is any solution to (\ref{7.49.1}) on a domain $D\subseteq\mathbb{R}^n$, then $f$ is real analytic in $D$ if and only if $\triangledown f$ does not vanish in $D$.\\
 \textit{(J. Lewis)}\\

\noindent\textbf{Update 7.49} No progress on this problem has been reported to us.\\

\noindent\textbf{Problem 7.50} Let $U$ be a connected open set in $\mathbb{R}^n$. Brelot and Choquet \cite{MR0062883} showed that the set of points on the boundary of $U$ which are accessible from the interior by (finite length) rectifiable paths supports harmonic measure. It is natural, in view of polygonal path connectedness of finely open sets, to ask if the same is true for finely open sets and the Keldysh measure.\\
\textit{(T. Lyons and B. \O ksendal)}\\

\noindent\textbf{Update 7.50} No progress on this problem has been reported to us.\\

\noindent\textbf{Problem 7.51} Let $D$ be a bounded strictly pseudo-convex domain in $\mathbb{C}^n, n>1$, with smooth boundary. Denote by $A^\infty(D)$ the class of functions analytic in $D$, continuous on $\overline{D}$, all of whose derivatives are continuous on $\overline{D}$. Let $E$ be a closed subset of the boundary of $D$ which is not a set of uniqueness for $A^\infty(D)$, i.e. there exists a function $f\not\equiv0$, which belongs to $A^\infty(D)$ such that $f$ vanishes exactly on $E$ and all of the derivatives of $f$ vanish on $E$. Is every closed subset of $E$ a set of non-uniqueness for $A^\infty(D)$? This is true in the case of the unit disc in $\mathbb{C}$.\\
\textit{(A. -M. Chollet)}\\

\noindent\textbf{Update 7.51} No progress on this problem has been reported to us.\\

\noindent\textbf{Problem 7.52} Let $K\subset\mathbb{C}^n$ be a compact set and let $P_0(K)$ be the set of all polynomials on $K$. The $P$-\textit{hull} of $K$, the \textit{polynomial convex hull} of $K$, is defined by
\[\text{$P$-hull } K=\{z\in\mathbb{C}^n:|p(z)|\leq\sup_k|P(z)|\text{ for all }p\in P_0(K)\}.\]
Let $P(K)$ be the uniform closure of $P_0(K)$ in $C(K)$, the continuous function on $K$. Let \v{S}ilov Bd $(P(K))$ denote the \v{S}ilov boundary of the uniform algebra $P(K)$. Determine all compact sets $K\subset\mathbb{C}^n (n>1)$ such that \v{S}ilov Bd $(P(K))=\text{Boundary }(\text{$P$-hull }K)$.\\
For $n=1$, every compact set $K\subset\mathbb{C}$ has this property. For $n>2$, examples of $K$ are compact sets and closed spheres.\\
\textit{(S. Kilambi)}\\

\noindent\textbf{Update 7.52} No progress on this problem has been reported to us.\\

\noindent\textbf{Problem 7.53} \begin{enumerate}[(i)]
\item (One-dimensional version.) Let $E$ be a compact set in $\mathbb{R}$; and, for
each $x \in E$, let $\delta_x > 0$ be given. Let $I_x = (x-\delta_x, x + \delta_x)$. For what values of $c$ can
one always find a disjoint collection of such intervals, $\{I_{x_j}\}$ say, such that
$\sum_j|I_{x_j}|\geq c|E|$? It is known that this is possible for $c = \frac{1}{2}$, but is impossible in general
for $c > \frac{2}{3}$.
\item ($n$-dimensional version.) Let $E$ be a compact set in $\mathbb{R}^n$; and let $K$ be an
open bounded symmetric convex set in $\mathbb{R}^n$. For each $x \in E$, let $\delta_x > 0$ be given; and
let $K_x = x + \delta_xK$, the dilation of $K$ by a factor $\delta_x$, centred at $x$. For what values of $c$
can one always find a disjoint collection of such sets, $\{K_x\}$ say, such that
$\sum_j|K_{x_j}|\geq c|E|$? If $c(K)$ denotes the best value, it is known that $2^{-n}\leq c(K)<1$.
\end{enumerate}
The facts in $(i)$ and $(ii)$, together with some sketchy information on $c(Q_n)$ (where
$Q_n$ is the $n$-cube) and $c(S_n)$ (where $S_n$ is the $n$-sphere), are given in \cite{MR0385047}, but there is no
information on the correct asymptotic behaviour of $c(Q_n)$ or $c(S_n)$. The problem has
applications to the best constants in results concerning the Hardy-Littlewood
maximal function.\\
\textit{(P. L. Walker)}\\

\noindent\textbf{Update 7.53} No progress on this problem has been reported to us.\\

\noindent\textbf{Problem 7.54} Let $\phi_t(z)= e^{tz}-1$, $\phi^1_t=\phi_t$ and $\phi^{n+1}_t=\phi_t\circ\phi^n_t$ for $n\geq1$; it follows that 
\[\phi^1_t(-1)=e^{-t}-1 = -t+\frac{t}{2!}-\ldots,\]
\[\phi^2_t(-1)=e^{t(e^{-t}-1)}-1 = -t^2+\ldots\text{ etc.}\]
Are the coefficients in these formal power series for $\{\phi^n_t(-1)\}^\infty_{n=1}$ uniformly bounded
by $1$ in modulus?\\
\textit{(P. J. Rippon)}\\

\noindent\textbf{Update 7.54} No progress on this problem has been reported to us.\\

\noindent\textbf{Problem 7.55} It is known that any quasi-conformal homeomorphism of
$B^n = \{x : x \in\mathbb{R}^n,|x|<1\}$ onto a Jordan domain $D$ in $\mathbb{R}^n$ can be extended to a
homeomorphism of $B^n$ onto $D$. If $\partial D$ is rectifiable (in the sense that $\Lambda^{n-1}(\partial D)<\infty$),
is $f|_{\partial B^n}$ absolutely continuous (in the sense that $\Lambda^{n-1}\big(f(E)\big)=0$ for every set $E$ in
$\partial B^n$ with $|E|=0$)? One can also ask the analogous question about $f^{-1}$.\\
When $n = 2$, the answer to both questions is `yes' for conformal mappings, but `no'
for quasi-conformal mappings. When $n = 3$, Gehring \cite{MR0382638} has proved that, if in
addition the function $f$ has a quasi-conformal extension to $\mathbb{R}^n$, then $f|_{\partial B^n}$ is
absolutely continuous; but, even in this special case, it is not known if $f^{-1}|_{\partial D}$ is
absolutely continuous.\\
\textit{(A. Baernstein II)}\\

\noindent\textbf{Update 7.55} No progress on this problem has been reported to us.\\

\noindent\textbf{Problem 7.56} Let $\Gamma$ be a closed Jordan curve in the extended plane, and suppose that
$\infty\in\Gamma$. Let $f_1, f_2$ map the upper and lower half-planes, respectively, onto the two
different domains in $\mathbb{C}$ bounded by $\Gamma$, with $f_1(\infty)=f_2(\infty)=\infty$. Then, if
$h=f^{-1}_2\circ f_2$, $h$ is a homeomorphism of $\mathbb{R}$ onto $\mathbb{R}$. It is known that $\Gamma$ is a quasi-circle
if and only if $h$ is quasi-symmetric (that is, there exists a constant $c$ such that
\[\frac{1}{c}\leq\frac{h(x+t)-h(x)}{h(x)-h(x-t)}\leq c\]
for all $x,t \in\mathbb{R}$). Can one characterise the function $h$ for general Jordan curves $\Gamma$?
In particular, can every function $h:\mathbb{R}\to\mathbb{R}$ be generated in this fashion?\\
\textit{(L. Bers; communicated by A. Baernstein II)}\\

\noindent\textbf{Update 7.56} A counter-example to the latter question has been given by Oikawa \cite{MR0125956}. A counter-example has also been given by Huber \cite{MR0425110}.\\

\noindent\textbf{Problem 7.57} For $n\geq2$, let
\[\mathbb{R}^n_+=\big\{(x_1,\ldots,x_n)\in\mathbb{R}^n:x_n>0\big\}\]
and
\[\mathbb{R}^n_-=\big\{(x_1,\ldots,x_n)\in\mathbb{R}^n:x_n<0\big\};\]
let $E$, $F$ be non-empty compact subsets of $\mathbb{R}^n_+$, $\mathbb{R}^n_-$ (respectively), and let $F^*$ denote
the symmetric image of $F$ in $\partial \mathbb{R}^n_+$. Denote by $\Delta(E,F)$, $\Delta(E,F^*)$ the families of all
curves in $\overline{\mathbb{R}^n}$ joining $E$ and $F$, $E$ and $F^*$ (respectively). Is it true that
\begin{equation}
\label{J7.57}
M\big(\Delta(E,F)\big)\leq M\big(\Delta(E,F^*)\big),
\end{equation}
where $M$ denotes the $n$-modulus of a curve family?\\
It is easy to show that strict inequality holds in (\ref{J7.57}) if $E$ and $F$ are balls. Also, if
$\Delta(E,F)$ is obtained from $\Delta(E, F^*)$ as a result of symmetrisation, then (\ref{J7.57}) holds. Note
too that it follows from the symmetry principle for the modulus that
$\frac{1}{2}M\big(\Delta(E,F)\big)\leq M\big(\Delta(E,F^*)\big)$, at least if $E\cap f^*=\emptyset$.\\
\textit{(M. Vuorinen)}\\

\noindent\textbf{Update 7.57} No progress on this problem has been reported to us.\\

\noindent\textbf{Problem 7.58} Let $E\subset[0,1]$ be a compact set on the positive real axis in $\mathbb{R}^2$; and let $E$
have positive conformal $2$-capacity, that is $M\big(\Delta(E, \partial B^2(2);\mathbb{R}^2)\big)> 0$ where
$\Delta(E,\partial B^2(2);\mathbb{R}^2)$ is the family of all curves joining $E$ to $\partial B^2(2)$ and $M(\Delta)$ is the
$2$-modulus of $\Delta$. Is it true that $M\big(\Delta(E, F; \mathbb{R}^2)\big) = \infty$, where $F = \mathbb{R}\setminus E$?\\
\textit{(A. A. Gon\^car; communicated by M. Vuorinen)}\\

\noindent\textbf{Update 7.58} No progress on this problem has been reported to us.\\

\noindent The following three problems (7.59, 7.60 and 7.61) are about polynomials in $n$ variables. We write $z = (z_1,\ldots,z_n)$ and $D = (D_1,\ldots,D_n)$ where $D_i$ denotes $\frac{\delta}{\delta z_i}$. Also $\mathcal{E}$ denotes the set of entire functions
in $\mathbb{C}^n$.\\

\noindent\textbf{Problem 7.59} Let $(P, Q)$ denote a pair of polynomials with the following property:
\begin{equation}
\label{J7.59}
\text{the map }f\mapsto P(D)(Qf)\text{ carries }\mathcal{E}\text{ bijectively onto }\mathcal{E}.
\end{equation}
If $(P, Q)$ has the property (\ref{J7.59}), is it necessarily true that $(Q, P)$ also has the property
(\ref{J7.59})? (There is little theoretical ground so far to support such a conjecture, but in all
examples the proposer has been able to check it is true.)\\
Note that it is fairly easy to show that, under the hypotheses above, the map
$F\mapsto Q(D)(PF)$ carries $\tilde{\mathcal{E}}$ bijectively onto $\tilde{\mathcal{E}}$, where $\tilde{\mathcal{E}} = \{F : F \in\mathcal{E}, F \text{ is of exponential type}\}$.\\
\textit{(H. S. Shapiro)}\\

\noindent\textbf{Update 7.59} Meril and Struppa \cite{MR806014} have shown that the answer to the first question is `no' in general.\\

\noindent\textbf{Problem 7.60} Let $P$ be a polynomial of degree $m$ in which the coefficient of $z^m_1$ is
non-zero, and let $Q(z) = z^m_1$.
\begin{enumerate}[(a)]
\item Does the pair $(P, Q)$ have the property (\ref{J7.59})?
\item Does the pair $(Q, P)$ have the property (\ref{J7.59})?
\end{enumerate}
The proposer can prove the conjectures in the case that
\[Q(z) = z^m_1 + \big(\text{polynomial in }(z_2,\ldots, z_n)\big).\]
Note that if $(a)$ were true, then it would follow that the non-characteristic Cauchy
problem with entire data on a hyperplane has a unique entire solution. (The
uniqueness follows from classical results; only the entirety is in question.)\\
\textit{(H. S. Shapiro)}\\

\noindent\textbf{Update 7.60} Meril and Struppa \cite{MR806014} have shown that the answer to $(a)$ and $(b)$ is `yes'.\\

\noindent\textbf{Problem 7.61} Is it true that, for any polynomial $P$ with complex coefficients, the
mapping $f\mapsto P^*(D)(Pf)$, where $P^*(z) = \overline{P(\overline{z})}$, is a bijection of $\mathcal{E}$ onto $\mathcal{E}$?\\
The proposer can prove that this is true when $P$ is a homogeneous polynomial;
and D. J. Newman told the proposer that he could prove the
injectivity half of the conjecture, but the proposer has seen no details of the proof.
The conjecture would follow if one could show that the partial differential equation
$P^*(D)(Pf) = z^\alpha$ where $P^*(z) = \overline{P(\overline{z})}$, has a solution $f$ that is entire and of
exponential type for every multi-index $\alpha$.\\
\textit{(H. S. Shapiro)}\\

\noindent\textbf{Update 7.61} No progress on this problem has been reported to us.\\

\noindent\textbf{Problem 7.62} Given a countable number of convergent series with positive decreasing
terms, one can find such a series converging more slowly than any of these. Without
making any assumption about the Continuum Hypothesis, can one associate with
every countable ordinal number $\alpha$ a convergent series $\sum^\infty_{n=1}x_{n,\alpha}$ with $0\leq x_{n+1,\alpha}\leq x_{n,\alpha}$
such that
\begin{enumerate}[(a)]
\item if $\alpha<\beta$, then $x_{n,\alpha}/x_{n,\beta}\to0$ as $n\to\infty$,\\
and
\item if $x_n > 0$ and $\sum^\infty_{n=1}x_n<\infty$, then there exists $\alpha$ such that $x_n/x_{n,\alpha}\to0$
as $n\to\infty$.
\end{enumerate}
[See also Problem 2.66.]\\
\textit{(A. Hinkkanen)}\\

\noindent\textbf{Update 7.62} No progress on this problem has been reported to us.\\

\noindent\textbf{Problem 7.63} For $n\geq 2$ and $\alpha > 0$, let
\[\widehat{T^\alpha_Rf}(\xi)=\Big(1-\frac{|\xi|^2}{R^2}\Big)^\alpha_+\hat{f}(\xi)\]
where $R > 0$ and $f\in\mathcal{S}(\mathbb{R}^n)$. Is it true, for all $f\in L^{2n/(n+1)}(\mathbb{R}^n)$ and all $\alpha > 0$, that
$T^\alpha_Rf(x)\to f(x)$ a.e. as $R\to\infty$? When $n = 2$ and $\{R_j\}$ is a lacunary sequence
tending to $\infty$, Carbery \cite{zbMATH03826463}, C\'ordoba and L\'opez-Melero \cite{MR638621} and also Igari \cite{igari1981} have shown that the
answer is `yes'. When $n\geq3$, the more `elementary' problem of norm convergence
remains unsolved. (See also \cite{zbMATH03834608}, \cite{MR0365016}, \cite{Francia}.)\\
\textit{(A. Carbery)}\\

\noindent\textbf{Update 7.63} Carbery, Rubio de Francia and Vega \cite{MR972135} have answered have answered the almost-everywhere convergence problem in the affirmative in all dimensions.\\

\noindent\textbf{Problem 7.64} Let $\Gamma$ be a Fuchsian group in $\mathbb{D}$, and let $i(z) \equiv z$. Is it true
that
\[\sum_{\gamma\in\Gamma}|\gamma'(0)|\geq\prod_{\gamma\in\Gamma,\,\gamma\neq i}|\gamma(0)|^2\,?\]
Note that this is an equivalent formulation of a problem connecting the Bergman
kernel function and the capacity function of a Riemann surface.\\
\textit{(N. Suita; communicated by Ch. Pommerenke)}\\

\noindent\textbf{Update 7.64} No progress on this problem has been reported to us.\\

\noindent\textbf{Problem 7.65} Let $\mathcal{W}$ be a hyperbolic Riemann surface, $G_\omega$ the Green's function with
pole $\omega\in\mathcal{W}$, and $\Gamma=\{[\gamma_n]\}$ the fundamental group. Let $\tilde{\Gamma}$ be the subgroup of
equivalence classes $[\gamma]$ such that, for every $\omega\in\mathcal{W}$, the harmonic conjugate $G^*_\omega$ of $G_\omega$ changes by an integral multiple of $2\pi$ on a representative path in $[\gamma]$. Certainly
$[\Gamma,\Gamma]\triangleleft\tilde{\Gamma}\triangleleft\Gamma$. Is it true that, if $[\Gamma, \Gamma] \neq\tilde{\Gamma}$, then $\mathcal{W}$ is necessarily of the form
$\mathcal{V}\setminus A$ where $\mathcal{V}$ is a hyperbolic surface and $A\neq\emptyset$ is a relatively closed subset of zero
logarithmic capacity?\\
The conjecture is true if $\tilde{\Gamma}=\Gamma$ and $\mathcal{V}$ is the unit disc; and it is false if
$[\Gamma,\Gamma] = \tilde{\Gamma}$. (The conjecture arises from work in \cite{MR674872}.)\\
\textit{(K. Stephenson)}\\

\noindent\textbf{Update 7.65} No progress on this problem has been reported to us.\\

\noindent\textbf{Problem 7.66} A continuous mapping $f:B^n\to\mathbb{R}^n$, where $B^n = \{x : x\in\mathbb{R}^n,|x|<1\}$, and $n\geq2$, is said to be \textit{proper} if $f^{-1}(K)$ is a compact subset of $B^n$ whenever $K$ is a
proper subset of $f(B^n)$. Let
\[B_f = \{z : z \in B^n, f\text{ is not a local homeomorphism at }z\}.\]
Is it true that, if
\begin{enumerate}[(i)]
\item $n\geq3$,
\item $f:B^n\to\mathbb{R}^n$ is proper and quasi-analytic,
\end{enumerate}
and
\begin{enumerate}[(iii)]
\item $B_f$ is compact
\end{enumerate}
then $f$ is necessarily injective?\\
Note that the mapping $f(z) = z^2$, where $z \in B^2$, shows that the conjecture is false
when $n = 2$. The conjecture is known to be true in the special case that $f(B^n) = B^n$,
$n\geq3$. The conjecture is a special case of a more general open problem in \cite{Vuorinen}.\\
\textit{(M. Vuorinen)}\\

\noindent\textbf{Update 7.66} No progress on this problem has been reported to us.\\

\noindent\textbf{Problem 7.67} Let $V$ be the zero set of some analytic function in a strictly pseudo-convex
domain $\Omega$ in $\mathbb{C}^2$. If $V$ has finite area inside $\Omega$, is it necessarily true that $V$ is
the zero set of some bounded analytic function on $\Omega$?\\
Bo Berndttson (no citation) has shown that the answer is `yes' when $\Omega$ is the ball; and easy
examples show that the answer is `no' for strictly pseudo-convex domains in $\mathbb{C}$ when
$n > 2$. Skoda (\cite{MR0385155}, \cite{MR0450620}) and independently Henkin \cite{MR0466634} have shown (under some
cohomology condition on $\Omega$) for all $n$ that $V$ is the zero set of a function of
Nevanlinna class only if $V$ satisfies the Blaschke condition.\\
\textit{(R. Zeinstra)}\\

\noindent\textbf{Update 7.67} No progress on this problem has been reported to us.\\

\noindent\textbf{Problem 7.68} For a Hadamard gap sequence $\{n_k\}^\infty_{k=1}$, $n_{k+1}/n_k\geq q>1$, is it true that the
measure in $[0,2\pi[$ of the set of those points $x$ for which
\[\liminf_{m\to\infty}\Big|\sum^m_{k=1}\cos(n_kx)-\xi\Big|=0,\hspace{1cm}\text{ for all }\xi\in\mathbb{R},\]
equals $2\pi$?\\ 
\textit{(T. Murai)}\\

\noindent\textbf{Update 7.68} No progress on this problem has been reported to us.\\

\noindent\textbf{Problem 7.69} Given an associative algebra $A$, with identity $1$ and countable basis, then
for a finite-dimensional subspace $V$ spanned by the vectors $\{e_j\}^k_{j=1}$ we have the
differential operator 
\[\sum^k_{j=1}e_j\frac{\partial}{\partial x_j}\]
acting on differentiable functions defined over domains in $V$ and taking their values
in $A$. We call a function $f: U\subseteq V\to A$ a \textit{left analytic function} if
\[\sum^k_{j=1}e_j\frac{\partial f}{\partial x_j}(x)=0\]
for all $x\in U$. Ryan \cite{MR879189} has shown that there is a generalised Cauchy integral
formula
\[f(x_0)=\int_{\partial M}G(x,x_0)\sum^k_{j=1}(-1)^je_j\,d\hat{x}_jf(x),\]
with real-analytic kernel $G(x,x_0)$, where $M$ is an arbitrary, compact real $n$-dimensional
manifold lying in $U$; and $x_0\in \mathring{M}$ if and only if there are elements
$\{p_j\}^k_{j=1}\subseteq A$ satisfying the relation
\[p_je_l+p_le_j=2\delta_{jl}.\]
\begin{enumerate}[(a)]
\item Is the result still valid if we only assume $G(x,x_0)$ to be a $C^1$ function?
\item What analogous result holds if we assume the algebra to be non-associative?
\end{enumerate}
\textit{(J. Ryan)}\\

\noindent\textbf{Update 7.69} No progress on this problem has been reported to us.\\

\noindent\textbf{Problem 7.70} Using arguments due to Ahlfors (see, for example \cite{MR863437}) any M\"obius
transformation in $\mathbb{R}^n$ can be written in the form $(ax + b)(cx+d)^{-1}$, where $x\in\mathbb{R}^n$ and
$a, b, c, d$ are elements of a Clifford algebra $A_n$ that satisfies certain constraints. It can
be shown that the linear differential equations whose solution spaces are conformally
invariant are of the type
\[D^kf_k((ax+b)(cx+d)^{-1})=0,\hspace{1cm} k\in\mathbb{N},\]
where $D$ is the Euclidean Dirac operator, and the associated conformal weight is
\[J_k(cx+d)=
\begin{cases}
(cx+d)\ast |cx+d|^{-n-1+k}&\text{for }k=2p-1,\\
|cx+d|^{-n+k}&\text{for }k=2p,
\end{cases}\]
where $\ast$ is the involution described in \cite{MR863437}.
\begin{enumerate}
\item What are the non-linear differential equations whose solution spaces are
conformally invariant?
\item Can their conformal weights also be expressed in terms of $cx + d$, and what
relationship do these solutions have to the linear conformally-invariant differential
equations?
\end{enumerate}
\textit{(J. Ryan)}\\

\noindent\textbf{Update 7.70} No progress on this problem has been reported to us.\\

\noindent\textbf{Problem 7.71} Given a domain of holomorphy $\Omega\subseteq\mathbb{C}^n$, $n\geq2$, what conditions are
required on $\Omega$ to admit an analytic function $p:\Omega\to\mathbb{C}$ which cannot be
analytically extended beyond the boundary of $\Omega$, and satisfies the complex
version 
\[\sum^n_{k=1}\frac{\partial^2}{\partial z^2_j}p(z)=0\]
of Laplace's equation?\\
\textit{(J. Ryan)}\\

\noindent\textbf{Update 7.71} No progress on this problem has been reported to us.\\

\noindent\textbf{Problem 7.72} Let $N$ denote the class of complex-valued $L^\infty$-functions $v$ on the unit disc
$U$ such that $\int_U v\phi\,dx\,dy=0$ whenever $\phi$ is analytic in $U$ with $\int_U|\phi(x+iy)|\,dx\,dy<\infty$.
The Cauchy principal value of
\[(Bv)(z)=\frac{-1}{\pi}\int_U\frac{v(\zeta)}{(z-\zeta)^2}\,d\xi\,d\eta,\hspace{1cm}\zeta=\xi+i\eta,\]
defines the Beurling transform $Bv$ of $v$.\\
Is it true that $Bv\in L^\infty$ and, furthermore, that
\begin{equation}
\label{K7.72}
\|Bv\|_\infty\leq C\|v\|_\infty,
\end{equation}
for some $C < \infty$, whenever $v\in N$. A weaker question is whether this holds for all
$v\in N\cap P$, where $P$ is the class of all polynomials in $z$ and $\overline{z}$. (The inequality (\ref{K7.72}) is true
at least for certain subclasses of $N\cap P$.)\\ 
\textit{(A. Hinkkanen)}\\

\noindent\textbf{Update 7.72} No progress on this problem has been reported to us.\\

\noindent\textbf{Problem 7.73} Let $D_1$, $D_2$ be domains in $\{|z|< R\}$, and let $\lambda_1(z)\,|dz|$ and $\lambda_2(z)\,|dz|$ be their hyperbolic metrics. What is the least number $A = A(R)$ such that the hyperbolic
metric $\lambda(z)\,|dz|$ of $D_1 \cap D_2$ satisfies the inequality
\[\lambda(z)< A(\lambda_1(z)+\lambda_2(z))?\]
\textit{(W. H. J. Fuchs)}\\

\noindent\textbf{Update 7.73} No progress on this problem has been reported to us.\\

\noindent\textbf{Problem 7.74} In their famous Acta paper \cite{MR1555303}, Hardy and Littlewood introduced the
celebrated Hardy-Littlewood maximal function in connection with complex function
theory. Since then it has proved an invaluable tool in real analysis. Here we ask some
questions about the dependence of constants on dimension.
Let $B$ be a convex compact symmetric body in $\mathbb{R}^n$, normalised to have Euclidean
volume $1$. Let the Hardy-Littlewood maximal functions be
\[Mf(x)=\sup_{t>0}\big(\frac{1}{t^n}\int_{tB}|f(x+y)|\,dy\Big)\]
and
\[\tilde{M}f(x)=\sup_{k\in Z}\big(\frac{1}{2^{kn}}\int_{2^kB}|f(x+y)|\,dy\Big).\]
\begin{enumerate}[(a)]
\item If $B$ is the Euclidean ball in $\mathbb{R}^n$, does there exist a constant $C$ such that
$\text{meas}\{x:\tilde{M}f(x)>\lambda\}\leq C\lambda^{-1}\|f\|_1$ for all $\lambda> 0$, with $C$ independent of $n$?
\item If so, what is the answer to the same question for $Mf$?
\item In the case that $n = 1$, a conjecture of F. Sonia and the proposer is that the
best constant in the inequality $\text{meas}\{x:Mf(x)>\lambda\}\leq C\lambda^{-1}\|f\|_1$
 is $C = \frac{3}{2}$. Prove this.
\item Let $1<p\leq\frac{3}{2}$. Can the best constant in the inequality
\[\|Mf\|_p\leq C_p\|f\|_p\]
be taken to be independent of $n$ and the body $B$? Even if $B = [-\frac{1}{2},\frac{1}{2}]^n$, can the constant
be chosen independent of $n$?
\end{enumerate}
The following relevant facts are known. For $(a)$ and $(b)$, the best known constants
are found in \cite{MR727348}. For $(c)$, Sonia and the proposer have shown that the answer is `yes'
if $\frac{3}{2}< p \leq\infty$; for $B$ with suitably-curved boundary, the answer is `yes' for $1< p\leq\infty$;
for the sphere, see \cite{MR699317}.\\ 
\textit{(A. Carbery)}\\

\noindent\textbf{Update 7.74(b)} Melas (see \cite{MR2011925}  and \cite{MR1973058})  investigates the exact value of the best possible constant $C$ for the weak-type $(1, 1)$ inequality for the one-dimensional centered Hardy-Littlewood maximal operator. In connection to this problem, the sharp bound for the weak-type $(1, 1)$ inequality for the $n$-dimensional Hardy operator is obtained by Zhao, Fu and Lu \cite{MR2972625}.\\
\textbf{(c)} Melas \cite{MR1897399}  studies the centered Hardy-Littlewood maximal operator acting on positive linear combinations of Dirac deltas, and uses this to obtain improvements in both the lower and upper bounds or the best constant $C$ in the $L^1 \rightarrow \text{weak } L^1$ inequality for this operator. A counter-example is given by Aldaz \cite{MR1606325}, and  Aldaz \cite{MR1745465} also shows that the general conjecture fails for every $n\geq2$, and also asymptotically, i.e. $\lim\inf c_n>\lim\left(\frac{1+2^{1/n}}{2}\right)^n=\sqrt{2}$.\\

\noindent\textbf{Problem 7.75} Prove or disprove the following statements about analytic capacity $\gamma$.
\begin{enumerate}[(a)]
\item If $E\subseteq\mathbb{C}$ is compact and $\phi$ is a $C^1$-diffeomorphism of $\mathbb{C}$ onto $\mathbb{C}$, then
$\gamma(E) = 0$ if and only if $\gamma(\phi(E))=0$. The statement is false if $\phi$ is a homeomorphism or
a quasi-conformal mapping.
\item If $E\subseteq\mathbb{C}$ is compact and $\phi\in \text{GL}(2,\mathbb{R})$, then $\gamma(E) = 0$ if and only if
$\gamma(\phi(E))=0$.
\item If $E$, $F$ are compact subsets of $\mathbb{C}$, then there exists a constant $K > 0$
(independent of the choice of $E$ and $F$) such that
\[\gamma(E\cup F)\leq K(\gamma(E)+\gamma(F)).\]
Perhaps one can take $K = 1$? See \cite{MR0350009} for related results.
\item If $K$ is a compact subset of $\mathbb{C}$ with $\gamma(K) = 0$, then
\[\gamma(E\setminus K) = \gamma(E)\]
for all compact subsets $E$ of $\mathbb{C}$. Here $\gamma(E\setminus K)$ means the inner capacity
\[\sup\{\gamma(L):L\subseteq E\setminus K, L\text{ compact}\}.\]
\end{enumerate}
An interesting special case would be that when $K$ is the `corner quarters square
Cantor set' (the so-called \textit{Garnett set}).\\ 
\textit{(A. G. O'Farrell)}\\

\noindent\textbf{Update 7.75} Part (c) is the same as Problem 7.39, which has been proved by Tolsa \cite{zbMATH02115780}. If $E\subset \mathbb{R}$ then $\gamma (E)=\tfrac{1}{4}\text{length}(E)$, see Pommerenke \cite{zbMATH03153281}.\\

\noindent\textbf{Problem 7.76} Prove or disprove the following statement: If $K$ is a compact subset of $\mathbb{C}$
with continuous analytic capacity $\alpha(K)= 0$, then
\[\alpha(E\setminus K)=\alpha(E)\]
for all compact subsets $E$ of $\mathbb{C}$.\\
An interesting special case is when $K$ is a $C^1$-arc. The case when $K$ is a
$C^{1+\varepsilon}$-arc has already been settled.\\
\textit{(A. G. O'Farrell)}\\

\noindent\textbf{Update 7.76} Note that the continuous analytic capacity is different from the analytic capacity used in Problem 7.75. No progress on this problem has been reported to us.\\

\noindent\textbf{Problem 7.77} We will say that $g(x)$ is a \textit{rearrangement} of $f(x)$ if
\[m\{x:g(x) < y\} = m\{x:f(x) < y\}\hspace{1cm}\text{ for all } y\in\mathbb{R},\]
where $m$ is Lebesgue measure and $f$ and $g$ are defined on some finite interval $I$. What
are those functions $f(x)$ for which $f'(x)$ is a rearrangement $f(x)$?\\
Obvious examples of such functions are $f(x)=ke^x$ on any interval and
$f(x) = k\sin x$ on $[0,\frac{1}{2}\pi]$. What others are there? How about other `differential
rearrangements' than $f'(z)\sim f(x)$?\\ 
\textit{(L. A. Rubel)}\\

\noindent\textbf{Update 7.77} No progress on this problem has been reported to us.\\

\noindent\textbf{Problem 7.78} Does there exist a sequence $\{z_n\}^\infty_1$ of distinct complex numbers such that
\[\sum\frac{1}{|z_n|}<\infty\hspace{1cm}\text{and}\hspace{1cm}\sum\frac{1}{z-z_n}\neq0,\]
for all $z\in\mathbb{C}$?\\
This has the following physical interpretation. If we imagine electrons (really unit-charged
wires perpendicular to the complex plane) placed at each point $z_n$, then these
generate a logarithmic potential given by $\sum\log|z-z_j|$. The gradient of this potential
is $\sum1/(z-z_n)$. Thus the question is whether such a field must always have an
equilibrium point - that is, a point where a free electron (or wire), once placed there,
would remain there.\\
Of course, the corresponding problem could be asked for $\mathbb{R}^n$, $n\geq3$ also.\\
\textit{(L. A. Rubel)}\\

\noindent\textbf{Update 7.78} For $n=2$ this was solved by Clunie, Eremenko and Rossi \cite{zbMATH00561242}.
They also have partial results for $n\geq 3$. Further generalisations
were obtained by Eremenko, Langley and Rossi \cite{zbMATH00569439}.\\

\noindent\textbf{Problem 7.79} Let $f$ be analytic on a domain $G$ in $\mathbb{C}$. We will say that a point $z_0\in G$
is a \textit{MacLane point} of $f$ if there exists some neighbourhood $N$ of $z_0$ such that the
restrictions to $N$ of the successive derivatives of $f$,
\[\{f^{(n)}|_N:n\in\mathbb{N}\},\]
form a normal family of functions on $N$. Let $M(f)$ denote the set of MacLane points
of $f$.\\ 
What can be said about the set $M(f)$, besides the fact that it is open? Must
$M(f)$ be connected? If $G$ is simply-connected, must $M(f)$ be simply-connected? Or
can $M(f)$ be an arbitrarily prescribed open subset of $G$?\\
Similar questions can be asked about functions meromorphic in $G$. Perhaps it is
more natural, also, to ask such questions about$\{f^{(n)}(z)/n!\}$ rather than simply $\{f^{(n)}(z)\}$?\\
Some relevant facts can be found in \cite{MR0053231} and \cite{MR0087741}.\\ 
\textit{(L. A. Rubel)}\\

\noindent\textbf{Update 7.79} No progress on this problem has been reported to us.\\

\noindent\textbf{Problem 7.80} Let $f$ be an inner function, with $f(0) = 0$; then $f$ induces an ergodic
(Lebesgue-) measure-preserving map of the circle onto itself. What is the entropy,
$h(f)$, of $f$? It is conjectured that $h(f) < \infty$ if and only if $f'$ belongs to the Nevanlinna
class; and that, in that case, then
\[h(f)=\frac{1}{2\pi}\int^{2\pi}_0\log\big|f'(e^{i\theta})\big|\,d\theta.\]
(See \cite{MR845869} and \cite{MR620121}.)\\ 
\textit{(J. L. Fernandez)}\\

\noindent\textbf{Update 7.80} No progress on this problem has been reported to us.\\

\noindent\textbf{Problem 7.81} Let $I$ denote the class of all inner functions. Then $I$, as a subset of $H^\infty$,
enjoys some of the properties that the collection of unimodular functions have as a
subset of $L^\infty_\mathbb{R}$ (the real-valued functions in $L^\infty(\Pi)$ \cite{MR942827}).\\
Is it true that if $\{\Lambda_n\} \subset (H^\infty)^*$ and if for each function $\phi\in I$ one has $|\Lambda_n(\phi)|\leq C(\phi)$,
then
\[\sup_n\|\Lambda_n\|_{(H^\infty)^*}<\infty?\]
This is known to hold if $\Lambda_n\in L^1\setminus H^1_0(\subset(H^\infty)^*)$ (\cite{MR959269}, \cite{MR916071}, \cite{MR942827}); the corresponding real-variable
result is also known \cite{MR942827}.\\ 
\textit{(J. L. Fernandez)}\\

\noindent\textbf{Update 7.81} No progress on this problem has been reported to us.\\

\noindent\textbf{Problem 7.82} Let $E\subset \mathbb{C}$, and the function $F:E\times \mathbb{D}\to\mathbb{C}$ satisfy the
following conditions:
\begin{enumerate}[(a)]
\item $F$ is injective on $E$, for all $z\in \mathbb{D}$;
\item $F$ is analytic in $z\in \mathbb{D}$, for each $w\in E$;
\item $F(z) = z$, when $w = 0$.
\end{enumerate}
Does there exist a function $G: \mathbb{C}\times \mathbb{D} \to \mathbb{C}$ that satisfies conditions $(a)$, $(b)$ and $(c)$,
and for which $G = F$ on $E$?\\
\textit{(D. Sullivan, W. Thurston, H. Royden; communicated by D. Hamilton)}\\

\noindent\textbf{Update 7.82} Slodkowski \cite{zbMATH00008459} has proved that the answer is `yes'.\\

\noindent\textbf{Problem 7.83} Let $G$ be a finitely-generated Fuchsian group of the first kind.
\begin{enumerate}[(a)]
\item If $F:\Pi \to \Pi$ is a $G$-invariant quasi-symmetric function, is $F$ totally singular?
\item Is the Teichmuller space $T_G$ dense in the space $S_G$ (Schwarzians of
$G$-invariant univalent functions)?
\end{enumerate}
\textit{(O. Lehto; L. Bers; I. Kra; communicated by D. Hamilton)}\\

\noindent\textbf{Update 7.83} No progress on this problem has been reported to us.\\

\section*{New Problems}

\noindent To appear.

%======================= CHAPTER 8 ====================================

\chapter{Spaces of Functions}

\underline{Preface by F. Holland (to appear)}\\

\noindent\textbf{Problem 8.1} Let $L = \{L\}^\infty_0$ be a non-negative increasing sequence such that
\[\sum^\infty_{k=0}L_kr^k<\infty\hspace{1cm}(0<r<1).\]
If $f(z)$ is analytic in $|z|< 1$, we will say that $f\in\mathcal{P}_L$ if there exists a constant $A$ such
that, for each integer $n\geq0$ and each polynomial $P_n$ of degree $n$,
\[\|P_n\ast f\|_\infty\leq AL_n\|P\|_\infty,\]
where $\ast$ denotes the Hadamard product. The infimum of such $A$ for a given $f\in\mathcal{P}_L$
defines a norm on $\mathcal{P}_L$ which then becomes a Banach space. A variety of linear
operators (such as subordination) have the property that they are norm-decreasing
on $\mathcal{P}_L$. This enables one to obtain sharp coefficient inequalities for subordinate
functions etc., once a function has been shown to lie in $\mathcal{P}_L$. The spaces $\mathcal{P}_L$ are
`large-growth' spaces; for example, the case $L_k = 1$ $(k\geq 0)$ is the space of Cauchy-Stieltjes transforms. Convex sequences $\{L_k\}^\infty_0$ (such as $L_k = k$) have the property
that
\[\sum^\infty_{k=0}L_kz^k\]
lies in the unit ball of $\mathcal{P}_L$.\\
Can one relate known spaces to these spaces, for example $H^p$, for $0 < p < 1$?\\
\textit{(T. Sheil-Small)}\\

\noindent\textbf{Update 8.1} No progress on this problem has been reported to us.\\
 
\noindent\textbf{Problem 8.2} Given functions $f_1,\ldots,f_N\in H^\infty$, let $I = I(f_1,\ldots,f_N)$ be the ideal of $H^\infty$
generated by $f_1,\ldots,f_N$ and let $J = J(f_1,\ldots,f_N)$ denote the set of all $g\in H^\infty$ for which
there exists a constant $C = C(g)\geq 0$ for which
\[|g(z)|\leq C[|f_1(z)|+\ldots+|f_N(z)|],\hspace{1cm}|z|<1.\]
$J$ is an ideal of $H^\infty$ which contains $I$. The corona theorem states that $I=H^\infty$ if
and only if $J = H^\infty$; in general $I\not\subseteq J$, and one seeks further relations between $I$ and
$J$ when these are proper ideals. In particular does there exist an absolute constant
$\kappa>0$ such that, if $g\in H^\infty$ and
\[|g(z)|^\kappa\leq C[|f_1(z)|+\ldots+|f_N(z)|],\hspace{1cm}|z|<1,\]
then necessarily $g\in I$? (If so, we must have $\kappa\geq2$.) As a special case, is it true that
$J^2\subset I$? This is true in appropriate algebras of functions defined in terms of faster
rates of growth as $|z|\to1$.\\
\textit{(J. J. Kelleher)}\\

\noindent\textbf{Update 8.2} Von Renteln points out that the problem is improperly stated. The problem should ask for a constant $k$, $k\geq2$ such that \linebreak\mbox{$|g(z)|\leq C\big(|f_1(z)|+\ldots+|f_N(z)|\big)$}, $|z|<1$, implies $g^k\in I(f_1,\ldots,f_N)$. Wolff \cite{WolffT} proved the result with $k=3$, and Treil \cite{MR1945294} showed the $k=3$ is best possible.\\

\noindent\textbf{Problem 8.3} For a bounded plane domain $D$, denote by $N(D)$ the class of functions analytic
on $D$ of bounded characteristic (i.e., all quotients of functions in $H^\infty(D)$ with nonvanishing
denominator); let $f_1,\ldots,f_N\in N(D)$ have no common zeros in $D$. Find
necessary and sufficient conditions on $f_1,\ldots,f_N$ in order that they generate the full
ring $N(D)$. Equivalently, if $g_1,\ldots,g_N\in H^\infty(D)$, when does the ideal generated by
$g_1,\ldots,g_N$ in $H^\infty(D)$, or in $N(D)$, contain a non-vanishing function?\\
For $D = \mathbb{D}$, for example, the zeros of $g_1,\ldots,g_N$ should not get too close
together as one approaches $\partial D$, and one would like to obtain a corona-type theorem
for this problem - i.e., a lower estimate for
\[|g_1(z)|+\ldots|g_N(z)|\]
in $\mathbb{D}$.\\
\textit{(J. J. Kelleher)}\\

\noindent\textbf{Update 8.3} For the unit disc, this problem is solved by Dahlberg, Kelleher and Taylor (not published). Earlier partial results are due to Mantel \cite{MR0457732} and von Renteln \cite{MR0486534}.\\

\noindent\textbf{Problem 8.4} Let $D$ be a simply-connected domain in the complex plane $\mathbb{C}$ and $A(D)$ the
ring of functions $f: D\to\mathbb{C}$ analytic in $D$. Bers (no citation) has shown that (for domains
of arbitrary connectivity) the algebraic structure of $A(D)$ determines the conformal
structure of $D$. Can $A(D)$ be the direct sum of two non-trivial subrings of itself?
(This would represent a generalisation of Taylor's Theorem.) Is there a more general
result for multiply-connected domains $D$?\\
\textit{(J. J. Kelleher)}\\

\noindent\textbf{Update 8.4} No progress on this problem has been reported to us.\\
  
\noindent\textbf{Problem 8.5} Let $G$ be a domain in $\mathbb{C}$, and $H(G)$ the ring of functions analytic on $G$. It is
known that, for two domains $G_1$ and $G_2$, $H(G_1)$ is isomorphic to $H(G_2)$ if and only
if $G_1$ and $G_2$ are conformally (or anticonformally) equivalent. What can be said
under only the hypothesis that $H(G_1)$ and $H(G_2)$ are elementarily equivalent in the
sense of model theory? For a large class of domains the corresponding problem has
been solved for $H_\mathbb{C}(G)$, the algebra of functions analytic on $G$, by Henson and Rubel
(see \cite{zbMATH03845559} and \cite{zbMATH03935013}).\\
\textit{(L. A. Rubel)}\\

\noindent\textbf{Update 8.5} Some real progress has been made by Becker, Hansen and Rubel \cite{MR584077}. The emphasis is on $H_\mathbb{C}(G)$. The problem as stated is open when $G_1$ and $G_2$ are just annuli.\\ 

\noindent\textbf{Problem 8.6} Let $A^p$, $p > 0$ be the Bergmann space of functions $f(z)$ analytic in $|z|< 1$
such that
\[\|f\|_p=\Huge({\int\int}_{|z|<1}|f(re^{i\theta})|^p\,r\,dr\,d\theta\Huge)^{1/p}<\infty;\]
clearly $H^p \subset A^p$. Horowitz \cite{MR0357747} has shown that if $f\in A^p$ and $f$ has zeros $\{z_k\}$ in $|z|< 1$
then
\begin{equation}
\label{C8.1}
\prod^n_{k=1}|z_k|^{-1}=O(n^{(1/p)+\varepsilon})\hspace{1cm}\text{ as }n\to\infty,
\end{equation}
for all $\varepsilon>0$; in (\ref{C8.1}) $\varepsilon$ cannot be replaced by $0$.\\
Recall that $\{z_k\}$ is a zero set for $H^p$ if and only if $\sum(1-|z_k|)<\infty$. Characterise the zero sets for $A^p$, or at least find some non-trivial converse to (\ref{C8.1}).\\
\textit{(P. L. Duren)}\\ 

\noindent\textbf{Update 8.6} For a finite set $F$ on $|z|=1$ construct the domain
\[G_F=\{|z|<1\}\setminus\cup_\nu R(I_\nu),\]
where $I_\nu$ are the complementary arcs of $F$ and
\[R(I)=\big\{z=re^{i\theta}:1-\frac{|I|}{2\pi}\leq r<1,e^{i\theta}\in\overline{I}\big\},\]
$|I|$ being the angular length of an arc $I$.\\
Set further
\[\kappa(F)=\frac{1}{2\pi}\sum_\nu\frac{|I_\nu|}{2\pi}\log\Big(\frac{2\pi}{|I_\nu|}+1\Big).\]
Then there exists an absolute constant $\lambda>1$ such that the two conditions
\begin{enumerate}[(a)]
\item \[\sup_F\Big\{\sum_{z_k\in G_F}(1-|z_k|)-\lambda p^{-1}\kappa(F)\Big\}<\infty,\]
\item \[\sup_F\Big\{\sum_{z_k\in G_F}(1-|z_k|)-\lambda^{-1}p^{-1}\kappa(F)\Big\}<\infty,\]
\end{enumerate}
are respectively necessary and sufficient for $z_k$ to be a zero set for $A^p$ ($\sup$ is taken over all finite $F$). From this it follows that 
\begin{enumerate}[(c)]
\item \[\sup_F\Big\{\kappa(F)^{-1}\sum_{z_k\in G_F}(1-|z_k|)\Big\}<\infty\]
\end{enumerate}
is necessary and sufficient for $\big\{z_k\big\}$ to be a zero set for some $A^p (p>0)$.\\
This is a modified version of a theorem of Korenblum \cite{MR0425124}.\\ 

\noindent\textbf{Problem 8.7} We use the notation of Problem 8.6. Let $\{z_n\}^\infty_1$ be a sequence of
points in $\mathbb{D}$ such that the kernel functions $k_n(z) = (1-z_nz)^{-2}$ do not span the
space $A^2$, and let $\{f_j\}^\infty_1$ be finite linear combinations of the $k_n$ which converge (in
norm) to some function $f$ in $A^2$. Prove that, if the sequence $\{f_j\}$ converges uniformly
to $0$ on some disc $\Delta$ in $\{|z| > 1\}$, then $f\equiv0$.\\
Similar problems can of course be stated for the functions $(1- z_nz)^{-1}$, and for
spaces other than $A^2$. This slightly `weird' problem arises in the theory of
generalised analytic continuation; it is known to be true if the closure of $\{z_n\}^\infty_1$ does
not contain all of $\{|z| = 1\}$; it is also known that the analogous result for $H^2$ is true.\\
\textit{(H. S. Shapiro)}\\

\noindent\textbf{Update 8.7} No progress on this problem has been reported to us.\\

\noindent\textbf{Problem 8.8} Suppose that $F$ is a relatively-closed subset of $\mathbb{D}$, and let
\[\|f\|_F=\sup\{|f(z)|;z\in F\}\]
for functions $f$ in the Bergman space $A^2$. Describe geometrically the set
\[\{z: |z|< 1, |f(z)|\leq\|f\|_F\text{ for all }f\in A^2\}.\]
\textit{(A. Stray)}\\

\noindent\textbf{Update 8.8} No progress on this problem has been reported to us.\\

\noindent\textbf{Problem 8.9} Determine
\[\|\Lambda\|=\sup\Big|{\int\int}_{|z|<1}f(z)\phi(z)\,d\sigma_z\Big|\]
over those $f\in A^1$ with $\|f\|_1\leq1$, where
\[\phi(z)=\text{sgn}(\text{Re } z).\]
It is shown by Reich and Strebel \cite{MR0361065} that $\|\Lambda\|<1$, and that there exists an extremal
function for the problem. See also Harrington and Ortel \cite{MR0409828}. Problems of this type are of
interest in connection with the theorem of Hamilton \cite{zbMATH03332697} in quasi-conformal mapping.\\
\textit{(K. Strebel; communicated by M. Ortel)}\\

\noindent\textbf{Update 8.9} Reich \cite{MR0499172} proved $\|\Lambda\|<1$. Another proof is due to Harrington and Ortel \cite{MR0425117} who also showed in \cite{MR0409828} that the supremum is actually attained. A different (analytic) proof that $\|\Lambda\|<1$ is given by Anderson \cite{MR569247}.\\

\noindent\textbf{Problem 8.10} If $f$ and $1/f$ belong to the Bergman space $A^2$, does it follow that $\mathcal{P}f$ is dense in $A^2$? Here $\mathcal{P}f$ denotes the set of all polynomial multiples of $f$, i.e.
\[\mathcal{P}f = \{pf: p\text{ a polynomial}\}.\]
More generally, if $f\in A^2$ and $|f(z)|\geq c(1-|z|)^a$ for some $a$, $c > 0$, then does it follow
that $\mathcal{P}f$ is dense in $A^2$? (For partial results in this direction, see \cite{MR0365150}).\\
\textit{(A. L. Shields)}\\

\noindent\textbf{Update 8.10} No progress on this problem has been reported to us.\\

\noindent\textbf{Problem 8.11} Let $g$ be a function in the space $D$ of functions analytic in $|z|< 1$ with
finite Dirichlet integral, i.e.,
\[{\int\int}_{|z|<1}|h'(z)|^2\,dx\,dy<\infty.\]
If $\mathcal{P}g$, as defined in Problem 8.10, is dense in $D$ and if, for some $f\in D$,
$|f(z)|\geq|g(z)|$ for all $|z|<1$, is it necessarily true that $\mathcal{P}f$ is dense in $D$?\\
The analogous result is true in $H^2$ and in $A^2$. Shields \cite{MR0361836} solved the special case
$g(z)\equiv1$.\\
\textit{(A. L. Shields)}\\

\noindent\textbf{Update 8.11} No progress on this problem has been reported to us.\\

\noindent\textbf{Problem 8.12} Let $A$ denote the set of functions continuous on $|z|= 1$ which extend
continuously to analytic functions on $|z|< 1$ (the disc algebra). Let $\tilde{A}$ denote the set
of functions of the form $f\circ\phi$ where $f$ ranges over $A$ and $\phi$ ranges over the set of
sense-preserving homeomorphisms of $|z| = 1$. Find a `good' characterisation of $\tilde{A}$.\\
Is there a function in $\tilde{A}$ which coincides with
\[\sum^\infty_{n=1}2^{-n}\exp(-i2^n\theta)\]
on some subset of $|z|= 1$ having positive measure? The proposer conjectures not. (This problem
is related to generalised analytic continuation.) \\
\textit{(H. S. Shapiro)}\\

\noindent\textbf{Update 8.12} No progress on this problem has been reported to us.\\

\noindent\textbf{Problem 8.13} Does there exist a singular measure in Zygmund's class $A^*$ all of whose
Fourier-Stieltjes coefficients are non-negative?\\
\textit{(F. Holland)}\\

\noindent\textbf{Update 8.13} No progress on this problem has been reported to us.\\

\noindent\textbf{Problem 8.14} Which functions in $L^\infty$ on the unit circle generate positive Hankel
operators?\\
\textit{(F. Holland)}\\

\noindent\textbf{Update 8.14} No progress on this problem has been reported to us.\\

\noindent\textbf{Problem 8.15} Characterise the Hankel operators on the Hardy space $H^2$ on the circle
that are of trace class. \\
\textit{(F. Holland)}\\

\noindent\textbf{Update 8.15} No progress on this problem has been reported to us.\\

\noindent\textbf{Problem 8.16} Miles \cite{Milesjoseph} and Rudin \cite{zbMATH03283536} have shown that in $\mathbb{C}^n$ a function analytic in the unit
polydisc and in the Hardy class $H^1$ may not be expressible as the product of two
functions in $H^2$, if $n\geq3$. Is this result also true for $n = 2$?\\ 
\textit{(J. G. Clunie)}\\

\noindent\textbf{Update 8.16} Rudin notes that the problem was settled for $n=2$, and hence for all $n\geq2$, by Rosay \cite{MR0377098}. The analogous problem is open (for all $n>1$) if the polydisc $U^n$ is replaced by the unit ball $B^n\subset\mathbb{C}^n$.\\

\noindent\textbf{Problem 8.17} In the ring of bounded analytic functions on the unit ball or the polydisc
in $n$ variables, is the intersection of two finitely-generated ideals again finitely
generated? This was proved for $n = 1$ by McVoy and Rubel \cite{MR0410393}.\\ 
\textit{(L. A. Rubel)}\\

\noindent\textbf{Update 8.17} No progress on this problem has been reported to us.\\

\noindent\textbf{Problem 8.18} For which simply-connected domains $D$ (with $0 \in D$) is it true that there is
a constant $K = K(D)$ such that
\begin{equation}
\label{J8.18}
{\int\int}_D|f|^2\,dx\,dy\leq K{\int\int}_D|f'|^2\,dx\,dy
\end{equation}
for all functions $f$ analytic on $D$ with $f(0) = 0$? This inequality (\ref{J8.18}) is known as the \textit{analytic Poincar\'e inequality}.\\ 
Courant and Hilbert \cite{MR0009069} have given a Jordan
domain for which (\ref{J8.18}) is false; and Hummel \cite{MR0084546} has given an example of a spiral
domain $D$ for which ${\int\int}_D|f|^2\,dx\,dy = \infty$ and ${\int\int}_D|f'|^2\,dx\,dy < \infty$.\\
\textit{(D. H. Hamilton)}\\

\noindent\textbf{Update 8.18} No progress on this problem has been reported to us.\\

\noindent\textbf{Problem 8.19} Let $\mu\geq1$ be a singular measure on the boundary of the unit disc $\mathbb{D}$; and
let $S_\mu$ be the corresponding singular inner function
\[S_\mu(z)=\exp\Big\{\int^{2\pi}_0\frac{z+e^{i\theta}}{z-e^{i\theta}}\,d\mu(\theta)\Big\},\hspace{1cm}z\in \mathbb{D}.\]
The function $S_\mu$ is said to be \textit{discrete} if $\mu$ is discrete, and \textit{continuous} if it has no
discrete part. If $S_\mu$ is discrete, does there exist some $\delta > 0$ such that $\|S_\mu-s_\nu\|_\infty>\delta$
for all continuous functions $S_\nu$?\\
The answer is `no' if `discrete' and `continuous' are interchanged \cite{MR0481015}.\\
\textit{(K. Stephenson)}\\

\noindent\textbf{Update 8.19} No progress on this problem has been reported to us.\\

\noindent\textbf{Problem 8.20} Let $f_1,\ldots,f_n$ and $g$ be $H^\infty$ functions on $\mathbb{D}$. If
$|g(z)| \leq |f_1(z)| + \ldots + |f_n(z)|$, are there necessarily functions $g_1,\ldots,g_n\in H^\infty$ such that
\begin{equation}
\label{J8.20}
g^2=f_1g_1+\ldots+f_ng_n?
\end{equation}
Wolff (see \cite{MR628971}, page 329) has shown that (\ref{J8.20}) holds with $g^2$ replaced by $g^3$, and Rao
(see \cite{MR0210910}; and \cite{MR628971}, Exercise 3, page 369) has shown that (\ref{J8.20}) is false with $g^2$ replaced by $g$. For
general background on the problem see \cite[Chp. VIII]{MR628971}.\\
\textit{(J. B. Garnett)}\\

\noindent\textbf{Update 8.20} Woytaszczyk notes that a paper by Tolokonnikov \cite{MR629839} is related to this problem.\\

\noindent\textbf{Problem 8.21} Let $w\geq0$ be an integrable function on the circle or on the line. Then $w$ is
said to belong to the \textit{class $A_2$} if
\[\sup_l\Big(\frac{1}{|l|}\int_lw\,dx\Big)\Big(\frac{1}{|l|}\int_l\frac{1}{w}\,dx\Big)<\infty\]
over all (appropriate) intervals $I$ ; and $w$ is said to satisfy the \textit{Helson-Szeg\"o condition}
if $\log w = u + \tilde{v}$ where $u\in L^\infty,\|v\|_\infty<\frac{1}{2}\pi$ and $\tilde{v}$ is the conjugate function of $v$. It is
known that the Helson-Szeg\"o condition is equivalent to $A_2$, since both conditions
are necessary and sufficient for
\begin{equation}
\label{J8.21}
\int_l|\tilde{f}|^2w\,dx\leq\text{ const.}\int_l|f|^2w\,dx.
\end{equation}
The problem is to prove that the Helson-Szeg\"o condition and $A_2$ are equivalent
without using (\ref{J8.21}). For information on the case $\|v\|_\infty<\pi$, see \cite{MR687639}.\\
\textit{(J. B. Garnett and P. W. Jones)}\\

\noindent\textbf{Update 8.21} No progress on this problem has been reported to us.\\

\noindent\textbf{Problem 8.22} Do there exist inner functions in any strictly pseudo-convex domain of
$\mathbb{C}^n$, $n\geq2$? Alexandrov (no citation) and independently L{\o}w (no citation) and Hakim and Sibony \cite{MR0457784} have
shown the existence of inner functions in the case of the ball.\\
\textit{(R. Zeinstra)}\\

\noindent\textbf{Update 8.22} There is a positive answer by L\o w (see \cite{MR2632653} and \cite{MR731340}) if the domain is smoothly bounded. Also, Aleksandrov \cite{MR745695} has constructed such inner functions if the domain has a $C^2$ boundary. \\

\noindent\textbf{Problem 8.23} In Pompeiu's formula, $f(z)$ and $\int_\Gamma f(\zeta)/(\zeta-z)\,d\zeta$ make sense for merely
continuous functions, but
\[\int\int_\Omega\frac{\partial f}{\partial\overline{\zeta}}\frac{1}{\zeta-z}\,d\zeta\,d\eta\]
needs at least some weak differentiability of $f$. It would be useful to extend the validity
of the formula, for example to cover the case of functions $f\in \text{Lip }\alpha$, $\alpha> 0$. This
prompts the following question. For which Borel sets $\Omega$ is the inequality
\[\Big|\left<\frac{\partial\chi_\Omega}{\partial\overline{z}},f\right>\Big|\leq C_\Omega\|f\|_{\text{Lip }\alpha}\]
valid? Putting it another way, for which Borel sets $\Omega$ does $\chi_\Omega$ act on the Besov space
$B^{\alpha-1}_{\infty,\infty}$? An interesting special case would be that when $\Omega$ is a `Swiss cheese'.\\
\textit{(A. G. O'Farrell)}\\

\noindent\textbf{Update 8.23} No progress on this problem has been reported to us.\\

\noindent\textbf{Problem 8.24} Let the sequence $\{M_k\}^\infty_0$ of positive numbers be such that
\[M_0=1\hspace{1cm}\text{ and }\hspace{1cm}\frac{M_{k+j}}{M_kM_j}\geq\begin{pmatrix}k+j\\j\end{pmatrix}.\]
Assume also the non-quasianalyticity condition that $\sum_k\big(M_k/M_{k+1}\big)<\infty$.\\
Consider those functions $f$ on a compact set $X \subset\mathbb{C}$ that are limits of rational
functions with poles of $X$ in the norm
\[g\mapsto\sum^\infty_{k=0}\frac{1}{M_k}\sup_{X}|g^{(k)}|;\]
these functions form a Banach algebra. Is $X$ its maximal ideal space? (The answer is
`yes' if $X$ is the unit disc or the unit interval.)\\ 
\textit{(A. G. O'Farrell)}\\

\noindent\textbf{Update 8.24} No progress on this problem has been reported to us.\\

\noindent\textbf{Problem 8.25} Let $\psi:S^1\to S^1$ be a direction-reversing homeomorphism, and let $A_\psi$ denote
the set of functions $f:S^1\to\mathbb{C}$ such that both $f$ and $f\circ\psi$ belong to the disc algebra.
When does $A_\psi$ contain only constant functions?\\
\textit{(A. G. O'Farrell)}\\

\noindent\textbf{Update 8.25} No progress on this problem has been reported to us.\\

\noindent\textbf{Problem 8.26} Let $\psi:S^1\to S^1$ be a homeomorphism. When is it true that $\text{Re }A = \text{Re }A \circ \psi$?
That is, when is each function $f$ in the real part of the disc algebra also of the form
$g\circ\psi$, for some function $g$ in the disc algebra?\\
O'Connell \cite{o'connell1973} has shown that it is necessary that $\psi$ is absolutely continuous, and
sufficient that $\psi$ is $C^{1+\varepsilon}$.\\
\textit{(A. G. O'Farrell)}\\

\noindent\textbf{Update 8.26} No progress on this problem has been reported to us.\\

\section*{New Problems}

\noindent To appear.

%======================= CHAPTER 9 ====================================

\chapter{Interpolation and Approximation}

\underline{Preface by J. L. Rovnyak (to appear)}\\

\noindent\textbf{Problem 9.1} A sequence $\{z_n\}^\infty_1$ in $|z|< 1$ is \textit{interpolating} for bounded analytic (harmonic)
functions if, for each bounded sequence $\{\alpha_n\}^\infty_1$ there exists a function $u(z)$ bounded and
analytic (harmonic) in $|z|< 1$ with $u(z_n) = \alpha_n (n\geq 1)$. It is known that a sequence is
interpolating for bounded harmonic functions if and only if it is interpolating for
bounded analytic functions (\cite{MR0141789}, \cite{MR628971}, \cite{varopoulos1977}); however all such proofs
require knowledge of conditions implying that a sequence interpolates for bounded
analytic functions. Find a simple proof of this equivalence which does not rely on
knowledge of such conditions.\\
\textit{(L. Zalcman)}\\

\noindent\textbf{Update 9.1} No progress on this problem has been reported to us.\\

\noindent\textbf{Problem 9.2} Let $f(z)\in H^\infty$, and let $\{z_n\}$ be a Blaschke sequence
\[\sum^\infty_{n=1}(\-|z_n|)<\infty.\]
\begin{enumerate}[(a)]
\item Does there always exist a Blaschke product $B(z)$ of norm not necessarily equal
to $1$, such that 
\begin{equation}
\label{Kstar9.2}
B(z_n)=f(z_n)
\end{equation}
for $n\geq1$? This is certainly the case if the Blaschke
sequence is uniformly separated, that is
\[\inf_n\prod_{m\neq n}\big|\frac{z_m-z_n}{1-\overline{z_m}z_n}\big|>0.\]
See, for example, Earl \cite{MR0419773}.
\item Is the unique $H^\infty$ function of minimal norm assuming the values $\{f(z_n)\}$ at $\{z_n\}$ a
constant multiple of a Blaschke product? The answer is `yes' if the sequence $\{z_n\}$ is
finite, see Earl \cite{MR0284588}. What happens if we know, in addition, that $f(z_n)\to0$ as $n\to\infty$? The
answer is again `yes' if $z_n\to1$ non-tangentially as $n\to\infty$.
\item One can ask questions $(a)$ and $(b)$ for sequences $\{z_n\}$ that are weakly separated,
viz
\[\inf_{m\neq n}\big|\frac{z_m-z_n}{1-\overline{z_m}z_n}\big|>0.\]
\item Also, one can ask questions $(a)$ and $(c)$ for inner functions instead of for
Blaschke products.
\end{enumerate}
\textit{(J. P. Earl and A. Stray)}\\

\noindent\textbf{Update 9.2} Stray \cite{MR940287} has proved that if there exist two distinct solutions $B(z)$ of (\ref{Kstar9.2}) which have minimum supremum norm, then one of them is a Blaschke product. Stray \cite{MR825925} has also shown that if $\big\{z_n\big\}$ is an interpolation sequence and $f(z_n)\to0$, then the solution of (\ref{Kstar9.2}) with minimum norm is unique and is a Blaschke product.\\

\noindent\textbf{Problem 9.3} \begin{enumerate} [(a)] 
\item Suppose that $f, f_1, f_2,\ldots,f_n\in H^\infty$, and
\[|f|\leq|f_1|+|f_2|+\ldots+|f_n|.\] 
Do there necessarily exist $h_1,h_2,\ldots,h_n\in H^\infty$ such that
\[f=f_1h_1+f_2h_2+\ldots+f_nh_n?\]
(If $|f_1|+|f_2|+\ldots+|f_n|\geq\delta>0$, this is the corona theorem.)
\item Suppose that $f_1,f_2\in H^\infty$. Do there necessarily exist $f\in H^\infty$ and $\delta>0$ such that
\[\delta(|f_1|+|f_2|)\leq|f|\leq|f_1|+|f_2|?\]
If the answer is `yes' is $f$ necessarily or possibly of the form $h_1f_1+h_2f_2$ for some
$h_1,h_2\in H^\infty$? (Notice that $(b)$ would imply $(a)$).
\end{enumerate}
\textit{(J. P. Earl)}\\

\noindent\textbf{Update 9.3(a)} The problem was originally posed by Rubel in Birtel's collections (see \cite[Problem 12, p. 347] {birtel}). A counter-example ($N=2$, $f_1$ and $f_2$ Blaschke products) was given by Rao \cite{MR0210910}. An easier counter-example is the following: Let $F_1(z)=1-z$, $F_2(z)=\exp[-(1+z)/(1-z)]$ and take $f_1=F_1^2$, $f_2=F_2^2$, $f=F_1F_2$. Then $|f|\leq|f_1|+|f_2|$. Since $F_1$ and $F_2$ are relatively prime, the existence of $h_1$, $h_2$ in $H^\infty$ with $f=f_1h_1+f_2h_2$ is equivalent to the existence of $g_1$, $g_2$ in $H^\infty$ with $1=F_1g_1+F_2g_2$. But this is impossible, since the right-hand side of the equation tends to zero on the real axis.\\

\noindent\textbf{Update 9.3(b)} The following is a counter-example: take $N=2$, $f_1=F_1$ and $f_2=F_2$ where $F_1$ and $F_2$ are defined as in Update 9.3(a). Assume that there exist $f$ in $H^\infty$ and $\delta>0$ such that $\delta(|f_1|+|f_2|)\leq|f|\leq |f_1+|f_2|$ holds. The left-hand side of the inequality implies $f\neq0$, $\delta|f_1|\leq|f|$ and $\delta|f_2|\leq|f|$, that is, $f$ is a divisor of $f_1$ and of $f_2$ (in the algebra $H^\infty$). Since $f_1$ is an outer function, so if $f$, and since $f_2$ is an inner function, so is $f$ (up to invertible elements). That is, $f$ has the factorisation $f=gh$ where $g$ is invertible in $H^\infty$ and $h$ is both outer and inner. Therefore $h$ is a non-zero constant function. This implies $|f|\geq\varepsilon>0$ for some positive $\varepsilon$, but this is a contradiction to $|f|\leq|f_1|+|f_2|$, since the right-hand side of the inequality tends to zero on the real axis.\\

\noindent\textbf{Problem 9.4} For each pair of $f,g\in H^\infty$, does there necessarily exist another pair of functions
$a,b\in H^\infty$ such that
\begin{equation}
\label{C9.1}
af+gb\neq0,\hspace{1cm}|z|<1\,?
\end{equation}
It is easy to see that a necessary condition for this is that $\log(|f|-|g|)$ have a
harmonic minorant. Is this condition also sufficient for (\ref{C9.1})? This problem is closely
related to Problem 9.3.\\
\textit{(B. A. Taylor; communicated by L. A. Rubel)}\\

\noindent\textbf{Update 9.4} This is a special case of Problem 8.3, and a positive solution has been given by Dahlberg, Kelleher and Taylor (not published).\\
  
\noindent\textbf{Problem 9.5} Let $K_1, K_2, K_3$ be disjoint closed sets in the extended complex plane, and
$C_1, C_2, C_3$ constants. Let $\rho_n(f)$ be the best rational approximation to the function
$f$ which equals $C_i$ on $K_i$ $(i = 1, 2, 3)$; i.e.
\[\rho_n(f)=\inf_{g\in R_n}\max_{z\in\cup_i K_i}|f(z)-g(z)|,\]
where $R_n$ is the class of rational functions $f$ order at most $n$. Find a geometric characterisation
of $\lim_{n\to\infty}\rho_n^{1/n}$. For the case of two sets, see Gonchar \cite{MR0231997}.\\
\textit{(T. Ganelius)}\\

\noindent\textbf{Update 9.5} No progress on this problem has been reported to us.\\ 

\noindent\textbf{Problem 9.6} Let $D$ be an open subset of the extended complex plane with non-empty
boundary $\partial D$, and let $F$ be a relatively-closed subset of $D$. Let $f$ be a function given
on $F$, and $\{f_n\}^\infty_1$ a sequence of functions analytic on $D$ such that $f_n\to f$ uniformly
on $F$. If $E\subset(\partial F_n\cap\partial D)$ and if $f$ extends continuously to $F\cup E$, can each $f_n$ be
extended continuously to $F\cup E$? The answer is `yes' if $D$ is the unit disc, or if $E$ is
compact.\\
\textit{(A. Stray)}\\ 

\noindent\textbf{Update 9.6} No progress on this problem has been reported to us.\\ 

\noindent\textbf{Problem 9.7} Let us call a closed set $E$ in $\mathbb{C}$ a \textit{weak Arakelian set} if, corresponding to each
function $g(z)$ continuous on $E$ and analytic in the interior of $E$, there exists an entire
function $g(z)$ such that, for any sequence $\{z_n\}^\infty_1$ in $E$, $|f(z_n)|\to\infty$ if and only if
$|g(z_n)|\to\infty$. Find a geometric characterisation of the weak Arakelian sets.\\
\textit{(L. A. Rubel)}\\

\noindent\textbf{Update 9.7} Goldstein \cite{goldstein1983} and Gauthier, Hengartner and Stray \cite{gauthier1989} have obtained necessary conditions and sufficient conditions for a closed subset of $C$ to be a weak Arakelian set, but have not obtained  conditions which are both necessary and sufficient.\\

\noindent\textbf{Problem 9.8} Let $\gamma$ be a Jordan arc in $\mathbb{C}^n$, $n\geq2$ such that the projections $\gamma_j$ on the complex coordinate planes $j=1,\ldots,n$ have area zero. Then $R(\gamma)=C(\gamma)$. Is it true that $P(\gamma)=C(\gamma)$? See Korevaar \cite{MR623473} and Wermer \cite{MR0394218}.\\
\textit{(J. Korevaar)}\\

\noindent\textbf{Update 9.8} No progress on this problem has been reported to us.\\

\noindent\textbf{Problem 9.9} Does the condition $\sum 1/p_n<\infty$ for positive integers $p_n$ guarantee that the sequences of powers $\{z^{p_n}\}$ fails to span $C(\gamma)$ for every Jordan arc $\gamma$? Korevaar and Dixon \cite{MR515123} have shown that for arcs of locally limited rotation (for example $C^1$ arcs), the condition
\[p_n\geq nL(n),\hspace{1cm}0<L(n)\uparrow,\hspace{1cm}\sum 1/nL(n)<\infty\]
assues a non-spanning sequence $\{z^{p_n}\}$.\\
\textit{(J. Korevaar)}\\

\noindent\textbf{Update 9.9} No progress on this problem has been reported to us.\\

\noindent\textbf{Problem 9.10} Let $F$ be a closed subest of $\mathbb{R}^n$, $n\geq2$. Call $F$ a set of harmonic approximation if every function continuous on $F$ and harmonic in the interior of $F$ can be uniformly approximated there by a harmonic function in $\mathbb{R}^n$. Give necessary and sufficient conditions that $F$ be a set of harmonic approximation. If $F$ is nowhere dense, \v{S}aginjan \cite{MR0293298} has done this. If $F$ is the closure of its interior, Gauthier, Ow and Goldstein \cite{MR576870} have given necessary conditions and sufficient conditions when $n=2$, but not necessary and sufficient conditions. This has applications to Rubel's Problem 9.7.\\
\textit{(M. Goldstein)}\\

\noindent\textbf{Update 9.10} No progress on this problem has been reported to us.\\

\noindent\textbf{Problem 9.11} Let $D$ be a planar domain. A sequence $\{z_j\}$ of points in $D$ is said to be an
\textit{interpolating sequence} if whenever $\{\alpha_j\}\in\ell^\infty$ there is a function $F\in H^\infty(D)$ such that
$F(z_j) = \alpha_j$, for all $j$.\\
Suppose that the sequence $\{z_j\}$ has the property that for each $j$ there is a function
$F_j \in H^\infty(D)$ such that $F_j(z_k)=0$ if $k\neq j$, $F_j(z_j)=1$ and $\|F\|_\infty\leq C$. Is $\{z_j\}$ necessarily
an interpolating sequence?\\
\textit{(P. W. Jones)}\\

\noindent\textbf{Update 9.11} No progress on this problem has been reported to us.\\

\noindent\textbf{Problem 9.12} Let $\Gamma\subset\mathbb{C}$ be a Jordan curve of logarithmic capacity $1$, and let $\phi$ be a
conformal map from the exterior of $\Gamma$ to the exterior of the unit circle such that
$\phi(\infty)=\infty$. We consider charge distributions on $\Gamma$ consisting of $n$ point charges $1/n$
at $n$th order Fekete points $z_1,\ldots,z_n$ on $\Gamma$, $n\in\mathbb{N}$. If $\Gamma$ is smooth enough, the
corresponding potentials
\[\frac{1}{n}\sum^n_{k=1}\log|z-z_k|\]
give approximations to $\log|\phi(z)|$ (outside $\Gamma$) and to $0$ (inside $\Gamma$) which are $O(1/n)$
away from $\Gamma$ (see \cite{MR0361030}, \cite{MR0361031}, \cite{MR705427}, \cite{MR0231990}).\\
Prove a similar result for the case where $\Gamma$ is a square. It does hold in the
degenerate case $\Gamma = [-2,2]$.\\ 
\textit{(J. Korevaar)}\\

\noindent\textbf{Update 9.12} No progress on this problem has been reported to us.\\

\noindent\textbf{Problem 9.13} Let $K$ be a compact subset of $\mathbb{R}^n$, $n\geq3$. For $\phi\in\mathcal{D}$, let $D(\phi)$ be a least-diameter disc containing $\text{spt }\phi$; let $d(\phi) = \text{diam }(\text{spt }\phi)$, and let
\[\|\phi\|_\ast=\|\phi\|_\infty+d(\phi)\cdot\|\triangledown\phi\|_\infty.\]
Are the following conditions equivalent for continuous functions $f:\mathbb{R}^n\to\mathbb{R}$?
\begin{enumerate}[(1)]
\item There exists a sequence $\{f_n\}^\infty_1$ of functions harmonic near $K$, such that
$f_n\to f$ uniformly on $K$.
\item There exists a function $n$ such that $n(\delta)$ decreases to $0$ as $\delta$ decreases to $0$, for
which
\[\Big|\int_{\mathbb{R}^n}f\delta\phi\,dx\Big|\leq\eta(d(\phi))\|\phi\|_\ast C(D(\phi)-X).\]

Here $C$ denotes the harmonic capacity of $\mathbb{R}^n$ obtained from the kernel $r^{-n+2}$.
\end{enumerate}
Note that $(1)$ implies $(2)$, and that $(2)$ implies $(1)$ if $f$ is a $C^2$-function. The
condition $(2)$ is formally analogous to one that occurs in rational approximation
theory.\\ 
\textit{(A. G. O'Farrell)}\\

\noindent\textbf{Update 9.13} No progress on this problem has been reported to us.\\

\noindent\textbf{Problem 9.14} Let $f$ be continuous on a compact subset $K$ of $\mathbb{C}$. If there exists a sequence
$\{f_n\}^\infty_1$ of functions analytic near $K$ for which $g_n\to f^2$ uniformly on $K$, does there
necessarily exist a sequence $\{h_n\}^\infty_1$ of functions analytic near $K$ for which $h_n\to f$ uniformly on $K$?\\
Paramanov has proved this under the stronger hypothesis that $f\in\text{Lip }(\frac{1}{2})$; it is also
true under the hypothesis that $f\in W^{1, p}$ $(p> 2)$.\\ 
\textit{(A. G. O'Farrell)}\\

\noindent\textbf{Update 9.14} No progress on this problem has been reported to us.\\

\noindent\textbf{Problem 9.15} Let $f_1$ and $f_2\in H^\infty(\text{unit disc}) = H^\infty$; and let the function $g\in H^\infty$ satisfy the inequality
\[|g(z)|\leq|f_1(z)|+|f_2(z)|.\]
Do there necessarily exist functions $g_1$ and $g_2$ in $H^\infty$ such that
\[g^2=f_1g_1+f_2g_2\,?\]
In other words, is it true that $g^2\in I(f_1, f_2)$ (the ideal generated by $f_1$ and $f_2$)?\\
Wolff has proved that $g^3\in I(f_1,f_2)$. Also Rao has given an example of
a function $g\notin I(f_1,f_2)$. For related results by Tolokonnikov, see \cite[p. 399]{MR734178}.\\ 
\textit{(J. Garnett)}\\

\noindent\textbf{Update 9.15} No progress on this problem has been reported to us.\\

\noindent\textbf{Problem 9.16} Let $\Gamma$ be a curve of the form
\[\{x+iA(x):-\infty<x<\infty\}\]
with
\[|A(x_1)-A(x_2)|\leq M|x_1-x_2|.\]
Let $E$ be a compact subset of $\Gamma$, $\Delta_1(t)> 0$, and let
\[\Omega=\mathbb{C}^*\setminus E,\hspace{1cm}\text{ where }\mathbb{C}^*=\mathbb{C}\cup\{\infty\}.\]
Prove the corona theorem for $\Omega$.\\
\textit{(J. Garnett)}\\

\noindent\textbf{Update 9.16} No progress on this problem has been reported to us.\\

\noindent\textbf{Problem 9.17} Let $K$ denote the $\frac{1}{3}$-Cantor set on $\mathbb{R}$; let $E = K\times K$, and let $\Omega=\mathbb{C}^*\setminus E$. Prove the corona theorem for $\Omega$.\\
\textit{(J. Garnett)}\\

\noindent\textbf{Update 9.17} No progress on this problem has been reported to us.\\

\section*{New Problems}

\noindent To appear.

\bibliography{MFO}

\begin{thebibliography}{100}

\bibitem{MR1799672}
M.~L. Agranovsky and Yakov Krasnov.
\newblock Quadratic divisors of harmonic polynomials in $\mathbb{R}^n$.
\newblock {\em J. Anal. Math.}, 82:379--395, 2000.

\bibitem{MR0412400}
D.~Aharonov and H.~S. Shapiro.
\newblock A minimal area problem in conformal mapping. ({A}bstract).
\newblock pages 1--5. London Math. Soc. Lecture Note Ser., No. 12, 1974.

\bibitem{AharonovShapiro}
D.~Aharonov and H.~S. Shapiro.
\newblock A minimal-area problem in conformal mapping - {P}reliminary {R}eport:
  {P}art {II}.
\newblock 1978.

\bibitem{MR0365150}
D.~Aharonov, H.~S. Shapiro, and A.~L. Shields.
\newblock Weakly invertible elements in the space of square-summable
  holomorphic functions.
\newblock {\em J. London Math. Soc. (2)}, 9:183--192, 1974/75.

\bibitem{zbMATH01369959}
D.~{Aharonov}, H.~S. {Shapiro}, and A.~Y. {Solynin}.
\newblock {A minimal area problem in conformal mapping.}
\newblock {\em {J. Anal. Math.}}, 78:157--176, 1999.

\bibitem{zbMATH01640257}
D.~{Aharonov}, H.~S. {Shapiro}, and A.~Y. {Solynin}.
\newblock {A minimal area problem in conformal mapping. II.}
\newblock {\em {J. Anal. Math.}}, 83:259--288, 2001.

\bibitem{Ahlfors1930}
L.~V. Ahlfors.
\newblock Untersuchungen zur {T}heorie der konformen {A}bbildung und der ganzen
  {F}unktionen.
\newblock {\em Acta Soc. Sci. Fenn.}, Nova Ser. I, no. 9, 1930.

\bibitem{MR1509465}
L.~V. Ahlfors.
\newblock \"{U}ber die {K}reise die von einer {R}iemannschen {F}l\"ache
  schlicht \"uberdeckt werden.
\newblock {\em Comment. Math. Helv.}, 5(1):28--38, 1933.

\bibitem{MR1501949}
L.~V. Ahlfors.
\newblock An extension of {S}chwarz's lemma.
\newblock {\em Trans. Amer. Math. Soc.}, 43(3):359--364, 1938.

\bibitem{MR863437}
L.~V. Ahlfors.
\newblock Clifford numbers and {M}\"obius transformations in {$\mathbb{R}^n$}.
\newblock In {\em Clifford algebras and their applications in mathematical
  physics ({C}anterbury, 1985)}, volume 183 of {\em NATO Adv. Sci. Inst. Ser. C
  Math. Phys. Sci.}, pages 167--175. Reidel, Dordrecht, 1986.

\bibitem{MR1545698}
L.~V. Ahlfors and H.~Grunsky.
\newblock \"{U}ber die {B}lochsche {K}onstante.
\newblock {\em Math. Z.}, 42(1):671--673, 1937.

\bibitem{MR990410}
H.~Aikawa.
\newblock Harmonic functions having no tangential limits.
\newblock {\em Proc. Amer. Math. Soc.}, 108(2):457--464, 1990.

\bibitem{MR1169019}
H.~Aikawa.
\newblock Integrability of superharmonic functions and subharmonic functions.
\newblock {\em Proc. Amer. Math. Soc.}, 120(1):109--117, 1994.

\bibitem{MR0199388}
W.~Al-Katifi.
\newblock On the asymptotic values and paths of certain integral and
  meromorphic functions.
\newblock {\em Proc. London Math. Soc. (3)}, 16:599--634, 1966.

\bibitem{MR1606325}
J.~M. Aldaz.
\newblock Remarks on the {H}ardy-{L}ittlewood maximal function.
\newblock {\em Proc. Roy. Soc. Edinburgh Sect. A}, 128(1):1--9, 1998.

\bibitem{MR1745465}
J.~M. Aldaz.
\newblock A remark on the centered {$n$}-dimensional {H}ardy-{L}ittlewood
  maximal function.
\newblock {\em Czechoslovak Math. J.}, 50(125)(1):103--112, 2000.

\bibitem{MR745695}
A.~B. Aleksandrov.
\newblock Inner functions on compact spaces.
\newblock {\em Funktsional. Anal. i Prilozhen.}, 18(2):1--13, 1984.

\bibitem{MR1503516}
J.~W. Alexander.
\newblock Functions which map the interior of the unit circle upon simple
  regions.
\newblock {\em Ann. of Math. (2)}, 17(1):12--22, 1915.

\bibitem{MR915430}
G.~R. Allan, A.~G. O'Farrell, and T.~J. Ransford.
\newblock A {T}auberian theorem arising in operator theory.
\newblock {\em Bull. London Math. Soc.}, 19(6):537--545, 1987.

\bibitem{MR513885}
A.~Ancona.
\newblock Principe de {H}arnack \`a la fronti\`ere et th\'eor\`eme de {F}atou
  pour un op\'erateur elliptique dans un domaine {L}ipschitzien.
\newblock {\em Ann. Inst. Fourier (Grenoble)}, 28(4):169--213, x, 1978.

\bibitem{MR0287295}
J.~M. Anderson.
\newblock Category theorems for certain {B}anach spaces of analytic functions.
\newblock {\em J. Reine Angew. Math.}, 249:83--91, 1971.

\bibitem{MR536388}
J.~M. Anderson.
\newblock Asymptotic values of meromorphic functions of smooth growth.
\newblock {\em Glasgow Math. J.}, 20(2):155--162, 1979.

\bibitem{MR569247}
J.~M. Anderson.
\newblock The extremum problem for analytic functions with finite area
  integral.
\newblock {\em Comment. Math. Helv.}, 55(1):87--96, 1980.

\bibitem{MR0199397}
J.~M. Anderson and J.~Clunie.
\newblock Slowly growing meromorphic functions.
\newblock {\em Comment. Math. Helv.}, 40:267--280, 1966.

\bibitem{MR0361090}
J.~M. Anderson, J.~Clunie, and Ch. Pommerenke.
\newblock On {B}loch functions and normal functions.
\newblock {\em J. Reine Angew. Math.}, 270:12--37, 1974.

\bibitem{MR979966}
J.~M. Anderson, J.~L. Fern\'andez, and A.~L. Shields.
\newblock Inner functions and cyclic vectors in the {B}loch space.
\newblock {\em Trans. Amer. Math. Soc.}, 323(1):429--448, 1991.

\bibitem{MR0206286}
N.~U. {A}rakeljan.
\newblock Entire functions of finite order with an infinite set of deficient
  values.
\newblock {\em Dokl. Akad. Nauk SSSR}, 170:999--1002, 1966.

\bibitem{Armitage}
D.~H. Armitage.
\newblock On the global integrability of superharmonic functions in balls.
\newblock {\em J. London Math. Soc. (2)}, 4:365--373, 1971.

\bibitem{MR540060}
D.~H. Armitage.
\newblock Spherical extrema of harmonic polynomials.
\newblock {\em J. London Math. Soc. (2)}, 19(3):451--456, 1979.

\bibitem{MR1277091}
D.~H. Armitage.
\newblock On solutions of elliptic equations that decay rapidly on paths.
\newblock {\em Proc. Amer. Math. Soc.}, 123(8):2421--2422, 1995.

\bibitem{MR813738}
D.~H. Armitage, T.~Bagby, and P.~M. Gauthier.
\newblock Note on the decay of solutions of elliptic equations.
\newblock {\em Bull. London Math. Soc.}, 17(6):554--556, 1985.

\bibitem{MR1274778}
D.~H. Armitage and M.~Goldstein.
\newblock Radial limiting behaviour of harmonic functions in cones.
\newblock {\em Complex Variables Theory Appl.}, 22(3-4):267--276, 1993.

\bibitem{zbMATH00714744}
K.~{Astala}.
\newblock {Area distortion of quasiconformal mappings.}
\newblock {\em {Acta Math.}}, 173(1):37--60, 1994.

\bibitem{MR861687}
K.~Astala and F.~W. Gehring.
\newblock Injectivity, the {BMO} norm and the universal {T}eichm\"uller space.
\newblock {\em J. Analyse Math.}, 46:16--57, 1986.

\bibitem{MR0126418}
F.~V. Atkinson.
\newblock On sums of powers of complex numbers.
\newblock {\em Acta Math. Acad. Sci. Hungar.}, 12:185--188, 1961.

\bibitem{MR0294640}
A.~Baernstein, II.
\newblock Proof of {E}drei's spread conjecture.
\newblock {\em Bull. Amer. Math. Soc.}, 78:277--278, 1972.

\bibitem{MR0374429}
A.~Baernstein, II.
\newblock Proof of {E}drei's spread conjecture.
\newblock {\em Proc. London Math. Soc. (3)}, 26:418--434, 1973.

\bibitem{MR0417406}
A.~Baernstein, II.
\newblock Integral means, univalent functions and circular symmetrization.
\newblock {\em Acta Math.}, 133:139--169, 1974.

\bibitem{MR0387590}
A.~Baernstein, II.
\newblock Some extremal problems for univalent functions, harmonic measures,
  and subharmonic functions.
\newblock pages 11--15. London Math. Soc. Lecture Note Ser., No. 12, 1974.

\bibitem{baernstein1976}
A.~Baernstein, II.
\newblock Univalence and bounded mean oscillation.
\newblock {\em Michigan Math. J.}, 23(3):217--223, 1976.

\bibitem{MR623463}
A.~Baernstein, II.
\newblock Analytic functions of bounded mean oscillation.
\newblock In {\em Aspects of contemporary complex analysis ({P}roc. {NATO}
  {A}dv. {S}tudy {I}nst., {U}niv. {D}urham, {D}urham, 1979)}, pages 3--36.
  Academic Press, London-New York, 1980.

\bibitem{MR846491}
A.~Baernstein, II.
\newblock Coefficients of univalent functions with restricted maximum modulus.
\newblock {\em Complex Variables Theory Appl.}, 5(2-4):225--236, 1986.

\bibitem{MR0437765}
A.~Baernstein, II and R.~Rochberg.
\newblock Means and coefficients of functions which omit a sequence of values.
\newblock {\em Math. Proc. Cambridge Philos. Soc.}, 81(1):47--57, 1977.

\bibitem{MR0114007}
I.~N. Baker.
\newblock Some entire functions with fixpoints of every order.
\newblock {\em J. Austral. Math. Soc.}, 1:203--209, 1959/1961.

\bibitem{Baker1964ae}
I.~N. Baker.
\newblock Entire functions with linearly distributed values.
\newblock {\em Mathematische Zeitschrift}, 86(4):263--267, Aug 1964.

\bibitem{MR0169989}
I.~N. Baker.
\newblock Fixpoints of polynomials and rational functions.
\newblock {\em J. London Math. Soc.}, 39:615--622, 1964.

\bibitem{MR0197725}
I.~N. Baker.
\newblock The distribution of fixpoints of entire functions.
\newblock {\em Proc. London Math. Soc. (3)}, 16:493--506, 1966.

\bibitem{MR0226009}
I.~N. Baker.
\newblock Repulsive fixpoints of entire functions.
\newblock {\em Math. Z.}, 104:252--256, 1968.

\bibitem{MR0264071}
I.~N. Baker.
\newblock Limit functions and sets of non-normality in iteration theory.
\newblock {\em Ann. Acad. Sci. Fenn. Ser. A I No.}, 467:11, 1970.

\bibitem{MR0419759}
I.~N. Baker.
\newblock An entire function which has wandering domains.
\newblock {\em J. Austral. Math. Soc. Ser. A}, 22(2):173--176, 1976.

\bibitem{MR865514}
A.~K. Bakhtin.
\newblock Extrema of linear functionals.
\newblock {\em Akad. Nauk Ukrain. SSR Inst. Mat. Preprint}, (25):8, 1986.

\bibitem{MR624232}
S.~B. Bank and R.~P. Kaufman.
\newblock On {B}riot-{B}ouquet differential equations and a question of {E}inar
  {H}ille.
\newblock {\em Math. Z.}, 177(4):549--559, 1981.

\bibitem{MR664047}
S.~B. Bank and I.~Laine.
\newblock On the oscillation theory of {$f^{\prime\prime}+Af=0$}\ where {$A$}\
  is entire.
\newblock {\em Trans. Amer. Math. Soc.}, 273(1):351--363, 1982.

\bibitem{MR0379816}
R.~Barnard and J.~L. Lewis.
\newblock Subordination theorems for some classes of starlike functions.
\newblock {\em Pacific J. Math.}, 56(2):333--366, 1975.

\bibitem{Barnardproblem}
R.~W. Barnard.
\newblock On the coefficient bounds on $f$ in ${S}$ when $f'$ is univalent.
\newblock {\em Notices Amer. Math. Soc.}, 26:A--69, 1979.

\bibitem{MR555840}
R.~W. Barnard and C.~Kellogg.
\newblock Applications of convolution operators to problems in univalent
  function theory.
\newblock {\em Michigan Math. J.}, 27(1):81--94, 1980.

\bibitem{MR3002600}
R.~W. Barnard, K.~Pearce, and A.~Y. Solynin.
\newblock Iceberg-type problems: estimating hidden parts of a continuum from
  the visible parts.
\newblock {\em Math. Nachr.}, 285(17-18):2042--2058, 2012.

\bibitem{MR694924}
R.~W. Barnard and T.~J. Suffridge.
\newblock On the simultaneous univalence of {$f$}\ and {$f^{\prime} $}.
\newblock {\em Michigan Math. J.}, 30(1):9--16, 1983.

\bibitem{PLMS:PLMS0445}
P.~D. Barry.
\newblock The minimum modulus of small integral and subharmonic functions.
\newblock {\em Proceedings of the London Mathematical Society},
  s3-12(1):445--495, 1962.

\bibitem{MR0156993}
P.~D. Barry.
\newblock On a theorem of {B}esicovitch.
\newblock {\em Quart. J. Math. Oxford Ser. (2)}, 14:293--302, 1963.

\bibitem{Barsegyan}
G.~A. Barsegyan.
\newblock On the relation between the behaviour of asymptotic values and
  {$a$}-points of meromorphic functions ({R}ussian).
\newblock {\em Akad. Nauk Armyan. SSR Dokl.}, 18(2):124--133, 1983.

\bibitem{MR507611}
K.~F. Barth, D.~A. Brannan, and W.~K. Hayman.
\newblock The growth of plane harmonic functions along an asymptotic path.
\newblock {\em Proc. London Math. Soc. (3)}, 37(2):363--384, 1978.

\bibitem{C}
K.~F. Barth, D.~A. Brannan, and W.~K. Hayman.
\newblock Research problems in complex analysis.
\newblock {\em Bull. London Math. Soc.}, 16(5):490--517, 1984.

\bibitem{MR660605}
K.~F. Barth and J.~G. Clunie.
\newblock A bounded analytic function in the unit disk with a level set
  component of infinite length.
\newblock {\em Proc. Amer. Math. Soc.}, 85(4):562--566, 1982.

\bibitem{MR0269834}
K.~F. Barth and W.~J. Schneider.
\newblock Entire functions mapping countable dense subsets of the reals onto
  each other monotonically.
\newblock {\em J. London Math. Soc. (2)}, 2:620--626, 1970.

\bibitem{MR0301195}
K.~F. Barth and W.~J. Schneider.
\newblock Entire functions mapping arbitrary countable dense sets and their
  complements onto each other.
\newblock {\em J. London Math. Soc. (2)}, 4:482--488, 1971/72.

\bibitem{MR0293089}
K.~F. Barth and W.~J. Schneider.
\newblock On a problem of {E}rd\"os concerning the zeros of the derivatives of
  an entire function.
\newblock {\em Proc. Amer. Math. Soc.}, 32:229--232, 1972.

\bibitem{MR3261731}
A.~Barton and L.~A. Ward.
\newblock A new class of harmonic measure distribution functions.
\newblock {\em J. Geom. Anal.}, 24(4):2035--2071, 2014.

\bibitem{MR0039804}
I.~E. Bazilevi\v{c}.
\newblock On distortion theorems and coefficients of univalent functions.
\newblock {\em Mat. Sbornik N.S.}, 28(70):147--164, 1951.

\bibitem{zbMATH01542658}
A.~F. {Beardon}.
\newblock {\em {Iteration of rational functions. Complex analytic dynamical
  systems. Paperback ed.}}, volume 132.
\newblock New York, NY: Springer, paperback ed. edition, 2000.

\bibitem{MR0299780}
J.~Becker.
\newblock L\"ownersche {D}ifferentialgleichung und quasikonform fortsetzbare
  schlichte {F}unktionen.
\newblock {\em J. Reine Angew. Math.}, 255:23--43, 1972.

\bibitem{MR623464}
J.~Becker.
\newblock Conformal mappings with quasiconformal extensions.
\newblock In {\em Aspects of contemporary complex analysis ({P}roc. {NATO}
  {A}dv. {S}tudy {I}nst., {U}niv. {D}urham, {D}urham, 1979)}, pages 37--77.
  Academic Press, London-New York, 1980.

\bibitem{MR584077}
J.~Becker, C.~W. Henson, and L.~A. {R}ubel.
\newblock First-order conformal invariants.
\newblock {\em Ann. of Math. (2)}, 112(1):123--178, 1980.

\bibitem{zbMATH03857435}
J.~{Becker} and Ch. {Pommerenke}.
\newblock {Schlichtheitskriterien und Jordangebiete.}
\newblock {\em {J. Reine Angew. Math.}}, 354:74--94, 1984.

\bibitem{MR951982}
J.~Becker and Ch. Pommerenke.
\newblock On the {H}ausdorff dimension of quasicircles.
\newblock {\em Ann. Acad. Sci. Fenn. Ser. A I Math.}, 12(2):329--333, 1987.

\bibitem{MR0256166}
M.~Behrens.
\newblock The corona conjecture for a class of infinitely connected domains.
\newblock {\em Bull. Amer. Math. Soc.}, 76:387--391, 1970.

\bibitem{MR2425327}
D.~Beliaev.
\newblock Integral means spectrum of random conformal snowflakes.
\newblock {\em Nonlinearity}, 21(7):1435--1442, 2008.

\bibitem{MR2680427}
D.~Beliaev and S.~Smirnov.
\newblock Random conformal snowflakes.
\newblock {\em Ann. of Math. (2)}, 172(1):597--615, 2010.

\bibitem{MR0355015}
E.~Beller and D.~J. Newman.
\newblock The minimum modulus of polynomials.
\newblock {\em Proc. Amer. Math. Soc.}, 45:463--465, 1974.

\bibitem{MR608327}
M.~Benedicks.
\newblock Positive harmonic functions vanishing on the boundary of certain
  domains in {$\mathbb{R}^{n}$}.
\newblock {\em Ark. Mat.}, 18(1):53--72, 1980.

\bibitem{MR1123803}
W.~Bergweiler.
\newblock Periodic points of entire functions: proof of a conjecture of
  {B}aker.
\newblock {\em Complex Variables Theory Appl.}, 17(1-2):57--72, 1991.

\bibitem{zbMATH00487133}
W.~{Bergweiler}.
\newblock {Iteration of meromorphic functions.}
\newblock {\em {Bull. Am. Math. Soc., New Ser.}}, 29(2):151--188, 1993.

\bibitem{MR1464835}
W.~Bergweiler.
\newblock Non-real periodic points of entire functions.
\newblock {\em Canad. Math. Bull.}, 40(3):271--275, 1997.

\bibitem{MR1325262}
W.~Bergweiler, J.~Clunie, and J.~Langley.
\newblock Proof of a conjecture of {B}aker concerning the distribution of
  fixpoints.
\newblock {\em Bull. London Math. Soc.}, 27(2):148--154, 1995.

\bibitem{MR1344897}
W.~Bergweiler and A.~Eremenko.
\newblock On the singularities of the inverse to a meromorphic function of
  finite order.
\newblock {\em Rev. Mat. Iberoamericana}, 11(2):355--373, 1995.

\bibitem{Bergweiler2014OnTB}
W.~Bergweiler and A.~Eremenko.
\newblock On the {B}ank-{L}aine conjecture.
\newblock 2014.

\bibitem{2015arXiv151005731B}
W.~{Bergweiler} and A.~{Eremenko}.
\newblock {Quasiconformal surgery and linear differential equations}.
\newblock {\em J. Analyse Math. To appear}, October 2015.

\bibitem{bergweiler_eremenko_hinkkanen_2018}
W.~Bergweiler, A.~Eremenko, and A.~Hinkkanen.
\newblock Entire functions with two radially distributed values.
\newblock {\em Mathematical Proceedings of the Cambridge Philosophical
  Society}, 165(1):93–108, 2018.

\bibitem{MR2024413}
W.~Bergweiler, A.~Eremenko, and J.~K. Langley.
\newblock Real entire functions of infinite order and a conjecture of {W}iman.
\newblock {\em Geom. Funct. Anal.}, 13(5):975--991, 2003.

\bibitem{zbMATH02190101}
W.~{Bergweiler}, A.~{Eremenko}, and J.~K. {Langley}.
\newblock {Zeros of differential polynomials in real meromorphic functions.}
\newblock {\em {Proc. Edinb. Math. Soc., II. Ser.}}, 48(2):279--293, 2005.

\bibitem{MR1782673}
F.~Berteloot and J.~Duval.
\newblock Une d\'emonstration directe de la densit\'e des cycles r\'epulsifs
  dans l'ensemble de {J}ulia.
\newblock In {\em Complex analysis and geometry ({P}aris, 1997)}, volume 188 of
  {\em Progr. Math.}, pages 221--222. Birkh\"auser, Basel, 2000.

\bibitem{MR1512383}
A.~S. Besicovitch.
\newblock On integral functions of order {$<1$}.
\newblock {\em Math. Ann.}, 97(1):677--695, 1927.

\bibitem{MR1877323}
Dimitrios Betsakos.
\newblock Geometric theorems and problems for harmonic measure.
\newblock {\em Rocky Mountain J. Math.}, 31(3):773--795, 2001.

\bibitem{MR0001370}
A.~Beurling.
\newblock Ensembles exceptionnels.
\newblock {\em Acta Math.}, 72:1--13, 1940.

\bibitem{MR0027954}
A.~Beurling.
\newblock On two problems concerning linear transformations in {H}ilbert space.
\newblock {\em Acta Math.}, 81:17, 1948.

\bibitem{MR0029980}
A.~Beurling.
\newblock Some theorems on boundedness of analytic functions.
\newblock {\em Duke Math. J.}, 16:355--359, 1949.

\bibitem{bieberbach1916koeffizienten}
L.~Bieberbach.
\newblock {\em {\"{U}}ber die {K}oeffizienten derjenigen {P}otenzreihen, welche
  eine schlichte {A}bbildung des {E}inheitskreises vermitteln}.
\newblock Reimer in Komm., 1916.

\bibitem{biernacki}
M.~Biernacki.
\newblock Sur les \'equations alg\'ebriques contenant des param\`etres
  arbitraires.
\newblock {\em Bull. del. Acad. Polon. Sci.}, Cl. III:562, 1927.

\bibitem{birtel}
F.~T. Birtel, editor.
\newblock {\em Function algebras}.
\newblock Proceedings of the {I}nternational {S}ymposium on {F}unction
  {A}lgebras, Tulane University, Chicago. Scott, Foresman, 1966.

\bibitem{Bishop2018}
C.~J. Bishop.
\newblock A transcendental {J}ulia set of dimension 1.
\newblock {\em Inventiones mathematicae}, 212(2):407--460, May 2018.

\bibitem{MR0133462}
E.~Bishop.
\newblock A general {R}udin-{C}arleson theorem.
\newblock {\em Proc. Amer. Math. Soc.}, 13:140--143, 1962.

\bibitem{MR0059365}
R.~P. Boas, Jr.
\newblock Asymptotic properties of functions of exponential type.
\newblock {\em Duke Math. J.}, 20:433--448, 1953.

\bibitem{MR0162918}
E.~Bombieri.
\newblock Sopra un teorema di {H}. {B}ohr e {G}. {R}icci sulle funzioni
  maggioranti delle serie di potenze.
\newblock {\em Boll. Un. Mat. Ital. (3)}, 17:276--282, 1962.

\bibitem{MR0185108}
E.~Bombieri.
\newblock On functions which are regular and univalent in a half-plane.
\newblock {\em Proc. London Math. Soc. (3)}, 14a:47--50, 1965.

\bibitem{bombieri2}
E.~L. Bombieri.
\newblock {\em Sull'integrazione approssimatat dell'equazione differenziale di
  {L}oewner e le sue applicazioni alla theoria delle funzioni univalenti}.
\newblock Universita di Milano, 1963.

\bibitem{zbMATH04173340}
M.~{Bonk}.
\newblock {On Bloch's constant.}
\newblock {\em {Proc. Am. Math. Soc.}}, 110(4):889--894, 1990.

\bibitem{MR0180776}
D.~Borwein.
\newblock On a class of convergent series of positive terms.
\newblock {\em J. London Math. Soc.}, 40:587--588, 1965.

\bibitem{zbMATH00741404}
P.~{Borwein}.
\newblock {The arc length of the lemniscate $\{\vert p(z)\vert =1\}$.}
\newblock {\em {Proc. Am. Math. Soc.}}, 123(3):797--799, 1995.

\bibitem{MR0344465}
V.~S. Bo\u\i\v{c}uk and A.~A. Gol'dberg.
\newblock On the three lines theorem.
\newblock {\em Mat. Zametki}, 15:45--53, 1974.

\bibitem{MR880292}
J.~Bourgain.
\newblock The metrical interpretation of superreflexivity in {B}anach spaces.
\newblock {\em Israel J. Math.}, 56(2):222--230, 1986.

\bibitem{MR0071620}
B.~V. Boyarski\u\i.
\newblock Homeomorphic solutions of {B}eltrami systems.
\newblock {\em Dokl. Akad. Nauk SSSR (N.S.)}, 102:661--664, 1955.

\bibitem{MR0412411}
D.~A. Brannan.
\newblock On coefficient problems for certain power series.
\newblock pages 17--27. London Math. Soc. Lecture Note Ser., No. 12, 1974.

\bibitem{MR0062883}
M.~Brelot and G.~Choquet.
\newblock Espaces et lignes de {G}reen.
\newblock {\em Ann. Inst. Fourier Grenoble}, 3:199--263 (1952), 1951.

\bibitem{MR0069968}
M.~Brelot and G.~Choquet.
\newblock Polyn\^omes harmoniques et polyharmoniques.
\newblock In {\em Second colloque sur les \'equations aux d\'eriv\'ees
  partielles, {B}ruxelles, 1954}, pages 45--66. G. Thone, Li\`ege; Masson \&
  Cie, Paris, 1955.

\bibitem{MR509942}
J.~E. Brennan.
\newblock The integrability of the derivative in conformal mapping.
\newblock {\em J. London Math. Soc. (2)}, 18(2):261--272, 1978.

\bibitem{MR766220}
J.~E. Brown.
\newblock A method for investigating geometric properties of support points and
  applications.
\newblock {\em Trans. Amer. Math. Soc.}, 287(1):285--291, 1985.

\bibitem{MR891142}
J.~E. Brown.
\newblock Level sets for functions convex in one direction.
\newblock {\em Proc. Amer. Math. Soc.}, 100(3):442--446, 1987.

\bibitem{MR1081693}
J.~E. Brown.
\newblock On the {S}endov conjecture for sixth degree polynomials.
\newblock {\em Proc. Amer. Math. Soc.}, 113(4):939--946, 1991.

\bibitem{Brown1973}
L.~Brown, B.~Schreiber, and B.~A. Taylor.
\newblock Spectral synthesis and the {P}ompeiu problem.
\newblock {\em Annales de l'institut Fourier}, 23(3):125--154, 1973.

\bibitem{MR1057941}
F.~Br\"uggemann.
\newblock On solutions of linear differential equations with real zeros; proof
  of a conjecture of {H}ellerstein and {R}ossi.
\newblock {\em Proc. Amer. Math. Soc.}, 113(2):371--379, 1991.

\bibitem{MR1159367}
F.~Br\"uggemann.
\newblock Proof of a conjecture of {F}rank and {L}angley concerning zeros of
  meromorphic functions and linear differential polynomials.
\newblock {\em Analysis}, 12(1-2):5--30, 1992.

\bibitem{MR577757}
D.~Bshouty.
\newblock A note on {H}adamard products of univalent functions.
\newblock {\em Proc. Amer. Math. Soc.}, 80(2):271--272, 1980.

\bibitem{MR0029984}
R.~C. Buck.
\newblock Integral valued entire functions.
\newblock {\em Duke Math. J.}, 15:879--891, 1948.

\bibitem{MR0259076}
J.~D. Buckholtz.
\newblock Zeros of partial sums of power series. {II}.
\newblock {\em Michigan Math. J.}, 17:5--14, 1970.

\bibitem{zbMATH06093940}
X.~{Buff} and A.~{Ch\'eritat}.
\newblock {Quadratic {J}ulia sets with positive area.}
\newblock {\em {Ann. Math. (2)}}, 176(2):673--746, 2012.

\bibitem{MR0466568}
A.~P. Calder\'on.
\newblock Cauchy integrals on {L}ipschitz curves and related operators.
\newblock {\em Proc. Nat. Acad. Sci. U.S.A.}, 74(4):1324--1327, 1977.

\bibitem{Camera}
G.~Camera.
\newblock Doctoral thesis.
\newblock {\em University of London}, 1977.

\bibitem{MR0299769}
D.~M. Campbell.
\newblock Majorization-subordination theorems for locally univalent functions.
\newblock {\em Bull. Amer. Math. Soc.}, 78:535--538, 1972.

\bibitem{MR0315120}
D.~M. Campbell.
\newblock Majorization-subordination theorems for locally univalent functions.
  {II}.
\newblock {\em Canad. J. Math.}, 25:420--425, 1973.

\bibitem{MR0349987}
D.~M. Campbell.
\newblock Majorization-subordination theorems for locally univalent functions.
  {III}.
\newblock {\em Trans. Amer. Math. Soc.}, 198:297--306, 1974.

\bibitem{Campbell}
D.~M. Campbell.
\newblock The limiting behaviour of $zf''(z)/f;(z)$ and two conjectures on
  univalent functions.
\newblock {\em Notices Amer. Math. Soc.}, 22:A--120, 1975.

\bibitem{MR2284052}
J.~Cantarella, J.~H.~G. Fu, R.~Kusner, J.~M. Sullivan, and N.~C. Wrinkle.
\newblock Criticality for the {G}ehring link problem.
\newblock {\em Geom. Topol.}, 10:2055--2116, 2006.

\bibitem{10.2307/24903359}
A.~Cant\'on, D.~Drasin, and A.~Granados.
\newblock Asymptotic values of meromorphic functions of finite order.
\newblock {\em Indiana University Mathematics Journal}, 59(3):1057--1095, 2010.

\bibitem{MR0064861}
C.~Carath\'eodory.
\newblock {\em Theory of functions of a complex variable. {V}ol. 2}.
\newblock Chelsea Publishing Company, New York, 1954.
\newblock Translated by F. Steinhardt.

\bibitem{zbMATH03826463}
A.~{Carbery}.
\newblock {The boundedness of the maximal Bochner-Riesz operator on $L\sp
  4$($R\sp 2$).}
\newblock {\em {Duke Math. J.}}, 50:409--416, 1983.

\bibitem{zbMATH03834608}
A.~{Carbery}.
\newblock {A weighted inequality for the maximal Bochner-Riesz operator on
  ${\mathbb R}\sp 2$.}
\newblock {\em {Trans. Am. Math. Soc.}}, 287:673--680, 1985.

\bibitem{MR972135}
A.~Carbery, J.~L. Rubio~de Francia, and L.~Vega.
\newblock Almost everywhere summability of {F}ourier integrals.
\newblock {\em J. London Math. Soc. (2)}, 38(3):513--524, 1988.

\bibitem{MR1555232}
T.~Carleman.
\newblock Extension d'un th\'eor\`eme de {L}iouville.
\newblock {\em Acta Math.}, 48(3-4):363--366, 1926.

\bibitem{MR0050011}
L.~Carleson.
\newblock Sets of uniqueness for functions regular in the unit circle.
\newblock {\em Acta Math.}, 87:325--345, 1952.

\bibitem{MR0141789}
L.~Carleson.
\newblock Interpolations by bounded analytic functions and the corona problem.
\newblock {\em Ann. of Math. (2)}, 76:547--559, 1962.

\bibitem{MR0470233}
L.~Carleson.
\newblock Asymptotic paths for subharmonic functions in {$\mathbb{R}^{n}$}.
\newblock {\em Ann. Acad. Sci. Fenn. Ser. A I Math.}, 2:35--39, 1976.

\bibitem{MR730079}
L.~Carleson.
\newblock On {$H^{\infty }$}\ in multiply connected domains.
\newblock In {\em Conference on harmonic analysis in honor of {A}ntoni
  {Z}ygmund, {V}ol. {I}, {II} ({C}hicago, {I}ll., 1981)}, Wadsworth Math. Ser.,
  pages 349--372. Wadsworth, Belmont, CA, 1983.

\bibitem{MR1162188}
L.~Carleson and P.~W. Jones.
\newblock On coefficient problems for univalent functions and conformal
  dimension.
\newblock {\em Duke Math. J.}, 66(2):169--206, 1992.

\bibitem{MR0025568}
F.~Carlson.
\newblock Sur les fonctions enti\`eres.
\newblock {\em Ark. Mat. Astr. Fys.}, 35A(14):18, 1948.

\bibitem{MR555223}
F.~W. Carroll.
\newblock A strongly annular function with countably many singular values.
\newblock {\em Math. Scand.}, 44(2):330--334, 1979.

\bibitem{MR1509288}
H.~Cartan.
\newblock Sur les syst\`emes de fonctions holomorphes \`a vari\'et\'es
  lin\'eaires lacunaires et leurs applications.
\newblock {\em Ann. Sci. \'Ecole Norm. Sup. (3)}, 45:255--346, 1928.

\bibitem{zbMATH03010843}
H.~{Cartan}.
\newblock {Sur les z\'eros des combinaisons lin\'eaires de $p$ fonctions
  holomorphes donn\'ees.}
\newblock {\em {Mathematica, Cluj}}, 7:5--31, 1933.

\bibitem{MR1512983}
M.~L. Cartwright.
\newblock Some inequalities in the theory of functions.
\newblock {\em Math. Ann.}, 111(1):98--118, 1935.

\bibitem{MR0130379}
M.~L. Cartwright and E.~F. Collingwood.
\newblock The radial limits of functions meromorphic in a circular disc.
\newblock {\em Math. Z.}, 76:404--410, 1961.

\bibitem{MR0315097}
I.~L. Chang.
\newblock On the zeros of power series with {H}adamard gaps-distribution in
  sectors.
\newblock {\em Trans. Amer. Math. Soc.}, 178:393--400, 1973.

\bibitem{doi:10.1112/blms/bds003}
J.~Chang.
\newblock On meromorphic functions whose first derivatives have finitely many
  zeros.
\newblock {\em Bulletin of the London Mathematical Society}, 44(4):703--715,
  2012.

\bibitem{MR0585642}
K.~H. Chang.
\newblock Asymptotic values of entire and meromorphic functions.
\newblock {\em Sci. Sinica}, 20(6):720--739, 1977.

\bibitem{MR0114930}
Z.~Charzy\'nski and M.~Schiffer.
\newblock A new proof of the {B}ieberbach conjecture for the fourth
  coefficient.
\newblock {\em Arch. Rational Mech. Anal.}, 5:187--193 (1960), 1960.

\bibitem{MR1360682}
H.~H. Chen and M.~L. Fang.
\newblock The value distribution of {$f^nf'$}.
\newblock {\em Sci. China Ser. A}, 38(7):789--798, 1995.

\bibitem{zbMATH00971867}
H.~H. {Chen} and P.~M. {Gauthier}.
\newblock {On Bloch's constant.}
\newblock {\em {J. Anal. Math.}}, 69:275--291, 1996.

\bibitem{MR0276445}
Y.~M. Chen and M.~C. Liu.
\newblock On {L}ittlewoods's conjectural inequalities.
\newblock {\em J. London Math. Soc. (2)}, 1:385--397, 1969.

\bibitem{MR0227423}
J.~A. Cima and P.~Colwell.
\newblock Blaschke quotients and normality.
\newblock {\em Proc. Amer. Math. Soc.}, 19:796--798, 1968.

\bibitem{MR0107715}
J.~Clunie.
\newblock On schlicht functions.
\newblock {\em Ann. of Math. (2)}, 69:511--519, 1959.

\bibitem{MR0143906}
J.~Clunie.
\newblock On integral and meromorphic functions.
\newblock {\em J. London Math. Soc.}, 37:17--27, 1962.

\bibitem{MR0213317}
J.~Clunie.
\newblock On a problem of {E}rd\"os.
\newblock {\em J. London Math. Soc.}, 42:133--136, 1967.

\bibitem{MR0214769}
J.~Clunie.
\newblock On a result of {H}ayman.
\newblock {\em J. London Math. Soc.}, 42:389--392, 1967.

\bibitem{MR0222262}
J.~Clunie and P.~Erd\"os.
\newblock On the partial sums of power series.
\newblock {\em Proc. Roy. Irish Acad. Sect. A}, 65:113--123 (1967), 1967.

\bibitem{zbMATH00561242}
J.~{Clunie}, A.~{Eremenko}, and J.~{Rossi}.
\newblock {On equilibrium points of logarithmic and Newtonian potentials.}
\newblock {\em {J. Lond. Math. Soc., II. Ser.}}, 47(2):309--320, 1993.

\bibitem{MR0170005}
J.~Clunie and W.~K. Hayman.
\newblock The maximum term of a power series.
\newblock {\em J. Analyse Math.}, 12:143--186, 1964.

\bibitem{MR0110814}
J.~Clunie and F.~R. Keogh.
\newblock On starlike and convex schlicht functions.
\newblock {\em J. London Math. Soc.}, 35:229--233, 1960.

\bibitem{MR0207974}
J.~Clunie and Ch. Pommerenke.
\newblock On the coefficients of univalent functions.
\newblock {\em Michigan Math. J.}, 14:71--78, 1967.

\bibitem{MR752388}
J.~Clunie and T.~Sheil-Small.
\newblock Harmonic univalent functions.
\newblock {\em Ann. Acad. Sci. Fenn. Ser. A I Math.}, 9:3--25, 1984.

\bibitem{MR623468}
J.~G. Clunie.
\newblock Some remarks on extreme points in function theory.
\newblock In {\em Aspects of contemporary complex analysis ({P}roc. {NATO}
  {A}dv. {S}tudy {I}nst., {U}niv. {D}urham, {D}urham, 1979)}, pages 137--146.
  Academic Press, London-New York, 1980.

\bibitem{MR0133397}
P.~J. Cohen.
\newblock On a conjecture of {L}ittlewood and idempotent measures.
\newblock {\em Amer. J. Math.}, 82:191--212, 1960.

\bibitem{MR687639}
R.~Coifman, P.~W. Jones, and J.~L. Rubio~de Francia.
\newblock Constructive decomposition of {BMO} functions and factorization of
  {$A_{p}$}\ weights.
\newblock {\em Proc. Amer. Math. Soc.}, 87(4):675--676, 1983.

\bibitem{MR0050010}
E.~F. Collingwood and M.~L. Cartwright.
\newblock Boundary theorems for a function meromorphic in the unit circle.
\newblock {\em Acta Math.}, 87:83--146, 1952.

\bibitem{MR0231999}
E.~F. Collingwood and A.~J. Lohwater.
\newblock {\em The theory of cluster sets}.
\newblock Cambridge Tracts in Mathematics and Mathematical Physics, No. 56.
  Cambridge University Press, Cambridge, 1966.

\bibitem{MR0365016}
A.~C\'ordoba.
\newblock A radial multiplier and a related {K}akeya maximal function.
\newblock {\em Bull. Amer. Math. Soc.}, 81:428--430, 1975.

\bibitem{MR638621}
A.~C\'ordoba and B.~L\'opez-Melero.
\newblock Spherical summation: a problem of {E}. {M}. {S}tein.
\newblock {\em Ann. Inst. Fourier (Grenoble)}, 31(3):x, 147--152, 1981.

\bibitem{MR0009069}
R.~Courant and D.~Hilbert.
\newblock {\em Methoden der {M}athematischen {P}hysik. {V}ols. {I}, {II}}.
\newblock Interscience Publishers, Inc., N.Y., 1943.

\bibitem{Cremer}
H.~Cremer.
\newblock Zur {Z}entrumproblem.
\newblock {\em Math. Ann.}, 98:151--163, 1928.

\bibitem{MR674281}
G.~Csordas, A.~J. Lohwater, and T.~Ramsey.
\newblock Lacunary series and the boundary behavior of {B}loch functions.
\newblock {\em Michigan Math. J.}, 29(3):281--288, 1982.

\bibitem{MR0322184}
B.~Dahlberg.
\newblock Mean values of subharmonic functions.
\newblock {\em Ark. Mat.}, 10:293--309, 1972.

\bibitem{MR0466593}
B.~E.~J. Dahlberg.
\newblock Estimates of harmonic measure.
\newblock {\em Arch. Rational Mech. Anal.}, 65(3):275--288, 1977.

\bibitem{MR0450558}
M.~Damodaran.
\newblock On the distribution of values of meromorphic functions of slow
  growth.
\newblock pages 17--21. Lecture Notes in Math., Vol. 599, 1977.

\bibitem{MR3525402}
A.~A. Danielyan.
\newblock {R}ubel's problem on bounded analytic functions.
\newblock {\em Ann. Acad. Sci. Fenn. Math.}, 41(2):813--816, 2016.

\bibitem{MR0124681}
H.~Davenport.
\newblock On a theorem of {P}. {J}. {C}ohen.
\newblock {\em Mathematika}, 7:93--97, 1960.

\bibitem{MR0350009}
A.~M. Davie.
\newblock Analytic capacity and approximation problems.
\newblock {\em Trans. Amer. Math. Soc.}, 171:409--444, 1972.

\bibitem{MR509941}
S.~T. Davies.
\newblock On the maximum term of an entire function without zeros.
\newblock {\em J. London Math. Soc. (2)}, 18(2):253--260, 1978.

\bibitem{MR735222}
B.~Davis and J.~L. Lewis.
\newblock Paths for subharmonic functions.
\newblock {\em Proc. London Math. Soc. (3)}, 48(3):401--427, 1984.

\bibitem{MR772434}
L.~de~Branges.
\newblock A proof of the {B}ieberbach conjecture.
\newblock {\em Acta Math.}, 154(1-2):137--152, 1985.

\bibitem{MR3162254}
J.~D\'egot.
\newblock Sendov conjecture for high degree polynomials.
\newblock {\em Proc. Amer. Math. Soc.}, 142(4):1337--1349, 2014.

\bibitem{MR0102672}
J.~Delsarte and J.~L. Lions.
\newblock Moyennes g\'en\'eralis\'ees.
\newblock {\em Comment. Math. Helv.}, 33:59--69, 1959.

\bibitem{MR533896}
J.~P. Demailly.
\newblock Fonctions holomorphes \`a croissance polynomiale sur la surface
  d'\'equation {$e^{x}+e^{y}=1$}.
\newblock {\em Bull. Sci. Math. (2)}, 103(2):179--191, 1979.

\bibitem{MR0075294}
A.~Denjoy.
\newblock L'allure asymptotique des fonctions enti\`eres d'ordre fini.
\newblock {\em C. R. Acad. Sci. Paris}, 242:213--218, 1956.

\bibitem{MR0167610}
R.~J. Distler.
\newblock The domain of univalence of certain classes of meromorphic functions.
\newblock {\em Proc. Amer. Math. Soc.}, 15:923--928, 1964.

\bibitem{Dive1}
P.~Dive.
\newblock Sur l'attraction des ellips\"oides homog\`enes.
\newblock {\em C.R. Acad. Sci. Paris}, 192:1443--1446, 1931.

\bibitem{Dive2}
P.~Dive.
\newblock Sur une propri\'et\'e exclusive des homono\"ides ellipsoidaux.
\newblock {\em C.R. Acad. Sci. Paris}, 193:141--142, 1931.

\bibitem{MR507458}
M.~Dixon and J.~Korevaar.
\newblock Nonspanning sets of powers on curves: analyticity theorem.
\newblock {\em Duke Math. J.}, 45(3):543--559, 1978.

\bibitem{MR728980}
A.~Douady.
\newblock Syst\`emes dynamiques holomorphes.
\newblock In {\em Bourbaki seminar, {V}ol. 1982/83}, volume 105 of {\em
  Ast\'erisque}, pages 39--63. Soc. Math. France, Paris, 1983.

\bibitem{MR936853}
A.~Douady.
\newblock Disques de {S}iegel et anneaux de {H}erman.
\newblock {\em Ast\'erisque}, (152-153):4, 151--172 (1988), 1987.
\newblock S\'eminaire Bourbaki, Vol. 1986/87.

\bibitem{zbMATH03455609}
D.~{Drasin}.
\newblock {A meromorphic function with assigned Nevanlinna deficiencies.}
\newblock {Symp. Complex Analysis, Canterbury 1973, 31-41 (1974).}, 1974.

\bibitem{Drasin1987}
D.~Drasin.
\newblock Proof of a conjecture of {F}. {N}evanlinna concerning functions which
  have deficiency sum two.
\newblock {\em Acta Mathematica}, 158(1):1--94, Jul 1987.

\bibitem{MR1405053}
D.~Drasin.
\newblock The minimum modulus of subharmonic functions of order one and a
  method of {F}ryntov.
\newblock {\em J. London Math. Soc. (2)}, 54(2):239--250, 1996.

\bibitem{zbMATH03353695}
D.~{Drasin} and A.~{Weitsman}.
\newblock {The growth of the {N}evanlinna proximity function and the
  logarithmic potential.}
\newblock {\em {Indiana Univ. Math. J.}}, 20:699--715, 1971.

\bibitem{MR746071}
V.~N. Dubinin.
\newblock Change of harmonic measure in symmetrization.
\newblock {\em Mat. Sb. (N.S.)}, 124(166)(2):272--279, 1984.

\bibitem{MR0021975}
D.~Dugu\'e.
\newblock Le d\'efaut au sens de {M}. {N}evanlinna d\'epend de l'origine
  choisie.
\newblock {\em C. R. Acad. Sci. Paris}, 225:555--556, 1947.

\bibitem{zbMATH02117216}
P.~{Duren}.
\newblock {\em {Harmonic mappings in the plane}}, volume 156.
\newblock Cambridge: Cambridge University Press, 2004.

\bibitem{MR861691}
P.~Duren and Y.~J. Leung.
\newblock Generalized support points of the set of univalent functions.
\newblock {\em J. Analyse Math.}, 46:94--108, 1986.

\bibitem{MR0176044}
P.~L. Duren.
\newblock On the {M}arx conjecture for starlike functions.
\newblock {\em Trans. Amer. Math. Soc.}, 118:331--337, 1965.

\bibitem{MR0185057}
P.~L. Duren.
\newblock Smoothness of functions generate. by {R}iesz products.
\newblock {\em Proc. Amer. Math. Soc.}, 16:1263--1268, 1965.

\bibitem{MR0268655}
P.~L. Duren.
\newblock {\em Theory of {$H^{p}$} spaces}.
\newblock Pure and Applied Mathematics, Vol. 38. Academic Press, New
  York-London, 1970.

\bibitem{MR0364623}
P.~L. Duren.
\newblock Estimation of coefficients of univalent functions by a {T}auberian
  remainder theorem.
\newblock {\em J. London Math. Soc. (2)}, 8:279--282, 1974.

\bibitem{MR540059}
P.~L. Duren.
\newblock Successive coefficients of univalent functions.
\newblock {\em J. London Math. Soc. (2)}, 19(3):448--450, 1979.

\bibitem{MR708494}
P.~L. Duren.
\newblock {\em Univalent functions}, volume 259 of {\em Grundlehren der
  Mathematischen Wissenschaften [Fundamental Principles of Mathematical
  Sciences]}.
\newblock Springer-Verlag, New York, 1983.

\bibitem{MR0199359}
P.~L. Duren, H.~S. Shapiro, and A.~L. Shields.
\newblock Singular measures and domains not of {S}mirnov type.
\newblock {\em Duke Math. J.}, 33:247--254, 1966.

\bibitem{MR0284588}
J.~P. Earl.
\newblock On the interpolation of bounded sequences by bounded functions.
\newblock {\em J. London Math. Soc. (2)}, 2:544--548, 1970.

\bibitem{MR0419773}
J.~P. Earl.
\newblock A note on bounded interpolation in the unit disc.
\newblock {\em J. London Math. Soc. (2)}, 13(3):419--423, 1976.

\bibitem{MR0180678}
A.~Edrei.
\newblock Sums of deficiencies of meromorphic functions.
\newblock {\em J. Analyse Math.}, 14:79--107, 1965.

\bibitem{MR0267098}
A.~Edrei.
\newblock A local form of the {P}hragm\'en-{L}indel\"of indicator.
\newblock {\em Mathematika}, 17:149--172, 1970.

\bibitem{10.2307/1993455}
A.~Edrei and W.~H.~J. Fuchs.
\newblock On the growth of meromorphic functions with several deficient values.
\newblock {\em Transactions of the American Mathematical Society},
  93(2):292--328, 1959.

\bibitem{EdreiAlbert;FuchsWolfgangH.J.1959Vdev}
A.~{Edrei} and W.~H.~J. {Fuchs}.
\newblock {Valeurs d\'eficientes et valeurs asymptotiques des fonctions
  m\'eromorphes.}
\newblock {\em {Comment. Math. Helv.}}, 33:258--295, 1959.

\bibitem{MR0477051}
A.~Edrei and W.~H.~J. Fuchs.
\newblock Asymptotic behavior of meromorphic functions with extremal spread.
  {I}.
\newblock {\em Ann. Acad. Sci. Fenn. Ser. A I Math.}, 2:67--111, 1976.

\bibitem{MR0585641}
A.~Edrei and W.~H.~J. Fuchs.
\newblock Asymptotic behavior of meromorphic functions with extremal spread.
  {II}.
\newblock {\em Ann. Acad. Sci. Fenn. Ser. A I Math.}, 3(1):141--168, 1977.

\bibitem{MR0122999}
A.~Edrei, W.~H.~J. Fuchs, and S.~Hellerstein.
\newblock Radial distribution and deficiencies of the values of a meromorphic
  function.
\newblock {\em Pacific J. Math.}, 11:135--151, 1961.

\bibitem{MR0087741}
A.~Edrei and G.~R. MacLane.
\newblock On the zeros of the derivatives of an entire function.
\newblock {\em Proc. Amer. Math. Soc.}, 8:702--706, 1957.

\bibitem{MR0285849}
L.~Ehrenpreis.
\newblock {\em Fourier analysis in several complex variables}.
\newblock Pure and Applied Mathematics, Vol. XVII. Wiley-Interscience
  Publishers A Division of J. Wiley \& Sons, New York-London-Sydney, 1970.

\bibitem{MR0204964}
\'A. Elbert.
\newblock A contribution to the problem of the convergence of the {F}ourier
  series of {$L^{2}$} integrable functions.
\newblock {\em Studia Sci. Math. Hungar.}, 1:147--151, 1966.

\bibitem{MR0209422}
\'A. Elbert.
\newblock \"{U}ber eine {V}ermutung von {E}rd\"os betreffs {P}olynome. {I}.
\newblock {\em Studia Sci. Math. Hungar.}, 1:119--128, 1966.

\bibitem{MR741643}
M.~M. Elhosh.
\newblock On successive coefficients of close-to-convex functions.
\newblock {\em Proc. Roy. Soc. Edinburgh Sect. A}, 96(1-2):47--49, 1984.

\bibitem{MR0003595}
P.~Erd\"os.
\newblock Note on some elementary properties of polynomials.
\newblock {\em Bull. Amer. Math. Soc.}, 46:954--958, 1940.

\bibitem{MR0014608}
P.~Erd\"os.
\newblock On a lemma of {L}ittlewood and {O}fford.
\newblock {\em Bull. Amer. Math. Soc.}, 51:898--902, 1945.

\bibitem{MR0168482}
P.~Erd\"os.
\newblock An interpolation problem associated with the continuum hypothesis.
\newblock {\em Michigan Math. J.}, 11:9--10, 1964.

\bibitem{MR0101311}
P.~Erd\"os, F.~Herzog, and G.~Piranian.
\newblock Metric properties of polynomials.
\newblock {\em J. Analyse Math.}, 6:125--148, 1958.

\bibitem{MR501765}
A.~Eremenko.
\newblock The growth of the {N}evanlinna proximity function.
\newblock {\em Sibirsk. Mat. Zh.}, 19(3):571--576, 717, 1978.

\bibitem{MR522410}
A.~Eremenko.
\newblock The set of asymptotic values of a finite order meromorphic function.
\newblock {\em Mat. Zametki}, 24(6):779--783, 893, 1978.

\bibitem{Eremenko1980}
A.~Eremenko.
\newblock Growth of entire and subharmonic functions on asymptotic curves.
\newblock {\em Siberian Mathematical Journal}, 21(5):673--683, Sep 1980.

\bibitem{MR592215}
A.~Eremenko.
\newblock The growth of entire and subharmonic functions on asymptotic curves.
\newblock {\em Sibirsk. Mat. Zh.}, 21(5):39--51, 189, 1980.

\bibitem{MR667974}
A.~Eremenko.
\newblock Meromorphic solutions of algebraic differential equations.
\newblock {\em Uspekhi Mat. Nauk}, 37(4(226)):53--82, 240, 1982.

\bibitem{MR686077}
A.~Eremenko.
\newblock Meromorphic solutions of equations of {B}riot-{B}ouquet type.
\newblock {\em Teor. Funktsi\u\i \ Funktsional. Anal. i Prilozhen.},
  (38):48--56, 127, 1982.

\bibitem{Eremenko1984}
A.~Eremenko.
\newblock Meromorphic solutions of first-order algebraic differential
  equations.
\newblock {\em Functional Analysis and Its Applications}, 18(3):246--248, Jul
  1984.

\bibitem{MR858915}
A.~Eremenko.
\newblock On the natural asymptotic curves of meromorphic functions.
\newblock {\em Complex Variables Theory Appl.}, 4(4):305--309, 1985.

\bibitem{MR853888}
A.~Eremenko.
\newblock Inverse problem of the theory of distribution of values for
  finite-order meromorphic functions.
\newblock {\em Sibirsk. Mat. Zh.}, 27(3):87--102, 223, 1986.

\bibitem{MR948815}
A.~Eremenko.
\newblock Entire functions bounded on the real axis.
\newblock {\em Dokl. Akad. Nauk SSSR}, 300(3):544--546, 1988.

\bibitem{MR1145577}
A.~Eremenko.
\newblock A counterexample to the {A}rakelyan conjecture.
\newblock {\em Bull. Amer. Math. Soc. (N.S.)}, 27(1):159--164, 1992.

\bibitem{MR1266090}
A.~Eremenko.
\newblock Meromorphic functions with small ramification.
\newblock {\em Indiana Univ. Math. J.}, 42(4):1193--1218, 1993.

\bibitem{erem2006ae}
A.~Eremenko.
\newblock A {M}arkov-type inequality for arbitrary plane continua.
\newblock 135, 07 2006.

\bibitem{MR2894426}
A.~Eremenko and A.~Gabrielov.
\newblock Singular perturbation of polynomial potentials with applications to
  {$PT$}-symmetric families.
\newblock {\em Mosc. Math. J.}, 11(3):473--503, 629--630, 2011.

\bibitem{zbMATH00868218}
A.~{Eremenko} and D.~H. {Hamilton}.
\newblock {On the area distortion by quasiconformal mappings.}
\newblock {\em {Proc. Am. Math. Soc.}}, 123(9):2793--2797, 1995.

\bibitem{zbMATH00569439}
A.~{Eremenko}, J.~{Langley}, and J.~{Rossi}.
\newblock {On the zeros of meromorphic functions of the form $f(z)=
  \sum\sb{k=1}\sp \infty {{a\sb k} \over {z-z\sb k}}$.}
\newblock {\em {J. Anal. Math.}}, 62:271--286, 1994.

\bibitem{eremlempert}
A~Eremenko and L.~Lempert.
\newblock An extremal problem for polynomials.
\newblock 122, 09 1994.

\bibitem{EremLiaoNg}
A.~Eremenko, L.~W. Liao, and T.~W. Ng.
\newblock Meromorphic solutions of higher order {B}riot-{B}ouquet differential
  equations.
\newblock 146, 07 2007.

\bibitem{MR918638}
A.~Eremenko and M.~Yu. Lyubich.
\newblock Examples of entire functions with pathological dynamics.
\newblock {\em J. London Math. Soc. (2)}, 36(3):458--468, 1987.

\bibitem{MR2210353}
A.~Eremenko and S.~Merenkov.
\newblock Nevanlinna functions with real zeros.
\newblock {\em Illinois J. Math.}, 49(4):1093--1110, 2005.

\bibitem{MR897006}
A.~Eremenko and M.~L. Sodin.
\newblock A proof of the conditional {L}ittlewood theorem on the distribution
  of the values of entire functions.
\newblock {\em Izv. Akad. Nauk SSSR Ser. Mat.}, 51(2):421--428, 448, 1987.

\bibitem{MR866467}
A.~Eremenko, M.~L. Sodin, and D.~F. Shia.
\newblock The minimum of the modulus of an entire function on a sequence of
  {P}\'olya peaks.
\newblock {\em Teor. Funktsi\u\i \ Funktsional. Anal. i Prilozhen.},
  (45):26--40, i, 1986.

\bibitem{eremhayman1999ae}
A.~E. {Eremenko} and W.~K. {Hayman}.
\newblock {On the length of lemniscates.}
\newblock In {\em {Paul Erd\H{o}s and his mathematics I. Based on the
  conference, Budapest, Hungary, July 4--11, 1999}}, pages 241--242. Berlin:
  Springer; Budapest: J\'anos Bolyai Mathematical Society, 2002.

\bibitem{MR1050784}
M.~Ess\'en.
\newblock A generalization of the {M}. {R}iesz theorem on conjugate functions
  and the {Z}ygmund {$L\log L$}-theorem to {$\mathbb{R}^d,\;d\ge 2$}.
\newblock {\em Ann. Acad. Sci. Fenn. Ser. A I Math.}, 15(1):83--105, 1990.

\bibitem{MR0450584}
M.~Ess\'en and D.~F. Shea.
\newblock Applications of {D}enjoy integral inequalities to growth problems for
  subharmonic and meromorphic functions.
\newblock pages 59--68. London Math. Soc. Lecture Note Ser., No. 12, 1974.

\bibitem{MR812321}
M.~Ess\'en, D.~F. Shea, and C.~S. Stanton.
\newblock A value-distribution criterion for the class {$L\,{\rm log}\,L$}, and
  some related questions.
\newblock {\em Ann. Inst. Fourier (Grenoble)}, 35(4):127--150, 1985.

\bibitem{MR1504787}
P.~Fatou.
\newblock Sur les \'equations fonctionnelles.
\newblock {\em Bull. Soc. Math. France}, 47:161--271, 1919.

\bibitem{MR1504797}
P.~Fatou.
\newblock Sur les \'equations fonctionnelles.
\newblock {\em Bull. Soc. Math. France}, 48:208--314, 1920.

\bibitem{MR1504792}
P.~Fatou.
\newblock Sur les \'equations fonctionnelles.
\newblock {\em Bull. Soc. Math. France}, 48:33--94, 1920.

\bibitem{MR1504820}
P.~Fatou.
\newblock Sur les fronti\`eres de certains domaines.
\newblock {\em Bull. Soc. Math. France}, 51:16--22, 1923.

\bibitem{MR1555220}
P.~Fatou.
\newblock Sur l'it\'eration des fonctions transcendantes {e}nti\`eres.
\newblock {\em Acta Math.}, 47(4):337--370, 1926.

\bibitem{Feng1976}
Jinfu Feng and T.~H. Macgregor.
\newblock Estimates on integral means of the derivatives of univalent
  functions.
\newblock {\em Journal d'Analyse Math{\'e}matique}, 29(1):203--231, Dec 1976.

\bibitem{MR601729}
P.~C. Fenton.
\newblock The minimum of small entire functions.
\newblock {\em Proc. Amer. Math. Soc.}, 81(4):557--561, 1981.

\bibitem{MR715310}
P.~C. Fenton.
\newblock Entire functions having asymptotic functions.
\newblock {\em Bull. Austral. Math. Soc.}, 27(3):321--328, 1983.

\bibitem{MR1651013}
P.~C. Fenton and J.~Rossi.
\newblock Cercles de remplissage for entire functions.
\newblock {\em Bull. London Math. Soc.}, 31(1):59--66, 1999.

\bibitem{MR769161}
J.~L. Fern\'andez.
\newblock On the growth and coefficients of analytic functions.
\newblock {\em Ann. of Math. (2)}, 120(3):505--516, 1984.

\bibitem{MR845869}
J.~L. Fern\'andez.
\newblock A note on entropy and inner functions.
\newblock {\em Israel J. Math.}, 53(2):158--162, 1986.

\bibitem{MR959269}
J.~L. Fern\'andez.
\newblock A boundedness theorem for {$L^1/H^1_0$}.
\newblock {\em Michigan Math. J.}, 35(2):227--231, 1988.

\bibitem{MR0335777}
C.~H. FitzGerald.
\newblock Quadratic inequalities and coefficient estimates for schlicht
  functions.
\newblock {\em Arch. Rational Mech. Anal.}, 46:356--368, 1972.

\bibitem{MR792819}
C.~H. FitzGerald and Ch. Pommerenke.
\newblock The de {B}ranges theorem on univalent functions.
\newblock {\em Trans. Amer. Math. Soc.}, 290(2):683--690, 1985.

\bibitem{MR768733}
C.~H. FitzGerald, B.~Rodin, and S.~E. Warschawski.
\newblock Estimates of the harmonic measure of a continuum in the unit disk.
\newblock {\em Trans. Amer. Math. Soc.}, 287(2):681--685, 1985.

\bibitem{MR1694289}
C.~H. FitzGerald and Frederick Weening.
\newblock Existence and uniqueness of rectilinear slit maps.
\newblock {\em Trans. Amer. Math. Soc.}, 352(2):485--513, 2000.

\bibitem{zbMATH06729355}
R.~{Fournier}.
\newblock {On a new proof and an extension of Jack's lemma.}
\newblock {\em {J. Appl. Anal.}}, 23(1):21--24, 2017.

\bibitem{MR0264160}
G.~Frank.
\newblock Picardsche {A}usnahmewerte bei {L}\"osungen linearer
  {D}ifferentialgleichungen.
\newblock {\em Manuscripta Math.}, 2:181--190, 1970.

\bibitem{MR0422615}
G.~Frank.
\newblock Eine {V}ermutung von {H}ayman \"uber {N}ullstellen meromorpher
  {F}unktionen.
\newblock {\em Math. Z.}, 149(1):29--36, 1976.

\bibitem{MR868452}
G.~Frank and S.~Hellerstein.
\newblock On the meromorphic solutions of nonhomogeneous linear differential
  equations with polynomial coefficients.
\newblock {\em Proc. London Math. Soc. (3)}, 53(3):407--428, 1986.

\bibitem{MR0430250}
G.~Frank, W.~Hennekemper, and G.~Polloczek.
\newblock \"{U}ber die {N}ullstellen meromorpher {F}unktionen und deren
  {A}bleitungen.
\newblock {\em Math. Ann.}, 225(2):145--154, 1977.

\bibitem{MR0012127}
O.~Frostman.
\newblock Sur les produits de {B}laschke.
\newblock {\em Kungl. Fysiografiska S\"allskapets i Lund F\"orhandlingar [Proc.
  Roy. Physiog. Soc. Lund]}, 12(15):169--182, 1942.

\bibitem{MR1687884}
A.~Fryntov.
\newblock On behavior of gap series on curves and a {$\cos\pi\lambda$}-type
  theorem.
\newblock {\em Complex Variables Theory Appl.}, 37(1-4):195--209, 1998.

\bibitem{zbMATH05580784}
A.~{Fryntov} and F.~{Nazarov}.
\newblock {New estimates for the length of the Erd\H os-Herzog-Piranian
  lemniscate.}
\newblock In {\em {Linear and complex analysis. Dedicated to V. P. Havin on the
  occasion of his 75th birthday}}, pages 49--60. Providence, RI: American
  Mathematical Society (AMS), 2009.

\bibitem{MR0159933}
W.~H.~J. Fuchs.
\newblock Proof of a conjecture of {G}. {P}\'olya concerning gap series.
\newblock {\em Illinois J. Math.}, 7:661--667, 1963.

\bibitem{zbMATH03739887}
W.~H.~J. {Fuchs}.
\newblock {A Phragmen-Lindel\"of theorem conjectured by D. J. Newman.}
\newblock {\em {Trans. Am. Math. Soc.}}, 267:285--293, 1981.

\bibitem{MR0379870}
B.~Fuglede.
\newblock Asymptotic paths for subharmonic functions.
\newblock {\em Math. Ann.}, 213:261--274, 1975.

\bibitem{MR0072959}
R.~M. Gabriel.
\newblock An extended principle of the maximum for harmonic functions in
  {$3$}-dimensions.
\newblock {\em J. London Math. Soc.}, 30:388--401, 1955.

\bibitem{zbMATH03834335}
D.~{Gaier}.
\newblock {\"Uber R\"aume konformer Selbstabbildungen ebener Gebiete.}
\newblock {\em {Math. Z.}}, 187:227--257, 1984.

\bibitem{MR0276742}
T.~W. Gamelin.
\newblock Localization of the corona problem.
\newblock {\em Pacific J. Math.}, 34:73--81, 1970.

\bibitem{MR0156980}
T.~Ganelius.
\newblock The zeros of the partial sums of power series.
\newblock {\em Duke Math. J.}, 30:533--540, 1963.

\bibitem{MR0060020}
P.~R. Garabedian and H.~L. Royden.
\newblock The one-quarter theorem for mean univalent functions.
\newblock {\em Ann. of Math. (2)}, 59:316--324, 1954.

\bibitem{MR0066457}
P.~R. Garabedian and M.~Schiffer.
\newblock A coefficient inequality for schlicht functions.
\newblock {\em Ann. of Math. (2)}, 61:116--136, 1955.

\bibitem{MR0070712}
P.~R. Garabedian and M.~Schiffer.
\newblock A proof of the {B}ieberbach conjecture for the fourth coefficient.
\newblock {\em J. Rational Mech. Anal.}, 4:427--465, 1955.

\bibitem{MR0454006}
J.~Garnett.
\newblock {\em Analytic capacity and measure}.
\newblock Lecture Notes in Mathematics, Vol. 297. Springer-Verlag, Berlin-New
  York, 1972.

\bibitem{MR628971}
J.~B. Garnett.
\newblock {\em Bounded analytic functions}, volume~96 of {\em Pure and Applied
  Mathematics}.
\newblock Academic Press, Inc. [Harcourt Brace Jovanovich, Publishers], New
  York-London, 1981.

\bibitem{MR721565}
J.~B. Garnett, F.~W. Gehring, and P.~W. Jones.
\newblock Conformally invariant length sums.
\newblock {\em Indiana Univ. Math. J.}, 32(6):809--829, 1983.

\bibitem{MR793236}
J.~B. Garnett and P.~W. Jones.
\newblock The corona theorem for {D}enjoy domains.
\newblock {\em Acta Math.}, 155(1-2):27--40, 1985.

\bibitem{zbMATH05296346}
J.~B. {Garnett} and D.~E. {Marshall}.
\newblock {\em {Harmonic measure. Paperback reprint of the hardback edition
  2005.}}, volume~2.
\newblock Cambridge: Cambridge University Press, paperback reprint of the
  hardback edition 2005 edition, 2008.

\bibitem{MR576870}
P.~M. Gauthier, M.~Goldstein, and W.~H. Ow.
\newblock Uniform approximation on unbounded sets by harmonic functions with
  logarithmic singularities.
\newblock {\em Trans. Amer. Math. Soc.}, 261(1):169--183, 1980.

\bibitem{gauthier1989}
P.M. Gauthier, W.~Hengartner, and A.~Stray.
\newblock A problem of {R}ubel concerning approximation on unbounded sets by
  entire functions.
\newblock {\em Rocky Mountain J. Math.}, 19(1):127--136, 03 1989.

\bibitem{MR0402038}
F.~W. Gehring.
\newblock The {$L^{p}$}-integrability of the partial derivatives of a
  quasiconformal mapping.
\newblock {\em Acta Math.}, 130:265--277, 1973.

\bibitem{MR0361008}
F.~W. Gehring.
\newblock The {H}ausdorff measure of sets which link in {E}uclidean space.
\newblock pages 159--167, 1974.

\bibitem{MR0382638}
F.~W. Gehring.
\newblock Lower dimensional absolute continuity properties of quasiconformal
  mappings.
\newblock {\em Math. Proc. Cambridge Philos. Soc.}, 78:81--93, 1975.

\bibitem{MR0499134}
F.~W. Gehring.
\newblock Spirals and the universal {T}eichm\"uller space.
\newblock {\em Acta Math.}, 141(1-2):99--113, 1978.

\bibitem{MR663670}
F.~W. Gehring, W.~K. Hayman, and A.~Hinkkanen.
\newblock Analytic functions satisfying {H}\"older conditions on the boundary.
\newblock {\em J. Approx. Theory}, 35(3):243--249, 1982.

\bibitem{MR749106}
F.~W. Gehring and Ch. Pommerenke.
\newblock On the {N}ehari univalence criterion and quasicircles.
\newblock {\em Comment. Math. Helv.}, 59(2):226--242, 1984.

\bibitem{MR0201635}
F.~W. Gehring and E.~Reich.
\newblock Area distortion under quasiconformal mappings.
\newblock {\em Ann. Acad. Sci. Fenn. Ser. A I No.}, 388:15, 1966.

\bibitem{MR0437767}
M.~A. Girnyk.
\newblock On the inverse problem of the theory of the distribution of values
  for functions that are analytic in the unit disc.
\newblock {\em Ukrain. Mat. \v{Z}.}, 29(1):32--39, 142, 1977.

\bibitem{zbMATH03908769}
D.~{Gnuschke} and Ch. {Pommerenke}.
\newblock {On the radial limits of functions with Hadamard gaps.}
\newblock {\em {Mich. Math. J.}}, 32:21--31, 1985.

\bibitem{MR1101984}
B.~Goddard.
\newblock Finite exponential series and {N}ewman polynomials.
\newblock {\em Proc. Amer. Math. Soc.}, 116(2):313--320, 1992.

\bibitem{MR0367172}
R.~M. Goel.
\newblock Functions starlike and convex of order {$\alpha $}.
\newblock {\em J. London Math. Soc. (2)}, 9:128--130, 1974/75.

\bibitem{MR0085396}
A.~A. Gol'dberg.
\newblock On one-valued integrals of differential equations of the first order.
\newblock {\em Ukrain. Mat. \v{Z}.}, 8:254--261, 1956.

\bibitem{MR0170006}
A.~A. Gol'dberg.
\newblock On the possible value of the lower order of an entire function with a
  finite deficient value.
\newblock {\em Dokl. Akad. Nauk SSSR}, 159:968--970, 1964.

\bibitem{MR0229827}
A.~A. Gol'dberg.
\newblock The distribution of the values of an entire function with respect to
  arguments.
\newblock {\em Acta Math. Acad. Sci. Hungar.}, 19:191--199, 1968.

\bibitem{goldbergunknown}
A.~A. {Gol'dberg}.
\newblock {\"{U}ber die {W}erteverteilung einer ganzen {F}unktion nach dem
  {A}rgument.}
\newblock {\em {Acta Math. Acad. Sci. Hung.}}, 19:191--199, 1968.

\bibitem{MR0320316}
A.~A. Gol'dberg.
\newblock The representation of a meromorphic function in the form of a
  quotient of entire functions.
\newblock {\em Izv. Vys\v{s}. U\v{c}ebn. Zaved. Matematika}, (10(125)):13--17,
  1972.

\bibitem{MR0335805}
A.~A. Gol'dberg.
\newblock The branched values of entire functions.
\newblock {\em Sibirsk. Mat. \v{Z}.}, 14:862--866, 911, 1973.

\bibitem{Gol'dberg1973}
A.~A. Gol'dberg.
\newblock On ramified values of entire functions.
\newblock {\em Siberian Mathematical Journal}, 14(4):599--602, Jul 1973.

\bibitem{MR0466551}
A.~A. Gol'dberg.
\newblock Counting functions of sequences of {$a$}-points for entire functions.
\newblock {\em Sibirsk. Mat. \v{Z}.}, 19(1):28--36, 236, 1978.

\bibitem{MR540239}
A.~A. Gol'dberg.
\newblock The minimum modulus of a meromorphic function of slow growth.
\newblock {\em Mat. Zametki}, 25(6):835--844, 956, 1979.

\bibitem{MR537357}
A.~A. Gol'dberg.
\newblock Sets on which the modulus of an entire function has a lower bound.
\newblock {\em Sibirsk. Mat. Zh.}, 20(3):512--518, 691, 1979.

\bibitem{MR760854}
A.~A. Gol'dberg.
\newblock Analytic functions mapping a disk on a disk.
\newblock {\em Izv. Vyssh. Uchebn. Zaved. Mat.}, (6):24--25, 1984.

\bibitem{GoldbergErem}
A.~A. Gol'dberg and A.~Eremenko.
\newblock On asymptotic curves of entire functions of finite order ({R}ussian).
\newblock {\em Mat. Sb. (N. S.)}, 79, 109 (151)(No. 4):555--581, 1982.

\bibitem{MR916442}
A.~A. Gol'dberg, A.~Eremenko, and M.~L. Sodin.
\newblock Exceptional values in the sense of {R}. {N}evanlinna and in the sense
  of {V}. {P}. {P}etrenko. {I}.
\newblock {\em Teor. Funktsi\u\i \ Funktsional. Anal. i Prilozhen.},
  (47):41--51, 1987.

\bibitem{MR0165105}
A.~A. Gol'dberg and B.~Ja. Levin.
\newblock Entire functions which are bounded on the real axis.
\newblock {\em Dokl. Akad. Nauk SSSR}, 157:19--21, 1964.

\bibitem{MR0280720}
A.~A. Gol'dberg and I.~V. Ostrovskii.
\newblock {\em Raspredelenie znacheni\u\i \ meromorfnykh funktsi\u\i}.
\newblock Izdat. ``Nauka'', Moscow, 1970.

\bibitem{MR883581}
A.~A. Gol'dberg and I.~V. Ostrovskii.
\newblock Indicators of entire absolutely monotone functions of finite order.
\newblock {\em Sibirsk. Mat. Zh.}, 27(6):33--49, 1986.

\bibitem{MR1054303}
A.~A. Gol'dberg and I.~V. Ostrovskii.
\newblock Indicators of finite-order entire functions that can be represented
  by {D}irichlet series.
\newblock {\em Dokl. Akad. Nauk Ukrain. SSR Ser. A}, (1):14--16, 85, 1990.

\bibitem{MR2435270}
A.~A. Gol'dberg and I.~V. Ostrovskii.
\newblock {\em Value distribution of meromorphic functions}, volume 236 of {\em
  Translations of Mathematical Monographs}.
\newblock American Mathematical Society, Providence, RI, 2008.
\newblock Translated from the 1970 {R}ussian original by M. Ostrovskii, With an
  appendix by A. Eremenko and J. K. Langley.

\bibitem{GoldbergTairova}
A.~A. Gol'dberg and V.~G. Tairova.
\newblock On integral functions with two finite completely multiple values
  ({R}ussian).
\newblock {\em Zap. Mech.-Mat-Fak. Harkov. Gos. Univ.}, XXIX:67--78, 1963.

\bibitem{goldstein1983}
M.~Goldstein.
\newblock An example of an {A}rakelian glove which is a weak {A}rakelian set.
\newblock {\em Illinois J. Math.}, 27(1):138--144, 03 1983.

\bibitem{MR679930}
M.~Goldstein, R.~R. Hall, T.~Sheil-Small, and H.~L. Smith.
\newblock Convexity preservation of inverse {E}uler operators and a problem of
  {S}. {M}iller.
\newblock {\em Bull. London Math. Soc.}, 14(6):537--541, 1982.

\bibitem{10.2307/2000876}
M.~Goldstein, W.~Haussman, and L.~Rogge.
\newblock On the mean value property of harmonic functions and best harmonic
  $l^1$-approximation.
\newblock {\em Transactions of the American Mathematical Society},
  305(2):505--515, 1988.

\bibitem{MR734993}
M.~Goldstein and J.~N. McDonald.
\newblock An extremal problem for nonnegative trigonometric polynomials.
\newblock {\em J. London Math. Soc. (2)}, 29(1):81--88, 1984.

\bibitem{Goluzin}
G.~M. Goluzin.
\newblock Einige {K}oeffizientenabsch\"tzungen f\"ur schilichte {F}unktionen.
\newblock {\em Rec. Math. [Mat. Sbornik] N. S.}, 3(2):321--30, 1938.

\bibitem{MR0043209}
G.~M. Goluzin.
\newblock On majorants of subordinate analytic functions. {I}.
\newblock {\em Mat. Sbornik N.S.}, 29(71):209--224, 1951.

\bibitem{MR0231997}
A.~A. Gon\v{c}ar.
\newblock Generalized analytic continuation.
\newblock {\em Mat. Sb. (N.S.)}, 76 (118):135--146, 1968.

\bibitem{MR0023910}
A.~W. Goodman.
\newblock On some determinants related to {$p$}-valent functions.
\newblock {\em Trans. Amer. Math. Soc.}, 63:175--192, 1948.

\bibitem{MR2937804}
A.~W. Goodman.
\newblock {\em On some determinants related to $p$-valent functions}.
\newblock ProQuest LLC, Ann Arbor, MI, 1948.
\newblock Doctoral Thesis--Columbia University.

\bibitem{MR0033886}
A.~W. Goodman.
\newblock On the {S}chwarz-{C}hristoffel transformation and {$p$}-valent
  functions.
\newblock {\em Trans. Amer. Math. Soc.}, 68:204--223, 1950.

\bibitem{MR0041920}
A.~W. Goodman.
\newblock Typically-real functions with assigned zeros.
\newblock {\em Proc. Amer. Math. Soc.}, 2:349--357, 1951.

\bibitem{MR0247057}
A.~W. Goodman.
\newblock The valence of certain means.
\newblock {\em J. Analyse Math.}, 22:355--361, 1969.

\bibitem{MR567984}
A.~W. Goodman.
\newblock Valence sequences.
\newblock {\em Proc. Amer. Math. Soc.}, 79(3):422--426, 1980.

\bibitem{MR760942}
A.~W. Goodman.
\newblock Convex functions of bounded type.
\newblock {\em Proc. Amer. Math. Soc.}, 92(4):541--546, 1984.

\bibitem{MR835886}
A.~W. Goodman.
\newblock More on convex functions of bounded type.
\newblock {\em Proc. Amer. Math. Soc.}, 97(2):303--306, 1986.

\bibitem{MR1005663}
A.~W. Goodman.
\newblock The complete multivalence of \v{C}akalov-{D}istler sums.
\newblock {\em Complex Variables Theory Appl.}, 11(2):87--93, 1989.

\bibitem{MR0040430}
A.~W. Goodman and M.~S. Robertson.
\newblock A class of multivalent functions.
\newblock {\em Trans. Amer. Math. Soc.}, 70:127--136, 1951.

\bibitem{zbMATH03634586}
A.~W. {Goodman} and E.~B. {Saff}.
\newblock {On univalent functions convex in one direction.}
\newblock {\em {Proc. Am. Math. Soc.}}, 73:183--187, 1979.

\bibitem{MR0269865}
N.~V. Govorov.
\newblock The estimation from below of a function that is subharmonic in a
  disk.
\newblock {\em Teor. Funkci\u\i \ Funkcional. Anal. i Prilo\v{z}en. Vyp.},
  6:130--150, 1968.

\bibitem{MR0247079}
N.~V. Govorov.
\newblock The {P}aley conjecture.
\newblock {\em Funkcional. Anal. i Prilo\v{z}en.}, 3(2):41--45, 1969.

\bibitem{MR1469316}
J.~Graczyk and G.~\'Swiatek.
\newblock Generic hyperbolicity in the logistic family.
\newblock {\em Ann. of Math. (2)}, 146(1):1--52, 1997.

\bibitem{MR0328061}
A.~Z. Grinshpan.
\newblock Logarithmic coefficients of functions of class {$S$}.
\newblock {\em Sibirsk. Mat. \v{Z}.}, 13:1145--1157, 1199, 1972.

\bibitem{grinspan1976}
A.~Z. Grinshpan.
\newblock Improve. bounds for the difference of the moduli of adjacent
  coefficients of univalent functions, in ``{S}ome questions in the modern
  theory of functions'' (in {R}ussian).
\newblock {\em Sib. Inst. Mat. Novosibirsk}, pages 41--45, 1976.

\bibitem{MR697766}
A.~Z. Grinshpan.
\newblock On the power stability for the {B}ieberbach inequality.
\newblock {\em Zap. Nauchn. Sem. Leningrad. Otdel. Mat. Inst. Steklov. (LOMI)},
  125:58--64, 1983.
\newblock Analytic number theory and the theory of functions, 5.

\bibitem{MR722895}
F.~Gross and C.~F. Osgood.
\newblock On the functional equation {$f^{n}+g^{n}=h^{n}$}\ and a new approach
  to a certain class of more general functional equations.
\newblock {\em Indian J. Math.}, 23(1-3):17--39, 1981.

\bibitem{MR1548975}
W.~Gro\ss.
\newblock \"{U}ber die {S}ingularit\"aten analytischer {F}unktionen.
\newblock {\em Monatsh. Math. Phys.}, 29(1):3--47, 1918.

\bibitem{Guckenheimer}
J.~Guckenheimer.
\newblock Endomorphisms of the {R}iemann sphere.
\newblock In {\em Proceedings Symposia in Pure Mathematics}, volume~15 of {\em
  Proceedings Symposia in Pure Mathematics}, pages 95--123. Amer. Math. Soc.,
  Providence, RI, 1970.

\bibitem{MR853961}
G.~G. Gundersen.
\newblock On the real zeros of solutions of {$f''+A(z)f=0$} where {$A(z)$} is
  entire.
\newblock {\em Ann. Acad. Sci. Fenn. Ser. A I Math.}, 11(2):275--294, 1986.

\bibitem{MR1660942}
G.~G. Gundersen.
\newblock Meromorphic solutions of {$f^6+g^6+h^6\equiv1$}.
\newblock {\em Analysis (Munich)}, 18(3):285--290, 1998.

\bibitem{MR1821651}
G.~G. Gundersen.
\newblock Meromorphic solutions of {$f^5+g^5+h^5\equiv 1$}.
\newblock {\em Complex Variables Theory Appl.}, 43(3-4):293--298, 2001.
\newblock The Chuang special issue.

\bibitem{MR1833253}
G.~G. Gundersen.
\newblock Solutions of {$f''+P(z)f=0$} that have almost all real zeros.
\newblock {\em Ann. Acad. Sci. Fenn. Math.}, 26(2):483--488, 2001.

\bibitem{zbMATH02113438}
G.~G. {Gundersen} and W.~K. {Hayman}.
\newblock {The strength of Cartan's version of Nevanlinna theory.}
\newblock {\em {Bull. Lond. Math. Soc.}}, 36(4):433--454, 2004.

\bibitem{MR0180696}
R.~C. Gunning and H.~Rossi.
\newblock {\em Analytic functions of several complex variables}.
\newblock Prentice-Hall, Inc., Englewood Cliffs, N.J., 1965.

\bibitem{zbMATH05900730}
P.~A. {Gunsul}.
\newblock {Value distribution for a class of small functions in the unit disk.}
\newblock {\em {Int. J. Math. Math. Sci.}}, 2011:24, 2011.

\bibitem{MR0457784}
M.~Hakim and N.~Sibony.
\newblock Quelques conditions pour l'existence de fonctions pics dans des
  domaines pseudoconvexes.
\newblock {\em Duke Math. J.}, 44(2):399--406, 1977.

\bibitem{MR0367545}
G.~Hal\'asz.
\newblock On a result of {S}alem and {Z}ygmund concerning random polynomials.
\newblock {\em Studia Sci. Math. Hungar.}, 8:369--377, 1973.

\bibitem{MR657497}
R.~R. Hall.
\newblock On a conjecture of {S}hapiro about trigonometric series.
\newblock {\em J. London Math. Soc. (2)}, 25(3):407--415, 1982.

\bibitem{MR848860}
R.~R. Hall.
\newblock On some theorems of {H}urwitz and {S}heil-{S}mall.
\newblock {\em Math. Proc. Cambridge Philos. Soc.}, 100(2):365--370, 1986.

\bibitem{zbMATH00020033}
R.~R. {Hall}, W.~K. {Hayman}, and A.~W. {Weitsman}.
\newblock {On asymmetry and capacity.}
\newblock {\em {J. Anal. Math.}}, 56:87--123, 1991.

\bibitem{MR0387585}
R.~R. Hall and J.~H. Williamson.
\newblock On a certain functional equation.
\newblock {\em J. London Math. Soc. (2)}, 12(2):133--136, 1975/76.

\bibitem{Hall}
T.~Hall.
\newblock Sur la mesure harmonique de certains ensembles.
\newblock {\em Ark. Mat. Astr. Fys.}, 25(28):8 pp., 1937.

\bibitem{DHHamilton}
D.~H. Hamilton.
\newblock Doctoral thesis.
\newblock {\em University of London}, 1980.

\bibitem{MR656110}
D.~H. Hamilton.
\newblock The extreme points of {$\Sigma $}.
\newblock {\em Proc. Amer. Math. Soc.}, 85(3):393--396, 1982.

\bibitem{MR684512}
D.~H. Hamilton.
\newblock The dispersion of the coefficients of univalent functions.
\newblock {\em Trans. Amer. Math. Soc.}, 276(1):323--333, 1983.

\bibitem{MR732121}
D.~H. Hamilton.
\newblock A sharp form of the {A}hlfors' distortion theorem, with applications.
\newblock {\em Trans. Amer. Math. Soc.}, 282(2):799--806, 1984.

\bibitem{MR915532}
D.~H. Hamilton.
\newblock Extremal boundary problems.
\newblock {\em Proc. London Math. Soc. (3)}, 56(1):101--113, 1988.

\bibitem{zbMATH03332697}
R.~S. {Hamilton}.
\newblock {Extremal quasiconformal mappings with prescribed boundary values.}
\newblock {\em {Trans. Am. Math. Soc.}}, 138:399--406, 1969.

\bibitem{MR0262512}
L.~J. Hansen.
\newblock Hardy classes and ranges of functions.
\newblock {\em Michigan Math. J.}, 17:235--248, 1970.

\bibitem{MR0460634}
L.~J. Hansen.
\newblock On the growth of entire functions bounded on large sets.
\newblock {\em Canad. J. Math.}, 29(6):1287--1291, 1977.

\bibitem{doi:10.1112/jlms/50.2.349}
W.~Hansen and N.~Nadirashvili.
\newblock Littlewood's one circle problem.
\newblock {\em Journal of the London Mathematical Society}, 50(2):349--360,
  1994.

\bibitem{MR1555303}
G.~H. Hardy and J.~E. Littlewood.
\newblock A maximal theorem with function-theoretic applications.
\newblock {\em Acta Math.}, 54(1):81--116, 1930.

\bibitem{MR0425117}
A.~Harrington and M.~Ortel.
\newblock The dilatation of an extremal quasi-conformal mapping.
\newblock {\em Duke Math. J.}, 43(3):533--544, 1976.

\bibitem{MR0409828}
A.~Harrington and M.~Ortel.
\newblock Two extremal problems.
\newblock {\em Trans. Amer. Math. Soc.}, 221(1):159--167, 1976.

\bibitem{MR0192059}
V.~P. Havin.
\newblock Boundary properties of integrals of {C}auchy type and of conjugate
  harmonic functions in regions with rectifiable boundary.
\newblock {\em Mat. Sb. (N.S.)}, 68 (110):499--517, 1965.

\bibitem{MR734178}
V.~P. Havin, S.~V. Hru\v{s}\v{c}\"ev, and N.~K. Nikol\v{s}ki\u\i, editors.
\newblock {\em Linear and complex analysis problem book}, volume 1043 of {\em
  Lecture Notes in Mathematics}.
\newblock Springer-Verlag, Berlin, 1984.
\newblock 199 research problems.

\bibitem{MR576895}
J.~Hawkes.
\newblock Probabilistic behaviour of some lacunary series.
\newblock {\em Z. Wahrsch. Verw. Gebiete}, 53(1):21--33, 1980.

\bibitem{hayman_1948}
W.~K. Hayman.
\newblock Some inequalities in the theory of functions.
\newblock {\em Mathematical Proceedings of the Cambridge Philosophical
  Society}, 44(2):159�178, 1948.

\bibitem{PLMS:PLMS0450}
W.~K. Hayman.
\newblock Inequalities in the theory of functions.
\newblock {\em Proceedings of the London Mathematical Society},
  s2-51(1):450--473, 1949.

\bibitem{MR0045210}
W.~K. Hayman.
\newblock Some applications of the transfinite diameter to the theory of
  functions.
\newblock {\em J. Analyse Math.}, 1:155--179, 1951.

\bibitem{MR0056083}
W.~K. Hayman.
\newblock The minimum modulus of large integral functions.
\newblock {\em Proc. London Math. Soc. (3)}, 2:469--512, 1952.

\bibitem{MR0054705}
W.~K. Hayman.
\newblock An integral function with a defective value that is neither
  asymptotic nor invariant under change of origin.
\newblock {\em J. London Math. Soc.}, 28:369--376, 1953.

\bibitem{MR0071536}
W.~K. Hayman.
\newblock The asymptotic behaviour of {$p$}-valent functions.
\newblock {\em Proc. London Math. Soc. (3)}, 5:257--284, 1955.

\bibitem{MR0070716}
W.~K. Hayman.
\newblock Uniformly normal families.
\newblock In {\em Lectures on functions of a complex variable}, pages 199--212.
  The University of Michigan Press, Ann Arbor, 1955.

\bibitem{MR0096810}
W.~K. Hayman.
\newblock Bounds for the large coefficients of univalent functions.
\newblock {\em Ann. Acad. Sci. Fenn. Ser. A.I, no.}, 250/13:13, 1958.

\bibitem{MR0108586}
W.~K. Hayman.
\newblock {\em Multivalent functions}.
\newblock Cambridge Tracts in Mathematics and Mathematical Physics, No. 48.
  Cambridge University Press, Cambridge, 1958.

\bibitem{MR0110807}
W.~K. Hayman.
\newblock Picard values of meromorphic functions and their derivatives.
\newblock {\em Ann. of Math. (2)}, 70:9--42, 1959.

\bibitem{MR0125223}
W.~K. Hayman.
\newblock On the growth of integral functions of asymptotic paths.
\newblock {\em J. Indian Math. Soc. (N.S.)}, 24:251--264 (1961), 1960.

\bibitem{MR0111839}
W.~K. Hayman.
\newblock Slowly growing integral and subharmonic functions.
\newblock {\em Comment. Math. Helv.}, 34:75--84, 1960.

\bibitem{M.F.}
W.~K. Hayman.
\newblock {\em Meromorphic functions}.
\newblock Oxford mathematical monographs. Clarendon Press, Oxford, 1964.

\bibitem{MR0168763}
W.~K. Hayman.
\newblock On the characteristic of functions meromorphic in the unit disk and
  of their integrals.
\newblock {\em Acta Math.}, 112:181--214, 1964.

\bibitem{MR0180679}
W.~K. Hayman.
\newblock On the characteristic of functions meromorphic in the plane and of
  their integrals.
\newblock {\em Proc. London Math. Soc. (3)}, 14a:93--128, 1965.

\bibitem{MR0213535}
W.~K. Hayman.
\newblock Mean {$p$}-valent functions with gaps.
\newblock {\em Colloq. Math.}, 16:1--21, 1967.

\bibitem{MR0252639}
W.~K. Hayman.
\newblock Note on {H}adamard's convexity theorem.
\newblock In {\em Entire {F}unctions and {R}elatE. {P}arts of {A}nalysis
  ({P}roc. {S}ympos. {P}ure {M}ath., {L}a {J}olla, {C}alif., 1966)}, pages
  210--213. Amer. Math. Soc., Providence, R.I., 1968.

\bibitem{MR0219715}
W.~K. Hayman.
\newblock On the second {H}ankel determinant of mean univalent functions.
\newblock {\em Proc. London Math. Soc. (3)}, 18:77--94, 1968.

\bibitem{MR0244487}
W.~K. Hayman.
\newblock On integral functions with distinct asymptotic values.
\newblock {\em Proc. Cambridge Philos. Soc.}, 66:301--315, 1969.

\bibitem{MR0274761}
W.~K. Hayman.
\newblock Some examples related to the cos {$\pi \rho $} theorem.
\newblock In {\em Mathematical {E}ssays {D}edicatE. to {A}. {J}. {M}acintyre},
  pages 149--170. Ohio Univ. Press, Athens, Ohio, 1970.

\bibitem{Hayman1972}
W.~K. Hayman.
\newblock On the {V}aliron deficiencies of integral functions of infinite
  order.
\newblock {\em Arkiv f{\"o}r Matematik}, 10(1):163--172, Dec 1972.

\bibitem{MR0316693}
W.~K. Hayman.
\newblock Differential inequalities and local valency.
\newblock {\em Pacific J. Math.}, 44:117--137, 1973.

\bibitem{Bprevious}
W.~K. Hayman.
\newblock Research problems in function theory: progress on the previous
  problems.
\newblock pages 143--154. London Math. Soc. Lecture Note Ser., No. 12, 1974.

\bibitem{Hayman1975}
W.~K. Hayman.
\newblock The minimum modulus of integral functions of order one.
\newblock {\em J. Analyse Math.}, 28:171--212, 1975.

\bibitem{MR0492265}
W.~K. Hayman.
\newblock On {I}versen's theorem for meromorphic functions with few poles.
\newblock {\em Acta Math.}, 141(1-2):115--145, 1978.

\bibitem{MR581802}
W.~K. Hayman.
\newblock On a conjecture of {L}ittlewood.
\newblock {\em J. Analyse Math.}, 36:75--95 (1980), 1979.

\bibitem{MR568639}
W.~K. Hayman.
\newblock The logarithmic derivative of multivalent functions.
\newblock {\em Michigan Math. J.}, 27(2):149--179, 1980.

\bibitem{MR1438606}
W.~K. Hayman.
\newblock The growth of solutions of algebraic differential equations.
\newblock {\em Atti Accad. Naz. Lincei Cl. Sci. Fis. Mat. Natur. Rend. Lincei
  (9) Mat. Appl.}, 7(2):67--73, 1996.

\bibitem{MR877651}
W.~K. Hayman and J.~A. Hummel.
\newblock Coefficients of powers of univalent functions.
\newblock {\em Complex Variables Theory Appl.}, 7(1-3):51--70, 1986.

\bibitem{MR0460672}
W.~K. Hayman and P.~B. Kennedy.
\newblock {\em Subharmonic functions. {V}ol. {I}}.
\newblock Academic Press [Harcourt Brace Jovanovich, Publishers], London-New
  York, 1976.
\newblock London Mathematical Society Monographs, No. 9.

\bibitem{MR820533}
W.~K. Hayman, D.~Kershaw, and T.~J. Lyons.
\newblock The best harmonic approximant to a continuous function.
\newblock In {\em Anniversary volume on approximation theory and functional
  analysis ({O}berwolfach, 1983)}, volume~65 of {\em Internat. Schriftenreihe
  Numer. Math.}, pages 317--327. Birkh\"auser, Basel, 1984.

\bibitem{MR820231}
W.~K. Hayman and B.~Kjellberg.
\newblock On the minimum of a subharmonic function on a connected set.
\newblock In {\em Studies in pure mathematics}, pages 291--322. Birkh\"auser,
  Basel, 1983.

\bibitem{MR1040924}
W.~K. Hayman and J.~Miles.
\newblock On the growth of a meromorphic function and its derivatives.
\newblock {\em Complex Variables Theory Appl.}, 12(1-4):245--260, 1989.

\bibitem{MR0486496}
W.~K. Hayman, S.~J. Patterson, and Ch. Pommerenke.
\newblock On the coefficients of certain automorphic functions.
\newblock {\em Math. Proc. Cambridge Philos. Soc.}, 82(3):357--367, 1977.

\bibitem{MR3155684}
W.~K. Hayman, T.~F. Tyler, and D.~J. White.
\newblock The {B}lumenthal conjecture.
\newblock In {\em Complex analysis and dynamical systems {V}}, volume 591 of
  {\em Contemp. Math.}, pages 149--157. Amer. Math. Soc., Providence, RI, 2013.

\bibitem{MR513412}
W.~K. Hayman and I.~Vincze.
\newblock A problem on entire functions.
\newblock In {\em Complex analysis and its applications ({R}ussian)}, pages
  591--594, 672. ``Nauka'', Moscow, 1978.

\bibitem{MR0361049}
W.~K. Hayman and A.~Weitsman.
\newblock On the coefficients and means of functions omitting values.
\newblock {\em Math. Proc. Cambridge Philos. Soc.}, 77:119--137, 1975.

\bibitem{MR639358}
W.~K. Hayman and J.~M.~G. Wu.
\newblock Level sets of univalent functions.
\newblock {\em Comment. Math. Helv.}, 56(3):366--403, 1981.

\bibitem{zbMATH00194706}
W.K. {Hayman}.
\newblock {\em {Subharmonic functions. Volume 2.}}
\newblock London etc.: Academic Press, Inc., 1989.

\bibitem{MR0219729}
L.~I. Hedberg.
\newblock Weighted mean square approximation in plane regions, and generators
  of an algebra of analytic functions.
\newblock {\em Ark. Mat.}, 5:541--552 (1965), 1965.

\bibitem{MR2130416}
H.~Hedenmalm and S.~Shimorin.
\newblock Weighted {B}ergman spaces and the integral means spectrum of
  conformal mappings.
\newblock {\em Duke Math. J.}, 127(2):341--393, 2005.

\bibitem{MR0130974}
M.~Heins.
\newblock A class of conformal metrics.
\newblock {\em Bull. Amer. Math. Soc.}, 67:475--478, 1961.

\bibitem{MR0285724}
S.~Hellerstein.
\newblock Some analytic varieties in the polydisc and the {M}\"untz-{S}zasz
  problem in several variables.
\newblock {\em Trans. Amer. Math. Soc.}, 158:285--292, 1971.

\bibitem{MR0168761}
S.~Hellerstein and J.~Korevaar.
\newblock The real values of an entire function.
\newblock {\em Bull. Amer. Math. Soc.}, 70:608--610, 1964.

\bibitem{MR847009}
S.~Hellerstein and J.~Rossi.
\newblock Zeros of meromorphic solutions of second order linear differential
  equations.
\newblock {\em Math. Z.}, 192(4):603--612, 1986.

\bibitem{MR0393486}
S.~Hellerstein and D.~F. Shea.
\newblock An extremal problem concerning entire functions with radially
  distributed zeros.
\newblock pages 81--87. London Math. Soc. Lecture Note Ser., No. 12, 1974.

\bibitem{MR0247087}
S.~Hellerstein and J.~Williamson.
\newblock Entire functions with negative zeros and a problem of {R}.
  {N}evanlinna.
\newblock {\em J. Analyse Math.}, 22:233--267, 1969.

\bibitem{MR0435393}
S.~Hellerstein and J.~Williamson.
\newblock Derivatives of entire functions and a question of {P}\'olya.
\newblock {\em Trans. Amer. Math. Soc.}, 227:227--249, 1977.

\bibitem{MR623485}
J.~A. Hempel.
\newblock Precise bounds in the theorems of {L}andau and {S}chottky.
\newblock In {\em Aspects of contemporary complex analysis ({P}roc. {NATO}
  {A}dv. {S}tudy {I}nst., {U}niv. {D}urham, {D}urham, 1979)}, pages 421--424.
  Academic Press, London-New York, 1980.

\bibitem{zbMATH03406294}
W.~{Hengartner} and G.~{Schober}.
\newblock {A remark on level curves for domains convex in one direction.}
\newblock {\em {Appl. Anal.}}, 3:101--106, 1973.

\bibitem{MR0466634}
G.~M. Henkin.
\newblock Solutions with bounds for the equations of {H}. {L}ewy and
  {P}oincar\'e-{L}elong. {C}onstruction of functions of {N}evanlinna class with
  given zeros in a strongly pseudoconvex domain.
\newblock {\em Dokl. Akad. Nauk SSSR}, 224(4):771--774, 1975.

\bibitem{zbMATH03845559}
C.~W. {Henson} and L.~A. {{R}ubel}.
\newblock {Some applications of Nevanlinna theory to mathematical logic:
  Identities of exponential functions.}
\newblock {\em {Trans. Am. Math. Soc.}}, 282:1--32, 1984.

\bibitem{zbMATH03935013}
C.~W. {Henson} and L.~A. {{R}ubel}.
\newblock {Correction to ``Some applications of Nevanlinna theory to
  mathematical logic: identities of exponential functions''.}
\newblock {\em {Trans. Am. Math. Soc.}}, 294:381, 1986.

\bibitem{MR538680}
M.~R. Herman.
\newblock Sur la conjugaison diff\'erentiable des diff\'eomorphismes du cercle
  \`a des rotations.
\newblock {\em Inst. Hautes \'Etudes Sci. Publ. Math.}, (49):5--233, 1979.

\bibitem{MR771920}
M.~R. Herman.
\newblock Exemples de fractions rationnelles ayant une orbite dense sur la
  sph\`ere de {R}iemann.
\newblock {\em Bull. Soc. Math. France}, 112(1):93--142, 1984.

\bibitem{MR796014}
M.~R. Herman.
\newblock Are there critical points on the boundaries of singular domains?
\newblock {\em Comm. Math. Phys.}, 99(4):593--612, 1985.

\bibitem{MR915567}
M.~R. Herman.
\newblock Recent results and some open questions on {S}iegel's linearization
  theorem of germs of complex analytic diffeomorphisms of {$\mathbb{C}^n$} near
  a fixed point.
\newblock In {\em V{III}th international congress on mathematical physics
  ({M}arseille, 1986)}, pages 138--184. World Sci. Publishing, Singapore, 1987.

\bibitem{MR0257359}
F.~Herzog and G.~Piranian.
\newblock The counting function for points of maximum modulus.
\newblock In {\em Entire {F}unctions and {R}elatE. {P}arts of {A}nalysis
  ({P}roc. {S}ympos. {P}ure {M}ath., {L}a {J}olla, {C}alif., 1966)}, pages
  240--243. Amer. Math. Soc., Providence, R. I., 1968.

\bibitem{MR0447549}
C.~G. Higginson.
\newblock The asymptotic {B}ieberbach conjecture for weakly {$p$}-valent
  functions.
\newblock {\em Proc. London Math. Soc. (3)}, 35(2):291--312, 1977.

\bibitem{MR0030000}
E.~Hille.
\newblock Remarks on a paper by {Z}eev {N}ehari.
\newblock {\em Bull. Amer. Math. Soc.}, 55:552--553, 1949.

\bibitem{MR0430651}
F.~Holland and J.~B. Twomey.
\newblock Fourier-{S}tieltjes series of measures in {Z}ygmund's class {$\Lambda
  \sb\ast $}.
\newblock {\em Proc. Roy. Irish Acad. Sect. A}, 76(26):289--299, 1976.

\bibitem{MR0296305}
R.~Hornblower.
\newblock A growth condition for the {M}ac{L}ane class {${\cal A}$}.
\newblock {\em Proc. London Math. Soc. (3)}, 23:371--384, 1971.

\bibitem{MR0364642}
R.~J.~M. Hornblower.
\newblock Subharmonic analogues of {M}ac{L}ane's classes.
\newblock {\em Ann. Polon. Math.}, 26:135--146, 1972.

\bibitem{MR0357747}
C.~Horowitz.
\newblock Zeros of functions in the {B}ergman spaces.
\newblock {\em Duke Math. J.}, 41:693--710, 1974.

\bibitem{MR0480979}
D.~Horowitz.
\newblock A further refinement for coefficient estimates of univalent
  functions.
\newblock {\em Proc. Amer. Math. Soc.}, 71(2):217--221, 1978.

\bibitem{MR1215974}
J.~H. Hubbard.
\newblock Local connectivity of {J}ulia sets and bifurcation loci: three
  theorems of {J}. -{C}. {Y}occoz.
\newblock In {\em Topological methods in modern mathematics ({S}tony {B}rook,
  {NY}, 1991)}, pages 467--511. Publish or Perish, Houston, TX, 1993.

\bibitem{MR0094452}
A.~Huber.
\newblock On subharmonic functions and differential geometry in the large.
\newblock {\em Comment. Math. Helv.}, 32:13--72, 1957.

\bibitem{MR0425110}
A.~Huber.
\newblock Isometrische und {K}onforme {V}erheftung.
\newblock {\em Comment. Math. Helv.}, 51(3):319--331, 1976.

\bibitem{MR916071}
S.~Hui.
\newblock An extension of a theorem of {J}. {F}ern\'andez.
\newblock {\em Bull. London Math. Soc.}, 20(1):34--36, 1988.

\bibitem{MR0084546}
J.~A. Hummel.
\newblock Counterexamples to the {P}oincar\'e inequality.
\newblock {\em Proc. Amer. Math. Soc.}, 8:207--210, 1957.

\bibitem{MR0304629}
J.~A. Hummel.
\newblock The {M}arx conjecture for starlike functions.
\newblock {\em Michigan Math. J.}, 19:257--266, 1972.

\bibitem{MR726986}
J.~S. Hwang.
\newblock A problem on automorphic functions and gap series.
\newblock {\em Bull. Inst. Math. Acad. Sinica}, 11(3):401--406, 1983.

\bibitem{Hyllengren1970}
A.~Hyllengren.
\newblock Valiron deficient values for meromorphic functions in the plane.
\newblock {\em Acta Mathematica}, 124(1):1--8, Jul 1970.

\bibitem{igari1981}
S.~Igari.
\newblock Decomposition theorem and lacunary convergence of riesz-bochner means
  of fourier transforms of two variables.
\newblock {\em Tohoku Math. J. (2)}, 33(3):413--419, 1981.

\bibitem{Iversen}
F.~Iversen.
\newblock Recherches sur les fonctions inverses des fonctions m\'eromorphes.
\newblock {\em The\,se de Helsingfors}, 1914.

\bibitem{MR0281897}
I.~S. Jack.
\newblock Functions starlike and convex of order {$\alpha $}.
\newblock {\em J. London Math. Soc. (2)}, 3:469--474, 1971.

\bibitem{MR0182732}
S.~Jaenisch.
\newblock Absch\"atzungen subharmonischer und ganzer {F}unktionen in der
  {U}mgebung einer {K}urve.
\newblock {\em Mitt. Math. Sem. Giessen Heft}, 65:iv+50, 1965.

\bibitem{MR709720}
M.~Jakob and A.~C. Offord.
\newblock The distribution of the values of a random power series in the unit
  disk.
\newblock {\em Proc. Roy. Soc. Edinburgh Sect. A}, 94(3-4):251--263, 1983.

\bibitem{MR0315110}
W.~Janowski.
\newblock Some extremal problems for certain families of analytic functions.
  {I}.
\newblock {\em Bull. Acad. Polon. Sci. S\'er. Sci. Math. Astronom. Phys.},
  21:17--25, 1973.

\bibitem{MR0096806}
J.~A. Jenkins.
\newblock {\em Univalent functions and conformal mapping}.
\newblock Ergebnisse der Mathematik und ihrer Grenzgebiete. Neue Folge, Heft
  18. Reihe: Moderne Funktionentheorie. Springer-Verlag,
  Berlin-G\"ottingen-Heidelberg, 1958.

\bibitem{MR922301}
J.~A. Jenkins.
\newblock Some estimates for harmonic measures.
\newblock In {\em Complex analysis, {I} ({C}ollege {P}ark, {M}d., 1985--86)},
  volume 1275 of {\em Lecture Notes in Math.}, pages 210--214. Springer,
  Berlin, 1987.

\bibitem{MR1162090}
J.~A. Jenkins.
\newblock Some estimates for harmonic measures. {III}.
\newblock {\em Proc. Amer. Math. Soc.}, 119(1):199--201, 1993.

\bibitem{MR0422606}
F.~John.
\newblock A criterion for univalency brought up to date.
\newblock {\em Comm. Pure Appl. Math.}, 29(3):293--295, 1976.

\bibitem{MR597809}
P.~W. Jones.
\newblock Estimates for the corona problem.
\newblock {\em J. Funct. Anal.}, 39(2):162--181, 1980.

\bibitem{zbMATH04091817}
P.~W. {Jones} and T.~H. {Wolff}.
\newblock {Hausdorff dimension of harmonic measures in the plane.}
\newblock {\em {Acta Math.}}, 161(1-2):131--144, 1988.

\bibitem{MR3225918}
I.~G. Kasmalkar.
\newblock On the {S}endov conjecture for a root close to the unit circle.
\newblock {\em Aust. J. Math. Anal. Appl.}, 11(1):Art. 4, 34, 2014.

\bibitem{MR0205864}
Gy. Katona.
\newblock On a conjecture of {E}rd\"os and a stronger form of {S}perner's
  theorem.
\newblock {\em Studia Sci. Math. Hungar}, 1:59--63, 1966.

\bibitem{MR1501602}
O.~D. Kellogg.
\newblock On the derivatives of harmonic functions on the boundary.
\newblock {\em Trans. Amer. Math. Soc.}, 33(2):486--510, 1931.

\bibitem{MR0293114}
J.~T. Kemper.
\newblock A boundary {H}arnack principle for {L}ipschitz domains and the
  principle of positive singularities.
\newblock {\em Comm. Pure Appl. Math.}, 25:247--255, 1972.

\bibitem{MR0083033}
P.~B. Kennedy.
\newblock A class of integral functions bounded on certain curves.
\newblock {\em Proc. London Math. Soc. (3)}, 6:518--547, 1956.

\bibitem{MR3051172}
D.~Khavinson, R.~Pereira, M.~Putinar, E.~B. Saff, and S.~Shimorin.
\newblock Borcea's variance conjectures on the critical points of polynomials.
\newblock In {\em Notions of positivity and the geometry of polynomials},
  Trends Math., pages 283--309. Birkh\"auser/Springer Basel AG, Basel, 2011.

\bibitem{Kierst}
S.~Kierst.
\newblock Sur l'ensemble des valeurs asymptotiques d'une fonction m\'eromorphe
  dans le cercle-unit\'e.
\newblock {\em Fund. Math.}, 27:226--233, 1936.

\bibitem{MR0159942}
B.~Kjellberg.
\newblock A theorem on the minimum modulus of entire functions.
\newblock {\em Math. Scand.}, 12:5--11, 1963.

\bibitem{Kjellberg}
B.~Kjellberg.
\newblock The convexity theorm of {H}adamard-{H}ayman.
\newblock In {\em {P}roceedings of the {S}ymposium in {M}athematics, {R}oyal
  {I}nstitute of {T}echnology, {S}tockholm, (June 1973)}, pages 87--114. 1973.

\bibitem{MR0184865}
D.~J. Kleitman.
\newblock On a lemma of {L}ittlewood and {O}fford on the distribution of
  certain sums.
\newblock {\em Math. Z.}, 90:251--259, 1965.

\bibitem{kobayashi1979}
T.~Kobayashi.
\newblock An entire function with linearly distributed values.
\newblock {\em Kodai Math. J.}, 2(1):54--81, 1979.

\bibitem{MR998240}
L.~K\"ohler.
\newblock Meromorphic functions sharing zeros and poles and also some of their
  derivatives sharing zeros.
\newblock {\em Complex Variables Theory Appl.}, 11(1):39--48, 1989.

\bibitem{zbMATH03806164}
S.~V. {Konyagin}.
\newblock {On the problem of Littlewood.}
\newblock {\em {Math. USSR, Izv.}}, 18:205--225, 1982.

\bibitem{MR0425124}
B.~Korenblum.
\newblock An extension of the {N}evanlinna theory.
\newblock {\em Acta Math.}, 135(3-4):187--219, 1975.

\bibitem{MR0361030}
J.~Korevaar.
\newblock Equilibrium distributions of electrons on roundish plane conductors.
  {I}.
\newblock {\em Nederl. Akad. Wetensch. Proc. Ser. A {\bf 77}=Indag. Math.},
  36:423--437, 1974.

\bibitem{MR0361031}
J.~Korevaar.
\newblock Equilibrium distributions of electrons on roundish plane conductors.
  {II}.
\newblock {\em Nederl. Akad. Wetensch. Proc. Ser. A {\bf 77}=Indag. Math.},
  36:438--456, 1974.

\bibitem{MR623473}
J.~Korevaar.
\newblock Polynomial and rational approximation in the complex domain.
\newblock In {\em Aspects of contemporary complex analysis ({P}roc. {NATO}
  {A}dv. {S}tudy {I}nst., {U}niv. {D}urham, {D}urham, 1979)}, pages 251--292.
  Academic Press, London-New York, 1980.

\bibitem{MR515123}
J.~{Korevaar} and M.~{Dixon}.
\newblock {Nonspanning sets of exponentials on curves.}
\newblock {\em {Acta Math. Acad. Sci. Hung.}}, 33:89--100, 1979.

\bibitem{MR705427}
J.~Korevaar and R.~A. Kortram.
\newblock Equilibrium distributions of electrons on smooth plane conductors.
\newblock {\em Nederl. Akad. Wetensch. Indag. Math.}, 45(2):203--219, 1983.

\bibitem{MR0164005}
J.~Korevaar and T.~L. McCoy.
\newblock Power series whose partial sums have few zeros in an angle.
\newblock {\em J. Math. Anal. Appl.}, 8:461--473, 1964.

\bibitem{MR0379810}
J.~Korevaar and P.~Pfluger.
\newblock Spanning sets of powers on wild {J}ordan curves.
\newblock {\em Nederl. Akad. Wetensch. Proc. Ser. A {\bf 77}=Indag. Math.},
  36:293--305, 1974.

\bibitem{zbMATH03989737}
J.~{Korevaar} and R.~{Zeinstra}.
\newblock {Transform\'ees de Laplace pour les courbes \`a pente born\'ee et un
  r\'esultat correspondant du type M\"untz-Szasz. (Laplace transforms along
  curves of bounded slope and a corresponding M\"untz-Sz\'asz type result).}
\newblock {\em {C. R. Acad. Sci., Paris, S\'er. I}}, 301:695--698, 1985.

\bibitem{MR572106}
T.~W. K\"orner.
\newblock On a polynomial of {B}yrnes.
\newblock {\em Bull. London Math. Soc.}, 12(3):219--224, 1980.

\bibitem{MR553047}
L.~J. Kotman.
\newblock An entire function with irregular growth and more than one deficient
  value.
\newblock In {\em Complex analysis {J}oensuu 1978 ({P}roc. {C}olloq., {U}niv.
  {J}oensuu, {J}oensuu, 1978)}, volume 747 of {\em Lecture Notes in Math.},
  pages 219--229. Springer, Berlin, 1979.

\bibitem{MR0185122}
T.~K\"ovari.
\newblock Asymptotic values of entire functions of finite order with density
  conditions.
\newblock {\em Acta Sci. Math. (Szeged)}, 26:233--237, 1965.

\bibitem{MR0176076}
T.~K\"ovari.
\newblock A gap-theorem for entire functions of infinite order.
\newblock {\em Michigan Math. J.}, 12:133--140, 1965.

\bibitem{MR0188438}
T.~K\"ovari.
\newblock On the asymptotic paths of entire functions with gap power series.
\newblock {\em J. Analyse Math.}, 15:281--286, 1965.

\bibitem{MR0193235}
T.~K\"ovari.
\newblock On the growth of entire functions of finite order with density
  conditions.
\newblock {\em Quart. J. Math. Oxford Ser. (2)}, 17:22--30, 1966.

\bibitem{MR0231990}
T.~K\"ovari and Ch. Pommerenke.
\newblock On the distribution of {F}ekete points.
\newblock {\em Mathematika}, 15:70--75, 1968.

\bibitem{MR551458}
G.~K. Kristiansen.
\newblock Some inequalities for algebraic and trigonometric polynomials.
\newblock {\em J. London Math. Soc. (2)}, 20(2):300--314, 1979.

\bibitem{MR924982}
S.~L. Krushkal.
\newblock The interpolating family of a univalent analytic function.
\newblock {\em Sibirsk. Mat. Zh.}, 28(5):88--94, 1987.

\bibitem{MR0335807}
V.~I. Krutin.
\newblock The size of the {N}evanlinna deficiencies of functions meromorphic in
  {$z<1$}.
\newblock {\em Izv. Akad. Nauk Armjan. SSR Ser. Mat.}, 8(5):347--358, 425,
  1973.

\bibitem{MR511294}
Y.~H. Ku.
\newblock Sur les familles normales de fonctions m\'eromorphes.
\newblock {\em Sci. Sinica}, 21(4):431--445, 1978.

\bibitem{MR0369683}
Y.~Kubota.
\newblock A coefficient inequality for certain meromorphic univalent functions.
\newblock {\em K\=odai Math. Sem. Rep.}, 26:85--94, 1974/75.

\bibitem{MR0437735}
Y.~Kubota.
\newblock Coefficients of meromorphic univalent functions.
\newblock {\em K\=odai Math. Sem. Rep.}, 28(2--3):253--261, 1976/77.

\bibitem{MR0069893}
H.~P. K\"unzi.
\newblock Konstruktion {R}iemannscher {F}l\"achen mit vorgegebener {O}rdnung
  der erzeugenden {F}unktionen.
\newblock {\em Math. Ann.}, 128:471--474, 1955.

\bibitem{MR532327}
W.~T. Lai.
\newblock The exact value of {H}ayman's constant in {L}andau's theorem.
\newblock {\em Sci. Sinica}, 22(2):129--134, 1979.

\bibitem{MR1200194}
J.~K. Langley.
\newblock On the deficiencies of composite entire functions.
\newblock {\em Proc. Edinburgh Math. Soc. (2)}, 36(1):151--164, 1993.

\bibitem{MR1241784}
J.~K. Langley.
\newblock Proof of a conjecture of {H}ayman concerning {$f$} and {$f''$}.
\newblock {\em J. London Math. Soc. (2)}, 48(3):500--514, 1993.

\bibitem{MR1374462}
J.~K. Langley.
\newblock Two results related to a question of {H}inkkanen.
\newblock {\em Kodai Math. J.}, 19(1):52--61, 1996.

\bibitem{MR0107003}
O.~Lehto.
\newblock The spherical derivative of meromorphic functions in the
  neighbourhood of an isolated singularity.
\newblock {\em Comment. Math. Helv.}, 33:196--205, 1959.

\bibitem{MR545203}
F.~D. Lesley.
\newblock On interior and conformal mappings of the disk.
\newblock {\em J. London Math. Soc. (2)}, 20(1):67--78, 1979.

\bibitem{MR0480980}
Y.~Leung.
\newblock Successive coefficients of starlike functions.
\newblock {\em Bull. London Math. Soc.}, 10(2):193--196, 1978.

\bibitem{MR534395}
Y.~Leung.
\newblock Robertson's conjecture on the coefficients of close-to-convex
  functions.
\newblock {\em Proc. Amer. Math. Soc.}, 76(1):89--94, 1979.

\bibitem{MR554395}
Y.~J. Leung.
\newblock Integral means of the derivatives of some univalent functions.
\newblock {\em Bull. London Math. Soc.}, 11(3):289--294, 1979.

\bibitem{MR0130979}
B.~Ja. Levin and I.~V. Ostrovskii.
\newblock The dependence of the growth of an entire function on the
  distribution of zeros of its derivatives.
\newblock {\em Sibirsk. Mat. \v{Z}.}, 1:427--455, 1960.

\bibitem{MR1637647}
G.~Levin and S.~van Strien.
\newblock Local connectivity of the {J}ulia set of real polynomials.
\newblock {\em Ann. of Math. (2)}, 147(3):471--541, 1998.

\bibitem{MR1545451}
V.~Levin.
\newblock Ein {B}eitrag zum {K}oeffizientenproblem der schlichten {F}unktionen.
\newblock {\em Math. Z.}, 38(1):306--311, 1934.

\bibitem{MR1501941}
N.~Levinson.
\newblock On the growth of analytic functions.
\newblock {\em Trans. Amer. Math. Soc.}, 43(2):240--257, 1938.

\bibitem{MR0206255}
M.~Lewin.
\newblock On a coefficient problem for bi-univalent functions.
\newblock {\em Proc. Amer. Math. Soc.}, 18:63--68, 1967.

\bibitem{MR735882}
J.~Lewis, J.~Rossi, and A.~Weitsman.
\newblock On the growth of subharmonic functions along paths.
\newblock {\em Ark. Mat.}, 22(1):109--119, 1984.

\bibitem{MR0499158}
J.~L. Lewis.
\newblock Note on the {N}evannlinna proximity function.
\newblock {\em Proc. Amer. Math. Soc.}, 69(1):129--134, 1978.

\bibitem{MR942832}
J.~L. Lewis and J.~M.~G. Wu.
\newblock On conjectures of {A}rakelyan and {L}ittlewood.
\newblock {\em J. Analyse Math.}, 50:259--283, 1988.

\bibitem{MR0437725}
C.~N. Linden.
\newblock The modulus of polynomials with zeros on the unit circle.
\newblock {\em Bull. London Math. Soc.}, 9(1):65--69, 1977.

\bibitem{MR0229836}
F.~S. Lisin.
\newblock A problem of mean square approximation related to the study of
  generators in the algebra {$l^{1}$}.
\newblock {\em Mat. Zametki}, 3:703--706, 1968.

\bibitem{MR1574417}
J.~E. Littlewood.
\newblock Mathematical {N}otes (4): {O}n a {T}heorem of {F}atou.
\newblock {\em J. London Math. Soc.}, S1-2(3):172, 1927.

\bibitem{Littlewood2}
J.~E. Littlewood.
\newblock On the coefficients of schlicht functions.
\newblock {\em Quart. J. Math.}, Oxford(9):14--20, 1938.

\bibitem{MR0012121}
J.~E. Littlewood.
\newblock {\em Lectures on the {T}heory of {F}unctions}.
\newblock Oxford University Press, 1944.

\bibitem{MR0049315}
J.~E. Littlewood.
\newblock On some conjectural inequalities, with applications to the theory of
  integral functions.
\newblock {\em J. London Math. Soc.}, 27:387--393, 1952.

\bibitem{MR0009656}
J.~E. Littlewood and A.~C. Offord.
\newblock On the number of real roots of a random algebraic equation. {III}.
\newblock {\em Rec. Math. [Mat. Sbornik] N.S.}, 12(54):277--286, 1943.

\bibitem{MR0243054}
A.~E. Livingston.
\newblock The coefficients of multivalent close-to-convex functions.
\newblock {\em Proc. Amer. Math. Soc.}, 21:545--552, 1969.

\bibitem{MR0072218}
A.~J. Lohwater, G.~Piranian, and W.~Rudin.
\newblock The derivative of a schlicht function.
\newblock {\em Math. Scand.}, 3:103--106, 1955.

\bibitem{MR0444948}
R.~R. London.
\newblock A note on {H}adamard's three circles theorem.
\newblock {\em Bull. London Math. Soc.}, 9(2):182--185, 1977.

\bibitem{MR2632653}
E.~L\o{}w.
\newblock {\em Inner functions and boundary values in ${H}^\infty(\Omega)$ and
  ${A}(\Omega)$ in smoothly bounded pseudoconvex domains}.
\newblock ProQuest LLC, Ann Arbor, MI, 1983.
\newblock Doctoral Thesis--Princeton University.

\bibitem{MR731340}
E.~L\o{}w.
\newblock Inner functions and boundary values in {$H^{\infty }(\Omega )$}\ and
  {$A(\Omega )$}\ in smoothly bounded pseudoconvex domains.
\newblock {\em Math. Z.}, 185(2):191--210, 1984.

\bibitem{Lowner1923}
K.~L\"owner.
\newblock Untersuchungen {\"{u}}ber schlichte konforme {A}bbildungen des
  {E}inheitskreises. i.
\newblock {\em Mathematische Annalen}, 89:103--121, 1923.

\bibitem{lucasnew}
K.~W. Lucas.
\newblock A two-point modulus bound for areally mean p-valent functions.
\newblock {\em Journal of the London Mathematical Society}, s1-43(1):487--494.

\bibitem{MR576642}
T.~J. Lyons.
\newblock Finely holomorphic functions.
\newblock {\em J. Funct. Anal.}, 37(1):1--18, 1980.

\bibitem{MR0576643}
T.~J. Lyons.
\newblock A theorem in fine potential theory and applications to finely
  holomorphic functions.
\newblock {\em J. Funct. Anal.}, 37(1):19--26, 1980.

\bibitem{MYuLyubich1986Tdor}
M.~Yu. Lyubich.
\newblock The dynamics of rational transforms: the topological picture.
\newblock {\em {R}ussian Mathematical Surveys}, 41(4):43--117, August 1986.

\bibitem{MR1459261}
M.~Yu. Lyubich.
\newblock Dynamics of quadratic polynomials. {I}, {II}.
\newblock {\em Acta Math.}, 178(2):185--247, 247--297, 1997.

\bibitem{MR1755441}
M.~Yu. Lyubich.
\newblock Dynamics of quadratic polynomials. {III}. {P}arapuzzle and {SBR}
  measures.
\newblock {\em Ast\'erisque}, (261):xii--xiii, 173--200, 2000.
\newblock G\'eom\'etrie complexe et syst\`emes dynamiques (Orsay, 1995).

\bibitem{MR807068}
R.~Ma\~n\'e.
\newblock On the instability of {H}erman rings.
\newblock {\em Invent. Math.}, 81(3):459--471, 1985.

\bibitem{MR0222272}
T.~H. MacGregor.
\newblock Certain integrals of univalent and convex functions.
\newblock {\em Math. Z.}, 103:48--54, 1968.

\bibitem{MR0236369}
T.~H. MacGregor.
\newblock The univalence of a linear combination of convex mappings.
\newblock {\em J. London Math. Soc.}, 44:210--212, 1969.

\bibitem{MR0367175}
T.~H. MacGregor.
\newblock A subordination for convex functions of order {$\alpha $}.
\newblock {\em J. London Math. Soc. (2)}, 9:530--536, 1974/75.

\bibitem{MR0050000}
A.~J. Macintyre.
\newblock Asymptotic paths of integral functions with gap power series.
\newblock {\em Proc. London Math. Soc. (3)}, 2:286--296, 1952.

\bibitem{MR0002823}
A.~J. Macintyre and W.~H.~J. Fuchs.
\newblock Inequalities for the logarithmic derivatives of a polynomial.
\newblock {\em J. London Math. Soc.}, 15:162--168, 1940.

\bibitem{MR0053231}
G.~R. MacLane.
\newblock Sequences of derivatives and normal families.
\newblock {\em J. Analyse Math.}, 2:72--87, 1952.

\bibitem{MR0148923}
G.~R. MacLane.
\newblock Asymptotic values of holomorphic functions.
\newblock {\em Rice Univ. Studies}, 49(1):83, 1963.

\bibitem{MR0274765}
G.~R. MacLane.
\newblock Exceptional values of {$f^{(n)}\,(z)$}, asymptotic values of
  {$f(z)$}, and linearly accessible asymptotic values.
\newblock In {\em Mathematical {E}ssays {D}edicatE. to {A}. {J}. {M}acintyre},
  pages 271--288. Ohio Univ. Press, Athens, Ohio, 1970.

\bibitem{MR0166418}
E.~Makai and P.~Tur\'an.
\newblock Hermite expansion and distribution of zeros of polynomials.
\newblock {\em Magyar Tud. Akad. Mat. Kutat\'o Int. K\"ozl.}, 8:157--163, 1963.

\bibitem{zbMATH03914673}
N.~G. {Makarov}.
\newblock {On the distortion of boundary sets under conformal mappings.}
\newblock {\em {Proc. Lond. Math. Soc. (3)}}, 51:369--384, 1985.

\bibitem{MR818450}
N.~G. Makarov.
\newblock A note on integral means of the derivative in conformal mapping.
\newblock {\em Proc. Amer. Math. Soc.}, 96(2):233--235, 1986.

\bibitem{MR1555170}
J.~Malmquist.
\newblock Sur les fonctions \`a un nombre fini de branches satisfaisant \`a une
  \'equation diff\'erentielle du premier ordre.
\newblock {\em Acta Math.}, 42(1):317--325, 1920.

\bibitem{MR0457732}
J.~H. Mantel.
\newblock Conditions for generating a nonvanishing bounded analytic function.
\newblock {\em Proc. Amer. Math. Soc.}, 66(1):62--64, 1977.

\bibitem{MR948195}
D.~E. Marshall and C.~Sundberg.
\newblock Harmonic measure and radial projection.
\newblock {\em Trans. Amer. Math. Soc.}, 316(1):81--95, 1989.

\bibitem{MR1512785}
A.~Marx.
\newblock Untersuchungen \"uber schlichte {A}bbildungen.
\newblock {\em Math. Ann.}, 107(1):40--67, 1933.

\bibitem{zbMATH03741921}
O.~C. {McGehee}, L.~{Pigno}, and B.~{Smith}.
\newblock {Hardy's inequality and the $L\sp 1$ norm of exponential sums.}
\newblock {\em {Ann. Math. (2)}}, 113:613--618, 1981.

\bibitem{MR0294614}
J.~E. McMillan.
\newblock On the boundary correspondence under conformal mapping.
\newblock {\em Duke Math. J.}, 37:725--739, 1970.

\bibitem{MR1312365}
C.~T. McMullen.
\newblock {\em Complex dynamics and renormalization}, volume 135 of {\em Annals
  of Mathematics Studies}.
\newblock Princeton University Press, Princeton, NJ, 1994.

\bibitem{MR0410393}
W.~S. McVoy and L.~A. {R}ubel.
\newblock Coherence of some rings of functions.
\newblock {\em J. Functional Analysis}, 21(1):76--87, 1976.

\bibitem{MR0249587}
A.~Meir and A.~Sharma.
\newblock On {I}lyeff's conjecture.
\newblock {\em Pacific J. Math.}, 31:459--467, 1969.

\bibitem{MR1897399}
A.~D. Melas.
\newblock On the centered {H}ardy-{L}ittlewood maximal operator.
\newblock {\em Trans. Amer. Math. Soc.}, 354(8):3263--3273, 2002.

\bibitem{MR1973058}
A.~D. Melas.
\newblock The best constant for the centered {H}ardy-{L}ittlewood maximal
  inequality.
\newblock {\em Ann. of Math. (2)}, 157(2):647--688, 2003.

\bibitem{MR2011925}
A.~D. Melas.
\newblock On a covering problem related to the centered {H}ardy-{L}ittlewood
  maximal inequality.
\newblock {\em Ark. Mat.}, 41(2):341--361, 2003.

\bibitem{MR806014}
A.~Meril and D.~C. Struppa.
\newblock Equivalence of {C}auchy problems for entire and exponential type
  functions.
\newblock {\em Bull. London Math. Soc.}, 17(5):469--473, 1985.

\bibitem{MR0417568}
N.~G. Meyers and A.~Elcrat.
\newblock Some results on regularity for solutions of non-linear elliptic
  systems and quasi-regular functions.
\newblock {\em Duke Math. J.}, 42:121--136, 1975.

\bibitem{MR0254245}
J.~Miles.
\newblock On the counting function for the {$a$}-values of a meromorphic
  function.
\newblock {\em Trans. Amer. Math. Soc.}, 147:203--222, 1970.

\bibitem{Milesjoseph}
J.~Miles.
\newblock A factorization theorem in ${H}^1({U}^3)$.
\newblock {\em Proc. Amer. Math. Soc.}, 52:319--322, 1975.

\bibitem{MR543739}
J.~Miles.
\newblock On entire functions of infinite order with radially distributed
  zeros.
\newblock {\em Pacific J. Math.}, 81(1):131--157, 1979.

\bibitem{MR698931}
J.~Miles.
\newblock Some examples of the dependence of the {N}evanlinna deficiency upon
  the choice of origin.
\newblock {\em Proc. London Math. Soc. (3)}, 47(1):145--176, 1983.

\bibitem{doi:10.1080/17476939808815143}
J.~Miles.
\newblock On a theorem of {H}ayman and {S}tewart.
\newblock {\em Complex Variables, Theory and Application: An International
  Journal}, 37(1-4):425--455, 1998.

\bibitem{MR0324041}
J.~Miles and D.~F. Shea.
\newblock An extremal problem in value-distribution theory.
\newblock {\em Quart. J. Math. Oxford Ser. (2)}, 24:377--383, 1973.

\bibitem{MR0492261}
J.~Miles and D.~Townsend.
\newblock Imaginary values of meromorphic functions.
\newblock {\em Indiana Univ. Math. J.}, 27(3):491--503, 1978.

\bibitem{MR829392}
J.~Miles and J.~Williamson.
\newblock A characterization of the exponential function.
\newblock {\em J. London Math. Soc. (2)}, 33(1):110--116, 1986.

\bibitem{MR0158056}
I.~M. Milin.
\newblock The area method in the theory of univalent functions.
\newblock {\em Dokl. Akad. Nauk SSSR}, 154:264--267, 1964.

\bibitem{MR0172991}
I.~M. Milin.
\newblock A bound for the coefficients of schlicht functions.
\newblock {\em Dokl. Akad. Nauk SSSR}, 160:769--771, 1965.

\bibitem{MR0231992}
I.~M. Milin.
\newblock Adjacent coefficients of univalent functions.
\newblock {\em Dokl. Akad. Nauk SSSR}, 180:1294--1297, 1968.

\bibitem{MR0369684}
I.~M. Milin.
\newblock {\em {O}dnolistnye funktsii i ortonormirovannye sistemy}.
\newblock Izdat. ``Nauka'', Moscow, 1971.

\bibitem{MR0450541}
I.~M. Milin.
\newblock Errata: ``\textit{{U}nivalent functions and orthonormal systems}''
  ({T}rans. {M}ath. {M}onographs, {V}ol. 49, {A}mer. {M}ath. {S}oc.,
  {P}rovidence, {R}.{I}., 1977).
\newblock {\em American Mathematical Society, Providence, R.I.,}, page~1, 1977.

\bibitem{MR0427620}
I.~M. Milin.
\newblock {\em Univalent functions and orthonormal systems}.
\newblock American Mathematical Society, Providence, R. I., 1977.
\newblock Translated from the {R}ussian, Translations of Mathematical
  Monographs, Vol. 49.

\bibitem{MR804182}
I.~M. Milin.
\newblock A problem for coefficients of {$p$}-multiply symmetric univalent
  functions.
\newblock {\em Mat. Zametki}, 38(1):66--73, 169, 1985.

\bibitem{MR506307}
S.~S. Miller and Petru~T. Mocanu.
\newblock Second-order differential inequalities in the complex plane.
\newblock {\em J. Math. Anal. Appl.}, 65(2):289--305, 1978.

\bibitem{zbMATH06257274}
J.~{Milnor}.
\newblock {\em {Dynamics in one complex variable. Reprint of the 2006 3rd
  edition published by Princeton University Press.}}, volume~62.
\newblock New Delhi: Hindustan Book Agency, reprint of the 2006 3rd edition
  published by princeton university press edition, 2012.

\bibitem{MR1505035}
C.~Miranda.
\newblock Sur un nouveau crit\`ere de normalit\'e pour les familles de
  fonctions holomorphes.
\newblock {\em Bull. Soc. Math. France}, 63:185--196, 1935.

\bibitem{MR627790}
M.~Misiurewicz.
\newblock On iterates of {$e^{z}$}.
\newblock {\em Ergodic Theory Dynamical Systems}, 1(1):103--106, 1981.

\bibitem{Moguodan}
G.~{Mo}.
\newblock {An improvement upon Cartan's theorem.}
\newblock {\em {Acta Math. Sin.}}, 25:287--296, 1982.

\bibitem{Molluzzo}
J.~Molluzzo.
\newblock Doctoral thesis.
\newblock {\em Yeshiva University}, 1972.

\bibitem{Montel}
P.~Montel.
\newblock Le\c cons sur les familles normales de fonctions analytiques et leur
  applications.
\newblock {\em Paris}, 1927.

\bibitem{MR0276471}
E.~Mues.
\newblock \"{U}ber eine {V}ermutung von {H}ayman.
\newblock {\em Math. Z.}, 119:11--20, 1971.

\bibitem{MR516609}
E.~Mues.
\newblock \"{U}ber ein {P}roblem von {H}ayman.
\newblock {\em Math. Z.}, 164(3):239--259, 1979.

\bibitem{zbMATH03682094}
T.~{Murai}.
\newblock {The value-distribution of lacunary series and a conjecture of
  Paley.}
\newblock {\em {Ann. Inst. Fourier}}, 31(1):135--156, 1981.

\bibitem{MR692343}
T.~Murai.
\newblock The boundary behaviour of {H}adamard lacunary series.
\newblock {\em Nagoya Math. J.}, 89:65--76, 1983.

\bibitem{MR723947}
T.~Murai.
\newblock The deficiency of entire functions with {F}ej\'er gaps.
\newblock {\em Ann. Inst. Fourier (Grenoble)}, 33(3):39--58, 1983.

\bibitem{Nazarov}
F.~Nazarov.
\newblock Growth of entire functions with sparse spectra.
\newblock {\em Unpublished, posted on
  http://users.math.msu.edu/users/fedja/prepr.html}.

\bibitem{zbMATH06373401}
F.~{Nazarov}, A.~{Nishry}, and M.~{Sodin}.
\newblock {Log-integrability of Rademacher Fourier series, with applications to
  random analytic functions.}
\newblock {\em {St. Petersbg. Math. J.}}, 25(3):467--494, 2014.

\bibitem{MR0083569}
Z.~Nehari.
\newblock On the coefficients of univalent functions.
\newblock {\em Proc. Amer. Math. Soc.}, 8:291--293, 1957.

\bibitem{FNevanlinna1929}
F.~Nevanlinna.
\newblock {\"{U}}ber eine {K}lasse meromorpher {F}unktionen.
\newblock {\em C. R. 7e Congr. Math. Scand., Oslo}, pages 81--3, 1929.

\bibitem{Newvanlinna1932}
R.~Nevanlinna.
\newblock {\"{U}}ber {R}iemannsche fl{\"a}chen mit endlich vielen
  {W}indungspunkten.
\newblock {\em Acta Mathematica}, 58(1):295--373, Dec 1932.

\bibitem{MR0344426}
R.~Nevanlinna.
\newblock {\em Eindeutige analytische {F}unktionen}.
\newblock Springer-Verlag, Berlin-New York, 1974.
\newblock Zweite Auflage, Reprint, Die Grundlehren der mathematischen
  Wissenschaften, Band 46.

\bibitem{MR0417418}
R.~Nevanlinna.
\newblock {\em Le th\'eor\`eme de {P}icard-{B}orel et la th\'eorie des
  fonctions m\'eromorphes}.
\newblock Chelsea Publishing Co., New York, 1974.
\newblock Reprinting of the 1929 original.

\bibitem{MR0170206}
D.~J. Newman.
\newblock Generators in {$l_{1}$}.
\newblock {\em Trans. Amer. Math. Soc.}, 113:393--396, 1964.

\bibitem{MR544895}
D.~J. Newman and M.~Slater.
\newblock Waring's problem for the ring of polynomials.
\newblock {\em J. Number Theory}, 11(4):477--487, 1979.

\bibitem{Watanabe}
X.~L. {Nguyen} and T.~{Watanabe}.
\newblock {A characterization of fine domains for a certain class of Markov
  processes with applications to Brelot harmonic spaces.}
\newblock {\em {Z. Wahrscheinlichkeitstheor. Verw. Geb.}}, 21:167--178, 1972.

\bibitem{MR620424}
P.~J. Nicholls and L.~R. Sons.
\newblock Automorphic functions with gap power series.
\newblock {\em Illinois J. Math.}, 25(3):383--389, 1981.

\bibitem{MR3008885}
D.~A. Nicks.
\newblock Wandering domains in quasiregular dynamics.
\newblock {\em Proc. Amer. Math. Soc.}, 141(4):1385--1392, 2013.

\bibitem{MR1545320}
W.~Nikliborc.
\newblock Eine {B}emerkung \"uber die {V}olumpotentiale {I}.
\newblock {\em Math. Z.}, 35(1):625--631, 1932.

\bibitem{Nikliborc}
W.~Nikliborc.
\newblock Eine {B}emerkung \"uber die {V}olumpotentiale {II}.
\newblock {\em Math. Z.}, 36(1):167--170, 1932.

\bibitem{NOONAN1972721}
J.~W Noonan.
\newblock Boundary behavior of functions with bounded boundary rotation.
\newblock {\em Journal of Mathematical Analysis and Applications}, 38(3):721 --
  734, 1972.

\bibitem{MR0306479}
J.~W. Noonan and D.~K. Thomas.
\newblock On the {H}ankel determinants of areally mean {$p$}-valent functions.
\newblock {\em Proc. London Math. Soc. (3)}, 25:503--524, 1972.

\bibitem{o'connell1973}
J.~M.~F. O'Connell.
\newblock Real parts of uniform algebras.
\newblock {\em Pacific J. Math.}, 46(1):235--247, 1973.

\bibitem{MR0307393}
A.~C. Offord.
\newblock The distribution of the values of a random function in the unit disk.
\newblock {\em Studia Math.}, 41:71--106, 1972.

\bibitem{MR0125956}
K.~Oikawa.
\newblock Welding of polygons and the type of {R}iemann surfaces.
\newblock {\em K\=odai Math. Sem. Rep.}, 13:37--52, 1961.

\bibitem{MR501158}
K.~Okabe.
\newblock On some open problems in the theory of quasiconformal mappings. {I}.
  {T}he {$L^{p}$}-integrability of the partial derivatives of a quasiconformal
  mapping.
\newblock {\em J. Hokkaido Univ. Ed. Sect. II A}, 28(2):65--69, 1977/78.

\bibitem{MR561991}
K.~Okabe.
\newblock On some open problems in the theory of quasiconformal mappings. {II}.
  {A}rea distortion under quasiconformal mappings.
\newblock {\em J. Hokkaido Univ. Ed. Sect. II A}, 30(1):5--12, 1979/80.

\bibitem{MR581203}
K.~Okabe.
\newblock On some open problems in the theory of quasiconformal mappings.
  {III}. {E}xtension of quasiconformal mappings from {$n$}\ to {$(n+1)$}\
  dimensions.
\newblock {\em J. Hokkaido Univ. Ed. Sect. II A}, 30(2):109--122, 1979/80.

\bibitem{MR568976}
B.~\O{}ksendal.
\newblock A {W}iener test for integrals of {B}rownian motion and the existence
  of smooth curves in nowhere dense sets.
\newblock {\em J. Funct. Anal.}, 36(1):72--87, 1980.

\bibitem{MR631668}
B.~\O{}ksendal.
\newblock Brownian motion and sets of harmonic measure zero.
\newblock {\em Pacific J. Math.}, 95(1):179--192, 1981.

\bibitem{Oliver}
P.~Oliver.
\newblock Doctoral thesis.
\newblock {\em London University}, 1975.

\bibitem{MR813642}
M.~Ortel and W.~Schneider.
\newblock Radial limits of functions of slow growth in the unit disk.
\newblock {\em Math. Scand.}, 56(2):287--310, 1985.

\bibitem{MR650781}
I.~B. Oshkin.
\newblock On a condition for the normality of families of holomorphic
  functions.
\newblock {\em Uspekhi Mat. Nauk}, 37(2(224)):221--222, 1982.

\bibitem{zbMATH03638660}
R.~{Osserman}.
\newblock {The isoperimetric inequality.}
\newblock {\em {Bull. Am. Math. Soc.}}, 84:1182--1238, 1978.

\bibitem{MR1749956}
M.~Overholt.
\newblock Sets of uniqueness for univalent functions.
\newblock {\em Canad. Math. Bull.}, 43(1):105--107, 2000.

\bibitem{MR0239073}
M.~Ozawa.
\newblock On the {B}ieberbach conjecture for the sixth coefficient.
\newblock {\em K\=odai Math. Sem. Rep.}, 21:97--128, 1969.

\bibitem{MR1057999}
X.~C. Pang.
\newblock Bloch's principle and normal criterion.
\newblock {\em Sci. China Ser. A}, 32(7):782--791, 1989.

\bibitem{zbMATH04159040}
X.~C. {Pang}.
\newblock {On normal criterion of meromorphic functions.}
\newblock {\em {Sci. China, Ser. A}}, 33(5):521--527, 1990.

\bibitem{MR0239072}
R.~N. Pederson.
\newblock A proof of the {B}ieberbach conjecture for the sixth coefficient.
\newblock {\em Arch. Rational Mech. Anal.}, 31:331--351, 1968/1969.

\bibitem{MR0311891}
R.~N. Pederson and M.~Schiffer.
\newblock A proof of the {B}ieberbach conjecture for the fifth coefficient.
\newblock {\em Arch. Rational Mech. Anal.}, 45:161--193, 1972.

\bibitem{MR1360570}
R.~P\'erez-Marco.
\newblock Sur une question de {D}ulac et {F}atou.
\newblock {\em C. R. Acad. Sci. Paris S\'er. I Math.}, 321(8):1045--1048, 1995.

\bibitem{MR1607557}
R.~P\'erez-Marco.
\newblock Fixed points and circle maps.
\newblock {\em Acta Math.}, 179(2):243--294, 1997.

\bibitem{Peterson}
K.~E. {Petersen}.
\newblock {\em {Brownian motion, Hardy spaces and bounded mean oscillation.}},
  volume~28.
\newblock Cambridge University Press, Cambridge. London Mathematical Society,
  London, 1977.

\bibitem{MR0257363}
V.~P. Petrenko.
\newblock Investigation of the structure of the set of positive deviations of
  meromorphic functions. {I}.
\newblock {\em Izv. Akad. Nauk SSSR Ser. Mat.}, 33:1330--1348, 1969.

\bibitem{MR0262473}
G.~Petruska.
\newblock A contribution to {B}loch's theorem.
\newblock {\em Ann. Univ. Sci. Budapest. E\"otv\"os Sect. Math.}, 12:39--42,
  1969.

\bibitem{MR0015488}
A.~Pfluger.
\newblock \"{U}ber ganze {F}unktionen ganzer {O}rdnung.
\newblock {\em Comment. Math. Helv.}, 18:177--203, 1946.

\bibitem{Pfluger1946}
A.~Pfluger.
\newblock {Z}ur {D}efektrelation ganzer {F}unktionen endlicher {O}rdnung.
\newblock {\em Commentarii Mathematici Helvetici}, 19(1):91--104, Dec 1946.

\bibitem{Picard1880}
E.~{Picard}.
\newblock {Sur une propri\'et\'e de certaines fonctions analogues aux fonctions
  alg\'ebriques.}
\newblock {\em {C. R. Acad. Sci., Paris}}, 89:1106--1108, 1880.

\bibitem{MR0111838}
G.~Piranian.
\newblock An entire function of restricted growth.
\newblock {\em Comment. Math. Helv.}, 33:322--324, 1959.

\bibitem{MR1581142}
A.~Plessner.
\newblock \"{U}ber das {V}erhalten analytischer {F}unktionen am {R}ande ihres
  {D}efinitionsbereiches.
\newblock {\em J. Reine Angew. Math.}, 158:219--227, 1927.

\bibitem{MR1544505}
G.~P\'olya.
\newblock \"{U}ber die {N}ullstellen sukzessiver {D}erivierten.
\newblock {\em Math. Z.}, 12(1):36--60, 1922.

\bibitem{MR0100753}
G.~P\'olya and I.~J. Schoenberg.
\newblock Remarks on de la {V}all\'ee {P}oussin means and convex conformal maps
  of the circle.
\newblock {\em Pacific J. Math.}, 8:295--334, 1958.

\bibitem{MR1492448}
G.~{P\'olya} and G.~{Szeg\"{o}}.
\newblock {\em {Problems and theorems in analysis II. Theory of functions,
  zeros, polynomials, determinants, number theory, geometry. Transl. from the
  German by C. E. Billigheimer. Reprint of the 1976 {E}nglish translation.}}
\newblock Berlin: Springer, reprint of the 1976 {E}nglish translation edition,
  1998.

\bibitem{MR0109875}
Ch. Pommerenke.
\newblock On some problems by {E}rd\"os, {H}erzog and {P}iranian.
\newblock {\em Michigan Math. J.}, 6:221--225, 1959.

\bibitem{MR0109208}
Ch. Pommerenke.
\newblock On the derivative of a polynomial.
\newblock {\em Michigan Math. J}, 6:373--375, 1959.

\bibitem{MR0120481}
Ch. Pommerenke.
\newblock Einige {S}\"atze \"uber die {K}apazit\"at ebener {M}engen.
\newblock {\em Math. Ann.}, 141:143--152, 1960.

\bibitem{zbMATH03153281}
Ch. {Pommerenke}.
\newblock {\"Uber die analytische Kapazit\"at.}
\newblock {\em {Arch. Math.}}, 11:270--277, 1960.

\bibitem{MR0151580}
Ch. Pommerenke.
\newblock On metric properties of complex polynomials.
\newblock {\em Michigan Math. J.}, 8:97--115, 1961.

\bibitem{MR0133448}
Ch. Pommerenke.
\newblock \"{U}ber die {M}ittelwerte und {K}oeffizienten multivalenter
  {F}unktionen.
\newblock {\em Math. Ann.}, 145:285--296, 1961/1962.

\bibitem{MR0158058}
Ch. Pommerenke.
\newblock \"{U}ber einige {K}lassen meromorpher schlichter {F}unktionen.
\newblock {\em Math. Z.}, 78:263--284, 1962.

\bibitem{MR0167611}
Ch. Pommerenke.
\newblock Lacunary power series and univalent functions.
\newblock {\em Michigan Math. J.}, 11:219--223, 1964.

\bibitem{MR0168772}
Ch. Pommerenke.
\newblock \"{U}ber die {F}aberschen {P}olynome schlichter {F}unktionen.
\newblock {\em Math. Z.}, 85:197--208, 1964.

\bibitem{MR0185105}
Ch. Pommerenke.
\newblock On the coefficients and {H}ankel determinants of univalent functions.
\newblock {\em J. London Math. Soc.}, 41:111--122, 1966.

\bibitem{MR0215976}
Ch. Pommerenke.
\newblock On the {H}ankel determinants of univalent functions.
\newblock {\em Mathematika}, 14:108--112, 1967.

\bibitem{MR0220920}
Ch. Pommerenke.
\newblock Relations between the coefficients of a univalent function.
\newblock {\em Invent. Math.}, 3:1--15, 1967.

\bibitem{MR0325976}
Ch. {Pommerenke}.
\newblock {Normal functions.}
\newblock {Proc. NRL Conf. classical Function Theory, Washington, D.C., 77-93
  (1970).}, 1970.

\bibitem{MR0284574}
Ch. Pommerenke.
\newblock On {B}loch functions.
\newblock {\em J. London Math. Soc. (2)}, 2:689--695, 1970.

\bibitem{MR0313496}
Ch. Pommerenke.
\newblock On the growth of the coefficients of analytic functions.
\newblock {\em J. London Math. Soc. (2)}, 5:624--628, 1972.

\bibitem{MR0357779}
Ch. Pommerenke.
\newblock On normal and automorphic functions.
\newblock {\em Michigan Math. J.}, 21:193--202, 1974.

\bibitem{MR0507768}
Ch. Pommerenke.
\newblock {\em Univalent functions}.
\newblock Vandenhoeck \& Ruprecht, G\"ottingen, 1975.
\newblock With a chapter on quadratic differentials by Gerd Jensen, Studia
  Mathematica/Mathematische Lehrb\"ucher, Band XXV.

\bibitem{MR620121}
Ch. Pommerenke.
\newblock On ergodic properties of inner functions.
\newblock {\em Math. Ann.}, 256(1):43--50, 1981.

\bibitem{MR813740}
Ch. Pommerenke.
\newblock On the integral means of the derivative of a univalent function.
  {II}.
\newblock {\em Bull. London Math. Soc.}, 17(6):565--570, 1985.

\bibitem{MR875238}
Ch. Pommerenke.
\newblock The growth of the derivative of a univalent function.
\newblock In {\em The {B}ieberbach conjecture ({W}est {L}afayette, {I}nd.,
  1985)}, volume~21 of {\em Math. Surveys Monogr.}, pages 143--152. Amer. Math.
  Soc., Providence, RI, 1986.

\bibitem{MR0008257}
H.~Rademacher.
\newblock On the {B}loch-{L}andau constant.
\newblock {\em Amer. J. Math.}, 65:387--390, 1943.

\bibitem{MR0210910}
K.~V. Rajeswara~Rao.
\newblock On a generalized corona problem.
\newblock {\em J. Analyse Math.}, 18:277--278, 1967.

\bibitem{MR0444974}
N.~V. Rao and D.~F. Shea.
\newblock Growth problems for subharmonic functions of finite order in space.
\newblock {\em Trans. Amer. Math. Soc.}, 230:347--370, 1977.

\bibitem{MR0499172}
E.~Reich.
\newblock An extremum problem for analytic functions with area norm.
\newblock {\em Ann. Acad. Sci. Fenn. Ser. A I Math.}, 2:429--445, 1976.

\bibitem{MR0361065}
E.~Reich and K.~Strebel.
\newblock Extremal quasiconformal mappings with given boundary values.
\newblock pages 375--391, 1974.

\bibitem{RempePHD}
L.~Rempe.
\newblock Doctoral thesis.
\newblock {\em Ch.-Albrechts-Universit\"at zu Kiel}, 2003.

\bibitem{MR2069014}
L.~Rempe.
\newblock On a question of {H}erman, {B}aker and {R}ippon concerning {S}iegel
  disks.
\newblock {\em Bull. London Math. Soc.}, 36(4):516--518, 2004.

\bibitem{MR0051308}
M.~S. Robertson.
\newblock A coefficient problem for functions regular in an annulus.
\newblock {\em Canadian J. Math.}, 4:407--423, 1952.

\bibitem{MR1545292}
W.~Rogosinski.
\newblock \"{U}ber positive harmonische {E}ntwicklungen und typisch-reelle
  {P}otenzreihen.
\newblock {\em Math. Z.}, 35(1):93--121, 1932.

\bibitem{MR692530}
M.~Roitman.
\newblock On roots of polynomials and of their derivatives.
\newblock {\em J. London Math. Soc. (2)}, 27(2):248--256, 1983.

\bibitem{MR0377098}
J.~P. Rosay.
\newblock Sur la non-factorisation des \'el\'ements de l'espace de {H}ardy
  {$H^{1}(U^{2})$}.
\newblock {\em Illinois J. Math.}, 19:479--482, 1975.

\bibitem{Rosenbloom}
P.~C. Rosenbloom.
\newblock On sequences of polynomials, especially sections of power series.
\newblock {\em Stanford University}, Doctoral Thesis, 1949.

\bibitem{MR831388}
J.~Rossi.
\newblock Second order differential equations with transcendental coefficients.
\newblock {\em Proc. Amer. Math. Soc.}, 97(1):61--66, 1986.

\bibitem{MR1304116}
J.~Rossi.
\newblock A sharp result concerning cercles de remplissage.
\newblock {\em Ann. Acad. Sci. Fenn. Ser. A I Math.}, 20(1):179--185, 1995.

\bibitem{MR1417690}
J.~Rossi and S.~Wang.
\newblock The radial oscillation of solutions to {ODE}'s in the complex domain.
\newblock {\em Proc. Edinburgh Math. Soc. (2)}, 39(3):473--483, 1996.

\bibitem{MR787658}
J.~Rossi and A.~Weitsman.
\newblock The growth of entire and harmonic functions along asymptotic paths.
\newblock {\em Comment. Math. Helv.}, 60(1):1--16, 1985.

\bibitem{MR0107008}
W.~C. Royster.
\newblock Coefficient problems for functions regular in an ellipse.
\newblock {\em Duke Math. J.}, 26:361--371, 1959.

\bibitem{MR0203036}
W.~C. Royster.
\newblock A {P}oisson integral formula for the ellipse and some applications.
\newblock {\em Proc. Amer. Math. Soc.}, 15:661--670, 1964.

\bibitem{MR0183866}
W.~C. Royster.
\newblock On the univalence of a certain integral.
\newblock {\em Michigan Math. J.}, 12:385--387, 1965.

\bibitem{MR0412395}
L.~A. {R}ubel, A.~L. Shields, and B.~A. Taylor.
\newblock Mergelyan sets and the modulus of continuity of analytic functions.
\newblock {\em J. Approximation Theory}, 15(1):23--40, 1975.

\bibitem{MR0237753}
Z.~Rubinstein.
\newblock On a problem of {I}lyeff.
\newblock {\em Pacific J. Math.}, 26:159--161, 1968.

\bibitem{Francia}
J.~L. Rubio De~Francia.
\newblock {\em Weighted norm inequalities and vector-valued inequalities
  ({P}roc. {C}onf. {H}armonic {A}nalysis, {M}inneapolis 1981)}.
\newblock Lecture Notes in Mathematics 908. Spring-Verlag, Berlin, 1981.

\bibitem{MR0079093}
W.~Rudin.
\newblock The radial variation of analytic functions.
\newblock {\em Duke Math. J.}, 22:235--242, 1955.

\bibitem{zbMATH03283536}
W.~{Rudin}.
\newblock {Function theory in polydiscs.}
\newblock {Mathematics Lecture Note Series. New York-Amsterdam: W.A. Benjamin,
  Inc., 188 p. (1969).}, 1969.

\bibitem{MR516096}
W.~Rudin.
\newblock Inner function images of radii.
\newblock {\em Math. Proc. Cambridge Philos. Soc.}, 85(2):357--360, 1979.

\bibitem{MR0318470}
S.~Ruscheweyh.
\newblock \"{U}ber die {F}altung schlichter {F}unktionen.
\newblock {\em Math. Z.}, 128:85--92, 1972.

\bibitem{MR0382626}
S.~Ruscheweyh.
\newblock Duality for {H}adamard products with applications to extremal
  problems for functions regular in the unit disc.
\newblock {\em Trans. Amer. Math. Soc.}, 210:63--74, 1975.

\bibitem{MR601721}
S.~Ruscheweyh.
\newblock Neighborhoods of univalent functions.
\newblock {\em Proc. Amer. Math. Soc.}, 81(4):521--527, 1981.

\bibitem{MR997971}
S.~Ruscheweyh and L.~C. Salinas.
\newblock On the preservation of direction-convexity and the {G}oodman-{S}aff
  conjecture.
\newblock {\em Ann. Acad. Sci. Fenn. Ser. A I Math.}, 14(1):63--73, 1989.

\bibitem{MR0447552}
S.~Ruscheweyh and T.~Sheil-Small.
\newblock Corrigendum: ``{H}adamard products of schlicht functions and the
  {P}\'olya-{S}choenberg conjecture'' ({C}omment. {M}ath. {H}elv. {\bf 48}
  (1973), 119--135).
\newblock {\em Comment. Math. Helv.}, 48:194, 1973.

\bibitem{MR0328051}
S.~Ruscheweyh and T.~Sheil-Small.
\newblock Hadamard products of {S}chlicht functions and the
  {P}\'olya-{S}choenberg conjecture.
\newblock {\em Comment. Math. Helv.}, 48:119--135, 1973.

\bibitem{MR0203034}
F.~B. Ryan.
\newblock The set of asymptotic values of a bounded holomorphic function.
\newblock {\em Duke Math. J.}, 33:477--484, 1966.

\bibitem{MR879189}
J.~Ryan.
\newblock Hypercomplex algebras, hypercomplex analysis and conformal
  invariance.
\newblock {\em Compositio Math.}, 61(1):61--80, 1987.

\bibitem{Sadi}
A.~{Sadi}.
\newblock {Some types of regularity for the Dirichlet problem.}
\newblock {\em {Nagoya Math. J.}}, 126:103--124, 1992.

\bibitem{MR0382609}
E.~B. Saff and T.~Sheil-Small.
\newblock Coefficient and integral mean estimates for algebraic and
  trigonometric polynomials with restricted zeros.
\newblock {\em J. London Math. Soc. (2)}, 9:16--22, 1974/75.

\bibitem{MR0414968}
E.~B. Saff and R.~S. Varga.
\newblock Zero-free parabolic regions for sequences of polynomials.
\newblock {\em SIAM J. Math. Anal.}, 7(3):344--357, 1976.

\bibitem{MR659788}
M.~Sakai.
\newblock Null quadrature domains.
\newblock {\em J. Analyse Math.}, 40:144--154 (1982), 1981.

\bibitem{MR0065679}
R.~Salem and A.~Zygmund.
\newblock Some properties of trigonometric series whose terms have random
  signs.
\newblock {\em Acta Math.}, 91:245--301, 1954.

\bibitem{MR2478263}
A.~Sauer.
\newblock How to detect {H}ayman directions.
\newblock {\em Comput. Methods Funct. Theory}, 9(1):57--64, 2009.

\bibitem{MR1505083}
M.~Schiffer.
\newblock Sur un probl\`eme d'extr\'emum de la repr\'esentation conforme.
\newblock {\em Bull. Soc. Math. France}, 66:48--55, 1938.

\bibitem{MR0486436}
G.~Schmeisser.
\newblock On {I}lieff's conjecture.
\newblock {\em Math. Z.}, 156(2):165--173, 1977.

\bibitem{MR630449}
G.~Schmieder.
\newblock {F}unktionen mit vorgeschriebenen {N}ull- und {V}erzweigungsstellen
  auf {R}iemannschen {F}l\"achen.
\newblock {\em Arch. Math. (Basel)}, 37(1):72--77, 1981.

\bibitem{MR1762263}
W.~Schwick.
\newblock A note on {Z}alcman's lemma.
\newblock {\em New Zealand J. Math.}, 29(1):71--72, 2000.

\bibitem{Selberg1928}
H.~L. Selberg.
\newblock {\"U}ber einige {E}igen schaften bei der {W}erteverteilung der
  meromorphen {F}unktionen endlicher {O}rdnung.
\newblock {\em Avh. Norske Vid.-Akad. Oslo 1(N.S)}, No. 7:1--17, 1928.

\bibitem{MR2722582}
B.~Sendov.
\newblock New conjectures in the {H}ausdorff geometry of polynomials.
\newblock {\em East J. Approx.}, 16(2):179--192, 2010.

\bibitem{MR0100094}
T.~Shah.
\newblock Goluzin's number {$(3-\surd 5)/2$} is the radius of superiority in
  subordination.
\newblock {\em Sci. Record (N.S.)}, 1:219--222, 1957.

\bibitem{MR942827}
H.~S. Shapiro.
\newblock A uniform boundedness principle concerning inner functions.
\newblock {\em J. Analyse Math.}, 50:183--188, 1988.

\bibitem{MR862041}
D.~F. Shea and L.~R. Sons.
\newblock Value distribution theory for meromorphic functions of slow growth in
  the disk.
\newblock {\em Houston J. Math.}, 12(2):249--266, 1986.

\bibitem{MR0320761}
T.~Sheil-Small.
\newblock On the convolution of analytic functions.
\newblock {\em J. Reine Angew. Math.}, 258:137--152, 1973.

\bibitem{MR623496}
T.~Sheil-Small.
\newblock Applications of the {H}adamard product.
\newblock In {\em Aspects of contemporary complex analysis ({P}roc. {NATO}
  {A}dv. {S}tudy {I}nst., {U}niv. {D}urham, {D}urham, 1979)}, pages 515--523.
  Academic Press, London-New York, 1980.

\bibitem{MR803610}
T.~Sheil-Small.
\newblock On the {F}ourier series of a finitely describe. convex curve and a
  conjecture of {H}. {S}. {S}hapiro.
\newblock {\em Math. Proc. Cambridge Philos. Soc.}, 98(3):513--527, 1985.

\bibitem{MR979605}
T.~Sheil-Small.
\newblock On the zeros of the derivatives of real entire functions and
  {W}iman's conjecture.
\newblock {\em Ann. of Math. (2)}, 129(1):179--193, 1989.

\bibitem{zbMATH04097742}
T.~{Sheil-Small}.
\newblock {On the zeros of the derivatives of real entire functions and Wiman's
  conjecture.}
\newblock {\em {Ann. Math. (2)}}, 129(1):179--193, 1989.

\bibitem{MR1083443}
T.~Sheil-Small.
\newblock Constants for planar harmonic mappings.
\newblock {\em J. London Math. Soc. (2)}, 42(2):237--248, 1990.

\bibitem{MR981503}
T.~Sheil-Small and E.~M. Silvia.
\newblock Neighborhoods of analytic functions.
\newblock {\em J. Analyse Math.}, 52:210--240, 1989.

\bibitem{MR0361836}
A.~L. Shields.
\newblock Cyclic vectors in some spaces of analytic functions.
\newblock {\em Proc. Roy. Irish Acad. Sect. A}, 74:293--296, 1974.
\newblock Spectral Theory Symposium (Trinity College, Dublin, 1974).

\bibitem{MR1285397}
M.~Shishikura.
\newblock The boundary of the {M}andelbrot set has {H}ausdorff dimension two.
\newblock {\em Ast\'erisque}, (222):7, 389--405, 1994.
\newblock Complex analytic methods in dynamical systems (Rio de Janeiro, 1992).

\bibitem{MR1626737}
M.~Shishikura.
\newblock The {H}ausdorff dimension of the boundary of the {M}andelbrot set and
  {J}ulia sets.
\newblock {\em Ann. of Math. (2)}, 147(2):225--267, 1998.

\bibitem{MR3320475}
S.~Sivasubramanian, R.~Sivakumar, S.~Kanas, and S.~A. Kim.
\newblock Verification of {B}rannan and {C}lunie's conjecture for certain
  subclasses of bi-univalent functions.
\newblock {\em Ann. Polon. Math.}, 113(3):295--304, 2015.

\bibitem{MR0385155}
H.~Skoda.
\newblock Z\'eros des fonctions de la classe de {N}evanlinna dans les ouverts
  strictement pseudoconvexes.
\newblock {\em C. R. Acad. Sci. Paris S\'er. A-B}, 280(24):Aii, A1677--A1680,
  1975.

\bibitem{MR0450620}
H.~Skoda.
\newblock Valeurs au bord pour les solutions de l'op\'erateur {$d''$}, et
  caract\'erisation des z\'eros des fonctions de la classe de {N}evanlinna.
\newblock {\em Bull. Soc. Math. France}, 104(3):225--299, 1976.

\bibitem{zbMATH00008459}
Z.~{Slodkowski}.
\newblock {Holomorphic motions and polynomial hulls.}
\newblock {\em {Proc. Am. Math. Soc.}}, 111(2):347--355, 1991.

\bibitem{MR2173369}
M.~A. Snipes and L.~A. W.
\newblock Realizing step functions as harmonic measure distributions of planar
  domains.
\newblock {\em Ann. Acad. Sci. Fenn. Math.}, 30(2):353--360, 2005.

\bibitem{MR575388}
G.~Somorjai.
\newblock On asymptotic functions.
\newblock {\em J. London Math. Soc. (2)}, 21(2):297--303, 1980.

\bibitem{MR0203030}
L.~R. Sons.
\newblock Value distribution and power series with moderate gaps.
\newblock {\em Michigan Math. J.}, 13:425--433, 1966.

\bibitem{MR889119}
L.~R. Sons.
\newblock Value distribution for unbounded functions in the unit disk.
\newblock {\em Complex Variables Theory Appl.}, 7(4):337--341, 1987.

\bibitem{MR0004903}
D.~C. Spencer.
\newblock On mean one-valent functions.
\newblock {\em Ann. of Math. (2)}, 42:614--633, 1941.

\bibitem{MR0041935}
G.~Springer.
\newblock The coefficient problem for schlicht mappings of the exterior of the
  unit circle.
\newblock {\em Trans. Amer. Math. Soc.}, 70:421--450, 1951.

\bibitem{MR699317}
E.~M. Stein.
\newblock Some results in harmonic analysis in {$\mathbb{R}^{n}$}, for
  {$n\rightarrow \infty $}.
\newblock {\em Bull. Amer. Math. Soc. (N.S.)}, 9(1):71--73, 1983.

\bibitem{MR727348}
E.~M. Stein and J.~O. Str\"omberg.
\newblock Behavior of maximal functions in {$\mathbb{R}^{n}$}\ for large {$n$}.
\newblock {\em Ark. Mat.}, 21(2):259--269, 1983.

\bibitem{MR595424}
N.~Steinmetz.
\newblock \"{U}ber das {A}nwachsen der {L}\"osungen homogener algebraischer
  {D}ifferentialgleichungen zweiter {O}rdnung.
\newblock {\em Manuscripta Math.}, 32(3-4):303--308, 1980.

\bibitem{MR686638}
N.~Steinmetz.
\newblock \"{U}ber die eindeutigen {L}\"osungen einer homogenen algebraischen
  {D}ifferentialgleichung zweiter {O}rdnung.
\newblock {\em Ann. Acad. Sci. Fenn. Ser. A I Math.}, 7(2):177--188, 1982.

\bibitem{MR928649}
N.~Steinmetz.
\newblock On the zeros of {$(f^{(p)}+a_{p-1}f^{(p-1)}+\cdots+a_0f)f$}.
\newblock {\em Analysis}, 7(3-4):375--389, 1987.

\bibitem{MR1081702}
N.~Steinmetz.
\newblock Linear differential equations with exceptional fundamental sets.
  {II}.
\newblock {\em Proc. Amer. Math. Soc.}, 117(2):355--358, 1993.

\bibitem{MR0481015}
K.~Stephenson.
\newblock Omitte. values of singular inner functions.
\newblock {\em Michigan Math. J.}, 25(1):91--100, 1978.

\bibitem{MR674872}
K.~Stephenson.
\newblock Analytic functions and hypergroups of function pairs.
\newblock {\em Indiana Univ. Math. J.}, 31(6):843--884, 1982.

\bibitem{MR825925}
A.~Stray.
\newblock Minimal interpolation by {B}laschke products.
\newblock {\em J. London Math. Soc. (2)}, 32(3):488--496, 1985.

\bibitem{MR940287}
A.~Stray.
\newblock Minimal interpolation by {B}laschke products. {II}.
\newblock {\em Bull. London Math. Soc.}, 20(4):329--332, 1988.

\bibitem{MR570689}
E.~W. Stredulinsky.
\newblock Higher integrability from reverse {H}\"older inequalities.
\newblock {\em Indiana Univ. Math. J.}, 29(3):407--413, 1980.

\bibitem{MR609659}
D.~Styer and D.~J. Wright.
\newblock Results on bi-univalent functions.
\newblock {\em Proc. Amer. Math. Soc.}, 82(2):243--248, 1981.

\bibitem{MR819553}
D.~Sullivan.
\newblock Quasiconformal homeomorphisms and dynamics. {I}. {S}olution of the
  {F}atou-{J}ulia problem on wandering domains.
\newblock {\em Ann. of Math. (2)}, 122(3):401--418, 1985.

\bibitem{MR0048559}
F.~Sunyer~Balaguer.
\newblock Sur la substitution d'une valeur exceptionnelle par une propri\'et\'e
  lacunaire.
\newblock {\em Acta Math.}, 87:17--31, 1952.

\bibitem{Talpur}
M.~N.~M. Talpur.
\newblock Doctoral thesis.
\newblock {\em London University}, 1967.

\bibitem{zbMATH03493031}
M.~N.~M. {Talpur}.
\newblock {A subharmonic analogue of Iversen's theorem.}
\newblock {\em {Proc. Lond. Math. Soc. (3)}}, 31:129--148, 1975.

\bibitem{MR0393514}
M.~N.~M. Talpur.
\newblock On the growth of subharmonic functions on asymptotic paths.
\newblock {\em Proc. London Math. Soc. (3)}, 32(2):193--198, 1976.

\bibitem{MR0357765}
D.~K. Thomas.
\newblock On the coefficients of meromorphic univalent functions.
\newblock {\em Proc. Amer. Math. Sco.}, 47:161--166, 1975.

\bibitem{MR0259075}
R.~Tijdeman.
\newblock {\em On the distribution of the values of certain functions}.
\newblock Doctoral dissertation, University of Amsterdam. Universiteit van
  Amsterdam, Amsterdam, 1969.

\bibitem{MR0286986}
R.~Tijdeman.
\newblock On the number of zeros of general exponential polynomials.
\newblock {\em Nederl. Akad. Wetensch. Proc. Ser. A 74 = Indag. Math.},
  33:1--7, 1971.

\bibitem{MR3155290}
E.~C. Titchmarsh.
\newblock {\em The theory of functions}.
\newblock Oxford University Press, Oxford, 1958.
\newblock Reprint of the second (1939) edition.

\bibitem{MR0176151}
E.~C. Titchmarsh.
\newblock {\em Eigenfunction expansions associated with second-order
  differential equations. {P}art {I}}.
\newblock Second Edition. Clarendon Press, Oxford, 1962.

\bibitem{MR1047599}
K.~Tohge.
\newblock On a problem of {H}inkkanen about {H}adamard products.
\newblock {\em Kodai Math. J.}, 13(1):101--120, 1990.

\bibitem{MR629839}
V.~A. Tolokonnikov.
\newblock Estimates in the {C}arleson corona theorem, ideals of the algebra
  {$H^{\infty }$}, a problem of {S}z.-{N}agy.
\newblock {\em Zap. Nauchn. Sem. Leningrad. Otdel. Mat. Inst. Steklov. (LOMI)},
  113:178--198, 267, 1981.
\newblock Investigations on linear operators and the theory of functions, XI.

\bibitem{zbMATH02115780}
X.~{Tolsa}.
\newblock {Painlev\'e's problem and the semiadditivity of analytic capacity.}
\newblock {\em {Acta Math.}}, 190(1):105--149, 2003.

\bibitem{MR0283202}
S.~Toppila.
\newblock Some remarks on exceptional values at {J}ulia lines.
\newblock {\em Ann. Acad. Sci. Fenn. Ser. A I No.}, 456:20, 1970.

\bibitem{MR0507887}
S.~Toppila.
\newblock On the counting function for the {$a$}-values of a meromorphic
  function.
\newblock {\em Ann. Acad. Sci. Fenn. Ser. A I Math.}, 2:565--572, 1976.

\bibitem{MR0585685}
S.~Toppila.
\newblock On {N}evanlinna's characteristic functions of entire functions and
  their derivatives.
\newblock {\em Ann. Acad. Sci. Fenn. Ser. A I Math.}, 3(1):131--134, 1977.

\bibitem{MR595173}
S.~Toppila.
\newblock On the length of asymptotic paths of entire functions of order zero.
\newblock {\em Ann. Acad. Sci. Fenn. Ser. A I Math.}, 5(1):13--15, 1980.

\bibitem{MR1945294}
S.~Treil.
\newblock Estimates in the corona theorem and ideals of {$H^\infty$}: a problem
  of {T}.\ {W}olff.
\newblock {\em J. Anal. Math.}, 87:481--495, 2002.
\newblock Dedicated to the memory of T. H. Wolff.

\bibitem{MR647809}
P.~Tukia and J.~V\"ais\"al\"a.
\newblock Quasiconformal extension from dimension {$n$}\ to {$n+1$}.
\newblock {\em Ann. of Math. (2)}, 115(2):331--348, 1982.

\bibitem{MR0016488}
Y.~Tumura.
\newblock Recherches sur la distribution des valeurs des fonctions analytiques.
\newblock {\em Jap. J. Math.}, 18:797--876, 1943.

\bibitem{MR0060548}
P.~Tur\'an.
\newblock {\em Eine neue {M}ethode in der {A}nalysis und deren {A}nwendungen}.
\newblock Akad\'emiai Kiad\'o, Budapest, 1953.

\bibitem{MR0219491}
P.~Tur\'an.
\newblock Problem 153.
\newblock {\em Mat. Lapok}, 17:215, 1966.

\bibitem{MR1662226}
T.~F. Tyler.
\newblock Maximum curves and isolated points of entire functions.
\newblock {\em Proc. Amer. Math. Soc.}, 128(9):2561--2568, 2000.

\bibitem{Valiron1925}
G.~Valiron.
\newblock Sur les valeurs asymptotiques de quelques fonctions m{\'e}romorphes.
\newblock {\em Rendiconti del Circolo Matematico di Palermo (1884-1940)},
  49(3):415--421, Oct 1925.

\bibitem{MR1556886}
G.~Valiron.
\newblock Sur une classe de fonctions enti\`eres admettant deux directions de
  {B}orel d'ordre divergent.
\newblock {\em Compositio Math.}, 1:193--206, 1935.

\bibitem{MR0021976}
G.~Valiron.
\newblock Valeurs exceptionnelles et valeurs d\'eficientes des fonctions
  m\'eromorphes.
\newblock {\em C. R. Acad. Sci. Paris}, 225:556--558, 1947.

\bibitem{MR0061658}
G.~Valiron.
\newblock {\em Fonctions analytiques}.
\newblock Presses Universitaires de France, Paris, 1954.

\bibitem{varopoulos1977}
N.~Th. Varopoulos.
\newblock {BMO} functions and the $\overline\partial$-equation.
\newblock {\em Pacific J. Math.}, 71(1):221--273, 1977.

\bibitem{MR0289800}
W.~A. Veech.
\newblock A converse to {G}auss' theorem.
\newblock {\em Bull. Amer. Math. Soc.}, 78:444--446, 1972.

\bibitem{MR0310269}
W.~A. Veech.
\newblock A zero-one law for a class of random walks and a converse to {G}auss'
  mean value theorem.
\newblock {\em Ann. of Math. (2)}, 97:189--216, 1973.

\bibitem{MR0486534}
M.~von Renteln.
\newblock Ideals in the {N}evanlinna class {$N$}.
\newblock {\em Mitt. Math. Sem. Giessen}, (Heft 123):57--65, 1977.
\newblock Dem Andenken an Karl Maruhn gewidmet.

\bibitem{MR0293298}
A.~A. \v{S}aginjan.
\newblock The uniform and tangent harmonic approximation of continuous
  functions on arbitrary sets.
\newblock {\em Mat. Zametki}, 9:131--142, 1971.

\bibitem{MR608334}
M.~Vuorinen.
\newblock On the existence of angular limits of {$n$}-dimensional
  quasiconformal mappings.
\newblock {\em Ark. Mat.}, 18(2):157--180, 1980.

\bibitem{Vuorinen}
M.~Vuorinen.
\newblock Queries: {N}o. 249.
\newblock {\em Notices Amer. Math. Soc.}, 28(7):607, 1981.

\bibitem{MR571725}
G.~Wagner.
\newblock On a problem of {E}rd\"os in {D}iophantine approximation.
\newblock {\em Bull. London Math. Soc.}, 12(2):81--88, 1980.

\bibitem{Wahlund1929}
A.~Wahlund.
\newblock {\"{U}}ber einen {Z}usammenhang zwischen dem {M}aximalbetrage der
  ganzen {F}unktion und seiner unteren {G}renze nach dem {J}ensensche theoreme.
\newblock {\em Ark. Mat.}, 21A(23):34pp, 1929.

\bibitem{MR1867264}
B.~L. Walden and L.~A. Ward.
\newblock Asymptotic behaviour of distributions of harmonic measure for planar
  domains.
\newblock {\em Complex Variables Theory Appl.}, 46(2):157--177, 2001.

\bibitem{MR0385047}
P.~L. Walker.
\newblock On rearranging maximal functions in {$\mathbb{R}^{n}$}.
\newblock {\em Proc. Edinburgh Math. Soc. (2)}, 19(4):363--369, 1974/75.

\bibitem{zbMATH01159053}
Y.~{Wang} and M.~{Fang}.
\newblock {Picard values and normal families of meromorphic functions with
  multiple zeros.}
\newblock {\em {Acta Math. Sin., New Ser.}}, 14(1):17--26, 1998.

\bibitem{MR0215980}
M.~Watson.
\newblock On functions that are bivalent in the unit circle.
\newblock {\em J. Analyse Math.}, 17:383--409, 1966.

\bibitem{MR0145052}
M.~Weiss and G.~Weiss.
\newblock On the {P}icard property of lacunary power series.
\newblock {\em Studia Math.}, 22:221--245, 1962/1963.

\bibitem{Weitsman1969}
A.~Weitsman.
\newblock Meromorphic functions with maximal deficiency sum and a conjecture of
  {F}. {N}evanlinna.
\newblock {\em Acta Mathematica}, 123(1):115, Dec 1969.

\bibitem{MR0387597}
A.~Weitsman.
\newblock A theorem on {N}evanlinna deficiencies.
\newblock {\em Acta Math.}, 128(1-2):41--52, 1972.

\bibitem{MR0387598}
A.~Weitsman.
\newblock Meromorphic functions with large sums of deficiencies.
\newblock pages 133--135. London Math. Soc. Lecture Note Ser., No. 12, 1974.

\bibitem{MR554396}
A.~Weitsman.
\newblock A symmetry property of the {P}oincar\'e metric.
\newblock {\em Bull. London Math. Soc.}, 11(3):295--299, 1979.

\bibitem{MR847954}
A.~Weitsman.
\newblock Symmetrization and the {P}oincar\'e metric.
\newblock {\em Ann. of Math. (2)}, 124(1):159--169, 1986.

\bibitem{Wermer}
J.~Wermer.
\newblock Non-rectifiable simple closed curve.
\newblock {\em Amer. Math. Monthly}, 64:372, 1957.

\bibitem{MR0394218}
J.~Wermer.
\newblock {\em Banach algebras and several complex variables}.
\newblock Springer-Verlag, New York-Heidelberg, second edition, 1976.
\newblock Graduate Texts in Mathematics, No. 35.

\bibitem{MR0165476}
G.~T. Whyburn.
\newblock {\em Topological analysis}.
\newblock Second, revised edition. Princeton Mathematical Series, No. 23.
  Princeton University Press, Princeton, N.J., 1964.

\bibitem{MR575386}
D.~R. Wilken and J.~Feng.
\newblock A remark on convex and starlike functions.
\newblock {\em J. London Math. Soc. (2)}, 21(2):287--290, 1980.

\bibitem{MR0261006}
J.~Winkler.
\newblock Zum {V}erzweigungsindex der {$a$}-{S}tellen ganzer und meromorpher
  {F}unktionen.
\newblock {\em Math. Z.}, 113:353--362, 1970.

\bibitem{MR0288268}
J.~Winkler.
\newblock \"{U}ber den {V}erzweigungsindex bei ganzen {F}unktionen.
\newblock {\em Manuscripta Math.}, 4:135--148, 1971.

\bibitem{MR0077620}
H.~Wittich.
\newblock {\em Neuere {U}ntersuchungen \"uber eindeutige analytische
  {F}unktionen}.
\newblock Ergebnisse der Mathematik und ihrer Grenzgebiete (N.F.), Heft 8.
  Springer-Verlag, Berlin-G\"ottingen-Heidelberg, 1955.

\bibitem{MR0197828}
H.~Wittich.
\newblock Zur {T}heorie linearer {D}ifferentialgleichungen im {K}omplexen.
\newblock {\em Ann. Acad. Sci. Fenn. Ser. A I No.}, 379:19, 1966.

\bibitem{MR0244490}
H.~Wittich.
\newblock {\em Neuere {U}ntersuchungen \"uber eindeutige analytische
  {F}unktionen}.
\newblock Zweite, korrigierte Auflage. Ergebnisse der Mathematik und ihrer
  Grenzgebiete, Band 8. Springer-Verlag, Berlin-New York, 1968.

\bibitem{MR0000288}
F.~Wolf.
\newblock An extension of the {P}hragm\'en-{L}indel\"of theorem.
\newblock {\em J. London Math. Soc.}, 14:208--216, 1939.

\bibitem{WolffT}
T.~H. Wolff.
\newblock A refinement of the corona theorem.
\newblock {\em Linear and complex analysis problem book}, 1043(Havin, V. P. and
  Hru\v{s}\v{c}\"ev, S. V. and Nikolv{s}ki\u\i , N. K. (eds.)), 1984.

\bibitem{zbMATH00799747}
T.~H. {Wolff}.
\newblock {Counterexamples with harmonic gradients in $\mathbb R\sp 3$.}
\newblock In {\em {Essays on Fourier analysis in honor of Elias M. Stein.
  Proceedings of the Princeton conference in harmonic analysis held at
  Princeton University, Princeton, NJ, May 13-17, 1991 in honor of Elias M.
  Stein's sixtieth birthday}}, pages 321--384. Princeton, NJ: Princeton Univ.
  Press, 1995.

\bibitem{MR1979771}
T.~H. Wolff.
\newblock Counterexamples to two variants of the {H}elson-{S}zeg\"o theorem.
\newblock {\em J. Anal. Math.}, 88:41--62, 2002.
\newblock Dedicated to the memory of Tom Wolff.

\bibitem{MR513884}
J.~M.~G. Wu.
\newblock Comparisons of kernel functions, boundary {H}arnack principle and
  relative {F}atou theorem on {L}ipschitz domains.
\newblock {\em Ann. Inst. Fourier (Grenoble)}, 28(4):147--167, vi, 1978.

\bibitem{MR825926}
J.~M.~G. Wu.
\newblock Length of paths for subharmonic functions.
\newblock {\em J. London Math. Soc. (2)}, 32(3):497--505, 1985.

\bibitem{MR2140307}
L.~Z. Y.
\newblock Meromorphic functions and also their first two derivatives have the
  same zeros.
\newblock {\em Ann. Acad. Sci. Fenn. Math.}, 30(1):205--218, 2005.

\bibitem{MR0507738}
S.~Yamashita.
\newblock On the {J}ohn constant.
\newblock {\em Math. Z.}, 161(2):185--188, 1978.

\bibitem{MR0457415}
C.~C. Yang.
\newblock A problem on polynomials.
\newblock {\em Rev. Roumaine Math. Pures Appl.}, 22(5):595--598, 1977.

\bibitem{MR653387}
L.~Yang.
\newblock Meromorphic functions and their derivatives.
\newblock {\em J. London Math. Soc. (2)}, 25(2):288--296, 1982.

\bibitem{MR576357}
L.~Yang and G.~H. Zhang.
\newblock Deficient values and asymptotic values of entire functions.
\newblock {\em Sci. Sinica}, (Special Issue II on Math.):190--203, 1979.

\bibitem{MR777374}
J.~C. Yoccoz.
\newblock Conjugaison diff\'erentiable des diff\'eomorphismes du cercle dont le
  nombre de rotation v\'erifie une condition diophantienne.
\newblock {\em Ann. Sci. \'Ecole Norm. Sup. (4)}, 17(3):333--359, 1984.

\bibitem{MR1367353}
J.~C. Yoccoz.
\newblock Th\'eor\`eme de {S}iegel, nombres de {B}runo et polyn\^omes
  quadratiques.
\newblock {\em Ast\'erisque}, (231):3--88, 1995.
\newblock Petits diviseurs en dimension $1$.

\bibitem{MR0348084}
L.~Zalcman.
\newblock Analyticity and the {P}ompeiu problem.
\newblock {\em Arch. Rational Mech. Anal.}, 47:237--254, 1972.

\bibitem{MR0335835}
L.~Zalcman.
\newblock Mean values and differential equations.
\newblock {\em Israel J. Math.}, 14:339--352, 1973.

\bibitem{MR0379852}
L.~Zalcman.
\newblock A heuristic principle in complex function theory.
\newblock {\em Amer. Math. Monthly}, 82(8):813--817, 1975.

\bibitem{MR3539851}
L.~Zalcman.
\newblock A tale of three theorems.
\newblock {\em Amer. Math. Monthly}, 123(7):643--656, 2016.

\bibitem{zbMATH04212730}
R.~{Zeinstra}.
\newblock {Zeros and regular growth of Laplace transforms along curves.}
\newblock {\em {J. Reine Angew. Math.}}, 424:1--15, 1992.

\bibitem{MR2972625}
F.~Zhao, Z.~Fu, and S.~Lu.
\newblock Endpoint estimates for {$n$}-dimensional {H}ardy operators and their
  commutators.
\newblock {\em Sci. China Math.}, 55(10):1977--1990, 2012.

\bibitem{MR1772425}
J.~H. Zheng.
\newblock Singularities and wandering domains in iteration of meromorphic
  functions.
\newblock {\em Illinois J. Math.}, 44(3):520--530, 2000.

\bibitem{MR1317004}
J.~H. Zhu.
\newblock Hayman direction of meromorphic functions.
\newblock {\em Kodai Math. J.}, 18(1):37--43, 1995.

\bibitem{MR517191}
I.~V. {Zhuravlev}.
\newblock {Some sufficient conditions for the quasi-conformal extension of
  analytic functions.}
\newblock {\em {Sov. Math., Dokl.}}, 19:1549--1552, 1978.

\bibitem{MR0291462}
V.~V. Zimogljad.
\newblock The order of growth of entire transcendental solutions of second
  order algebraic differential equations.
\newblock {\em Mat. Sb. (N.S.)}, 85(127):286--302, 1971.

\end{thebibliography}
\bibliographystyle{plain}

\vspace{1cm}
\begin{table}[H]
\setlength\extrarowheight{2pt}
\centering
\centerline{Table 1: \textit{Abbreviations for reference documents}}
\vspace{0.5cm}
\begin{tabularx}{\textwidth}{|c|C|}
\hline
A     & W. K. Hayman, Research problems in function theory, Athlone Press, London, 1967. \\ \hline
B& J. G. Clunie and W. K. Hayman (eds.), Symposium on complex analysis, Canterbury, 1973, LMS Lect. Notes ser. 12, Cambridge Univ. Press, 1974, 143-180.\\\hline
C& J. M. Anderson, K. F. Barth, D. A. Brannan and W. K. Hayman, Research problems in complex analysis, Bull. LMS, 9 (1977), 129-162.\\\hline
D, E, F, G, H& D. M. Campbell, J. G. Clunie and W. K. Hayman, Aspects of contemporary complex analysis (eds: D. A. Brannan and J. G. Clunie, Acad. Press 1980), 527-572.\\\hline
I, J& K. F. Barth, D. A. Brannan and W. K. Hayman, Research problems in complex analysis, Bull. LMS, 16 (1984), 490-517.\\\hline
K, L& D. A. Brannan and W. K. Hayman, Research problems in complex analysis, Bull. LMS, 21 (1989), 1-35.\\\hline
M& A. Eremenko, Progress in Entire and Meromorphic Functions, available from his homepage.\\ \hline
N&W. K. Hayman and E. F. Lingham, Research problems in function theory (This edition) (TBC: add details)\\ \hline
\end{tabularx}
\end{table}

\vspace{1cm}
\begin{table}[H]
\centering
\centerline{Table 2: \textit{List of problems proposed}}
\vspace{0.5cm}
\begin{tabular}{|cccccccc|}
\hline
Topic & A    & B     & C     & H     & I     & L  &N   \\ \hline
1     & 1-23 & 24-29 & 30    & -     & 31-36 & 37-43 &TBC\\ \hline
2     & 1-32 & 33-46 & 47-57 & 58-64 & 65-68 & 69-90&TBC \\ \hline
3     & 1-10& 11-15 & 16-18 & 19-20 & 21-30 & 31-35&TBC\\ \hline
4     & 1-21 & 22-24 & 25-27 & 28-31 & - &  -&TBC\\ \hline
5     & 1-21 & 22-37 & 38-58 & 59-66 & 67-70 & 71-79&TBC \\ \hline
6     & 1-26 & 27-31 & 32-63 & 64-82 & 83-95 & 96-117&TBC \\ \hline
7     & 1-8 & 9-22 & 23-32 & 33-52 & 53-67 & 68-83 &TBC\\ \hline
8     & - & - & 1-17 & - & 18-22 & 23-26&TBC\\ \hline
9     & - & - & 1-7 & 8-10 & 11 & 12-16 &TBC\\ \hline
\end{tabular}
\end{table}

\vspace{1cm}
\centerline{Table 3: \textit{Comments on problems}}
\begin{longtable}[c]{|cccccccccc|}

\hline
Problem & 1          & 2          & 3    & 4       & 5           & 6          & 7    & 8 & 9                     \\ 
\hline

\endhead
\hline
\endfoot
\endlastfoot

1              & B          & B          & C, G & D       & B, C        & B, K       &      &   &                       \\
2              & B          & B          & B, F & B, D. F &             & B, C, J, K &      & D & K                     \\
3              & B, J       & D, J       & B    &         & B, C        &            &      & D & D                     \\
4              & B, F       & B, F       &      &         & C, G        & C          &      &   & D                     \\
5              &            & B, F       & B    & B, D, K & K           & B, E, K    & B    & D &                       \\
6              & E, K, M    & C, M       & D    & C, D    & K           &            & B    & D &                       \\
7              & B          & D, J       & B    & J       & C, G, J, K, & E, K       & B, E &   & D                     \\
8              & B          & B, C, E, G &      &         & B           & B, E       &      &   &                       \\
9              & B          & E, M       &      & B       & E           & B          &      & D &                       \\
10             & B, K       & D, J       & B, K &         &             & B          & D &   &                       \\
11             & B, K       & M          & D    &         & J           & B          & J    &   &                       \\
12             & D          & M          & D    & B       & B, J        & E          & D    &   &                       \\
13             & B          & J          & C, G &         &             & B, E       &      &   &                       \\
14             & B          & C, G       &      & C, D    & B, E        & B, E       & D    &   &                       \\
15             & B          &            &      &         & D, K        & B, E       &      &   &                       \\
16             & C, M       & B, C       &      &         & B           & B          & D, J & D &                       \\
17             & B          & B, E, F, J & D    & B, F    &             &            & J    &   &                       \\
18             & B, C, G, M & C          &      & D, K    &             &            & D, F &   &                       \\
19             & D, M       & J, K, M    &      &         & J           &            &      &   &                       \\
20             & D, E, M    & M          &      & B, F    &             &            & D    & K &                       \\
21             & C, K       & B, M       &      &         &             &            &      &   &                       \\
22             & B          & B, J       &      &         &             &            & D    & K &                       \\
23             & B, J, K    & M          & K    &         & D           & E          &      &   &                       \\
24             & C          & C          &      &         & D           &            & J    &   &                       \\
25             & C, G       & B, K, M    & K    &         &             &            &      &   &                       \\
26             &            & J, M       & K    & J       &             & J          &      &   &                       \\
27             & D          &            & K    &         & J, K        &            &      &   &                       \\
28             & C, D, J    & B, E, F    &      &         & D, F        &            &      &   &                       \\
29             & K          & B, F       &      & J       &             &            &      &   &                       \\
30             &            & B          &      &         & D, F        & C, D, K    & D    &   &                       \\
31             &            & B          &      &         & D           & F          &      &   &                       \\
32             &            & C, G, K    &      &         &             &            &      &   &                       \\
33             & M          &            &      &         & C, D, G     &            &      &   &                       \\
34             & K          & M          &      &         &             & D          &      &   &   \\
35             & M          & M          &      &         & C           & D          &      &   &   \\
36             &            &            &      &         & D           & K          &      &   &   \\
37             & M          & J          &      &         & D           & J          &      &   &   \\
38             &            & F, M       &      &         & D           & J          &      &   &   \\
39             &            &            &      &         &             & K          &      &   &   \\
40             &            & C, J       &      &         & K           &            & K    &   &   \\
41             &            & J          &      &         & K           & D          & K    &   &   \\
42             & M          & D, J       &      &         &             & K          &      &   &   \\
43             &            & D          &      &         &             & D          & J    &   &   \\
44             &            & C, G       &      &         &             & D          &      &   &   \\
45             &            &            &      &         &             & D          & K    &   &   \\
46             &            &            &      &         &             & D          & J    &   &   \\
47             &            &            &      &         & J           & D          &      &   &   \\
49             &            &            &      &         & D, J        &            &      &   &   \\
50             &            &            &      &         & D           &            &      &   &   \\
52             &            & M          &      &         &             &            &      &   &   \\
54             &            &            &      &         & D           &            &      &   &   \\
55             &            & D          &      &         &             &            &      &   &   \\
56             &            &            &      &         &             &            & K    &   &   \\
57             &            & D          &      &         &             & D          &      &   &                       \\
58             &            & J          &      &         & G, J        &            &      &   &                       \\
59             &            &            &      &         &             &            & K    &   &                       \\
60             &            &            &      &         &             & D          & K    &   &                       \\
61             &            &            &      &         & J           &            &      &   &                       \\
62             &            & J          &      &         &             &            &      &   &                       \\
63             &            & J          &      &         &             &            &      &   &                       \\
64             &            & M          &      &         &             &            &      &   &                       \\
65             &            & M          &      &         & J           &            &      &   &                       \\
66             &            & K          &      &         & K           & J          &      &   &                       \\
67             &            &            &      &         &             & J          &      &   &                       \\
69             &            & M          &      &         &             & J          &      &   &                       \\
70             &            & M          &      &         &             & K          &      &   &                       \\
71             &            & M          &      &         &             &            &      &   &                       \\
72             &            & M          &      &         &             &            &      &   &                       \\
73             &            &            &      &         &             & J          &      &   &                       \\
75             &            & M          &      &         &             &            &      &   &                       \\
78             &            & M          &      &         &             &            &      &   &                       \\
79             &            & M          &      &         &             & K          &      &   &                       \\
83             &            & M          &      &         &             &            &      &   &                       \\
84             &            &            &      &         &             & K          &      &   &                       \\
86             &            & M          &      &         &             &            &      &   &                       \\
87             &            & M          &      &         &             & K          &      &   &                       \\
88             &            & M          &      &         &             &            &      &   &                       \\
90             &            &            &      &         &             & K          &      &   &                       \\
94             &            &            &      &         &             & K          &      &   &                       \\ \hline
\end{longtable}

\end{document}